\documentclass[preprint,11pt]{elsarticle}
\usepackage[bookmarks=true,colorlinks=true,linkcolor=blue]{hyperref}

\pdfoutput=1

\biboptions{sort&compress}

\usepackage{color}
\definecolor{jwbGreen}{rgb}{0, .6, 0}
\definecolor{darkBlue}{rgb}{.1, 0, .6}

\usepackage{soul}

\usepackage[usenames,dvipsnames,svgnames,table]{xcolor}

\usepackage{amssymb,amsmath}
\usepackage{amsthm}

\usepackage{calc}

\usepackage{calc}
\usepackage[margin=1.in]{geometry}

\renewcommand{\url}[1]{}
\newcommand{\citeCount}[1]{}
\newcommand{\f}[2]{\frac{#1}{#2}}
\newtheorem{theorem}{Theorem}

\newtheorem{lemma}{Lemma}

\newtheorem*{AMPInterfaceCondition}{AMP Interface Condition}

\newcommand{\esup}{{(e)}}
\newcommand{\psup}{{(p)}}

\newcommand{\dt}{\Delta t}

\newcommand{\dn}{\Delta n}

\newcommand{\bogus}[1]{{}}

\newenvironment{myIndent}%
 {\list{}{\leftmargin=0.1in\rightmargin=0.1in}\item[]}%
  {\endlist}

\newenvironment{myProcedure}[1]%
{\noindent\textbf{Procedure~}{#1}\begin{myIndent}\em}
{\end{myIndent}}

\newcommand{\eqdef}{\overset{{\rm def}}{=}}

\newcommand{\qvb}{\qv_b}
\newcommand{\cma}{\,,\,}

\newcommand{\av}{\mathbf{ a}}
\newcommand{\bv}{\mathbf{ b}}

\newcommand{\fv}{\mathbf{ f}}
\newcommand{\gv}{\mathbf{ g}}

\newcommand{\iv}{\mathbf{ i}}

\newcommand{\nv}{\mathbf{ n}}

\newcommand{\qv}{\mathbf{ q}}
\newcommand{\rv}{\mathbf{ r}}

\newcommand{\tv}{\mathbf{ t}}

\newcommand{\vv}{\mathbf{ v}}
\newcommand{\wv}{\mathbf{ w}}
\newcommand{\xv}{\mathbf{ x}}
\newcommand{\yv}{\mathbf{ y}}

\newcommand{\Dv}{\mathbf{ D}}

\newcommand{\Fv}{\mathbf{ F}}
\newcommand{\Gv}{\mathbf{ G}}

\newcommand{\Iv}{\mathbf{ I}}

\newcommand{\Tv}{\mathbf{ T}}

\newcommand{\Wv}{\mathbf{ W}}

\newcommand{\half}{{1\over2}}

\newcommand{\Real}{{\mathbb R}}

\newcommand{\zerov}{\mathbf{0}}

\newcommand{\Dc}{{\mathcal D}}

\newcommand{\Fc}{{\mathcal F}}
\newcommand{\Gc}{{\mathcal G}}
\newcommand{\Gs}{{\mathcal G}}

\newcommand{\omegav}{\boldsymbol{\omega}}
\newcommand{\tauv}{\boldsymbol{\tau}}

\newcommand{\sigmav}{\boldsymbol{\sigma}}

\newcommand{\grad}{\nabla}

\newcommand{\tableFont}{\scriptsize}

\newcommand{\num}[2]{#1e#2} 

\clearpage

\newcommand{\ms}{\bar{m}}
\newcommand{\rhos}{\rho_b}
\newcommand{\rhob}{\rho_b}

\newcommand{\vs}{\bar{v}}

\newcommand{\ds}{\Delta s}

\newcommand{\OmegaF}{\Omega}

\newcommand{\OmegaB}{{\Omega_b}}
\newcommand{\partialB}{\Gamma_b}
\newcommand{\GammaB}{\Gamma_b}

\newcommand{\strutt}{\rule{0pt}{9pt}}

\newcommand{\nd}{{n_d}} 

\newcommand{\Gcrd}{\Gc_{\rm{rd}}}

\newcommand\tnv{{\mathbf{t}}}

\renewcommand\tv{\tnv}

\newcommand{\mrb}{m_{b}}
\newcommand{\xb}{x_b}
\newcommand{\xvb}{\xv_b}
\newcommand{\vvb}{\vv_b}
\newcommand{\avb}{\av_b}
\newcommand{\bvb}{\bv_b}
\newcommand{\xvcm}{\xv_b}
\newcommand{\vvcm}{\vv_b}
\newcommand{\avcm}{\av_b}
\newcommand{\dotxvcm}{\dot\xv_b}
\newcommand{\dotvvcm}{\dot\vv_b}

\newcommand{\omegavb}{\omegav_b}
\newcommand{\Eb}{E_b}
\newcommand{\fvbe}{\fv_e}
\newcommand{\gvbe}{\gv_e}

\newcommand{\mb}{{m_b}}
\newcommand{\mass}{{m_b}}

\newcommand{\adc}{\beta_d}

\newcommand{\Da}{\Dc}
\newcommand{\Dvv}{\Da^{v v}}
\newcommand{\Dvw}{\Da^{v \omega}}
\newcommand{\Dwv}{\Da^{\omega v}}
\newcommand{\Dww}{\Da^{\omega \omega}}
\newcommand{\Dvva}{\tilde\Da^{v v}}
\newcommand{\Dvwa}{\tilde\Da^{v \omega}}
\newcommand{\Dwva}{\tilde\Da^{\omega v}}
\newcommand{\Dwwa}{\tilde\Da^{\omega \omega}}
\newcommand{\Ib}{\Iv_b}

\newcommand{\dArea}{dS}

\newcommand{\ampRB}{AMP-RB}

\newcommand{\alphas}{\bar{\alpha}}

\newcommand{\ADtensor}{\Dv}
\newcommand{\ADtensora}{\tilde\Dv}

\newcommand{\Wvv}{\Wv^{v}}
\newcommand{\Wvw}{\Wv^{\omega}}

\newlength{\ycbTop}
\newlength{\ycbMid}%

\newcommand{\bodyStepIComment}{Preliminary body evolution step}

\def\aa{r_1}
\def\bb{r_2}

\newcommand{\tpRB}{{TP-RB}}

\newcommand{\fap}{\fv}

\newcommand{\solidWidth}{\pistonWidth}
\newcommand{\solidHeight}{\pistonHeight}
\newcommand{\solidArea}{A_b}
\newcommand{\rampFunction}{{\cal R}}

\newcommand{\channelWidth}{{x_c}}

\newcommand{\channelBottom}{{y_0}}
\newcommand{\channelTop}{{y_1}}
\newcommand{\diskRadius}{{R_b}}

\newcommand{\errFormat}[1]{$E_j^{(#1)}$}

\newcommand{\hj}{h^{(j)}}

\newcommand{\channelRight}{L}
\newcommand{\channelHeight}{H}
\newcommand{\pistonWidth}{W_b}
\newcommand{\pistonHeight}{H_b}
\newcommand{\Gcp}{\Gc_p}
\newcommand{\Amplitude}{\alpha_b}

\usepackage{tikz}

\newlength{\tfwidth}
\newlength{\tfheight}
\newlength{\tfxa}
\newlength{\tfxb}
\newlength{\tfya}
\newlength{\tfyb}
\newcommand{\trimFigWithBox}[6]{%
\setlength\fboxsep{0pt}%
\setlength\fboxrule{1.0pt}
\fbox{\includegraphics[width=#2, clip, trim=#3 #4 #5 #6]{#1}}%
}
\newcommand{\trimFigNoBox}[6]{%
\setlength\fboxsep{1pt}
\setlength\fboxrule{0.0pt}
\fbox{\includegraphics[width=#2, clip, trim=#3 #4 #5 #6]{#1}}%
}
\newcommand{\trimFigHeightWithBox}[6]{%
\setlength\fboxsep{0pt}%
\setlength\fboxrule{1.0pt}
\fbox{\includegraphics[height=#2, clip, trim=#3 #4 #5 #6]{#1}}%
}
\newcommand{\trimFigHeightNoBox}[6]{%
\setlength\fboxsep{1pt}
\setlength\fboxrule{0.0pt}
\fbox{\includegraphics[height=#2, clip, trim=#3 #4 #5 #6]{#1}}%
}

\newcommand{\trimFig}[6]{%
\setlength{\tfwidth}{(#2+#2*\real{#3})+#2*\real{#4}}
\setlength{\tfheight}{(#2+#2*\real{#5})+#2*\real{#6}}%
\setlength{\tfxa}{\tfwidth*\real{#3}}%
\setlength{\tfxb}{\tfwidth*\real{#4}}%
\setlength{\tfya}{\tfheight*\real{#5}}%
\setlength{\tfyb}{\tfheight*\real{#6}}%
\trimFigNoBox{#1}{#2}{\tfxa}{\tfya}{\tfxb}{\tfyb}%
}

\newsavebox\figBox

\newcommand{\trimw}[6]{%
\sbox\figBox{\includegraphics{#1}}
\setlength{\tfwidth}{\the\wd\figBox}
\setlength{\tfheight}{\the\ht\figBox}
\setlength{\tfxa}{\tfwidth*\real{#3}}%
\setlength{\tfxb}{\tfwidth*\real{#4}}%
\setlength{\tfya}{\tfheight*\real{#5}}%
\setlength{\tfyb}{\tfheight*\real{#6}}%
\trimFigNoBox{#1}{#2}{\tfxa}{\tfya}{\tfxb}{\tfyb}%
}

\newcommand{\trimwb}[6]{%

\sbox\figBox{\includegraphics{#1}}
\setlength{\tfwidth}{\the\wd\figBox}
\setlength{\tfheight}{\the\ht\figBox}
\setlength{\tfxa}{\tfwidth*\real{#3}}%
\setlength{\tfxb}{\tfwidth*\real{#4}}%
\setlength{\tfya}{\tfheight*\real{#5}}%
\setlength{\tfyb}{\tfheight*\real{#6}}%
\trimFigWithBox{#1}{#2}{\tfxa}{\tfya}{\tfxb}{\tfyb}%
}

\newcommand{\trimh}[6]{%
\sbox\figBox{\includegraphics{#1}}
\setlength{\tfwidth}{\the\wd\figBox}
\setlength{\tfheight}{\the\ht\figBox}
\setlength{\tfxa}{\tfwidth*\real{#3}}%
\setlength{\tfxb}{\tfwidth*\real{#4}}%
\setlength{\tfya}{\tfheight*\real{#5}}%
\setlength{\tfyb}{\tfheight*\real{#6}}%
\trimFigHeightNoBox{#1}{#2}{\tfxa}{\tfya}{\tfxb}{\tfyb}%
}

\newcommand{\trimhb}[6]{%

\sbox\figBox{\includegraphics{#1}}
\setlength{\tfwidth}{\the\wd\figBox}
\setlength{\tfheight}{\the\ht\figBox}
\setlength{\tfxa}{\tfwidth*\real{#3}}%
\setlength{\tfxb}{\tfwidth*\real{#4}}%
\setlength{\tfya}{\tfheight*\real{#5}}%
\setlength{\tfyb}{\tfheight*\real{#6}}%
\trimFigHeightWithBox{#1}{#2}{\tfxa}{\tfya}{\tfxb}{\tfyb}%
}
\usepackage{algorithm} 
\usepackage{algpseudocode} 

\begin{document}

\small

\begin{frontmatter}
\title{
A stable partitioned FSI algorithm for rigid bodies and incompressible flow.
       Part II: General formulation}

\author[rpi]{J.~W.~Banks\fnref{DOEThanks,PECASEThanks}}
\ead{banksj3@rpi.edu}

\author[rpi]{W.~D.~Henshaw\corref{cor1}\fnref{DOEThanks,NSFgrantNew}}
\ead{henshw@rpi.edu}

\author[rpi]{D.~W.~Schwendeman\fnref{DOEThanks,NSFgrantNew}}
\ead{schwed@rpi.edu}

\author[rpi]{Qi Tang\fnref{QiThanks}}
\ead{tangq3@rpi.edu}

\address[rpi]{Department of Mathematical Sciences, Rensselaer Polytechnic Institute, Troy, NY 12180, USA.}

\cortext[cor1]{Department of Mathematical Sciences, Rensselaer Polytechnic Institute, 110 8th Street, Troy, NY 12180, USA.}

\fntext[QiThanks]{Research supported by the Eliza Ricketts Postdoctoral Fellowship.}

\fntext[DOEThanks]{This work was supported by contracts from the U.S. Department of Energy ASCR Applied Math Program.}

\fntext[NSFgrantNew]{Research supported by the National Science Foundation under grant DMS-1519934.}

\fntext[PECASEThanks]{Research supported by a U.S. Presidential Early Career Award for Scientists and Engineers.}

\begin{abstract}
A stable partitioned algorithm is developed for fluid-structure interaction (FSI) problems involving
viscous incompressible flow and rigid bodies. This {\em added-mass partitioned} (AMP) algorithm
remains stable, without sub-iterations, for light and even zero mass rigid bodies when added-mass
and viscous added-damping effects are large.  The scheme is based on a generalized Robin interface
condition for the fluid pressure that includes terms involving the linear acceleration and angular
acceleration of the rigid body.  Added mass effects are handled in the Robin condition by inclusion
of a boundary integral term that depends on the pressure.  Added-damping effects due to the viscous
shear forces on the body are treated by inclusion of added-damping tensors that are derived through
a linearization of the integrals defining the force and torque.  Added-damping effects may be
important at low Reynolds number, or, for example, in the case of a rotating cylinder or rotating
sphere when the rotational moments of inertia are small.  In this second part of a two-part series,
the general formulation of the AMP scheme is presented including the form of the AMP interface
conditions and added-damping tensors for general geometries.  A fully second-order accurate
implementation of the AMP scheme is developed in two dimensions based on a fractional-step method for the
incompressible Navier-Stokes equations using finite difference methods and overlapping grids to
handle the moving geometry. The numerical scheme is verified on a number of difficult benchmark
problems.

\end{abstract}


\end{frontmatter}

\clearpage
\tableofcontents

\clearpage
\section{Introduction} \label{sec:intro}

We describe a new numerical approach for fluid-structure interaction (FSI) problems involving the motion of rigid bodies in an incompressible fluid.  The approach, referred to as the {\ampRB}~scheme, is a partitioned algorithm in which the equations for the fluid and the rigid bodies are handled using separate solvers.
This is in contrast to monolithic schemes where the whole system of equations are solved simultaneously at each time step.  A significant challenge for partitioned schemes is stability, especially for light bodies (or even zero-mass bodies) when the effects of added mass and added damping are important.\footnote{A brief physical explanation of added-mass and added-damping effects is given in the Introduction of Part~I.}  In addition, it can be difficult to achieve second-order accuracy, or higher, for partitioned time-stepping schemes due to errors caused by the numerical treatment of the matching conditions at the interface between the fluid and the rigid bodies.

The {\ampRB}~scheme is based on a fractional-step approach for the fluid in which the velocity is advanced in one stage, and the pressure is determined in a second stage~\cite{ICNS,splitStep2003}.  The viscous terms in the stress tensor are handled implicitly so that the velocity can be advanced with a larger stable time step.  The key ingredients of the~\ampRB~scheme are contained in the added-mass partitioned (AMP) interface conditions which couple the equations of motion of the rigid body to a compatibility boundary condition for the pressure on the surface of the body.  These conditions are derived at a continuous level by matching the {\em acceleration} of the body to that of the fluid.  As a result, the integration of the equations of motion of the rigid-body are coupled strongly to the update of the fluid pressure;  this ensures the proper balances of forces at the interface thereby suppressing instabilities due to added-mass effects.  Suppressing instabilities due to added-damping is more subtle.  The fluid forces on the body depend on the viscous shear stresses which, in turn, implicitly depend on the velocity of the body. This implicit dependence of the fluid forces on the body velocity is explicitly exposed and, after some simplifying approximations, is expressed in terms of {\em added-damping tensors} which are incorporated into the AMP interface conditions as a means to overcome added-damping instabilities.

This paper is the second of a two-part series of papers in which the {\ampRB}~scheme is developed
and analyzed.  The work in Part~I~\cite{rbinsmp2016r} introduced the scheme and applied it to
various model problems.  A stability analysis of the {\ampRB}~scheme was performed, and it was shown
that the new scheme remains stable, without sub-time-step iterations, even for light or zero-mass
rigid bodies when added-mass and added-damping effects are large.  In this paper, we extend the
formulation of Part~I to general three-dimensional 
configurations\footnote{Even though the numerical results are restricted to two dimensions, we believe that
the derivation and description of the added-damping tensors is more clearly presented in three dimensions.}. 
 An important new feature in this
extension is the generalization of the added-damping tensors that are incorporated into the
{\ampRB}~scheme to treat the effects of added damping.  Exact formulas are derived for the
added-damping tensors which involve solutions to variational problems given by two discrete vector
Helmholtz equations.  Approximate added-damping tensors, convenient for use in {\ampRB}~scheme, are
obtained using the results from the model problem analysis in Part~I.  These approximate tensors are
readily evaluated at the initial time using surface integrals over a given body, which can be
computed either analytically or numerically.  The scheme is implemented in two dimensions for arbitrary body motions
using moving overlapping grids~\cite{mog2006,Koblitz2016}.  Numerical examples in two dimensions are
provided that verify the accuracy and stability of the scheme for some challenging problems
involving light and zero-mass rigid bodies.

The FSI regime involving incompressible flows and moving rigid bodies is of great practical and scientific interest.
The reader is referred to the introduction of Part I for a discussion
of the literature, particularly as related to the issue of added-mass instabilities for partitioned schemes.
Here we note the long history of using composite overlapping (overset, Chimera) grids for simulating rigid
bodies in fluids, going back to the early work on  aircraft store-separation by Dougherty and Kuan~\cite{Dougherty}.
Recent use of overlapping grids for simulating aircraft, rotorcraft, wind-turbines, 
rockets, spacecraft, ships and underwater vehicles can be found at the overset grid symposium website~\cite{oversetGridSymposiumWebSite}.

The remainder of the paper is organized as follows.
The governing equations are given in Section~\ref{sec:governing}.
The analytic forms of the added-damping tensors are introduced in Section~\ref{sec:addedDampingTensors} and
these are subsequently used in Section~\ref{sec:AMPDPinterfaceConditions} to define the {\ampRB} interface conditions.
The {\ampRB} time-stepping algorithm, which uses a fractional-step scheme for the velocity-pressure form
of the Navier-Stokes equations is outlined in Section~\ref{sec:algorithm}.
Approximations to the added-damping tensors are developed in Section~\ref{sec:generalAddedDamping}.
Section~\ref{sec:numericalApproach} gives a brief outline of the moving overlapping grid approach.
Numerical results are given in Section~\ref{sec:numericalResults} and concluding remarks are made in Section~\ref{sec:conclusions}. 
\ref{sec:shearStressLemma} provides
the derivation of a formula for the shear stress used in determining the added-damping tensors, 
while 
\ref{sec:addedDampingTensorExamples} provides the explicit form of the 
added-damping tensors for some different shaped rigid bodies (e.g., rectangle, disk, sphere)

\section{Governing equations} \label{sec:governing}

We consider the fluid-structure coupling of an incompressible fluid and one or more rigid bodies 
as illustrated in Figure~\ref{fig:twoFallingBodies},
although the subsequent discussion will consider a single rigid body for simplicity.
The fluid occupies the domain $\xv\in\OmegaF(t)$ while the rigid body lies in the domain
$\xv\in \OmegaB(t)$, where $\xv$ is position and $t$ is time.  The coupling of the fluid and
body occurs along the interface $\GammaB(t)=\bar\OmegaF(t)\cap {\bar\Omega}_b(t)$. It is
assumed that the fluid is governed by the incompressible Navier-Stokes equations, which in an
Eulerian frame are given by
\begin{alignat}{3}
  &  \frac{\partial\vv}{\partial t} + (\vv\cdot\grad)\vv 
                 =  \frac{1}{\rho} \grad\cdot\sigmav  , \qquad&& \xv\in\OmegaF(t) ,  \label{eq:fluidMomentum}  \\
  & \grad\cdot\vv =0,  \quad&& \xv\in\OmegaF(t) ,  \label{eq:fluidDiv3d}
\end{alignat}
where $\rho$ is the (constant) fluid density and $\vv=\vv(\xv,t)$ is the fluid velocity.  The fluid stress tensor, $\sigmav=\sigmav(\xv,t)$, is given by 
\begin{equation}
  \sigmav = -p \Iv + \tauv,\qquad \tauv = \mu \left[ \grad\vv + (\grad\vv)^T \right],
  \label{eq:fluidStress}
\end{equation}
where $p=p(\xv,t)$ is the pressure, $\Iv$ is the identity tensor, $\tauv$ is the viscous stress tensor, and $\mu$ is the (constant) fluid viscosity. For future reference, the components of a vector such as $\vv$ will be denoted by $v_m$, $m=1,2,3$ (i.e.~$\vv=[v_1, v_2, v_3]^T$), while components of a tensor such as $\sigmav$ will be denoted by $\sigma_{mn}$, $m,n=1,2,3$.  The velocity-divergence form of the equations given by~\eqref{eq:fluidMomentum} and~\eqref{eq:fluidDiv3d} require appropriate initial and boundary conditions, as well as conditions on $\GammaB(t)$ where the behaviour of the fluid is coupled to that of the solid (as discussed below).

%
%
{
\newcommand{\figWidth}{5.25cm}
\newcommand{\trimfig}[2]{\trimFig{#1}{#2}{.275}{.5}{.23}{.23}}
\begin{figure}[htb]
\begin{center}
\resizebox{14cm}{!}{
\begin{tikzpicture}[scale=1]
  \useasboundingbox (0.0,.95) rectangle (16.5,7.);  
  \draw(0.0, 0) node[anchor=south west,xshift=-4pt,yshift=+0pt] {\trimfig{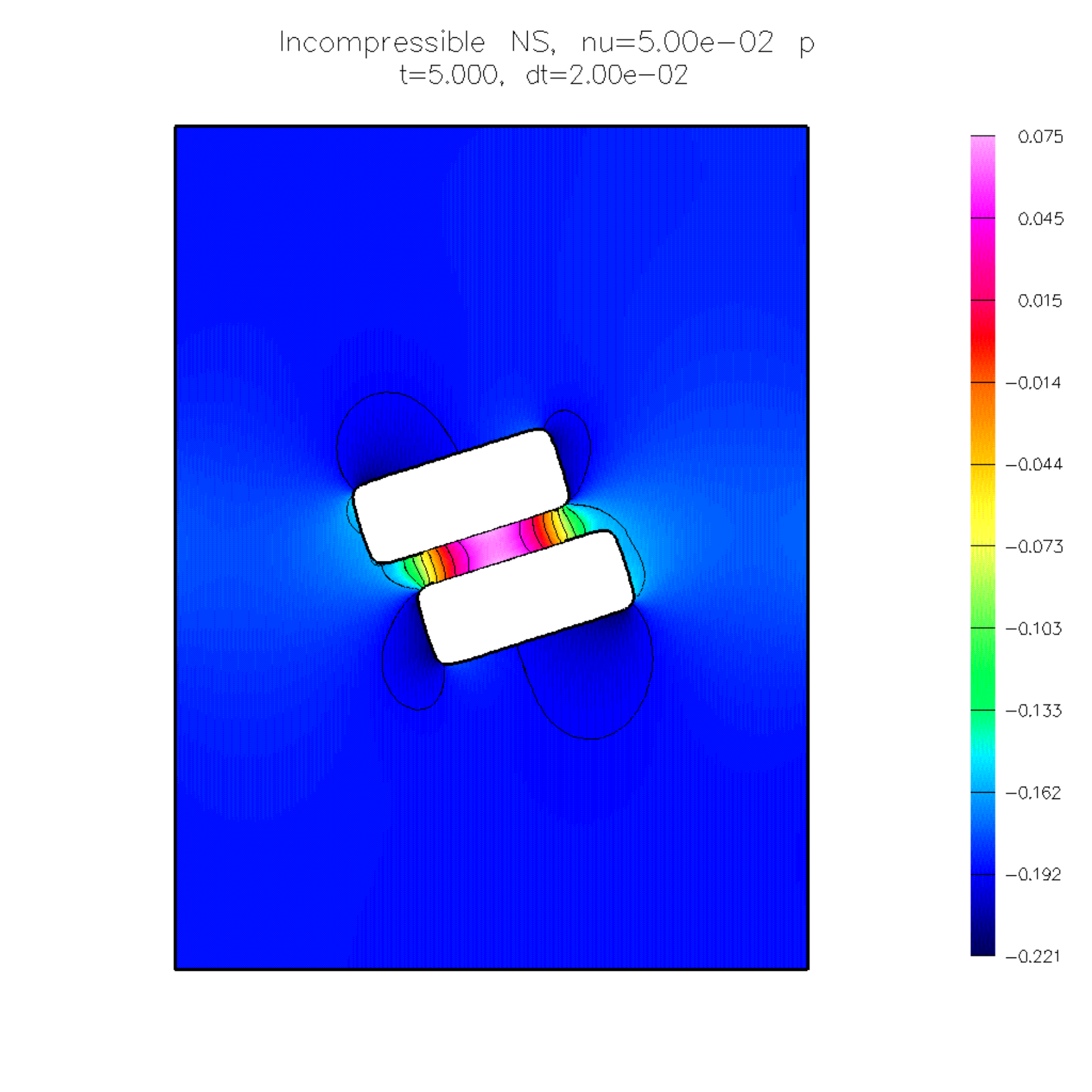}{\figWidth}};
  \draw(5.4, 0) node[anchor=south west,xshift=-4pt,yshift=+0pt] {\trimfig{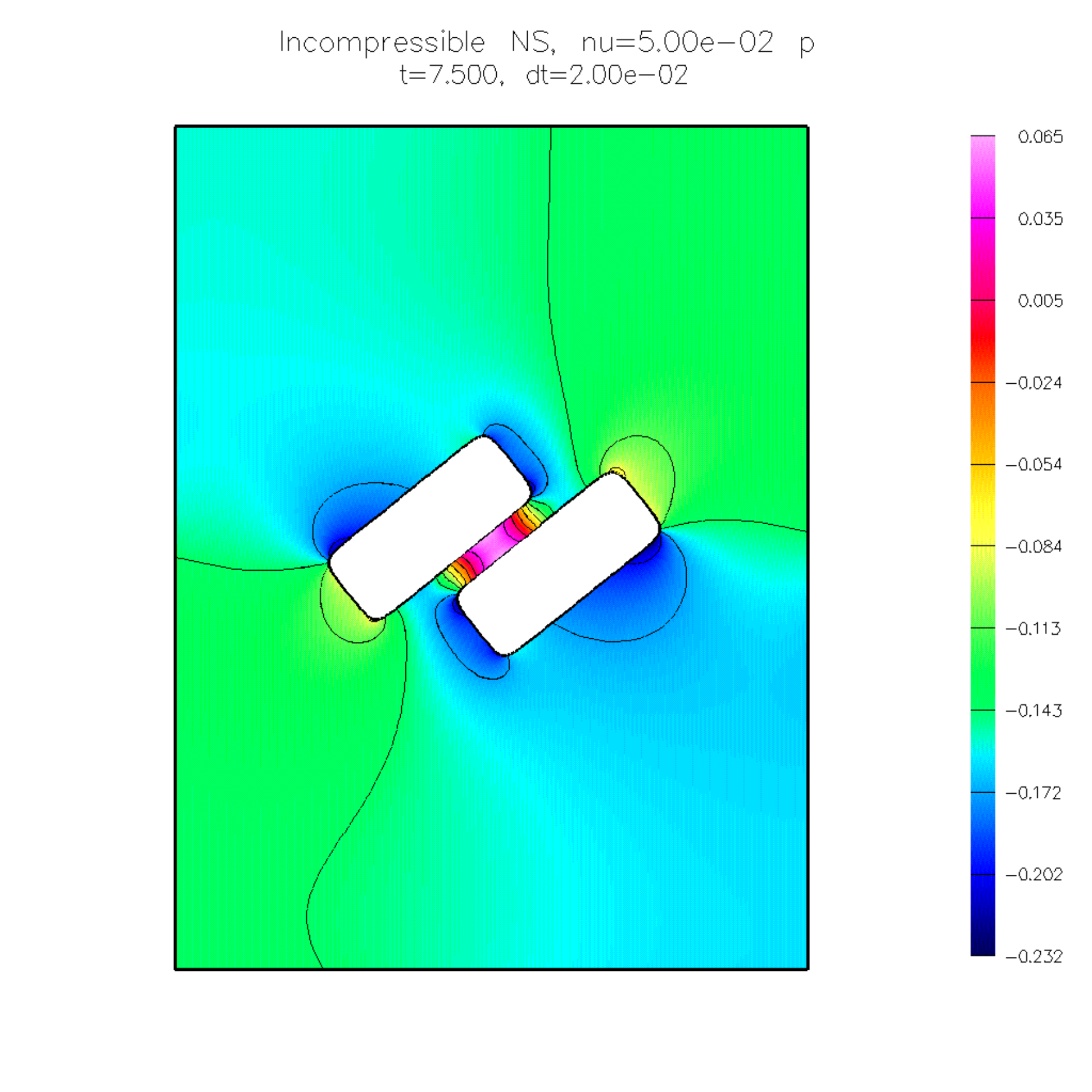}{\figWidth}};
  \draw(10.8,0) node[anchor=south west,xshift=-4pt,yshift=+0pt] {\trimfig{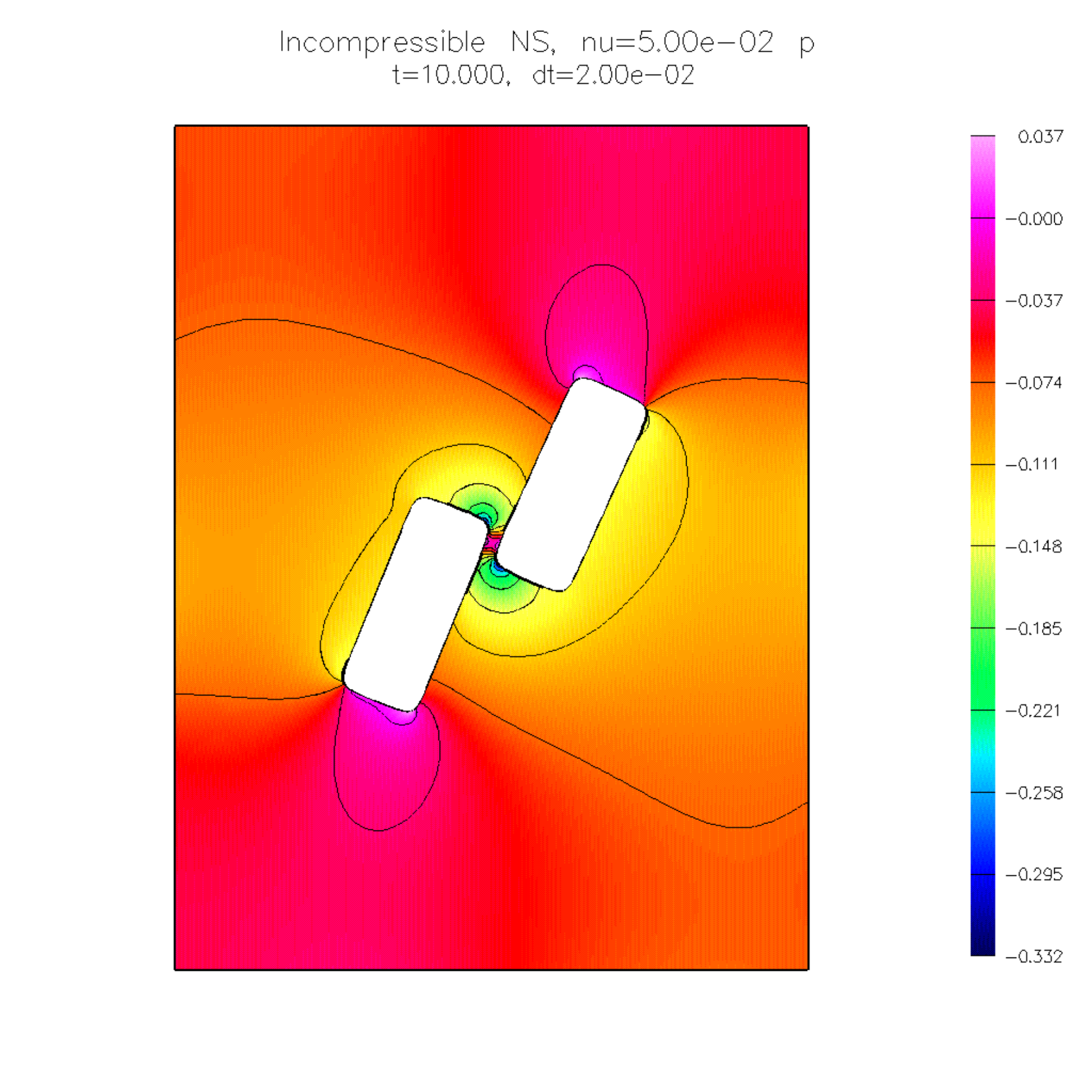}{\figWidth}};
%
\end{tikzpicture}
}
\end{center}
  \caption{Two rectangular bodies moving through an incompressible fluid under the influence of gravity showing
contours of the fluid pressure. Time increases
from left to right. The
heavier upper-body moves downward while the lighter lower-body rises.}
  \label{fig:twoFallingBodies}
\end{figure}
}

An elliptic equation for the fluid pressure can be derived from~\eqref{eq:fluidMomentum}--\eqref{eq:fluidStress} as
\begin{alignat}{3}
  &  \Delta p &= - \rho \grad\vv:(\grad\vv)^T ,  \qquad&& \xv\in\OmegaF(t) , \label{eq:fluidPressure}
\end{alignat}
where
\begin{align*}
\grad\vv:(\grad\vv)^T 
      ~ \eqdef \sum_{i=1}^{3}\sum_{j=1}^{3} \frac{\partial v_i}{\partial x_j}\frac{\partial v_j}{\partial x_i} .
\end{align*}
In the velocity-pressure form of the incompressible Navier-Stokes equations, the Poisson equation~\eqref{eq:fluidPressure} is used in place of~\eqref{eq:fluidDiv3d}. This form requires an additional boundary condition, and a natural choice is $\grad\cdot\vv=0$ for $\xv\in\partial\OmegaF(t)$, see~\cite{ICNS} for example.

The rigid body is defined in terms of the translation and rotation of the centre of mass. In discussing 
the rigid body motion we use the following notation: 
\begin{alignat*}{3}
  & \mrb\in\Real              ~&&: \text{mass of the rigid body},  \\
  & \rhob\in\Real              ~&&: \text{density of the rigid body},  \\
  & \Ib(t)\in\Real^{3\times 3}  ~&&: \text{moment of inertia matrix} , \\
  & \xvb(t)\in\Real^3         ~&&: \text{position of the centre of mass},\\ 
  & \vvb(t)\in\Real^3         ~&&: \text{velocity of the centre of mass},\\ 
  & \omegavb(t)\in\Real^3     ~&&: \text{angular velocity},\\ 
  & \Eb(t)\in\Real^{3\times 3} ~&&: \text{matrix with columns being the principle axes of inertia}, \\
  & \avb(t)\in\Real^3         ~&&: \text{linear acceleration of the centre of mass},\\ 
  & \bvb(t)\in\Real^3         ~&&: \text{angular acceleration of the centre of mass}, \\
  & \fvbe(t)\in\Real^3        ~&&: \text{external body force (e.g.~gravity)}, \\
  & \gvbe(t)\in\Real^3        ~&&: \text{external torque}. 
\end{alignat*}
The equations of motion for the rigid body are then given by the Newton-Euler equations,
\begin{align}
    \mrb \avb & = \int_{\GammaB} (-p \nv + \tauv\nv) \, \dArea + \fvbe  ,   \label{eq:linearAccelerationEquation} \\
    \Ib \bvb &= - \omegavb\times \Ib \omegavb + \int_{\GammaB} (\rv-\xvb)\times( -p\nv + \tauv\nv  ) \, \dArea + \gvbe,
                 \label{eq:angularAccelerationEquation}\\
   \dotvvcm &= \avb , \label{eq:centerOfMassVelocityEquation} \\
   \dotxvcm &= \vvcm , \label{eq:centerOfMassPositionEquation} \\
    \dot\omegav_b &= \bvb , \label{eq:angularVelocityEquation} \\
   \dot \Eb &= \omegavb\times \Eb ,  \label{eq:axesOfInertiaEquation}
\end{align}
where $\rv$ denotes a point on the surface of the body and 
 $\nv=\nv(\rv,t)$ is the outward unit normal to the body.
These equations of motion require
initial conditions for $\xvb(0)$, $\vvb(0)$, $\omegavb(0)$, and $\Eb(0)$. 

The motion of a point $\rv(t)$ on the surface of the body is given by a translation together with a rotation about the
initial centre of mass,
\begin{align}
    \rv(t) &= \xvcm(t) + R(t) (\rv(0)-\xvcm(0)) , \label{eq:bodySurface}
\end{align}
where $R(t)$ is the rotation matrix given by
\begin{align}
   R(t) &= \Eb(t) \Eb^{T}(0).    \label{eq:rotationMatrix}
\end{align}
The velocity and acceleration of that point 
are then given by
\begin{alignat}{3}
 &  \dot{\rv}(t) = \vvb(t) +\omegavb\times(\rv(t)-\xvb(t)), \quad&& \rv\in\GammaB, \label{eq:BodyPointVelocity}\\
 & \ddot{\rv}(t) = \avb + \bvb\times(\rv-\xvb)  + \omegavb\times\big[ \omegavb\times(\rv(t)-\xvcm(t))\big], \quad&& \rv\in\GammaB.\label{eq:BodyPointAcceleration} 
\end{alignat}
On the interface between the fluid and solid, the fluid velocity must match the solid velocity. 
If $\rv=\rv(t)$ denotes a point on the surface of the body $\GammaB$, 
then the fluid velocity $\vv(\rv(t),t)$ at the surface satisfies
\begin{align}
 &  \vv(\rv(t),t)  = \dot{\rv}(t), \qquad \rv\in\GammaB, \label{eq:RBsurfaceVelocity}
\end{align}
where $\dot{\rv}(t)$ is given in~\eqref{eq:BodyPointVelocity}.

\newcommand{\Fvt}{\Fc_\mu}
\newcommand{\Gvt}{\Gc_\mu}
\section{Added-damping tensors}\label{sec:addedDampingTensors}

Viscous shear stress on the surface of the rigid body generates forces and torques that 
are coupled to the motion of the rigid body. If this coupling is not properly treated,
then a numerical scheme may become unstable due to an added-damping instability.  Such an instability
arises, for example, when an over-estimate of the shear stress leads to an over-correction of
the rigid-body velocity, which in turn leads to an even larger over-estimate of the stress, and so on.
In this section, we introduce exact formulas for the added-damping tensors that
account for the coupling of the force and torque due to viscous shear stress and the motion of the rigid body.  These
tensors are then used in the
description of the~\ampRB~interface conditions given next in Section~\ref{sec:AMPDPinterfaceConditions}.
For the numerical implementation of the~\ampRB~interface conditions, we only require approximations of added-damping tensors, and these approximations are discussed in Section~\ref{sec:generalAddedDamping}.

Added-damping effects, which arise from the forces and torques on the body due the viscous shear stress, 
are represented by the terms involving $\tauv\nv$ in the surface integrals~\eqref{eq:linearAccelerationEquation}
and~\eqref{eq:angularAccelerationEquation}
for the linear and
angular accelerations of the rigid body, respectively.
These two contributions to the accelerations of the body are given by
\begin{equation}
   \Fvt(\vv,\vvb,\omegavb) = \int_{\GammaB} \tauv\nv\,\dArea, \qquad 
   \Gvt(\vv,\vvb,\omegavb) = \int_{\GammaB}  (\rv-\xvb)\times (\tauv\nv)\,\dArea,
\label{eq:surfaceIntegrals}
\end{equation}
where it is indicated that these integrals are functions of $\vv$, but also implicitly of the body velocities, $\vvb$ and $\omegavb$, as well. 
To expose this dependence consider 
linear approximations of $\Fvt$ and $\Gvt$ about some predicted states, $\vv^p$, $\vvb^p$, and $\omegavb^p$, which are assumed to be known.  These approximations have the form
\begin{align}
   \Fvt(\vv,\vvb,\omegavb) \approx \Fvt(\vv^p,\vvb^p,\omegavb^p) - \Dvv\, (\vvb-\vvb^p) - \Dvw\,(\omegavb-\omegavb^p) ,  \label{eq:addedDampingFdws} \\
   \Gvt(\vv,\vvb,\omegavb) \approx \Gvt(\vv^p,\vvb^p,\omegavb^p) - \Dwv\,(\vvb-\vvb^p) - \Dww\,(\omegavb-\omegavb^p), \label{eq:addedDampingGdws}
\end{align}
where the {\em added-damping tensors}, $\Dvv$, $\Dvw$, $\Dwv$, and $\Dww$, are $3\times 3$ matrices (assuming three-dimensional flow) defined by
\begin{alignat}{3}
&  \Dvv  \eqdef -\frac{\partial\Fvt}{\partial\vvb}  , \qquad&& 
  \Dvw  \eqdef -\frac{\partial\Fvt}{\partial\omegavb} , \label{eq:dampingTensorsF} \\
&  \Dwv  \eqdef -\frac{\partial\Gvt}{\partial\vvb}  , \qquad&& 
  \Dww  \eqdef -\frac{\partial\Gvt}{\partial\omegavb}   . \label{eq:dampingTensorsG}
\end{alignat}
The approximations in~\eqref{eq:addedDampingFdws} and~\eqref{eq:addedDampingGdws} are sufficient to
expose the dependence of 
the applied forces on the velocities $\vvb$ and $\omegavb$.  Note that the added-damping
tensors introduced here are an extension of the scalar {\em added-damping coefficients}, $\Dc\sp{u}$
and $\Dc\sp{\omega}$, obtained in Part~I~\cite{rbinsmp2016r} for two model FSI problems.

Exact formulas for the added-damping tensors can be obtained by first noting that the shear stress on the surface of the rigid body, $\tauv\nv$, appearing in the integrals in~\eqref{eq:surfaceIntegrals} can be expressed in terms of the normal derivative of the fluid velocity on the surface and the angular velocity of the rigid body.  This form is given in the following theorem.

\begin{theorem} \label{theorem:shearStressRigidBody}
The shear stress of a viscous incompressible fluid with velocity $\vv(\xv,t)$ and dynamic viscosity $\mu$ on the surface of a rigid body with angular velocity $\omegavb(t)$ is given by 
\begin{align}
   \tauv\nv 
   &= \mu \left\{ (\Iv-\nv\nv^T)\frac{\partial\vv}{\partial n} + \nv \times \omegavb \right \},
   \qquad \xv \in \GammaB,  \label{eq:viscousTractiondws}
\end{align} 
where $\Iv$ is the $3\times3$ identity matrix.
\end{theorem}
\noindent The proof of this theorem is given in \ref{sec:shearStressLemma}.

Using~\eqref{eq:viscousTractiondws} to eliminate $\tauv\nv$ in the integrals in~\eqref{eq:surfaceIntegrals} gives
\begin{align*}
\Fvt(\vv,\vvb,\omegavb) &= \mu\int_{\GammaB}\left\{(\Iv-\nv \nv^T) \frac{\partial\vv}{\partial n} 
                               + \nv \times \omegavb \right\} \,\dArea, \\
\Gvt(\vv,\vvb,\omegavb) &= \mu\int_{\GammaB}  (\rv-\xvb)\times \left\{(\Iv-\nv \nv^T) \frac{\partial\vv}{\partial n} 
                               + \nv \times \omegavb \right\}\,\dArea,
\end{align*}
Define a skew symmetric matrix $[\xv]_\times \in \Real^{3\times3}$,
\begin{align*}
	[\xv]_\times \eqdef
	\begin{bmatrix}
		0 & -x_3 & x_2 \\
		x_3 & 0 & -x_1 \\
		-x_2 & x_1 &  0 
	\end{bmatrix},
\end{align*}
so that the cross product can be rewritten as matrix multiplication, 
\begin{align*}
	\xv \times \yv = [\xv]_\times \yv.
\end{align*}
Therefore, the added-damping tensors in~\eqref{eq:dampingTensorsF} and~\eqref{eq:dampingTensorsG} can be expressed as
\begin{align}
\Dvv &= -\mu\int_{\GammaB} \left\{(\Iv-\nv \nv^T) \frac{\partial}{\partial n}\left(\frac{\partial\vv}{\partial\vvb}\right)\right\} \,\dArea, \label{eq:DvvExact} \\
\Dvw &= -\mu\int_{\GammaB} \left\{ (\Iv-\nv \nv^T)\frac{\partial}{\partial n}\left(\frac{\partial\vv}{\partial\omegavb}\right) + [\nv]_\times \right\} \,\dArea, \label{eq:DvwExact} \\
\Dwv &= -\mu\int_{\GammaB} [\rv-\xvb]_\times \left\{ (\Iv-\nv \nv^T) \frac{\partial}{\partial n}\left(\frac{\partial\vv}{\partial\vvb}\right)\right\} \,\dArea, \label{eq:DwvExact} \\
\Dww &= -\mu\int_{\GammaB} [\rv-\xvb]_\times \left\{ (\Iv-\nv \nv^T) \frac{\partial}{\partial n}\left(\frac{\partial\vv}{\partial\omegavb}\right) + [\nv]_\times \right\} \,\dArea, \label{eq:DwwExact} 
\end{align}
where the derivatives, $\partial_{\vvb}\vv$ and~$\partial_{\omegavb}\vv$, solve certain variations problems as discussed in Section~\ref{sec:generalAddedDamping}.
For the typical case, when the rigid body is fully immersed in the fluid,
the terms involving $[\nv]_{\times}$ in the added-damping tensors~\eqref{eq:DvwExact} and~\eqref{eq:DwwExact} can be evaluated using Gauss's theorem. 
Assuming that the body surface $\GammaB$ is sufficiently smooth 
so that the Gauss divergence theorem holds (this is true, roughly, if the surface is piecewise smooth with
a finite number of edges and corners~\cite{WidderAdvancedCalculus}), then 
\begin{align*}
	\int_{\GammaB} n_i \, \dArea = 0, \quad\text{and}\quad
         \int_{\GammaB} x_i n_j \, \dArea = 
	\int_{\OmegaB} \frac{\partial x_i}{\partial x_j} \, d\xv = \delta_{ij} V_b,
\end{align*}
where $V_b$ is the volume of the rigid body.
Therefore, two of the added-damping tensors can be simplified as
\begin{align}
\Dvw &= -\mu\int_{\GammaB} \left\{ (\Iv-\nv \nv^T)\frac{\partial}{\partial n}\left(\frac{\partial\vv}{\partial\omegavb}\right) \right\} \,\dArea, \label{eq:DvwExactNew} \\
\Dww &= 2 \mu V_b\, \Iv -\mu\int_{\GammaB} [\rv-\xvb]_\times \left\{ (\Iv-\nv \nv^T) \frac{\partial}{\partial n}\left(\frac{\partial\vv}{\partial\omegavb}\right)\right\} \,\dArea. \label{eq:DwwExactNew} 
\end{align}

The four added-damping tensors can be collected in a larger composite tensor, $\ADtensor(t) \in\Real^{6\times 6}$, having the form
\begin{align}
  \ADtensor \eqdef   \begin{bmatrix}
    \Dvv & \Dvw  \\
    \Dwv & \Dww
  \end{bmatrix} .  \label{eq:AddedDampingTensor} 
\end{align}
Note that the added-damping tensors change with time as the body rotates according to the formula
\begin{align}
   \Da^{\alpha\beta}(t) = R(t)\, \Da^{\alpha\beta}(0) \, R^T(t) ,  \label{eq:addedDampingInTime}
\end{align}
where $R(t)$ is the rotation matrix given in~\eqref{eq:rotationMatrix}.
Thus, the integrals defining the entries in~$\ADtensor(t)$ can be precomputed at $t=0$, and then~\eqref{eq:addedDampingInTime} can be used to determine~$\ADtensor(t)$ at later times.
We note that the added-damping tensors in two dimensions can be easily derived 
from the formulae in three dimensions 
by, for example, considering cylindrical bodies of unit depth in the $z$-direction (with no top or bottom) in which case 
the surface integrals will reduce to line integrals. The resulting expressions can then be 
restricted to two-dimensional motions.

Motivated by the analysis in~\cite{rbinsmp2016r}, it is convenient in the~\ampRB~algorithm to express the
linearizations in~\eqref{eq:addedDampingFdws} and~\eqref{eq:addedDampingGdws} in terms of the body
accelerations $\avb$ and $\bvb$. This is accomplished by including a factor of the discrete time-step $\dt$ and
an additional added-damping parameter $\adc$, 
\begin{align}
   \Fvt(\vv,\vvb,\omegavb) \approx \Fvt(\vv^p,\vvb^p,\omegavb^p) - \adc\, \dt \, \Dvv\, (\avb-\avb^p) 
                 - \adc\, \dt \,\Dvw\,(\bvb-\bvb^p) ,  \label{eq:addedDampingFdot} \\
   \Gvt(\vv,\vvb,\omegavb) \approx \Gvt(\vv^p,\vvb^p,\omegavb^p) - \adc\, \dt \,\Dwv\,(\avb-\avb^p) 
                 - \adc\, \dt \,\Dww\,(\bvb-\bvb^p). \label{eq:addedDampingGdot}
\end{align}
As discussed at length in
Part~I~\cite{rbinsmp2016r}, the choice $\adc=1$ lies in a range of values for $\adc$ that yield
stable schemes even for massless bodies. As a result that choice is made here for all computations.
Finally, note that while the formulas in~\eqref{eq:DvvExact}--\eqref{eq:DwwExact} are exact,
suitable approximations for the~\ampRB~algorithm can be obtained by considering discrete variational
problems for the derivatives~$\partial_{\vvb}\vv$ and~$\partial_{\omegavb}\vv$ appearing in the four
integrals.  These approximations are discussed in Section~\ref{sec:generalAddedDamping}.

\section{The added-mass and added-damping interface conditions}
\label{sec:AMPDPinterfaceConditions} 

In this section, the interface conditions used in the~\ampRB~scheme to account for
added-mass and added-damping effects are derived at a continuous level. These conditions are 
generalizations of the conditions given in Part~I~\cite{rbinsmp2016r} for some simplified model problems.
A discretization of the interface conditions obtained here is given later in the description of
the~\ampRB~scheme.

Added-mass effects arise from the contribution of the pressure to the force and torque on the rigid body.
To account for these effects and to obtain the proper force-balance for light rigid-bodies, the 
fluid and solid {\em accelerations} are matched at the surface of the rigid body.  Using~\eqref{eq:BodyPointAcceleration}, this implies
\[
D_t\vv =  \avb + \bvb\times(\rv-\xvb)  + \omegavb\times\big[ \omegavb\times(\rv(t)-\xvcm(t))\big], \quad \rv\in\GammaB,
\]
where $\rv(t)$ denotes a point on the surface of the body and $D_t=\partial_t+\dot{\rv}\cdot\grad$ denotes the total derivative following a point on the surface. Using the momentum equation~\eqref{eq:fluidMomentum} to eliminate $D_t\vv$
leads to the vector {\em acceleration compatibility condition} given by
\[
-\frac{1}{\rho}\Big[ \grad p -\mu\Delta\vv \Big] 
      =  \avb + \bvb\times(\rv-\xvb)  + \omegavb\times\big[ \omegavb\times(\rv(t)-\xvcm(t))\big], \quad \rv\in\GammaB.
\]
The normal component of the vector condition can be written in the form
\begin{align}
    \partial_n p
              & = -\rho\nv^T \Big( \avb + \bvb\times(\rv-\xvb)  + \omegavb\times\big[ \omegavb\times(\rv(t)-\xvcm(t))\big]  \Big)  
                + \mu \nv^T\Delta \vv ,  \quad \rv\in\GammaB, \label{eq:pressureBC} 
\end{align}
where $\nv=\nv(\rv,t)$ is the outward normal to the body.  

Added-damping effects arise from the influence of the motion of the rigid body on the integrals in~\eqref{eq:surfaceIntegrals} involving the shear stress on the body as discussed in the previous section.  The
linearizations in~\eqref{eq:addedDampingFdot}--\eqref{eq:addedDampingGdot} reveal this influence in terms
of the accelerations of the body, and these can be added to the interface conditions in~\eqref{eq:linearAccelerationEquation} and~\eqref{eq:angularAccelerationEquation}.  
Combining these conditions with the acceleration
compatibility condition in~\eqref{eq:pressureBC} leads to the primary~\ampRB~interface
condition.  
\begin{AMPInterfaceCondition} The \ampRB~interface conditions on the surface of a rigid body $\rv\in\GammaB$ for the pressure equation~\eqref{eq:fluidPressure} are 
\begin{align}
    \partial_n p + \rho\nv^T \Big( \avb + \bvb\times(\rv-\xvb)   \Big)  
               &= -\rho\nv^T \Big( \omegavb\times\big[ \omegavb\times(\rv(t)-\xvcm(t))\big]  \Big)  
                + \mu \nv^T\Delta \vv , \label{eq:AMPDPpressure} \\
 \left\{ 
\begin{bmatrix} \mrb\, I_{3\times 3} & 0 \smallskip\\  0 & \Ib  \end{bmatrix}
+ \adc \, \dt\, \ADtensor 
\right\} 
\begin{bmatrix} \avb \smallskip\\ \bvb  \end{bmatrix}
+
\Fv(p) 
& = 
-  \begin{bmatrix} 0 \smallskip\\ \omegavb\times \Ib \omegavb\end{bmatrix} + 
         \Gv(\vv) + \adc \, \dt\, \ADtensor \begin{bmatrix} \avb \smallskip\\ \bvb  \end{bmatrix},  \label{eq:AMPDPrigidBodyAcceleration}
\end{align}
where
\begin{equation}
\Fv(p) \eqdef \begin{bmatrix}
                 -\displaystyle \int_{\partialB} p\, \nv \, \dArea \\ 
                 -\displaystyle \int_{\partialB} (\rv-\xvb)\times (p\, \nv ) \, \dArea 
             \end{bmatrix} ,
\qquad
\Gv(\vv) \eqdef 
\begin{bmatrix}
                \displaystyle \Fvt \, + \fvbe \\
                \displaystyle \Gvt \, + \gvbe
              \end{bmatrix} \,=\,  
\begin{bmatrix}
                \displaystyle \int_{\partialB} \tauv\nv \, \dArea + \fvbe \\
                \displaystyle \int_{\partialB} (\rv-\xvb)\times (\tauv\nv) \, \dArea + \gvbe
              \end{bmatrix}.  
       \label{eq:Pdef}
\end{equation}
Here, $\dt$ is the time-step used in the~\ampRB~time-stepping scheme and $\adc$ is the added-damping parameter typically set to one. 
\end{AMPInterfaceCondition}

In the~\ampRB~time-stepping scheme, predicted values for $\vv$, $\avb$, $\bvb$,
$\omegavb$ and $\xvb$
are used to evaluate the right-hand sides of the AMP interface conditions in~\eqref{eq:AMPDPpressure} and~\eqref{eq:AMPDPrigidBodyAcceleration}. 
The fluid pressure $p$ and body accelerations, $\avb$ and $\bvb$, are
then updated by solving the Poisson equation in~\eqref{eq:fluidPressure} together with the AMP interface
conditions.
Note that the same term involving $\ADtensor$ with components $\Da^{\alpha\beta}$ appears on both sides of~\eqref{eq:AMPDPrigidBodyAcceleration}.  Thus, the components of this interface condition reduce to the conditions in~\eqref{eq:linearAccelerationEquation} and~\eqref{eq:angularAccelerationEquation}, so that the
interface condition is exact for {\em any}
choice of the added-damping tensor $\ADtensor$, including the approximations obtained later in Section~\ref{sec:generalAddedDamping}.
The purpose of including the term involving $\ADtensor$ in~\eqref{eq:AMPDPrigidBodyAcceleration} is to 
cancel the leading contributions to the added-damping embedded in the integrals~$\Fvt$ and~$\Gvt$ in $\Gv$. 
These contributions now appear on the left-hand side of~\eqref{eq:AMPDPrigidBodyAcceleration} in the calculation of the body accelerations.

\section{The \ampRB~time-stepping algorithm} \label{sec:algorithm}

We now describe a second-order accurate version of the~\ampRB~time-stepping algorithm.  In this scheme,
the fluid variables are computed in a fractional-step 
manner with the fluid velocity advanced in one stage, followed by an update of the fluid pressure in a subsequent stage.
A key element of the algorithm is the~\ampRB~interface conditions, which are used to couple the update of the pressure with the calculation of the accelerations of the rigid body.  These interface conditions are designed to suppress instabilities due to added-mass and added-damping effects as discussed previously.

Before describing the time-stepping algorithm in detail, let us first introduce
some notation. Let $\xv_\iv$ denote the grid-point coordinates on a discrete mesh, where $\iv=(i_1,i_2,i_3)$ is a multi-index,
and let $t^n=n\dt$, $n=0,1,2,3,\ldots$ denote the discrete times in terms
of the time-step $\dt$. 
The {\em composite grid} covering the fluid domain, $\Omega(t)$, at time $t=t^n$ is denoted by $\Gc^n$.  This grid is composed of a collection of component grids, some of which are fitted to the surface of the rigid body and move in time, as described in Section~\ref{sec:numericalApproach}.
Let $\vv_\iv^n \approx \vv(\xv_\iv,t^n)$ and $p_{\iv}^n \approx p(\xv_{\iv},t^n)$ denote grid functions for
the fluid velocity and pressure, respectively,
and let $\xvb^n\approx\xvb(t^n)$, $\vvb^n\approx \vvb(t^n)$, etc., denote time-discrete approximations
of the rigid-body variables.
Let $\grad_h$ and $\Delta_h$ denote some appropriate discrete approximations to the gradient and
Laplacian operators, respectively; the precise form of these approximations is not important for the
present discussion.
Let $\OmegaF_h$ denote the set of indices $\iv$ in the interior of the fluid grid,
$\Gamma_h$ denote the set of indices $\iv$ on the interface, and
\begin{align}
  \fap_{\iv}^n &\eqdef \big((\vv_{\iv}^n-\wv_{\iv}^n\,)\cdot\grad_h\big) \vv_{\iv}^n + \grad_h  p_\iv^{n}, \label{eq:advPressure}
\end{align}
denote the advection and pressure terms in the momentum equations, where the {\em grid velocity}, $\wv_{\iv}^n$,
is included from the transformation to a moving coordinate system.
Finally, denote the state of the rigid body as a vector given by 
\[
  \qvb^n \eqdef [\avb^n \cma \bvb^n \cma \vvb^n \cma \omegavb^n \cma \xvb^n \cma \Eb^n ]. 
\]

\begin{algorithm}\caption{\rm Added-mass partitioned (\ampRB) scheme}\small
\[
\begin{array}{l}
\hbox{// \textsl{\bodyStepIComment}}\smallskip\\
1.\quad (\, \qvb^\esup \cma \Gc^\esup \,)=\hbox{\bf predictAndExtrapRigidBody}(\, \qvb^n \cma \qvb^{n-1} \,) 
  \bigskip\\
\hbox{// \textsl{Prediction steps}}\smallskip\\
2.\quad (\,\vv_\iv^\psup\,)=\hbox{\bf advanceFluidVelocity}(\,\qvb^\esup \cma 2\fv_\iv^n - \fv_\iv^{n-1} \cma \fv_\iv^n \cma \vv_{\iv}^n \cma \Gc^\esup\,) 
  \smallskip\\
3.\quad (\, p_\iv^\psup \cma \avb^\psup \cma \bvb^\psup )=  
                 \hbox{\bf updatePressureAndBodyAcceleration}(\,\qvb^\esup \cma \vv_\iv^\psup\cma \Gc^\esup\,)
   \smallskip\\
4.\quad (\qvb^\psup \cma \Gc^\psup)=\hbox{\bf advanceBodyGivenAcceleration}(\, \avb^\psup \cma \bvb^\psup \cma \qvb^n \,)
  \bigskip\\
\hbox{// \textsl{Correction steps}}\smallskip\\
5.\quad (\,\vv_\iv^{n+1}\,)=\hbox{\bf advanceFluidVelocity}(\,\qvb^\psup \cma \fv_\iv^\psup \cma \fv_\iv^n \cma \vv_{\iv}^n \cma \Gc^\psup\,) 
  \smallskip\\
6.\quad (\, p_\iv^{n+1} \cma \avb^{n+1} \cma \bvb^{n+1} )=
              \hbox{\bf updatePressureAndBodyAcceleration}(\,\qvb^\psup \cma \vv_\iv^{n+1} \cma \Gc^\psup\,)
   \smallskip\\
7.\quad (\qvb^{n+1} \cma \Gc^{n+1})=\hbox{\bf advanceBodyGivenAcceleration}(\, \avb^{n+1} \cma \bvb^{n+1} \cma \qvb^n \,)
  \bigskip\\
\hbox{// \textsl{Fluid-velocity correction step (optional)}}\smallskip\\
8.\quad (\,\vv_\iv^{n+1}\,)=\hbox{\bf advanceFluidVelocity}(\,\qvb^{n+1} \cma \fv_\iv^{n+1} \cma \fv_\iv^n \cma \vv_{\iv}^n\cma \Gc^{n+1}\,) 
\end{array}
\]
\label{alg:ampRB}
\end{algorithm}

The~\ampRB~time-stepping scheme described in Algorithm~\ref{alg:ampRB} is a predictor-corrector-type fractional-step scheme for the fluid velocity and pressure that incorporates the AMP interface conditions coupling the motion of the rigid body.  
The algorithm involves a set of four procedures, which are given below, and is a generalization of the~\ampRB~scheme discussed in detail in Part~I~\cite{rbinsmp2016r} for a simplified model problem.
The key steps in the algorithm are Steps~3 and~6 where the pressure and body accelerations are computed
with the aim to suppress instabilities caused by added-mass
and added-damping effects.
The correction steps, Steps~5--7, can be repeated any number of times (or omitted completely), but we usually apply these steps once as this
increases the stability region of the scheme 
to allow a larger advection time-step. 
The velocity correction in Step~8 is optional, although including this step 
improves the stability of the scheme for added-damping effects, as noted in Part~I, and may be needed when these effects are large.

\bigskip
\begin{myProcedure}{$(\, \qvb \cma \Gc \,)=\hbox{\bf predictAndExtrapRigidBody}(\, \qvb^n \cma \qvb^{n-1} \,)$}
\noindent Predict the accelerations of the rigid body at $t^{n+1}$ using 
linear extrapolation in time, 
\begin{align*}
&    \avb = 2 \avb^n - \avb^{n-1} , \qquad 
    \bvb = 2 \bvb^n - \bvb^{n-1} .
\end{align*}
Predict the primary rigid-body degrees-of-freedom at $t^{n+1}$ using 
a leap-frog scheme in time, 
\begin{alignat*}{3}
&    \xvb = \xvb^{n-1} + 2\dt\, \vvb^{n}  , \qquad&& 
    \vvb = \vvb^{n-1} + 2\dt\, \avb^{n} , \\
&    \omegavb = \omegavb^{n-1} + 2\dt\, \bvb^{n} , \qquad&& 
    \Eb = \Eb^{n-1} + 2\dt\, \omegavb^{n}\times\Eb^{n} ,
\end{alignat*}
and then update the state of the rigid body at $t^{n+1}$ given by
\[
\qvb=[\avb \cma \bvb \cma \vvb \cma \omegavb \cma \xvb \cma \Eb]. 
\]
The predicted moving grid at the new time, $\Gc$,
can now be generated given the 
position of the body surface $\rv(t)$ at $t=t^{n+1}$ determined from~\eqref{eq:bodySurface} and~\eqref{eq:rotationMatrix} using
$\xvb$ and $\Eb$.
\label{proc:erb}
\end{myProcedure}

\medskip
\begin{myProcedure}{$(\,\vv_\iv\,)=\hbox{\bf advanceFluidVelocity}(\,\qvb^* \cma\fap_{\iv}^*\cma \fap_{\iv}^n\cma \vv_{\iv}^n\cma \Gc^*\,) $}
\noindent Advance the fluid velocity with a semi-implicit scheme using the given values $\fap_{\iv}^*$ and $\fap_{\iv}^n$ for the advection and 
pressure gradient terms~\eqref{eq:advPressure}, and setting the fluid velocity on the surface of the rigid body 
to match that of the body using the given state $\qvb^*$,
\begin{alignat*}{3}
&   \rho\frac{\vv_{\iv} -\vv^n_{\iv}}{\dt} + \frac{1}{2}\fap_{\iv}^* + \frac{1}{2} \fap_{\iv}^n 
         = \frac{\mu}{2}\Big( \Delta_h \vv_{\iv} + \Delta_h \vv_\iv^n \Big)  , \qquad&&\iv\in\OmegaF_h, \\
&   \vv_\iv = \vvb^{*} + \omegavb^*\times( \rv_\iv^*-\xvb^*),\quad \grad_h\cdot\vv_\iv=0, \quad \text{Extrapolate ghost: } \tv_m^T\vv_\iv, \qquad&&\iv\in\Gamma_h. 
\end{alignat*}
Here the ghost point values for the velocity are determined by the
divergence boundary condition together with extrapolation of the tangential components. Appropriate conditions on the velocity are
specified on other boundaries, e.g.~no-slip conditions on $\partial\Omega\backslash\GammaB$.
\end{myProcedure}

\medskip
\begin{myProcedure}{$(\, p_\iv\cma \avb\cma \bvb )=\hbox{\bf updatePressureAndBodyAcceleration}(\,\qvb^* \cma \vv_\iv^* \cma \Gc^*\,)$}
\noindent Update the fluid pressure and rigid-body accelerations by solving 
\begin{align*}
	&    \Delta_h p_\iv = -\grad_h\vv^{*}_\iv : {(\grad_h\vv^{*}_\iv)^T} , \qquad \iv \in \Omega_h\cup\Gamma_h, 
\end{align*}
together with the interface conditions\footnote{The term $\Delta\vv$ in the boundary 
condition~\eqref{eq:pressureInterfaceEquation} has been replaced 
by $-\grad\times\grad\times\vv$ to avoid a viscous time-step restriction, see for example,~\cite{Petersson00}.}
\begin{align}
   \nv_\iv^T\grad_h p_\iv + \rho{\nv_\iv^T} \Big( \avb + \bvb\times(\rv_\iv^{*}-\xvb^* )   \Big)  
               &= -\rho\nv_\iv^T \Big( \omegavb^* \times\big[ \omegavb^*\times(\rv^*_\iv-\xvcm^*)\big]  \Big)  \nonumber \\
               &\qquad\qquad\qquad - \mu \nv_\iv^T\, (\grad_h\times\grad_h\times \vv_\iv^*) , \qquad\iv\in\Gamma_h, \label{eq:pressureInterfaceEquation} \\
 \left\{ 
\begin{bmatrix} \mrb\, I_{3\times 3} & 0 \\  0 & \Ib  \end{bmatrix}
+ \dt\,\adc\,\ADtensor^*
 \right\} 
\begin{bmatrix} \avb \\ \bvb  \end{bmatrix}
+
\Fv(p_\iv) 
& = 
-  \begin{bmatrix} 0 \\ \omegavb^*\times \Ib \omegavb^*\end{bmatrix} + 
         \Gv(\vv_\iv^*) + \dt\, \adc\,
      \ADtensor^* \begin{bmatrix} \avb^* \\ \bvb^*  \end{bmatrix} .
    \label{eq:RBaccelerationInterface}
\end{align}
Here $\ADtensor^*$ is the added damping tensor~\eqref{eq:AddedDampingTensor},
computed using~\eqref{eq:addedDampingInTime}, and $\Fv$ and $\Gv$ 
are defined in~\eqref{eq:Pdef}.  Appropriate boundary conditions for pressure are
specified on the other boundaries.
\end{myProcedure}

\medskip
\begin{myProcedure}{$(\qvb \cma \Gc)=\hbox{\bf advanceBodyGivenAcceleration}(\, \avb^* \cma \bvb^*\cma \qvb^n \,)$}
\noindent Use the trapezoidal rule to
correct the  positions and velocities of the rigid body using the given accelerations $\avb^*$ and $\bvb^*$,
\begin{alignat*}{3}
&   \vvb = \vvb^{n} + \frac{\dt}{2}\big( \avb^{*} + \avb^{n}\big)   , \qquad&& 
    \omegavb = \omegavb^{n} + \frac{\dt}{2}\big( \bvb^{*} + \bvb^{n}\big) , \\
&    \xvb = \xvb^{n} + \frac{\dt}{2}\big( \vvb + \vvb^{n}\big)  , \qquad&& 
    \Eb = \Eb^{n} + \frac{\dt}{2}\big( \omegavb\times\Eb^{*} + \omegavb^{n}\times\Eb^{n} \big)  .
\end{alignat*}
and then update the state of the rigid body at $t^{n+1}$ given by
\[
\qvb=[\avb^* \cma \bvb^* \cma \vvb \cma \omegavb \cma \xvb \cma \Eb]. 
\]
The moving grid, $\Gc$, is corrected using the current position of the body.
\label{proc:abga}
\end{myProcedure}

\newcommand{\stageI}{PB}
\newcommand{\stageII}{PV}
\newcommand{\stageIII}{P}
\newcommand{\stageIV}{CB}
\newcommand{\stageV}{CV}

\section{Discrete and approximate added-damping tensors}\label{sec:generalAddedDamping}

In this section we return to the formulas for the added-damping tensors in~\eqref{eq:DvvExact}--\eqref{eq:DwwExact} and describe an approach to obtain discrete approximations that can be used in the~\ampRB~time-stepping scheme.  Following the discussion in Part~I~\cite{rbinsmp2016r}, we begin with a fully-coupled (monolithic) discretization of the governing equations.  The discrete equations for the fluid in velocity-pressure form are
\begin{alignat}{3}
  \rho \left[ \frac{\vv_\iv^{n+1} - \vv_\iv^n}{\dt} \right] + \frac{3}{2} \fv_\iv^{n} - \half \fv_\iv^{n-1} &= 
      \mu\Big[ \alpha \Delta_h \vv_\iv^{n+1} + (1-\alpha)  \Delta_h \vv_\iv^{n}\Big], \quad&& \iv\in\Omega_h, \label{eq:AD1g}  \\
	   \Delta_h p_\iv^{n+1} &= -\grad_h\vv_{\iv}^{n+1} : {(\grad_h\vv_{\iv}^{n+1}) ^T} ,   \quad&&\iv\in\Omega_h ,\label{eq:AD4g}
\end{alignat}
where $\fv_{\iv}^n$ is given in~\eqref{eq:advPressure}.  Boundary conditions for the discrete velocity, $\vv_\iv\sp{n+1}$, and pressure, $p_\iv\sp{n+1}$, on the surface of the rigid body at $t\sp{n+1}$ are
\begin{alignat}{3}
 \vv_\iv^{n+1} &= \vvcm^{n+1} + \omegavb^{n+1}\times(\rv_\iv^{n+1}-\xvb^{n+1}),  \quad&&\iv\in\Gamma_h , \label{eq:AD2g} \\
 ({\nv_\iv} \cdot\grad_h) p_\iv^{n+1} &= {\nv_\iv^T}\Big( \avb^{n+1} + \bvb^{n+1}\times(\rv_\iv^{n+1}-\xvb^{n+1})  \nonumber \\
 & \qquad + \omegavb^{n+1}\times(\omegavb^{n+1}\times (\rv_\iv^{n+1}-\xvb^{n+1})) + \nu\Delta_h \vv_\iv^{n+1} \Big) , \quad
             &&\iv\in\Gamma_h,  \label{eq:AD5g}
\end{alignat}
These boundary conditions involve the discrete quantities, $\xvb^{n+1}$, $\vvcm^{n+1}$, etc., of the body which are determined by the following discretization of the Euler-Newton equations:
\begin{alignat}{3}
 &  \mass \avcm^{n+1} = \sum_{\iv\in\Gamma_h} \Big\{ {\nv_\iv^T} 
         \big(-p_\iv^{n+1}\Iv +  \Tv_\iv^{n+1} \big) \,{\Delta S}_\iv \Big\}~+\fvbe(t^{n+1}),    \label{eq:AD3Dab} \\
 &  \Ib \bvb^{n+1} + \omegavb^{n+1}\times \Ib \omegavb^{n+1} = 
        \sum_{\iv\in\Gamma_h} \Big\{(\rv_\iv^{n+1}-\xvb^{n+1})\times (- p_\iv^{n+1} \nv_\iv + \Tv_\iv^{n+1} ) \, {\Delta S}_\iv \Big\}
             ~+ \gvbe(t^{n+1}), \label{eq:AD7g} \\ 
 &  \frac{\vvcm^{n+1} - \vvcm^n}{\dt} =  \alphas \avcm^{n+1}  + (1-\alphas) \avcm^{n} \label{eq:AD3Dbb}  , \\
 &  \frac{\omegavb^{n+1} - \omegavb^n}{\dt} =  \alphas \bvb^{n+1}  + (1-\alphas) \bvb^{n} \label{eq:AD9g}  , \\
 &  \frac{\xvb^{n+1} - \xvb^n}{\dt} =  \alphas \vvb^{n+1}  + (1-\alphas) \vvb^{n} \label{eq:AD8g} , \\
 &  \frac{\Eb^{n+1} - \Eb^n}{\dt} =  \alphas \bigl(\omegavb^{n+1}\times\Eb^{n+1}\bigr)  + (1-\alphas) \bigl(\omegavb^{n}\times\Eb^{n}\bigr) \label{eq:AD10g} ,
\end{alignat}
where $\Tv_{\iv}^{n+1}$ is a discrete approximation of the shear stress, $\tauv\nv$, on the surface of the body, $\iv\in\Gamma_h$.  The discretization of the fluid also requires boundary conditions for $\iv\in\partial\Omega_h\backslash\Gamma_h$, but the choice of these boundary conditions is not important for the present discussion.  The parameters $\alpha$ in~\eqref{eq:AD1g} and $\alphas$ in~\eqref{eq:AD3Dbb}--\eqref{eq:AD10g} are time-stepping parameters, and typical values are taken to be $1/2$ corresponding to the trapezoidal rule.  The symbol, ${\Delta S}_\iv$, in~\eqref{eq:AD3Dab} and~\eqref{eq:AD7g} denotes the area weights in the quadrature-rule approximation of the surface integrals on the body. 

Discrete approximations for the derivatives, $\partial_{\vvb}\vv$ and~$\partial_{\omegavb}\vv$, that
appear in the formulas for the added-damping tensors in~\eqref{eq:DvvExact}--\eqref{eq:DwwExact} can
be obtained by solving variational problems derived from the discretization of the governing
equations.  The key equations among the coupled set of equations are~\eqref{eq:AD1g}
and~\eqref{eq:AD2g} for~$\vv_\iv\sp{n+1}$.  Let $\Wvv_\iv\in\Real^{3\times3}$ and $\Wvw_\iv\in\Real^{3\times3}$ given by
\[
\Wvv_\iv \eqdef \frac{\partial{\vv_\iv^{n+1}}}{\partial\vvcm^{n+1}},  \qquad
\Wvw_\iv \eqdef \frac{\partial{\vv_\iv^{n+1}}}{\partial\omegavb^{n+1}},
\]
denote discrete approximations of the derivatives, $\partial_{\vvb}\vv$ and~$\partial_{\omegavb}\vv$.  Taking the derivative of the equations with respect to $\vvb^{n+1}$ (with the other discrete quantities for the body held fixed) shows that $\Wvv_\iv$ satisfies a discrete Helmholtz problem given by 
\begin{alignat}{3} 
 & \frac{\rho}{\dt} \Wvv_\iv =   \mu \alpha \Delta_h \Wvv_\iv, \qquad&& \iv\in\Omega_h, \label{eq:AD1h}  \\
 &  \Wvv_\iv = \Iv ,  && \iv\in\Gamma_h , \label{eq:AD2h}
\end{alignat}
where $\Iv$ is the $3\times3$ identity and with homogeneous boundary conditions for $\Wvv_\iv$ on $\iv\in\partial\Omega_h\backslash\Gamma_h$ of the same form as that imposed for~\eqref{eq:AD1g}.  Similarly, taking the derivative with respect to $\omegavb^{n+1}$ gives
\begin{alignat}{3} 
 & \frac{\rho}{\dt} \Wvw_\iv =   \mu \alpha \Delta_h \Wvw_\iv, \qquad&& \iv\in\Omega_h, \label{eq:AD1i}  \\
 &  \Wvw_\iv = -[\rv_\iv^{n+1}-\xv_b^{n+1}]_\times ,  && \iv\in\Gamma_h , \label{eq:AD2i}
\end{alignat} 
with homogeneous boundary conditions for $\Wvw_\iv$ on $\iv\in\partial\Omega_h\backslash\Gamma_h$ as before.  
Assuming that $\Wvv_\iv$ and $\Wvw_\iv$ can be found, discrete forms corresponding 
to~\eqref{eq:DvvExact}, \eqref{eq:DwvExact}, \eqref{eq:DvwExactNew}, and \eqref{eq:DwwExactNew} are
\begin{alignat}{3} 
   &  \Dvv_h \eqdef - \mu \sum_{\iv\in\Gamma_h} \left\{(\Iv - \nv_\iv \nv_\iv^T) D_{nh}\Wvv_\iv\right\} \,{\Delta S}_\iv, \label{eq:DvvDiscrete} \\
   &  \Dvw_h \eqdef - \mu \sum_{\iv\in\Gamma_h}  \left\{ (\Iv - \nv_\iv \nv_\iv^T) D_{nh}\Wvw_\iv\right\} \,{\Delta S}_\iv , \label{eq:DvwDiscrete} \\
   &  \Dwv_h \eqdef - \mu \sum_{\iv\in\Gamma_h} [\rv_\iv^{n+1}-\xv_b^{n+1}]_\times \left\{ (\Iv - \nv_\iv \nv_\iv^T) D_{nh}\Wvv_\iv\right\} \,{\Delta S}_\iv , \label{eq:DvwvDiscrete} \\
   &  \Dww_h \eqdef 2\mu V_b \, \Iv - \mu \sum_{\iv\in\Gamma_h}  [\rv_\iv^{n+1}-\xv_b^{n+1}]_\times \left\{ (\Iv - \nv_\iv \nv_\iv^T) D_{nh}\Wvw_\iv \right\} \,{\Delta S}_\iv , \label{eq:DwwDiscrete}
\end{alignat}
where  $D_{nh}$ is a discrete approximation of the normal derivative on the surface of the rigid body.

While it is possible to evaluate the added-damping tensors in~\eqref{eq:DvvDiscrete}--\eqref{eq:DwwDiscrete} by solving the discrete Helmholtz problems for $\Wvv_\iv$ and $\Wvw_\iv$, this is considered to be too much effort when only approximations to these tensors are needed for the~\ampRB~time-stepping scheme.
Following the approximations made in the model-problem analysis in Part~I~\cite{rbinsmp2016r}, it is postulated that the normal derivatives of $\Wvv_\iv$ and $\Wvw_\iv$ on the surface of the rigid body can be approximated by 
\begin{equation}
D_{nh} \Wvv_\iv \approx -\frac{1}{\dn_\iv}\Iv,  \qquad 
D_{nh} \Wvw_\iv \approx \frac{1}{\dn_\iv}  [\rv_\iv^{n+1}-\xv_b^{n+1}]_\times, \qquad \iv\in\Gamma_h, \label{eq:Wapprox}
\end{equation}
where $\dn_\iv$ is the {\em added-damping length-scale} parameter given by
\begin{equation}
\dn_\iv \eqdef \frac{\Delta s_{n,\iv}}{1-e\sp{-\delta_\iv}},\qquad \delta_\iv \eqdef \frac{\Delta s_{n,\iv}}{\sqrt{\alpha\nu\dt}}, \qquad \iv\in\Gamma_h.
\label{eq:dn}
\end{equation}
Here, $\Delta s_{n,\iv}$ is the mesh spacing in the normal direction, $\alpha$ is the time-stepping parameter in~\eqref{eq:AD1g}, and $\nu=\mu/\rho$ is the kinematic viscosity.  As noted in Part~I, $\dn_\iv$ varies with $\delta_\iv$, which is the ratio of the mesh spacing in the normal direction to a {\em viscous} grid spacing, and it takes the limiting values
\[
  \dn_\iv \sim \begin{cases}
                 \sqrt{\alpha \nu \dt} & \text{for $\delta_\iv\rightarrow 0$} , \\
                 \Delta s_{n,\iv}      & \text{for $\delta_\iv\rightarrow \infty$} .
           \end{cases}
\] 
Using the approximations in~\eqref{eq:Wapprox} with $\dn_\iv$ given in~\eqref{eq:dn} leads to the
simplified added-damping tensors
\begin{alignat}{3} 
  & \Dvva_h \eqdef \mu  \sum_{\iv\in\Gamma_h} \frac{1}{\dn_\iv} \, (\Iv - \nv_\iv \nv_\iv^T) \,{\Delta S}_\iv , \label{eq:ADvv} \\
  & \Dvwa_h \eqdef \mu \sum_{\iv\in\Gamma_h} \frac{1}{\dn_\iv} \, (\Iv - \nv_\iv \nv_\iv^T)  [\rv_\iv^{n+1}-\xv_b^{n+1}]_\times^T \,{\Delta S}_\iv , \label{eq:ADvw} \\
  & \Dwva_h \eqdef \mu \sum_{\iv\in\Gamma_h} \frac{1}{\dn_\iv} \,  [\rv_\iv^{n+1}-\xv_b^{n+1}]_\times  (\Iv - \nv_\iv \nv_\iv^T)   \,{\Delta S}_\iv , \label{eq:ADwv} \\
  & \Dwwa_h \eqdef \mu \sum_{\iv\in\Gamma_h} \frac{1}{\dn_\iv} \, [\rv_\iv^{n+1}-\xv_b^{n+1}]_\times (\Iv - \nv_\iv \nv_\iv^T)  [\rv_\iv^{n+1}-\xv_b^{n+1}]_\times^T \,{\Delta S}_\iv .\label{eq:ADww}
\end{alignat}
The term involving the volume of body has been dropped from the approximations in~\eqref{eq:ADww} 
since this term will generally be small compared to the term that was kept.

The approximate composite added-damping tensor given by
\begin{equation}
  \ADtensora \eqdef   \begin{bmatrix}
    \Dvva_h & \Dvwa_h  \\
    \Dwva_h & \Dwwa_h
  \end{bmatrix} ,  \label{eq:AddedDampingTensorApprox} 
\end{equation}
is used in the~\ampRB~time-stepping scheme.  We note that $\ADtensora$ is symmetric and positive semi-definite, and thus it provides a damping contribution to the calculation of the accelerations of the rigid body in~\eqref{eq:RBaccelerationInterface}.
We note also that the added-damping tensors in $\ADtensora$ need only be computed for the problem configuration at the initial time $t\sp{n}=0$, since the tensors at later times can be computed based on the rotation of the rigid body according to~\eqref{eq:addedDampingInTime}.

\section{Numerical approach using moving composite grids} \label{sec:numericalApproach}

Our numerical approach for the solution of the equations governing an
FSI initial-boundary-value problem is based on the use of moving (and possibly deforming)
composite grids. This flexible approach enables the use of efficient structured grids
for complex geometry and provides smooth, high-quality grids for moving grid 
problems even as the geometry undergoes large changes.  
We have previously applied this approach to the coupling of 
compressible flows and rigid bodies~\cite{mog2006,lrb2013}, compressible flows
with deforming bodies in~\cite{fsi2012,flunsi2016}, and incompressible flows
with rigid bodies~\cite{Koblitz2016} and deforming beams~\cite{beamins2016}. 

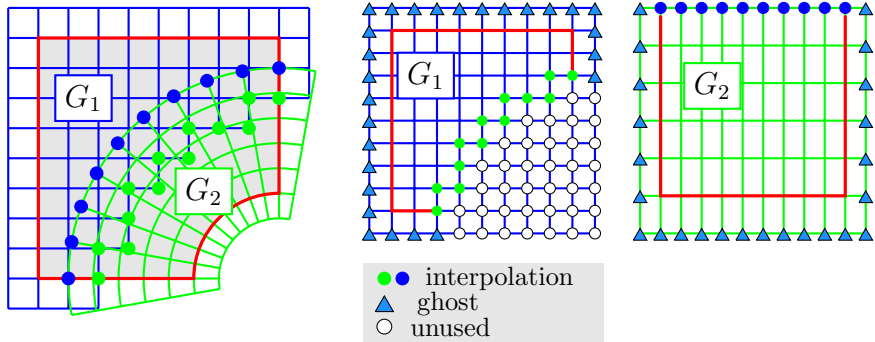
\begin{figure}[hbt]
\begin{center}
\begin{tikzpicture}[scale=.8]
\useasboundingbox (.75,1.2) rectangle (15.5,6.2);  
%
\begin{scope}[xshift=1cm,yshift=1cm]
\fill[black!10!white,xshift=.5cm,yshift=.5cm] (0,0) -- (2.583333,0) arc (180:90:1.416667) -- (4.,4.) -- (0,4.) -- (0,0);
\draw[-,thick,blue,yshift=.0 cm] 
   \foreach \x/\y in {1.5/0,1.5/.5,2/1,2/1.5,2.5/2,3/2.5,4/3,5/3.5,5/4,5/4.5,5/5}{ (0,\y) -- (\x,\y) }
   \foreach \x/\y in {0/0,.5/0,1/0,1.5/0,2/1,2.5/2,3/2.5,3.5/3,4/3,4.5/3.5,5/3.5}{ (\x,\y) -- (\x,5) };
  \begin{scope}[xshift=4.5cm,yshift=0.5cm]
    \draw[thick,green] \foreach \r in {1.000000,1.416667,1.833333,2.250000,2.666667,3.083333,3.500000}{ (0,\r) arc (90:190:\r)  (0,\r) arc (90:80:\r) };
    \draw[thick,green]
     (0.173648,0.984808)  -- (0.607769,3.446827)
     (0.000000,1.000000)  -- (0.000000,3.500000)
     (-0.173648,0.984808) -- (-0.607769,3.446827)
     (-0.342020,0.939693) -- (-1.197071,3.288924)
     (-0.500000,0.866025) -- (-1.750000,3.031089)
     (-0.642788,0.766044) -- (-2.249757,2.681156)
     (-0.766044,0.642788) -- (-2.681156,2.249757)
     (-0.866025,0.500000) -- (-3.031089,1.750000)
     (-0.939693,0.342020) -- (-3.288924,1.197071)
     (-0.984808,0.173648) -- (-3.446827,0.607769)
     (-1.000000,0.000000) -- (-3.500000,0.000000)
     (-0.984808,-0.173648) -- (-3.446827,-0.607769);
  \end{scope}
  \draw[very thick,red,xshift=.5cm,yshift=.5cm] (0,0) -- (2.583333,0) arc (180:90:1.416667) -- (4.,4.) -- (0,4.) -- (0,0);
%
   \filldraw[green] (1.5,.5)  circle (3pt)
                 (1.5,1 )  circle (3pt)
                 (2  ,1 )  circle (3pt)
                 (2 ,1.5)  circle (3pt)
                 (2 , 2 )  circle (3pt)
                 (2.5,2 )  circle (3pt)
                 (2.5,2.5) circle (3pt)
                 (3 , 2.5) circle (3pt)
                 (3  ,3 )  circle (3pt)
                 (3.5,3 )  circle (3pt)
                 (4  ,3. ) circle (3pt)
                 (4  ,3.5) circle (3pt)
                 (4.5,3.5) circle (3pt);
%
  \begin{scope}[xshift=4.5cm,yshift=0.5cm]
      \filldraw[blue]
       (0.000000,3.500000)    circle (3pt)
       (-0.607769,3.446827)   circle (3pt)
       (-1.197071,3.288924)  circle (3pt) 
       (-1.750000,3.031089)  circle (3pt) 
       (-2.249757,2.681156)  circle (3pt) 
       (-2.681156,2.249757)  circle (3pt) 
       (-3.031089,1.750000)  circle (3pt) 
       (-3.288924,1.197071)  circle (3pt) 
       (-3.446827,0.607769)  circle (3pt) 
       (-3.500000,0.000000)  circle (3pt);
  \end{scope}
   \draw (1.25,3.5) node[thick,draw=blue,fill=white] {\large$G_1$};
   \draw (3.25,1.95) node[thick,draw=green,fill=white] {\large$G_2$};
\end{scope}
%
\definecolor{ghostColour}{named}{DodgerBlue}
\newcommand{\mytrix}{(\x,-.15) -- ++(.3,0) -- ++(-.15,.26) -- (\x,-.15)}
\newcommand{\mytriy}{(-.15,\y) -- ++(.3,0) -- ++(-.15,.26) -- (-.15,\y)}
\begin{scope}[xshift=7cm,yshift=2.25cm,scale=.75]
\draw[-,thick,blue,yshift=.0 cm] 
   \foreach \x in {0,.5,...,5}{ (\x,0) -- (\x,5) }
   \foreach \y in {0,.5,...,5}{ (0,\y) -- (5,\y) };
  \draw[very thick,red,xshift=.5cm,yshift=.5cm] (1.,0) -- (.0,0) -- (.0,4.) -- (4.,4.) -- (4.,3.);
   \filldraw[green] (1.5,.5)  circle (3pt)
                 (1.5,1 )  circle (3pt)
                 (2  ,1 )  circle (3pt)
                 (2 ,1.5)  circle (3pt)
                 (2 , 2 )  circle (3pt)
                 (2.5,2 )  circle (3pt)
                 (2.5,2.5) circle (3pt)
                 (3 , 2.5) circle (3pt)
                 (3  ,3 )  circle (3pt)
                 (3.5,3 )  circle (3pt)
                 (4  ,3. ) circle (3pt)
                 (4  ,3.5) circle (3pt)
                 (4.5,3.5) circle (3pt);
  \filldraw[fill=white,draw=black]  \foreach \x in {2,2.5,...,5}{ (\x,.0) circle (3.5pt) };
  \filldraw[fill=white,draw=black]  \foreach \x in {2,2.5,...,5}{ (\x,.5) circle (3.5pt) };
  \filldraw[fill=white,draw=black]  \foreach \x in {2.5,3,...,5}{ (\x,1.) circle (3.5pt) };
  \filldraw[fill=white,draw=black]  \foreach \x in {2.5,3,...,5}{ (\x,1.5) circle (3.5pt) };
  \filldraw[fill=white,draw=black]  \foreach \x in {3,3.5,...,5}{ (\x,2.0) circle (3.5pt) };
  \filldraw[fill=white,draw=black]  \foreach \x in {3.5,4,...,5}{ (\x,2.5) circle (3.5pt) };
  \filldraw[fill=white,draw=black]  \foreach \x in {4.5,5}      { (\x,3.0) circle (3.5pt) };
  \draw[fill=ghostColour,xshift=-.15cm,yshift=0cm]  \foreach \x in {.5,1.,1.5}{ \mytrix };  
  \draw[fill=ghostColour,xshift=-.15cm,yshift=5cm]  \foreach \x in {.5,1.,...,5}{ \mytrix };  
  \draw[fill=ghostColour,xshift=0cm,yshift=-.15cm]  \foreach \y in {0,.5,...,5}{ \mytriy };
  \draw[fill=ghostColour,xshift=5cm,yshift=-.15cm]  \foreach \y in {3.5,4,4.5}{ \mytriy };
   \draw (1.25,3.5) node[thick,draw=blue,fill=white] {\large$G_1$};
\end{scope}
\begin{scope}[xshift=11.5cm,yshift=2.25cm,scale=.75]
\draw[-,thick,green,yshift=.0 cm] 
   \foreach \x in {0,.454545,...,5}{ (\x,0) -- (\x,5) }
   \foreach \y in {0,.833333,...,5}{ (0,\y) -- (5,\y) };
 \draw[very thick,red,xshift=.454545cm,yshift=.833333cm] (0.,4) -- (.0,0) -- (4.0909,0.) -- (4.0909,4);
 \filldraw[blue]  \foreach \x in {.454545,.909090,...,4.545454}{ (\x,5) circle (3.5pt) };
 \draw[fill=ghostColour,xshift=-.15cm]  \foreach \x in {.454545,.909090,...,4.545454}{ \mytrix };
 \draw[fill=ghostColour,yshift=-.15cm]  \foreach \y in {0,.833333,...,5}{ \mytriy };
 \draw[fill=ghostColour,xshift=5cm,yshift=-.15cm]  \foreach \y in {0,.833333,...,5}{ \mytriy };
\end{scope}
\begin{scope}[xshift=7cm,yshift=.7cm]
  \fill[black!10!white,xshift=-.1cm,yshift=-.25cm] (0,0) -- (4,0) -- (4.,1.3) -- (0,1.3) -- (0,0);
  \filldraw[green,xshift=.0cm,yshift=.8cm] (.25,.0)  circle (3pt);
  \filldraw[blue,xshift=.3cm,yshift=.8cm] (.25,.0)  circle (3pt);
  \draw[xshift=.0cm,yshift=.8cm] (.5,0) node[anchor=west,xshift=6] {\small interpolation};
  \draw[fill=ghostColour,xshift=.0cm,yshift=.4cm] (.35,0) \foreach \x in {.1}{ \mytrix } node[anchor=west,xshift=12,yshift=3] {\small ghost};
  \draw[fill=white,draw=black,xshift=.0cm,yshift=.0cm] (.25,0) circle (3.5pt) node[anchor=west,xshift=6] {\small unused};
\end{scope}
\begin{scope}[xshift=11.5cm,yshift=2.25cm,scale=.75]
   \draw (1.6,3.27) node[thick,draw=green,fill=white] {\large$G_2$};
\end{scope}
\end{tikzpicture}
\end{center}
\caption{Left: an overlapping grid consisting of two
structured curvilinear component grids, $\xv=G_1(\rv)$ and $\xv=G_2(\rv)$. Middle and right: 
component grids for the square and annular grids in the unit square parameter space $\rv$. Grid
 points are classified as discretization points, interpolation points or unused points. Ghost points
 are used to apply boundary conditions.
}   
 \label{fig:overlappingGridCartoon}
\end{figure}

An overlapping grid, $\Gs$, consists of a set of structured component grids, $\{G_g\}$,
$g=1,\ldots,{\mathcal N}$, that covers
the fluid domain, $\Omega(t)$, and overlap where
the component grids meet, see Figure~\ref{fig:overlappingGridCartoon}. 
Discrete solutions defined on different grids are matched by interpolation at the {\em interpolation points}; 
for the second-order accurate
computations performed here, a tensor-product Lagrange interpolation formula is used (i.e. quadratic
interpolation using a three-point stencil in each direction).
Typically, boundary-fitted curvilinear grids are used near the
boundaries (rigid-body surfaces and external boundaries), while one or more
background Cartesian grids are used to handle the bulk of the fluid domain. 
Each component grid is a logically rectangular, curvilinear grid in $\nd$~space
dimensions, and is defined by a smooth
mapping from parameter space~$\rv$ (the unit square or cube) to physical
space~$\xv$,
\[
  \xv = \Gv(\rv,t),\qquad \rv\in[0,1]^\nd,\qquad \xv\in\Real^\nd.
\]
For the present FSI problem, there are usually one or more background grids which are static, while boundary-fitted grids attached
to the boundary of the rigid body move over time.
The Navier-Stokes equations for the fluid are transformed to the unit-square reference coordinates, $\rv$, 
using the chain rule. On moving grids, the equations are transformed to a moving coordinate
system which introduces the grid velocity into the advection terms as indicated
in~\eqref{eq:advPressure}. The resulting equations are then discretized using standard 
finite-difference approximations for the derivatives with respect to~$\rv$.  Second-order accurate
approximations are used in the present implementation of the time-stepping scheme.
For more details on the discretization approach see~\cite{ICNS,mog2006,beamins2016}. 

\paragraph{Time-step determination} The discrete time-step $\dt$ is generally determined by a CFL-type
stability condition based
on the advection terms in the fluid momentum equations.  The viscous terms are treated
implicitly so that there is no stability constraint on~$\dt$ arising from these terms.
A maximum value for the time-step is also enforced in case the advection terms
are small, and this maximum value is often taken to be proportional to the grid spacing (in non-dimensional variables).
An additional consideration when choosing the time-step for the~\ampRB~scheme is that
since the entries in the added-damping tensors depend on $\dt$ (through the factor $\sqrt{\nu\dt}$), it is helpful if $\dt$ varies
slowly as the simulation progresses. 
It has been found that if $\dt$ jumps by a significant fraction, e.g. by
a factor of two, then for difficult problems with very light bodies, the body acceleration may also
experience a small jump (as can can be seen, for example, in Figure~\ref{fig:risingDropCurves}). This jump is local in time and the acceleration
recovers after a few time-steps but this can be avoided by either enforcing that $\dt$ varies smoothly, or
by adjusting the form of the added-damping coefficients so that these vary smoothly in time.

\section{Numerical Results} \label{sec:numericalResults}

Numerical results are now presented that demonstrate the stability and 
accuracy properties of the {\ampRB}~scheme as implemented using moving 
overlapping grids\footnote{The examples presented in this section are available with the Overture software at~{\tt overtureFramework.org}.}.
The first problem involves the one-dimensional motion of a piston in a rectangular fluid chamber
while the second problem involves the rotation of a solid disk in an annular fluid chamber.  These problems are
motivated by two of the model problems discussed in Part~I~\cite{rbinsmp2016r}.
The first problem highlights added-mass effects while the second problem isolates
the effects of added-damping.  Unlike the model problems considered in Part~I, the
present problems consider finite amplitude translations and rotations, and overlapping
grids are used to handle the moving geometry.  Exact solutions are available for
both problems, and these solution can be used to assess the stability and
accuracy of the general {\ampRB}~scheme. 

Subsequently, four other challenging problems are considered to illustrate the behaviour of the
scheme.  
In the first case, a zero-mass cylindrical body moves within a fluid-filled channel. This
problem is designed to provide a clean benchmark problem for evaluating the stability and
accuracy of the {\ampRB}~scheme when both added-mass and added-damping effects are important. 
Convergence rates are computed from a self-convergence grid refinement study.
The next problem considers a standard test problem examined in the literature consisting 
of a moderately-heavy cylindrical-body that falls, under gravity, in a fluid
channel. This problem is used 
to check the results of the {\ampRB}~scheme with results from other schemes.
The third of these four problems involves a light rectangular-shaped
body rising under buoyancy forces to the top of a closed fluid chamber.
This problem demonstrates
the need for added-mass corrections, and also the need for the added-damping corrections
of the~\ampRB~scheme for the case of a non-cylindrically-shaped body. 
For non-cylindrical bodies, one might posit (incorrectly, it turns out)
that any finite-size added-mass corrections would be sufficient to stabilize the scheme, without the need
for added-damping corrections, since added-mass effects are proportional to changes in the
acceleration of the body while added-damping effects are proportional to changes in the velocity; the
latter generally being $O(\Delta t)$ smaller than the former.
Despite this intuition, it is found in practice that added-damping effects must be properly treated to
maintain stability of the scheme.
In the last problem, the interaction between two rectangular-shaped bodies, one rising and
one falling, in a rectangular fluid chamber is simulated to demonstrate the performance
of the scheme with multiple bodies in close proximity.

\subsection{One-dimensional motion of a piston and a rectangular fluid chamber} \label{sec:piston2d}

Consider the horizontal motion of a two-dimensional rigid-body (a {\em piston}) located at one
end of a rectangular channel of fluid as shown in Figure~\ref{fig:pistonGrids}.  The problem is
posed with initial conditions and boundary conditions so that an exact solution exists which only
depends on the horizontal coordinate~$x$ and on time~$t$.  The pressure in the fluid
provides a force on the rigid-body, which implies the presence of added-mass effects,
whereas there are no forces due
to viscous shear for this problem so that added-damping effects are negligible.  Thus, this FSI problem provides a
good test of the added-mass properties of the~\ampRB~scheme. 
Note that although the exact solution for this problem does not depend on~$y$, the numerical
solution of the full incompressible Navier-Stokes equations in two dimensions is computed
on general moving grids to verify the overall approach.

{
\newcommand{\rbxa}{0} \newcommand{\rbxb}{4}
\newcommand{\rbya}{0} \newcommand{\rbyb}{4}
\newcommand{\plotRigidBody}{
\fill[fill=red!20,draw=red,line width=2pt] 
        (\rbxa,\rbya) -- (\rbxb,\rbya) -- (\rbxb,\rbyb) -- (\rbxa,\rbyb) -- cycle;
}
\newcommand{\figWidth}{4.22cm}
\newcommand{\trimfig}[2]{\trimh{#1}{#2}{.33}{.025}{.28}{.28}}
\newcommand{\figWidtha}{7.0cm}
\newcommand{\trimfiga}[2]{\trimFig{#1}{#2}{.0}{.0}{.0}{.0}}
\begin{figure}[htb]
\begin{center}
\resizebox{7cm}{!}{
\begin{tikzpicture}[scale=1]
  \useasboundingbox (-.2,.85) rectangle (11,9.2);  
  \begin{scope}[yshift=5cm]
    \draw(4.0,0.0) node[anchor=south west,xshift=-15pt,yshift=-8pt] {\trimfig{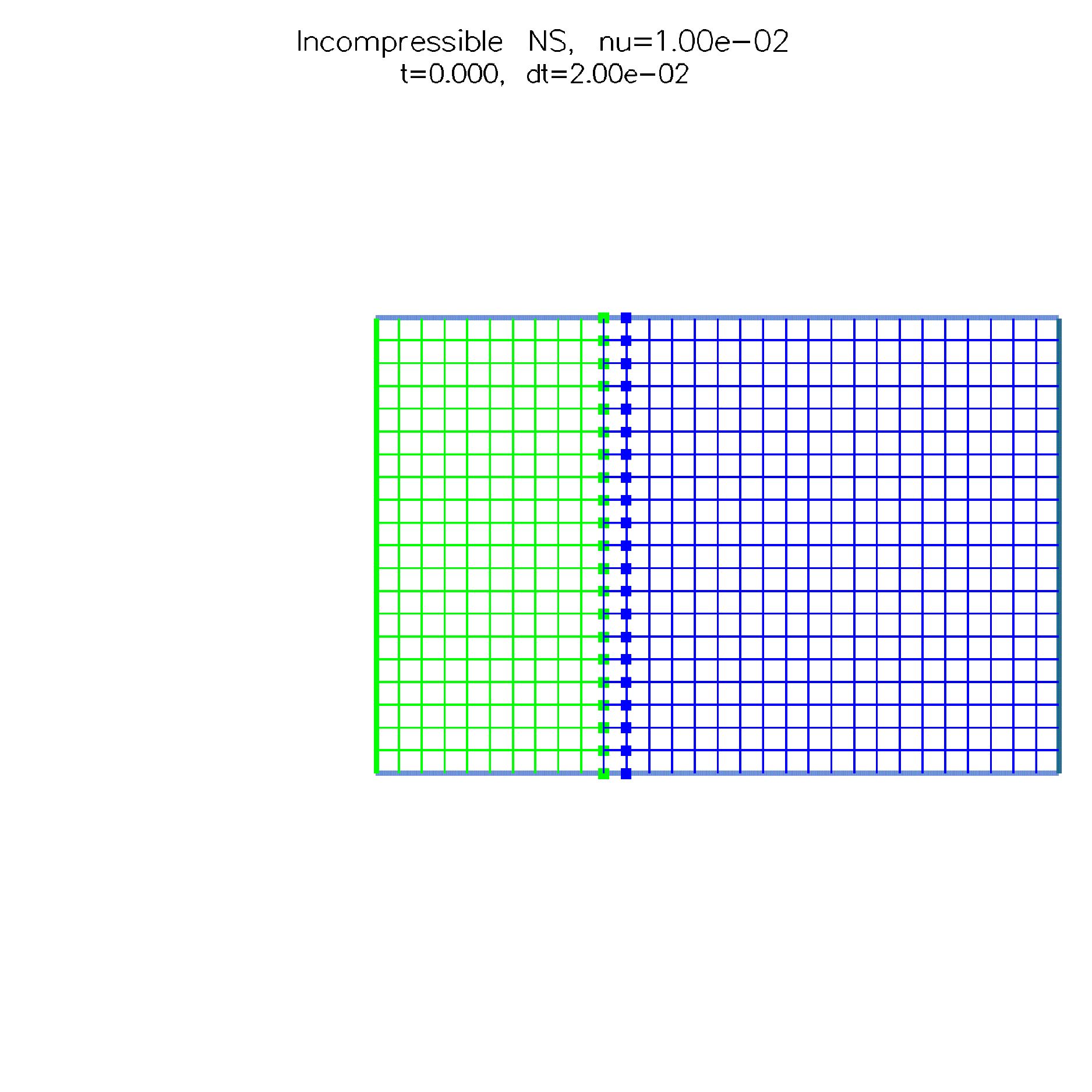}{\figWidth}};
    \begin{scope}[xshift=0cm,yshift=0cm]
      \plotRigidBody
    \end{scope}
    \draw[<->,thick,yshift=1cm] (0,0) -- (4,0); 
     \draw (3,1) node[anchor=north] {$\pistonWidth$};
     \draw[<->,thick,xshift=2cm] (0,0) -- (0,4); 
     \draw (2,2) node[anchor=west] {$\pistonHeight$};
  %
    \begin{scope}[xshift=2cm,yshift=-4pt]
     \draw[->,thick,yshift=0pt] (-2.25,0) -- (9,0) node [anchor=north] {$x$}; 
     \draw[-,thick,xshift=4pt,yshift=0pt ] (8,.1) -- (8.,-.1) node [anchor=north] {$\channelRight$}; 
     \draw[-,thick,xshift=2pt,yshift=0pt ] (2,.1) -- (2.,-.1) node [anchor=north] {$x_I(t)$}; 
     \draw[-,thick,xshift=0pt,yshift=0pt ] (0,.1) -- (0.,-.1) node [anchor=north] {$0$}; 
    \end{scope}
  \end{scope}
%
  \begin{scope}[yshift=0cm]
  \draw(4.0, 0.0) node[anchor=south west,xshift=-15pt,yshift=-8pt] {\trimfig{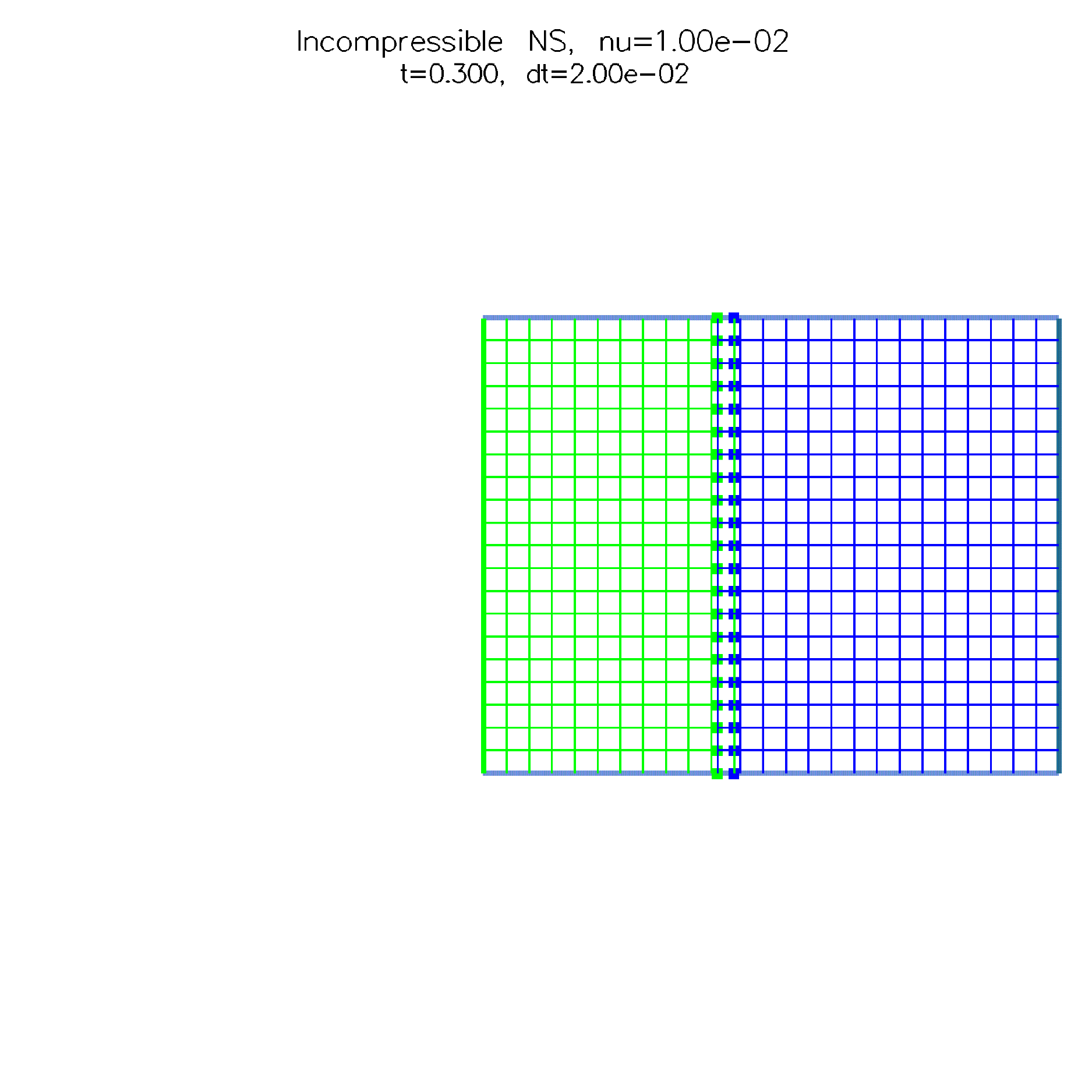}{\figWidth}};
    \begin{scope}[xshift=.95cm,yshift=.0cm]
      \plotRigidBody
      \draw[->,line width=2.pt,red,yshift=0pt] (1.5,2) -- (2.5,2);
    \end{scope}
  \end{scope}
\end{tikzpicture}
}
\end{center}
  \caption{Composite grids for a piston adjacent to a fluid channel. Grids are shown at two times.
    The green fluid grid adjacent to the piston moves over time. The blue background grid remains fixed.
}
  \label{fig:pistonGrids}
\end{figure}
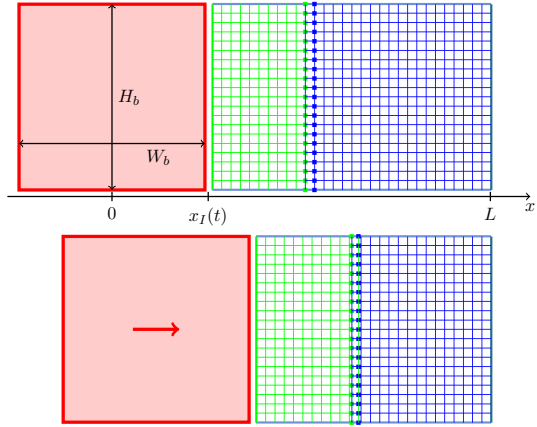
}

For this piston problem, the fluid occupies the rectangular channel
$\OmegaF(t)=[x_I(t),\channelRight]\times[0,\channelHeight]$, where $x_I(t)$ is the position of the
interface between the fluid and solid, $\channelRight$ is the position of the right end of the channel, and $\channelHeight$ is the height of the channel. The rigid body, located at the left end of the fluid channel,
has mass~$\mb$, width~$\pistonWidth$ and height~$\pistonHeight=\channelHeight$.  
The interface condition on the fluid at the
piston face is the velocity matching condition in~\eqref{eq:RBsurfaceVelocity}.  The boundary conditions on the top and bottom of the fluid channel
are slip-walls (equivalent to symmetry conditions), while the condition on the fluid at
$x=\channelRight$ is taken to be
\[
    p(\channelRight,y,t) = p_\channelRight(t), \qquad v_2(\channelRight,y,t) = 0,\qquad y\in[0,\channelHeight],\quad t>0,
\]
where $p_\channelRight(t)$ is a given applied pressure.  The solution for the body velocity $\vvb=(v_{1,b},v_{2,b})\sp{T}$, fluid velocity $\vv=(v_{1},v_{2})\sp{T}$ and fluid pressure $p$ is given by
\begin{align}
v_{1,b}(t) &= \int_0\sp{t}{-Hp_\channelRight(\tau)\,d\tau\over \mb+M_a(\tau)}+v_{1,b}(0), \label{eq:horizontalVelocity} \\
p(x,t) &= p_\channelRight(t) + \left({\channelRight-x\over\channelRight-x_I(t)}\right){-p_\channelRight(t)\over \mb/M_a(t)+1}, \label{eq:channelPressure} \\
v_1(t) &= v_{1,b}(t),
\end{align}
and $v_2=v_{2,b}=0$.  Here, the horizontal position of the body centre and interface are given by
\begin{equation}
x_b(t) = \int_0\sp{t} v_{1,b}(\tau)\,d\tau + x_b(0), \qquad x_I(t)=x_b(t)+\pistonWidth/2,
\label{eq:horizontalPosition}
\end{equation}
and the {\em added-mass} is determined analytically as
\begin{equation}
M_a(t)=\rho\channelHeight\bigl(\channelRight-x_I(t)\bigr).
\label{eq:addedMass}
\end{equation}
The solution in~\eqref{eq:horizontalVelocity}--\eqref{eq:addedMass} is an extension of the one for
the added-mass model problem in Part~I~\cite{rbinsmp2016r}.  Here, however, there is no
small-amplitude linearization about a fixed interface position, $x_I(t)=0$, as was done in Part~I.
While the solution can be determined from a given
applied pressure and initial states, $x_b(0)$ and $v_{1,b}(0)$, it is simpler to choose $x_b(t)$ and
$v_{1,b}(t)=\dot x_b(t)$, and then back out $p_\channelRight(t)$ and $p(x,t)$
using~\eqref{eq:horizontalVelocity} and~\eqref{eq:channelPressure}, respectively.  Following this
approach, we set
\[
x_b(t) = \Amplitude\sin(2 \pi t),\qquad \Amplitude=1/4,
\]
to determine solutions for various choices of the parameters in the verification tests below.

{
\newcommand{\figWidtha}{7.0cm}
\newcommand{\trimfiga}[2]{\trimh{#1}{#2}{.0}{.0}{.0}{.0}}
\begin{figure}[htb]
\begin{center}
\resizebox{12cm}{!}{
\begin{tikzpicture}[scale=1]
  \useasboundingbox (-.2,1.2) rectangle (17.,7.);  
  \draw(0,0) node[anchor=south west,xshift=-4pt,yshift=+0pt] {\trimfiga{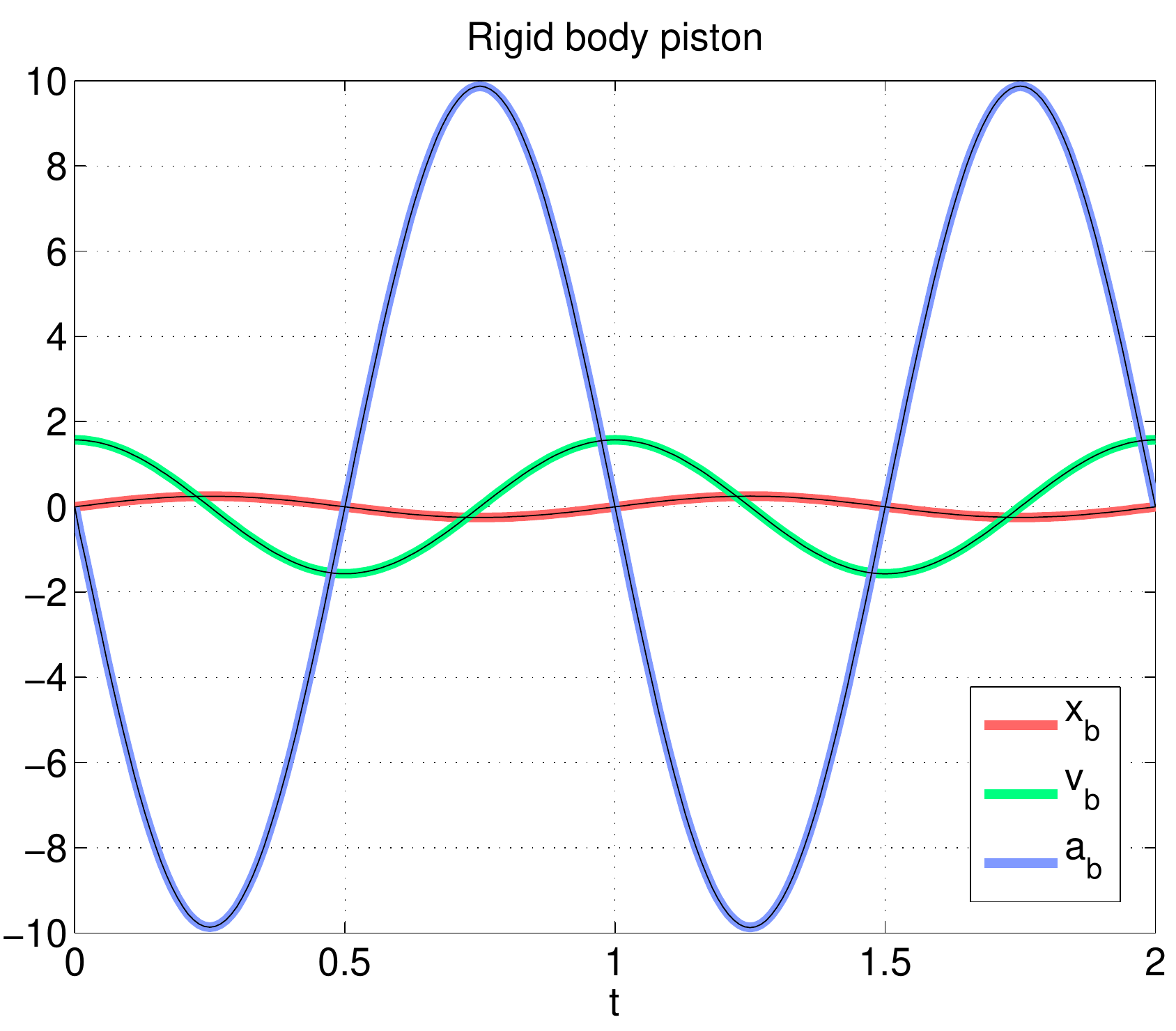}{\figWidtha}};
  \draw(8.5,0) node[anchor=south west,xshift=-4pt,yshift=+0pt] {\trimfiga{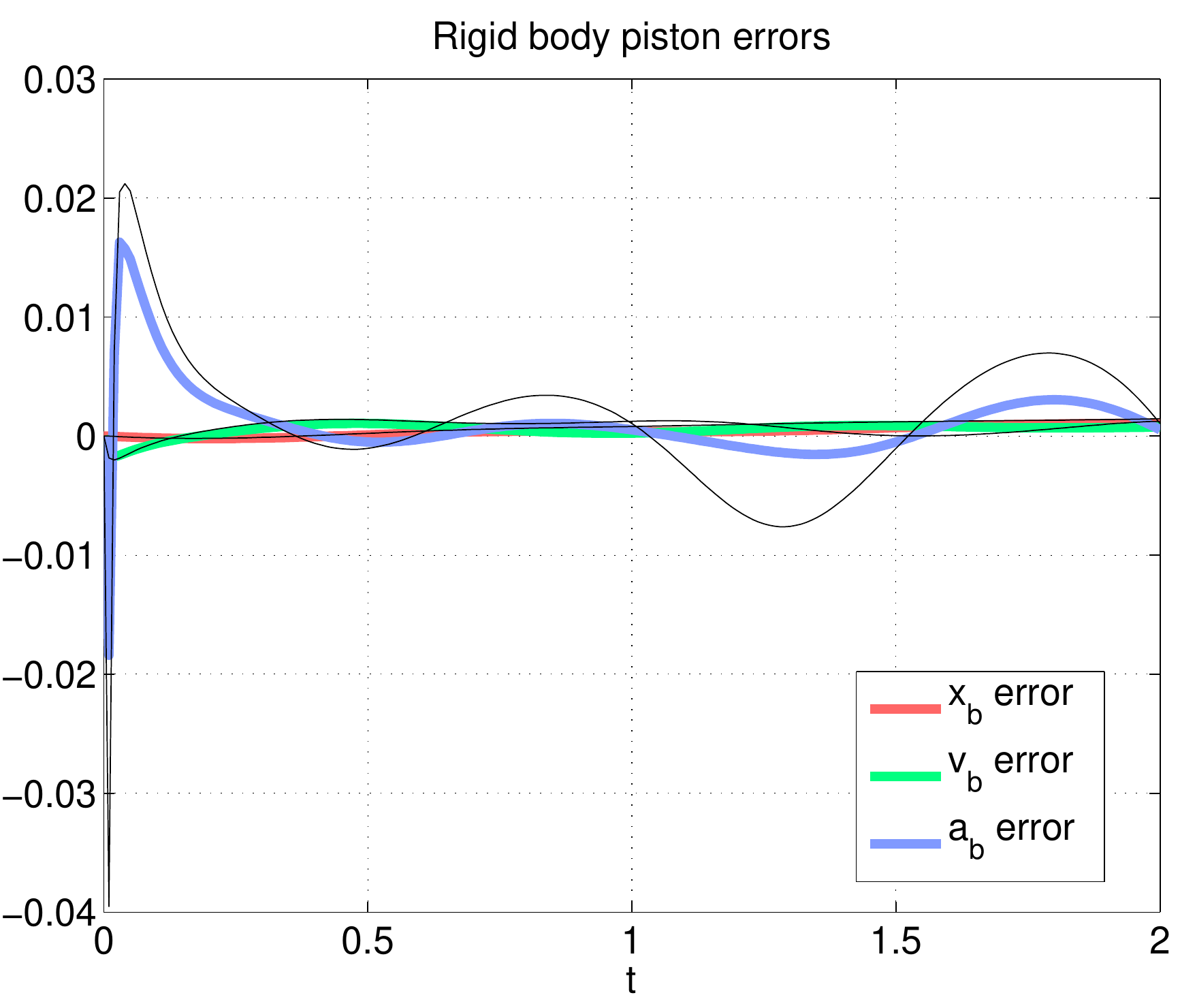}{\figWidtha}};
%
\end{tikzpicture}
}
\end{center}
  \caption{Computed piston motion and errors for $\rhos=1$ (coloured lines) and $\rhos=.01$ (black lines). Results for $\Gcp^{(4)}$.
   }
  \label{fig:pistonCurves}
\end{figure}
}

Numerical solutions are computed for the piston problem using the composite grid, denoted by
$\Gcp^{(j)}$, with resolution factor~$j$. As illustrated in
Figure~\ref{fig:pistonGrids}, the grid consists of two component grids.  The first is a body-fitted
rectangular grid (green in the figure) of fixed width $1/2$ and height $\channelHeight=1$,
which is attached to the right face of the rigid body.  This moving grid overlaps with a static
Cartesian background grid (blue in the figure) covering the region,
$[-3/4,\channelRight]\times[0,\channelHeight]$ with $\channelRight=3/2$, which covers the remaining
portion of the fluid channel.  The grid spacing for both component grids is uniform and chosen to be
$\hj=1/(10 j)$ in each direction.  The width and height of the rigid body are taken to be
$\pistonWidth=\pistonHeight=1$ so that $\mb=\rho_b$, where $\rho_b$ is the density of the solid.
The density of the fluid is taken to be $\rho=1$ and its viscosity is $\mu=0.1$.  In our numerical
tests, the density of the fluid is held fixed, while the density of the solid is varied to study the
effects of added mass.

{ 
\newcommand{\errp}{$E_j^{p}$}
\newcommand{\erruv}{$E_j^{\vv}$}
\newcommand{\errxA}{$E_j^{\xvb}$}
\newcommand{\errvA}{$E_j^{\vvb}$}
\newcommand{\erraA}{$E_j^{\avb}$}
\newcommand{\noOrder}{--}
%
\newcommand{\pistonTableI}{%
\begin{tabular}{|c|c|c|c|c|c|c|c|c|c|c|} \hline 
  \multicolumn{11}{|c|}{Piston motion, $\rhob=10$} \\ \hline 
\strutt$\hj$&     \errp     &  r   &     \erruv     &  r   &     \errxA     &  r   &     \errvA     &  r   &     \erraA     &  r    \\[3pt] \hline 
  1/10   & \num{1.3}{-3} &      & \num{9.9}{-3} &      & \num{3.7}{-2} &      & \num{9.9}{-3} &      & \num{1.2}{-4} &     \\ \hline
  1/20   & \num{4.4}{-5} & \noOrder & \num{2.1}{-3} & 4.6  & \num{9.1}{-3} & 4.1  & \num{2.1}{-3} & 4.6  & \num{4.3}{-6} & \noOrder \\ \hline
  1/40   & \num{5.4}{-8} & \noOrder  & \num{5.3}{-4} & 4.1  & \num{2.3}{-3} & 3.9  & \num{5.3}{-4} & 4.1  & \num{4.6}{-9} & \noOrder \\ \hline
  1/80   & \num{4.7}{-10} &\noOrder   & \num{1.3}{-4} & 4.0  & \num{6.0}{-4} & 3.9  & \num{1.3}{-4} & 4.0  & \num{9.8}{-13} & \noOrder \\ \hline
  rate   &    \noOrder       &      &    2.07       &      &    1.98       &      &    2.07       &      &    \noOrder       &      \\ \hline
\end{tabular}
}
\newcommand{\pistonTableII}{%
\begin{tabular}{|c|c|c|c|c|c|c|c|c|c|c|} \hline 
  \multicolumn{11}{|c|}{Piston motion, $\rhob=1$} \\ \hline 
\strutt$\hj$&     \errp     &  r   &     \erruv     &  r   &     \errxA     &  r   &     \errvA     &  r   &     \erraA     &  r    \\[3pt] \hline 
  1/10   & \num{1.3}{-4} &      & \num{4.8}{-2} &      & \num{4.6}{-2} &      & \num{4.8}{-2} &      & \num{1.1}{-4} &     \\ \hline
  1/20   & \num{1.1}{-6} & \noOrder  & \num{9.9}{-3} & 4.8  & \num{1.1}{-2} & 4.3  & \num{9.9}{-3} & 4.8  & \num{5.7}{-7} & \noOrder \\ \hline
  1/40   & \num{8.7}{-9} & \noOrder  & \num{2.4}{-3} & 4.1  & \num{2.7}{-3} & 4.0  & \num{2.4}{-3} & 4.1  & \num{1.1}{-9} & \noOrder \\ \hline
  1/80   & \num{7.8}{-12} &\noOrder  & \num{6.1}{-4} & 4.0  & \num{6.9}{-4} & 3.9  & \num{6.1}{-4} & 4.0  & \num{7.5}{-13} &\noOrder \\ \hline
  rate   &    \noOrder       &      &    2.09       &      &    2.02       &      &    2.09       &      &    \noOrder       &      \\ \hline
\end{tabular}
}
\newcommand{\pistonTableIII}{%
\begin{tabular}{|c|c|c|c|c|c|c|c|c|c|c|} \hline 
  \multicolumn{11}{|c|}{Piston motion, $\rhob=.001$} \\ \hline 
\strutt$\hj$&     \errp     &  r   &     \erruv     &  r   &     \errxA     &  r   &     \errvA     &  r   &     \erraA     &  r    \\[3pt] \hline

  1/10   & \num{6.4}{-5} &      & \num{8.3}{-2} &      & \num{5.5}{-2} &      & \num{8.3}{-2} &      & \num{5.1}{-4} &     \\ \hline
  1/20   & \num{3.3}{-6} & \noOrder & \num{1.7}{-2} & 4.9  & \num{1.3}{-2} & 4.4  & \num{1.7}{-2} & 4.8  & \num{3.1}{-5} &  \noOrder \\ \hline
  1/40   & \num{2.5}{-8} & \noOrder  & \num{4.2}{-3} & 4.1  & \num{3.1}{-3} & 4.0  & \num{4.2}{-3} & 4.1  & \num{2.8}{-8} & \noOrder  \\ \hline
  1/80   & \num{4.2}{-10} &\noOrder  & \num{1.0}{-3} & 4.0  & \num{8.0}{-4} & 3.9  & \num{1.0}{-3} & 4.0  & \num{1.5}{-12} & \noOrder \\ \hline
  rate   &    \noOrder       &      &    2.10       &      &    2.03       &      &    2.10       &      &    \noOrder       &      \\ \hline
\end{tabular}
}
{
\begin{table}[hbt]\tableFont 
\begin{center}
  \pistonTableI \\
\bigskip
  \pistonTableII \\
\bigskip
  \pistonTableIII
\caption{Piston motion. Maximum errors and estimated convergence rates at $t=1$ computed using the~\ampRB~scheme
   for a heavy, $\rhob=10$, medium, $\rhob=1$, and very light, $\rhob=0.001$, moving piston. 
    The column labeled "r" provides the ratio of the errors at the current grid spacing to that on the next coarser grid.
 }
\label{table:piston}
\end{center}
\end{table}
}
} 

Figure~\ref{fig:pistonCurves} shows a time history of the position, velocity and acceleration of the
rigid body, together with the errors, as computed with the~\ampRB~scheme using the grid $\Gcp^{(4)}$
and $\dt=10^{-2}$. Results are shown for a medium-light body with $\rhob=1$ (coloured curves) and a very
light body $\rhob=0.01$ (black curves).  The~\ampRB~scheme remains stable for both cases, and that
the solution curves for $x_b(t)$, $v_{1,b}(t)$ and $a_{1,b}(t)$ for each value of $\rhob$ are in
good agreement with the exact solution.
The observed stability behaviour of the~\ampRB~scheme agrees with the analysis in Part~I for an added-mass model problem (linearized about a fixed interface position).  For the model problem, it was found that the~\ampRB~scheme is
stable for {\rm any} ratio of the mass of the body to that of the added mass of the fluid.  An analysis of a
traditional-partitioned (TP) scheme (with no sub-time-step iterations) for the same model problem found that the
TP~scheme is stable if and only if $\mb>M_a$.  For the present piston problem, the interface position varies with time so that the added mass given in~\eqref{eq:addedMass} varies with time as well.  However, a condition
for stability of the TP~scheme can be estimated as $\mb>\max(M_a(t))=5/4$.  As a check, we find that the TP~scheme is
unstable for calculations using $\Gcp^{(4)}$ and $\dt=10^{-2}$
if $\rhob =\mb \approx 1.5$ and smaller.  This result
is in reasonable agreement with the estimate.

To assess the accuracy of the~\ampRB~scheme, errors are computed on a sequence of grids of increasing resolution. 
The time-step is chosen as $\dt_j=1/(10 j)$ for $\Gcp^{(j)}$ and the equations are integrated to $t_{{\rm final}}=1$. 
Table~\ref{table:piston} presents the results of this grid refinement study for a heavy, medium and light body. 
The tables show the max-norm errors
and estimated convergence rates computed by a least-squares fit to the 
logarithm of the error. 
The max-norm error of component $q$ for a calculation using grid $\Gcp^{(j)}$ is denoted by $E_j^{q}$. 
For vector quantities such as $\vv$, the corresponding error $E_j^{\vv}$ is the maximum taken over all components of the vector.
The results show that the scheme is second-order accurate in the fluid velocity, body velocity and body position.
The fluid pressure and rigid-body acceleration converge more rapidly for this simple problem 
since the second-order accurate approximation to the pressure equation is exact for a solution with linear variation in~$x$.

\subsection{Rotating disk in an annular fluid chamber} \label{sec:rotatingDisk}

In this section the rotating disk-in-a-disk problem described in Part~I~\cite{rbinsmp2016r} is
solved using moving overlapping grids to verify the accuracy and stability of the~\ampRB~scheme.
The geometry of the two-dimensional problem is shown in Figure~\ref{fig:diskInDiskGrid}.  An annular
fluid chamber is bounded by a rotating solid disk of radius~$r=\aa$ and a static wall at $r=\bb$.
The disk has uniform density $\rhob$ so
that its moment of inertia (for rotation in the plane) is $I_b=\rhob\pi\aa\sp4/2$.  An exact
solution for the case when the rigid body rotates about its centre is given in Part~I.  The
solution for the fluid velocity and pressure is given in terms of the polar coordinates $(r,\theta)$
by
\[
\vv(r,\theta,t) = v_\theta(r,t) \begin{bmatrix} -\sin\theta \\ \cos\theta \end{bmatrix},\qquad
p(r,t) = \rho\int_{\aa}^{r}  \frac{[v_\theta(s,t)]^2}{s} \, ds + p_0,
\]
where $p_0$ is a constant and
\[
v_\theta(r,t) =  \Amplitude \left[
           \frac{J_1(\lambda r) Y_1(\lambda \bb) - J_1(\lambda \bb) Y_1(\lambda r)}
                {J_1(\lambda \aa) Y_1(\lambda \bb) - J_1(\lambda \bb) Y_1(\lambda \aa)}\right]
           \,e^{-\lambda^2 \nu t},\]
is the circumferential component of the velocity.  Here, $J_1$ and $Y_1$ are Bessel functions of order one, $\nu=\mu/\rho$ is the kinematic viscosity of the fluid, $\Amplitude$ is the velocity at $r=\aa$ and $t=0$, and $\lambda$ is a constant determined by the matching conditions at the interface with the rigid body (see Part~I).  The angular velocity of the rigid body is given by
\[
\omega_b(t)={v_\theta(\aa,t)\over \aa}={\Amplitude e^{-\lambda^2 \nu t}\over\aa},
\]
where $\omega_b$ corresponds to the third component of the angular velocity vector $\omegavb$.
Note that the solution exhibits viscous shear at the surface of the rigid body so that this problem provides a good test of the stability of the~\ampRB~scheme when added-damping effects are important.  Added-mass effects are negligible for this problem since the fluid pressure is uniform along the surface of the rigid body.

{
\newcommand{\drawContour}[7]{%
\begin{scope}[#1]
\draw(0.0,0) node[anchor=south west,xshift=-4pt,yshift=+0pt] {\trimfiga{fig/#2}{\figWidtha}};
  \draw(.5,.5) node[draw,fill=white,anchor=west,xshift=2pt,yshift=1pt] {\scriptsize #3};
\begin{scope}[xshift=15pt,yshift=-3pt]
  \draw (\xcb,\ycb) node[anchor=south west,xshift=0.cm,yshift=.5cm,rotate=-90] {\trimfigcb{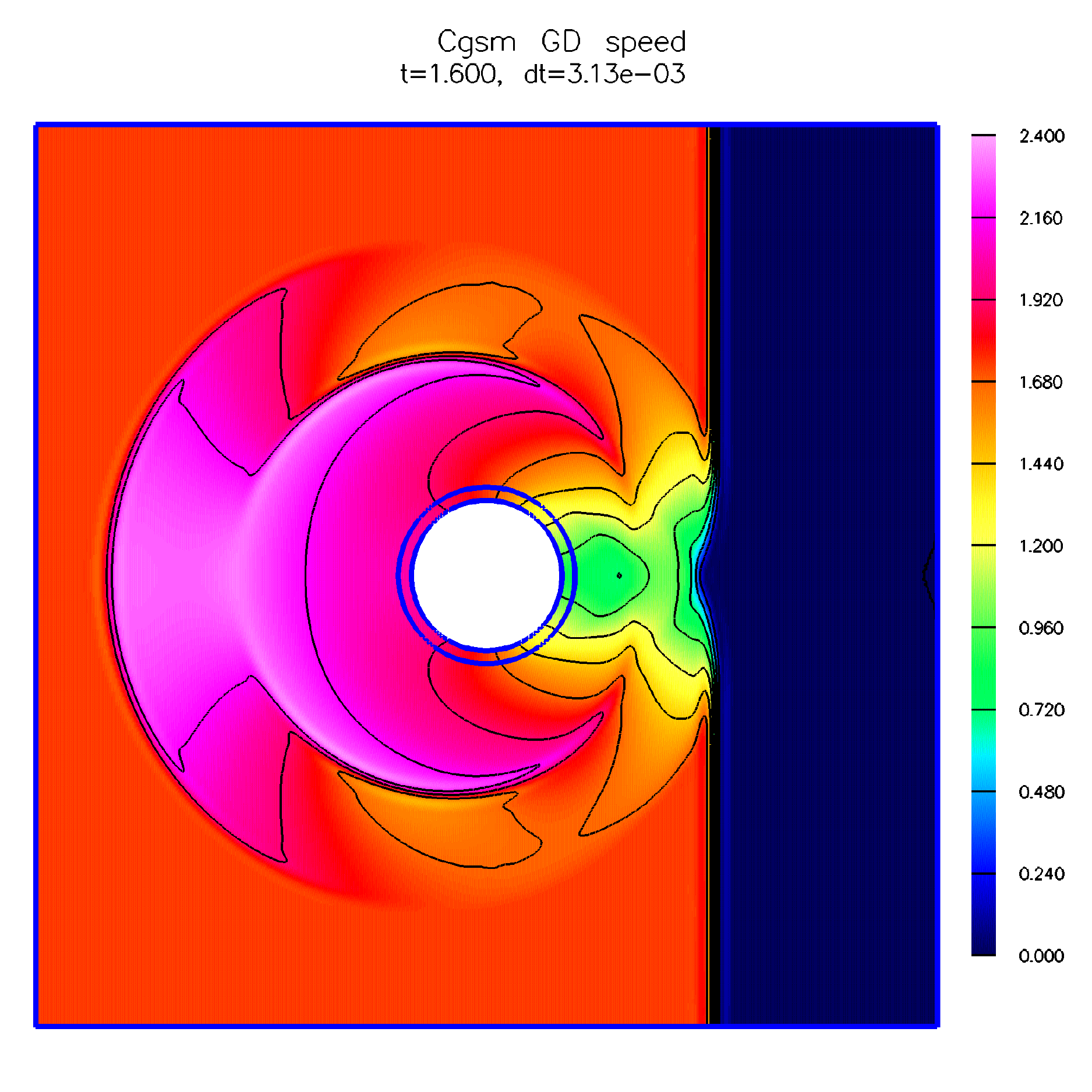}{\cbWidth}{\cbHeight}};
  \draw (.8,0) node[anchor=north,xshift=+3pt,yshift=+2pt] {\scriptsize $#6$};
  \draw (4.4,0) node[anchor=north,xshift=+0pt,yshift=+2pt] {\scriptsize $#7$};
\end{scope}
\end{scope}
}
\newcommand{\cbWidth}{.2cm}
\newcommand{\cbHeight}{4cm}
\newcommand{\xcb}{.5cm}
\newcommand{\ycb}{-.2cm}
\setlength{\ycbTop}{\ycb+\cbHeight}
\setlength{\ycbMid}{\ycb+\cbHeight*\real{.5}}
\newcommand{\trimfigcb}[3]{\includegraphics[width=#2, height=#3, clip, trim=17cm 2.35cm 1.65cm 2.35cm]{#1}}
%
%
\def\rad{1.27}
\newcommand{\plotDisk}{
\fill[fill=red!20,draw=red,line width=2pt] 
      plot[samples=100, domain=0.:360] ( {\rad*cos(\x)} , {\rad*sin(\x)} ) -- cycle ;
}
\newcommand{\figWidth}{5.5cm}
\newcommand{\trimfig}[2]{\trimh{#1}{#2}{.02}{.02}{.02}{.02}}
\newcommand{\figWidtha}{5.5cm}
\newcommand{\trimfiga}[2]{\trimh{#1}{#2}{.015}{.14}{.07}{.09}}
\begin{figure}[htb]
\begin{center}
\resizebox{14cm}{!}{
\begin{tikzpicture}[scale=1]
  \useasboundingbox (0.0,0.4) rectangle (16.4,5.8);  
  \draw(-.2,0) node[anchor=south west,xshift=-4pt,yshift=+0pt] {\trimfig{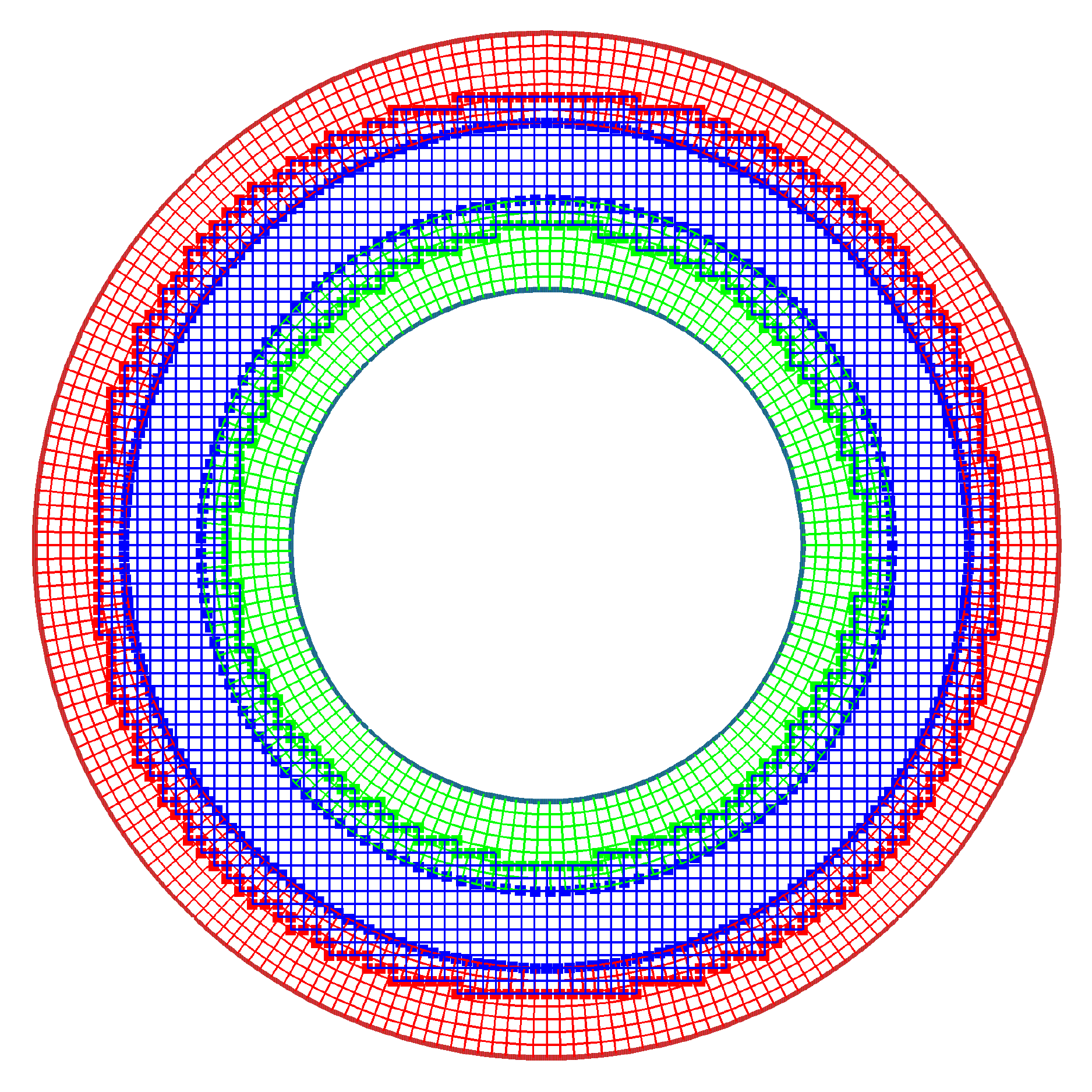}{\figWidth}};
 \drawContour{xshift=5.6cm,yshift=0cm}{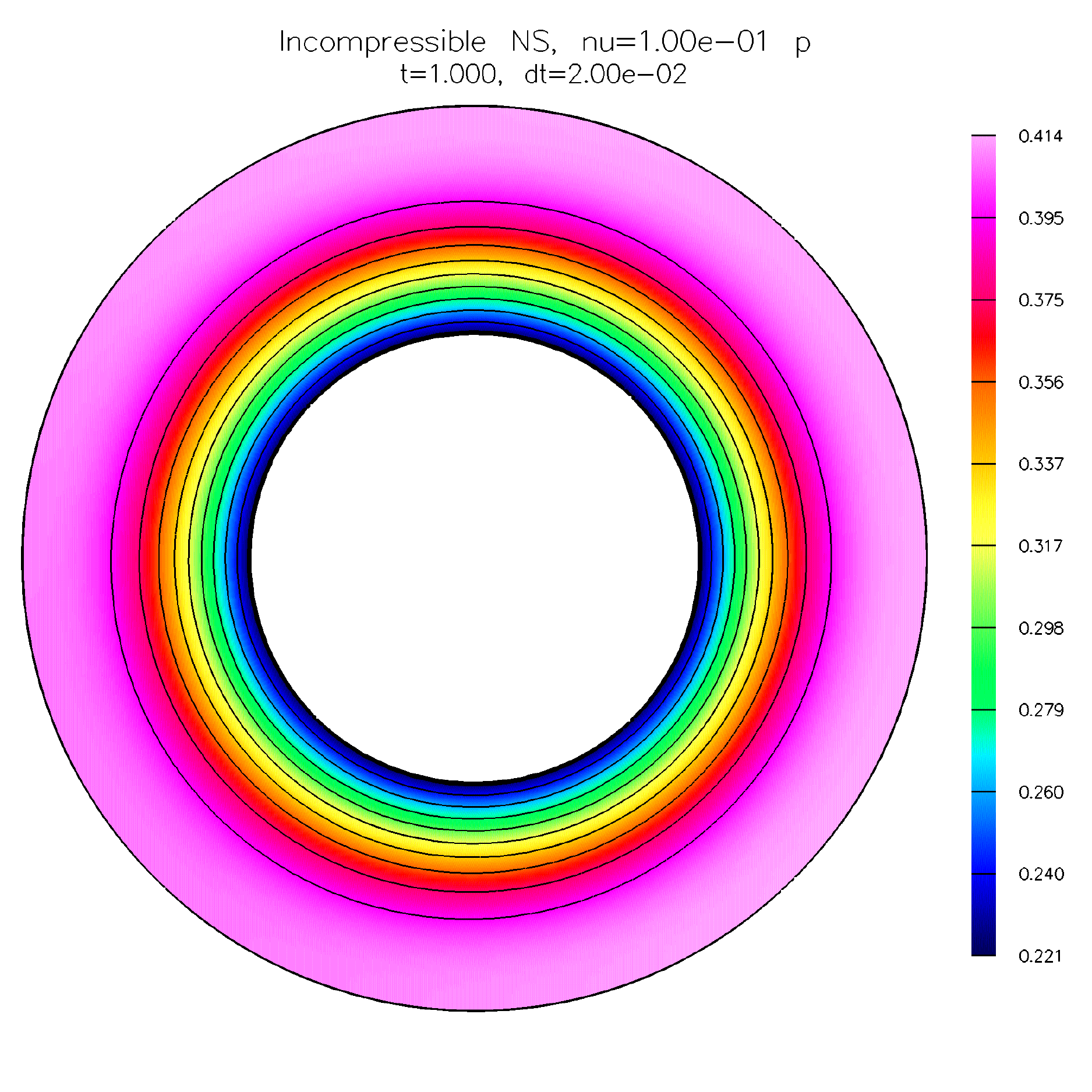}{$p$}{$p$}{$t=2.0$}{.221}{.414};
 \drawContour{xshift=11.2cm,yshift=0}{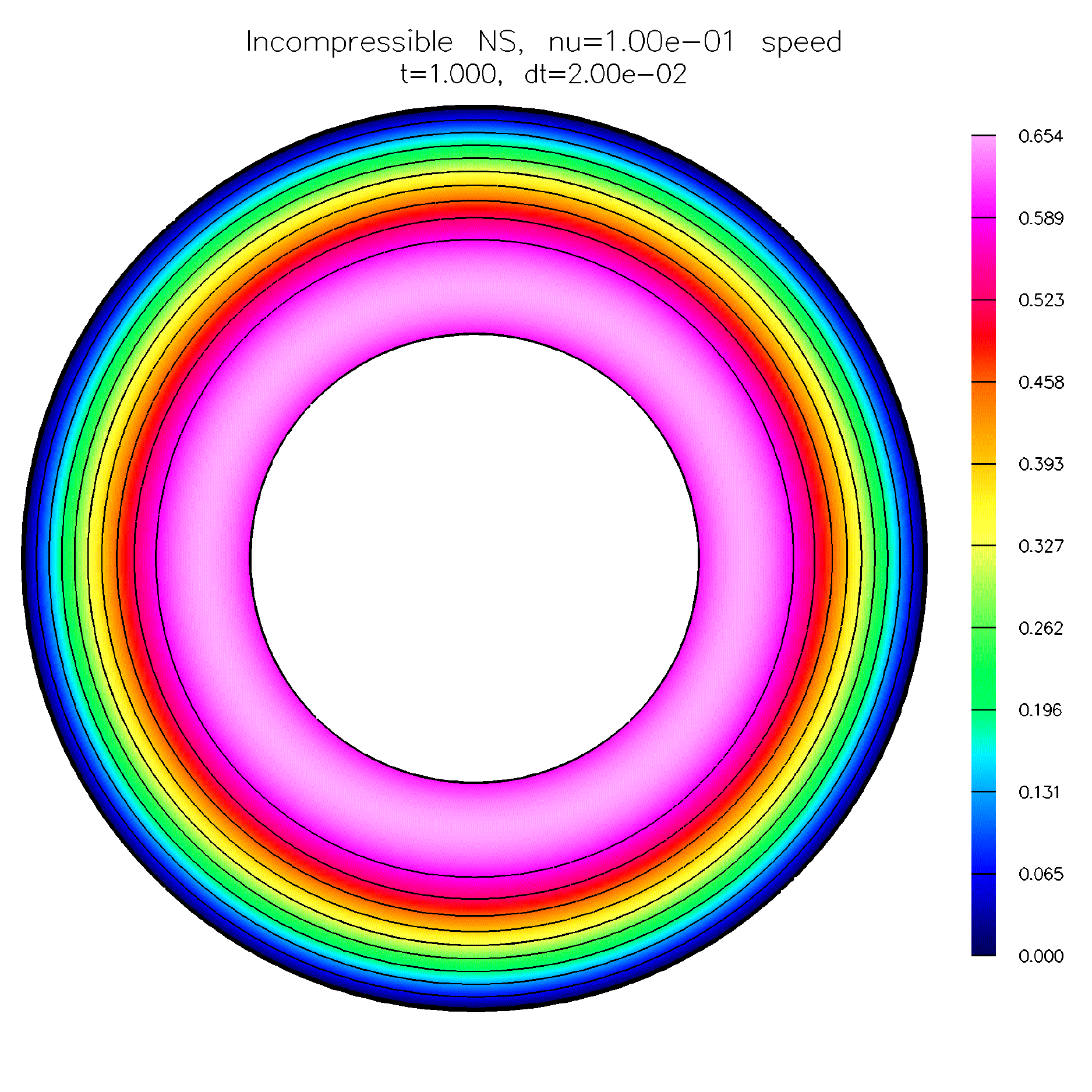}{$\vert\vv\vert$}{$p$}{$t=4.0$}{0.0}{.654};
%
  \begin{scope}[xshift=8.76cm,yshift=2.91cm]
    \plotDisk
    \draw[very thick,->,red] (\rad*.5-.1,0) arc (0:320:.5cm); %
  \end{scope}
\end{tikzpicture}
}
\end{center}
  \caption{Left: Coarse composite-grid $\Gcrd^{(2)}$ for a solid disk immersed in a disk of an incompressible fluid. 
    The green fluid grid adjacent to the
    rigid disk (solid red) moves in time. The red and blue fluid grids remain fixed.
    Computed pressure (middle) and speed (right) at $t=1.0$ on grid $\Gcrd^{(4)}$. 
   }
  \label{fig:diskInDiskGrid}
\end{figure}
}

Numerical solutions are computed using a composite grid for the annular fluid domain, denoted by
$\Gcrd^{(j)}$, as shown in~Figure~\ref{fig:diskInDiskGrid}.  The composite grid consists of two
annular boundary-fitted component grids and one background Cartesian grid (blue).  The inner annular
grid (green) is attached to, and moves with, the surface of the rigid body at $\aa=1$, while the
outer annular grid (red) is fixed to the static outer boundary of the fluid domain at $\bb=2$.  Both
annular grids have radial extent equal to $0.3$, and the representative grid spacing is $\hj=1/(10
j)$ for all component grids in each of their coordinate directions.  The initial conditions are
taken from the exact solution (with $p_0=0$), no-slip boundary conditions are applied at $r=\bb$,
and interface matching conditions are applied at the surface of the rigid body, $r=\aa$.  The fluid
density and viscosity are taken to be $\rho=1$ and $\mu=0.1$ for all calculations in this
section. Note that the exact solution contains only rotations of the body, while the numerical
solution allows both rotational and translational motions of the body.

Figure~\ref{fig:diskInDiskGrid} shows contours of the pressure, $p$, and speed, $|\vv|$, at $t=1.0$
for the case of a light body with $\rhob=0.01$ and for $\Amplitude=1$.  The solution is computed
using the~\ampRB~scheme with velocity correction on the composite grid $\Gcrd^{(4)}$ with $\adc=1$
and $\dt$ determined by the advection time-step restriction. 
The pressure is seen to depend only on the radius~$r$ and the motion of the body is
dominated by rotation, in agreement with the exact solution.  Even though the body is light there
are no instabilities due to added-damping effects, as predicted by the
stability analysis in Part~I.

{
\newcommand{\figWidth}{7.cm}
\newcommand{\trimfig}[2]{\trimFig{#1}{#2}{.0}{.0}{.0}{.0}}
\begin{figure}[htb]
\begin{center}
\resizebox{12cm}{!}{
\begin{tikzpicture}[scale=1]
  \useasboundingbox (0.25,1) rectangle (15.,5.8);  
  \draw(0,0) node[anchor=south west,xshift=-4pt,yshift=+0pt] {\trimfig{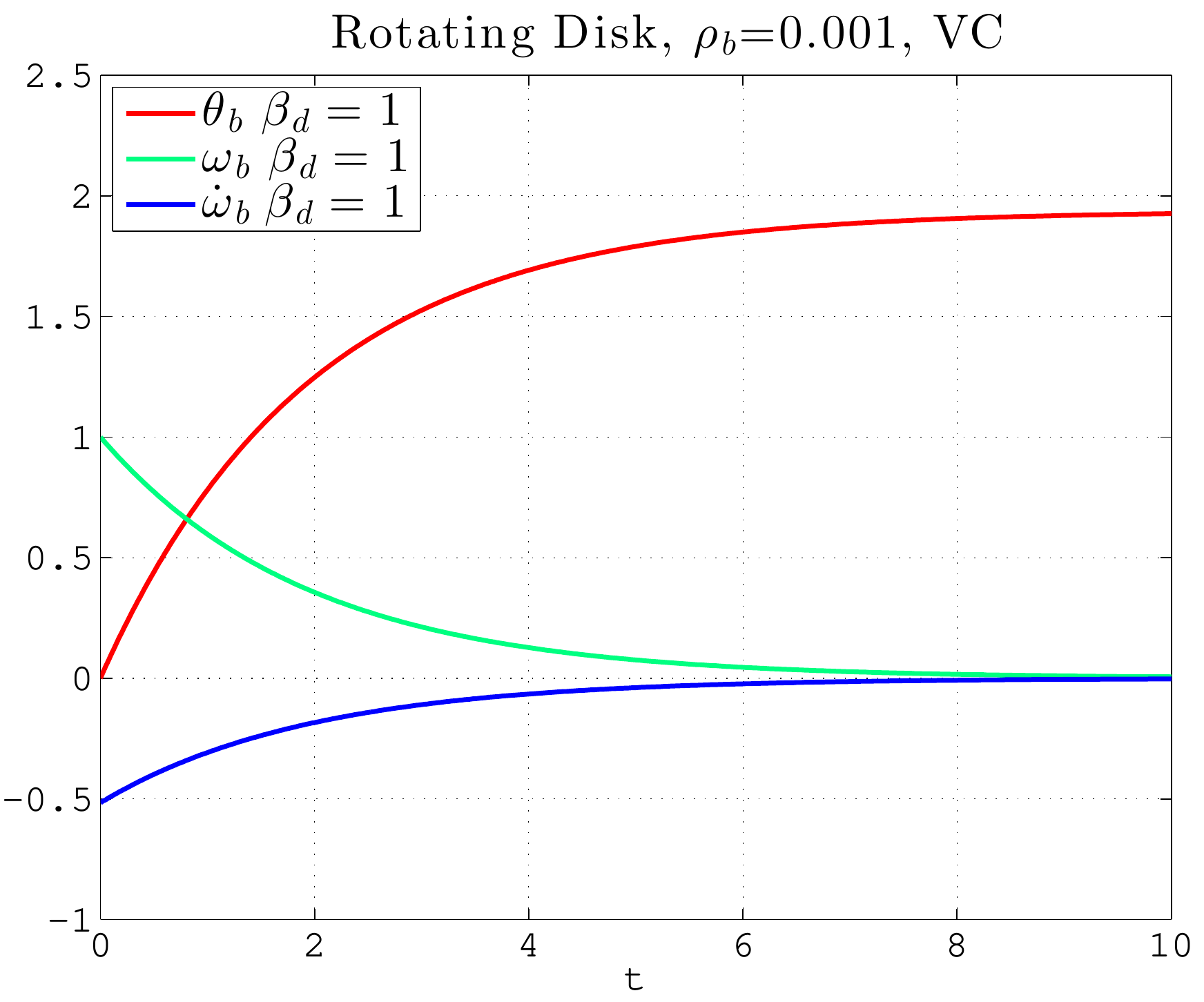}{\figWidth}};
  \draw(7.5,0) node[anchor=south west,xshift=-4pt,yshift=+0pt] {\trimfig{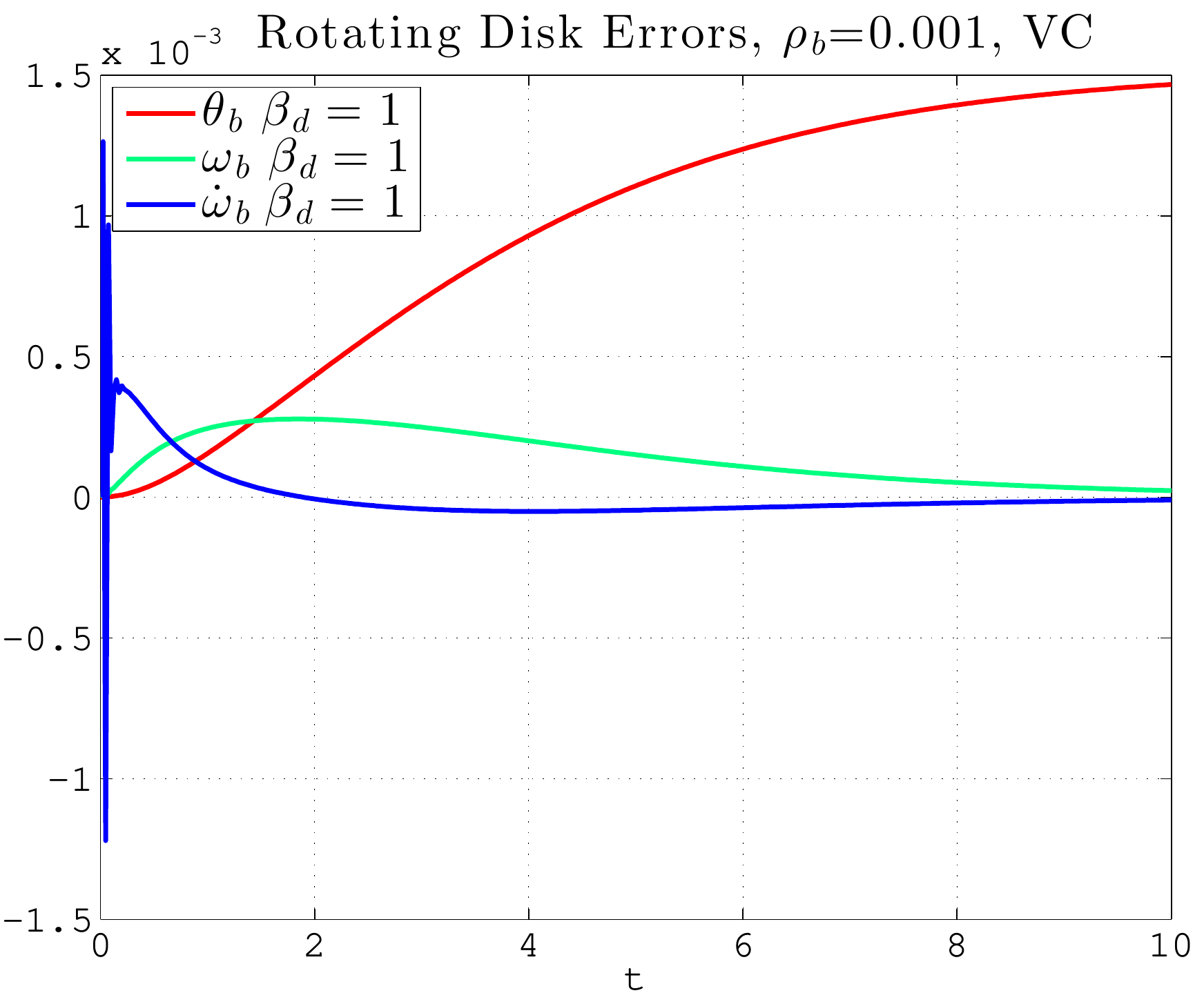}{\figWidth}};
\end{tikzpicture}
}
\end{center}
  \caption{Rotating disk. Rigid-body rotation-angle, angular velocity and angular acceleration for the rotation of a 
   very light rigid disk, $\rhob=0.001$, using grid $\Gcrd^{(4)}$ and $\dt=.05$. Left solutions. Right: errors.
    }
  \label{fig:shearDiskMotionCurves}
\end{figure}
}

Figure~\ref{fig:shearDiskMotionCurves} shows the time history of the rotation angle $\theta_b$,
velocity $\omega_b$, and acceleration $\dot\omega_b$, of the rigid body, from a calculation of very
light rigid-body with $\rhob=10^{-3}$. A plot of the errors in these quantities is also shown.  The
solution is computed using the~\ampRB~scheme with velocity correction on the grid $\Gcrd^{(4)}$
with~$\dt=.05$ and $\adc=1$.  
The inner disk, with
angular velocity equal to~$\Amplitude=1$ initially, rotates counter-clockwise and slows 
due to the viscous shear stress from the fluid.  
Even for these light bodies, the numerical
  solution shows no signs of instabilities.  However, the error in angular acceleration
  $\dot\omega_b$ is seen to possess a few small start-up oscillations near $t=0$, which are expected
  for these difficult simulations with small $\rhob$.

{ 
\newcommand{\errp}{$E_j^{p}$}
\newcommand{\erruv}{$E_j^{\vv}$}
\newcommand{\errwC}{$E_j^{\omega_b}$}
\newcommand{\errwtC}{$E_j^{\dot\omega_b}$}

\newcommand{\rotatingDiskTableI}{%
\begin{tabular}{|c|c|c|c|c|c|c|c|c|} \hline 
  \multicolumn{9}{|c|}{Rotating disk, $\rhob=10$} \\ \hline 
\strutt  $\hj$  &     \errp     &  r   &     \erruv     &  r   &     \errwC     &  r   &     \errwtC     &  r    \\[3pt] \hline 
  1/10   & \num{6.9}{-5} &      & \num{1.2}{-4} &      & \num{3.0}{-5} &      & \num{2.9}{-5} &     \\ \hline
  1/20   & \num{1.6}{-5} & 4.3  & \num{3.0}{-5} & 4.1  & \num{8.5}{-6} & 3.6  & \num{7.9}{-6} & 3.6 \\ \hline
  1/40   & \num{4.0}{-6} & 4.0  & \num{7.4}{-6} & 4.0  & \num{2.8}{-6} & 3.0  & \num{2.6}{-6} & 3.1 \\ \hline
  1/80   & \num{8.9}{-7} & 4.5  & \num{1.9}{-6} & 4.0  & \num{7.8}{-7} & 3.6  & \num{7.1}{-7} & 3.7 \\ \hline
  rate   &    2.08       &      &    2.01       &      &    1.74       &      &    1.76       &      \\ \hline
\end{tabular}
}
\newcommand{\rotatingDiskTableII}{%
\begin{tabular}{|c|c|c|c|c|c|c|c|c|} \hline 
  \multicolumn{9}{|c|}{Rotating disk, $\rhob=1$} \\ \hline 
\strutt  $\hj$  &     \errp     &  r   &     \erruv     &  r   &     \errwC     &  r   &     \errwtC     &  r    \\[3pt] \hline 
  1/10   & \num{1.2}{-5} &      & \num{4.8}{-5} &      & \num{5.1}{-6} &      & \num{1.8}{-5} &     \\ \hline
  1/20   & \num{2.5}{-6} & 4.7  & \num{9.5}{-6} & 5.1  & \num{2.6}{-6} & 1.9  & \num{5.0}{-6} & 3.5 \\ \hline
  1/40   & \num{5.9}{-7} & 4.2  & \num{2.5}{-6} & 3.8  & \num{6.2}{-7} & 4.2  & \num{1.3}{-6} & 3.8 \\ \hline
  1/80   & \num{1.6}{-7} & 3.6  & \num{6.2}{-7} & 4.0  & \num{1.7}{-7} & 3.6  & \num{3.4}{-7} & 3.8 \\ \hline
  rate   &    2.06       &      &    2.08       &      &    1.67       &      &    1.89       &      \\ \hline
\end{tabular}
}
\newcommand{\rotatingDiskTableIV}{%
\begin{tabular}{|c|c|c|c|c|c|c|c|c|} \hline 
  \multicolumn{9}{|c|}{Rotating disk, $\rhob=0$} \\ \hline 
\strutt  $\hj$  &     \errp     &  r   &     \erruv     &  r   &     \errwC     &  r   &     \errwtC     &  r    \\[3pt] \hline 
  1/10   & \num{3.0}{-5} &      & \num{1.0}{-4} &      & \num{5.5}{-5} &      & \num{4.5}{-5} &     \\ \hline 
  1/20   & \num{6.2}{-6} & 4.8  & \num{2.5}{-5} & 4.1  & \num{1.7}{-5} & 3.3  & \num{1.1}{-5} & 4.1 \\ \hline 
  1/40   & \num{1.5}{-6} & 4.1  & \num{6.4}{-6} & 3.9  & \num{5.0}{-6} & 3.4  & \num{2.9}{-6} & 3.8 \\ \hline 
  1/80   & \num{3.4}{-7} & 4.4  & \num{1.6}{-6} & 3.9  & \num{1.4}{-6} & 3.7  & \num{7.4}{-7} & 3.9 \\ \hline 
  rate   &    2.13       &      &    2.00       &      &    1.78       &      &    1.97       &      \\ \hline
\end{tabular}
}
{
\begin{table}[hbt]\tableFont 
\begin{center}
  \rotatingDiskTableI \\
\medskip
  \rotatingDiskTableII \\
\medskip
  \rotatingDiskTableIV
\caption{Rotating disk. Maximum errors and estimated convergence rates at $t=1$ computed using the~\ampRB~scheme
   with velocity correction for a heavy ($\rhob=10$), medium ($\rhob=1$) and zero mass ($\rhob=0$) rotating disk.
  }
\label{table:rotatingDisk}
\end{center}
\end{table}
}
} 

Results of a grid convergence study are given in Table~\ref{table:rotatingDisk}. 
The max-norm errors at $t=1$ are computed by comparing numerical solutions to the
exact solution for the cases of $\rhob=10$, $1$ and zero, all using $\Amplitude=0.1$.
Convergence rates are also given in the tables and these are estimated using a
least-squares fit to the log of the  errors.  The~\ampRB~scheme with velocity correction and $\adc=1$ is
used for all of the results.
The errors are seen to converge at rates close to second-order accuracy for all cases, including
the case of a zero-mass rigid body where added-damping effects are strongest.

\subsection{Zero mass disk in a counter-flow} \label{sec:cylDrop}
Consider a buoyant and very light rigid disk that is allowed to move within a rectangular fluid domain
as shown in Figure~\ref{fig:cylDropGrid}.
The flow in the channel and the motion of the solid disk are initiated from
rest by smoothly turning on a pressure gradient (by imposing a time-dependent pressure at inflow and fixing the
pressure at outflow) and smoothly turning on a gravitational body force.

{
\newcommand{\drawContour}[7]{%
\begin{scope}[#1]
\draw(0.0,0) node[anchor=south west,xshift=-4pt,yshift=+0pt] {\trimfig{fig/#2}{\figWidth}};
  \draw(1.5,1.5) node[draw,fill=white,anchor=west,xshift=4pt,yshift=0pt] {\scriptsize #3};
  \draw(2.2,1.5) node[draw,fill=white,anchor=west,xshift=2pt,yshift=0pt] {\scriptsize #5};
\begin{scope}[xshift=28pt,yshift=18pt]
  \draw (\xcb,\ycb) node[anchor=south west,xshift=0.cm,yshift=.5cm,rotate=-90] {\trimfigcb{fig/colourBarLines}{\cbWidth}{\cbHeight}};
  \draw (.8,0) node[anchor=north,xshift=+3pt,yshift=+2pt] {\scriptsize $#6$};
  \draw (4.4,0) node[anchor=north,xshift=+0pt,yshift=+2pt] {\scriptsize $#7$};
\end{scope}
\end{scope}
}
\newcommand{\cbWidth}{.2cm}
\newcommand{\cbHeight}{4cm}
\newcommand{\xcb}{.5cm}
\newcommand{\ycb}{-.2cm}
\setlength{\ycbTop}{\ycb+\cbHeight}
\setlength{\ycbMid}{\ycb+\cbHeight*\real{.5}}
\newcommand{\trimfigcb}[3]{\includegraphics[width=#2, height=#3, clip, trim=17cm 2.35cm 1.65cm 2.35cm]{#1}}
\def\rad{.5}
\newcommand{\plotDisk}{
\fill[fill=red!20,draw=red,line width=2pt] 
      plot[samples=100, domain=0.:360] ( {\rad*cos(\x)} , {\rad*sin(\x)} ) -- cycle ;
}
\newcommand{\figWidth}{9.25cm}
\newcommand{\trimfig}[2]{\trimh{#1}{#2}{.2}{.25}{.025}{.1}}
\newcommand{\trimfiga}[2]{\trimh{#1}{#2}{.27}{.27}{.025}{.1}}
\newcommand{\labelsize}{\small}
\begin{figure}[htb]
\begin{center}
\resizebox{12cm}{!}{
\begin{tikzpicture}[scale=1]
  \useasboundingbox (-.5,1) rectangle (15.5,9.75);  
  \begin{scope}[xshift=2.5cm,yshift=7.2cm]
      \fill[fill=blue!15] (-2.,-6) -- (2.,-6) -- (2,2) -- (-2,2) -- cycle; 
      \draw[->,very thick,blue] (-1,2.5) -- (-1,1.75); 
      \draw[->,very thick,blue] ( 1,2.5) -- ( 1,1.75); 
      \draw[->,thick] (-2.5,0) -- (2.5,0) node [anchor=north] {$x$}; 
      \draw[->,thick] (0,-6.5) -- (0,2.5) node [anchor=east] {$y$}; 
      \draw[-,very thick] (-2,-6) -- (-2,2);
      \draw[-,very thick] (2,-6) -- (2,2);
    \begin{scope}[xshift=1cm,yshift=0cm]  
      \plotDisk
    \end{scope}
    \fill[black] (1,.0) circle (2 pt); 
    \draw(1,-.5) node[anchor=north] {\labelsize$\frac{1}{2}$}; 
    \draw(-2,-6) node[anchor=north] {\labelsize$-1$};
    \draw(2,-6) node[anchor=north] {\labelsize$1$};
  \end{scope}
  \draw(4.7,0) node[anchor=south west,xshift=-4pt,yshift=+0pt] {\trimfig{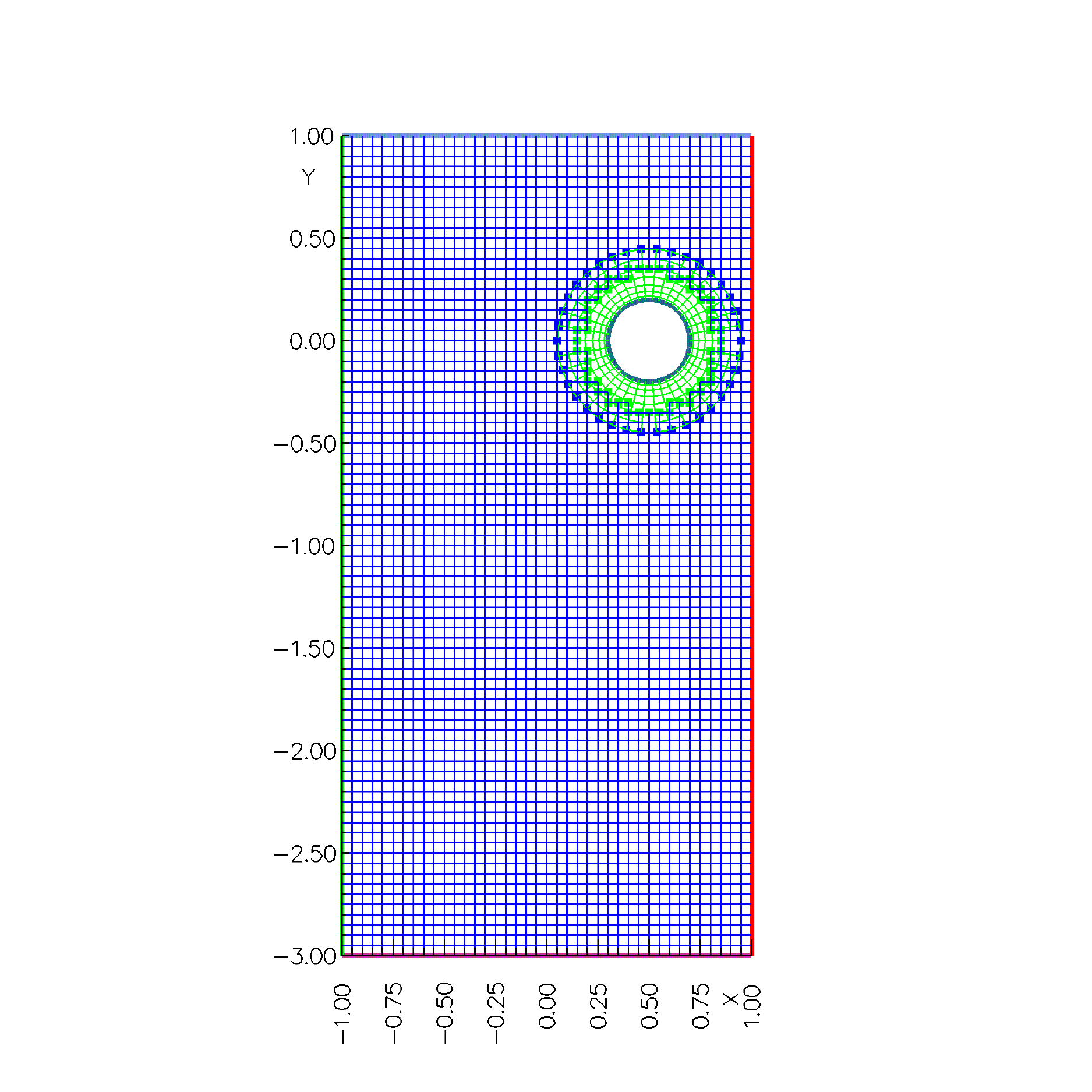}{\figWidth}};
  \drawContour{xshift=10.0cm,yshift=0cm}{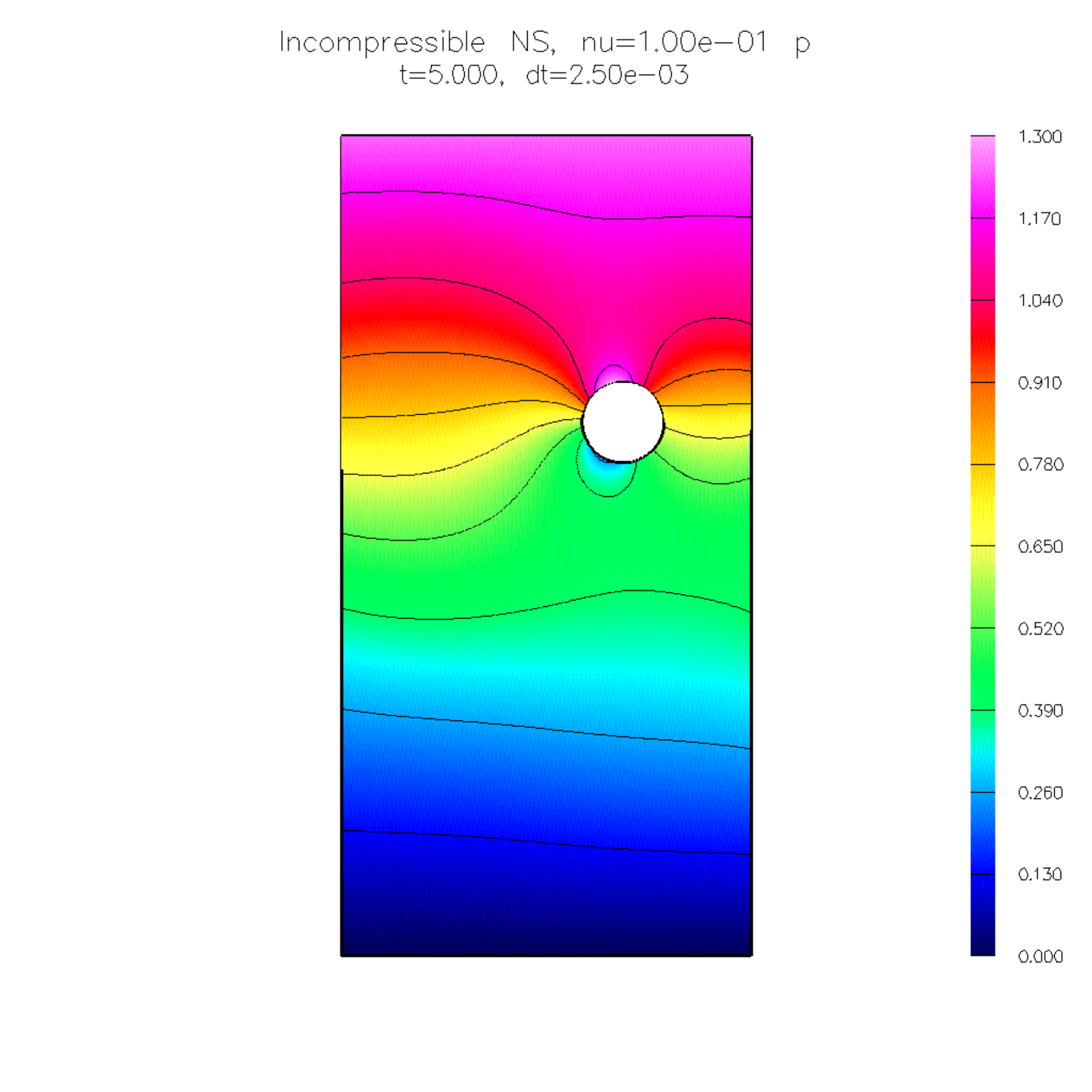}{$p$}{$p$}{$t=5$}{0}{1.3};
\end{tikzpicture}
}
\end{center}
  \caption{Rising disk in a counter-flow. Left: geometrical configuration at $t=0$.
  Middle: Composite grid $\Gcrd^{(2)}$ at $t=0$ (coarse grid). Right: contours of the pressure at $t=5$ using grid $\Gcrd^{(8)}$.
}
  \label{fig:cylDropGrid}
\end{figure}
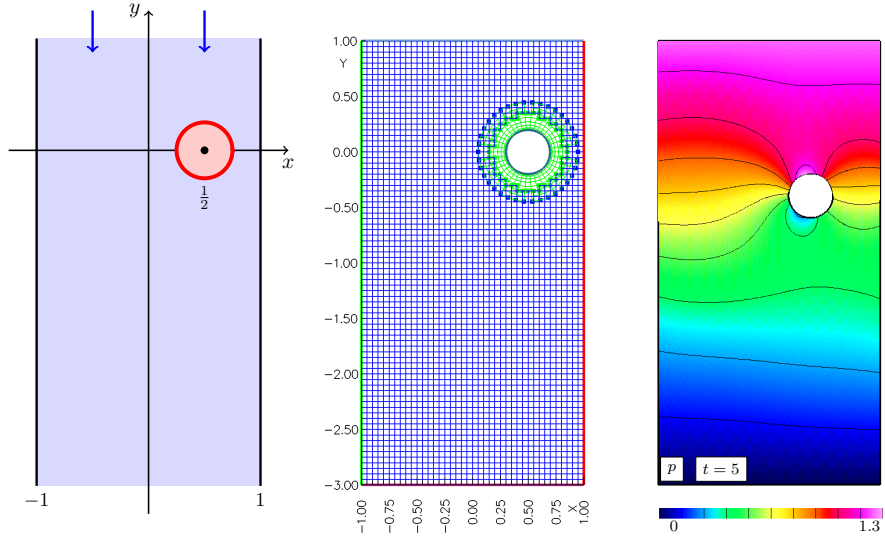
}

The geometry of the problem is shown in Figure~\ref{fig:cylDropGrid}.  The fluid channel
covers the rectangular domain $[-\channelWidth,\channelWidth]\times[\channelBottom,\channelTop]$,
where $\channelWidth=1$, $\channelBottom=-3$ and $\channelTop=1$.  
The solid disk, with radius $\diskRadius=0.2$, is initially centred at $(x,y)=(0.5,0)$. 
The fluid density and viscosity are taken to be $\rho=1$ and $\mu=1$, respectively, and the
density of the solid disk is zero, $\rhob=0$.  
At the top boundary of the channel the pressure is specified, the tangential component of the
velocity is set to zero, and the normal component of the velocity is extrapolated on the ghost points.
The specified pressure at the top ramps smoothly from zero to a value $P_{{\rm top}}=1.25$ according to 
\[
    p(x,\channelTop,t) = P_{{\rm top}} \, \rampFunction(t), 
\]
where $\rampFunction(t)$ is a {\em ramp} function given by
\begin{align}
    \rampFunction(t) =          \begin{cases}
            (35+(-84+(70- 20 t) t) t) t^4,\quad     & \text{for $0\le t \le 1$}, \\
                1,   & \text{for $t>1$}.
               \end{cases}   \label{eq:ramp}
\end{align}
The ramp function satisfies
$\rampFunction=\rampFunction\sp\prime=\rampFunction\sp{\prime\prime}=\rampFunction\sp{\prime\prime\prime}=0$
at $t=0$ and has three continuous derivatives at $t=1$.
No-slip boundary conditions are taken on the left and right sides of the channel, and on
the surface of the solid disk.  
An outflow boundary condition is taken at the bottom boundary of the channel: the
pressure and the normal derivative of the normal component of velocity are both set to zero, and
the tangential component of the velocity is extrapolated to the ghost points.
The body force due to gravity is smoothly turned on according to 
\[
\fvbe(t)=\pi\diskRadius\sp2(\rho_b-\rho) \, \gv \, \rampFunction(t),
\]
where the acceleration due to gravity is taken to be $\gv=[0,\,-4]^T$ for the present simulation.

{
\newcommand{\drawContour}[7]{%
\begin{scope}[#1]
\draw(0.0,0) node[anchor=south west,xshift=-4pt,yshift=+0pt] {\trimfig{fig/#2}{\figWidth}};
\draw(1,6.6) node[draw,fill=white,anchor=west,xshift=2pt,yshift=0pt] {\scriptsize #5};
\begin{scope}[xshift=8pt,yshift=4pt]
  \draw (\xcb,\ycb) node[anchor=south west,xshift=0.cm,yshift=.5cm,rotate=-90] {\trimfigcb{fig/colourBarLines}{\cbWidth}{\cbHeight}};
  \draw (.8,0) node[anchor=north,xshift=+3pt,yshift=+2pt] {\scriptsize $#6$};
  \draw (3.7,0) node[anchor=north,xshift=+0pt,yshift=+2pt] {\scriptsize $#7$};
\end{scope}
\end{scope}
}
\newcommand{\cbWidth}{.2cm}
\newcommand{\cbHeight}{3.3cm}
\newcommand{\xcb}{.5cm}
\newcommand{\ycb}{-.2cm}
\setlength{\ycbTop}{\ycb+\cbHeight}
\setlength{\ycbMid}{\ycb+\cbHeight*\real{.5}}
\newcommand{\trimfigcb}[3]{\includegraphics[width=#2, height=#3, clip, trim=17cm 2.35cm 1.65cm 2.35cm]{#1}}
%
\newcommand{\figWidth}{4.25cm}
\newcommand{\trimfig}[2]{\trimw{#1}{#2}{.245}{.25}{.075}{.08}}
\begin{figure}[htb]
\begin{center}
\resizebox{15cm}{!}{
\begin{tikzpicture}[scale=1]
  \useasboundingbox (1.0,.5) rectangle (16.,7.);  
%
\begin{scope}[yshift=0cm]
 \drawContour{xshift= 0.0cm,yshift=0.0cm}{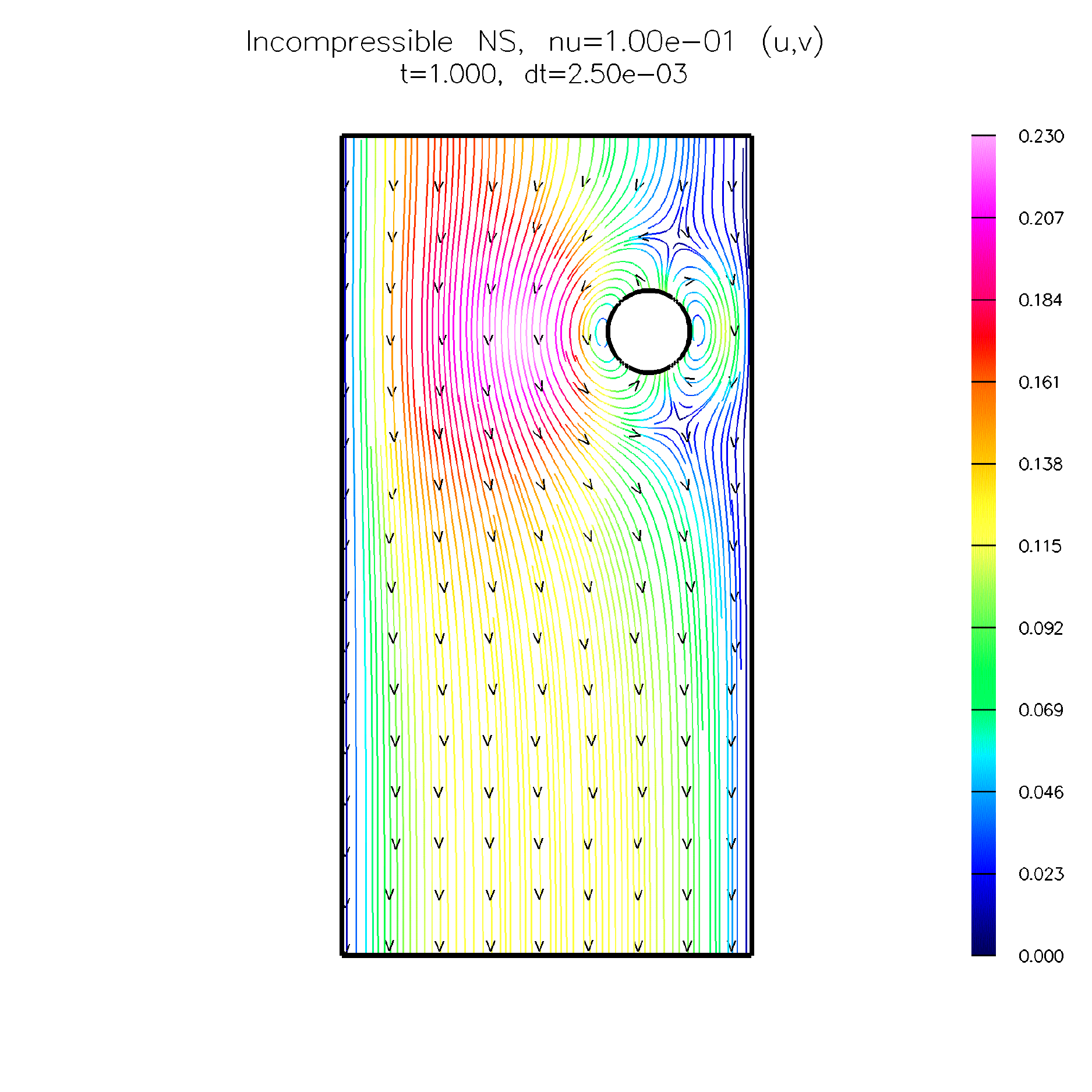}{$p$}{$p$}{$t=1$}{0.0}{0.23};
 \drawContour{xshift= 4.0cm,yshift=0.0cm}{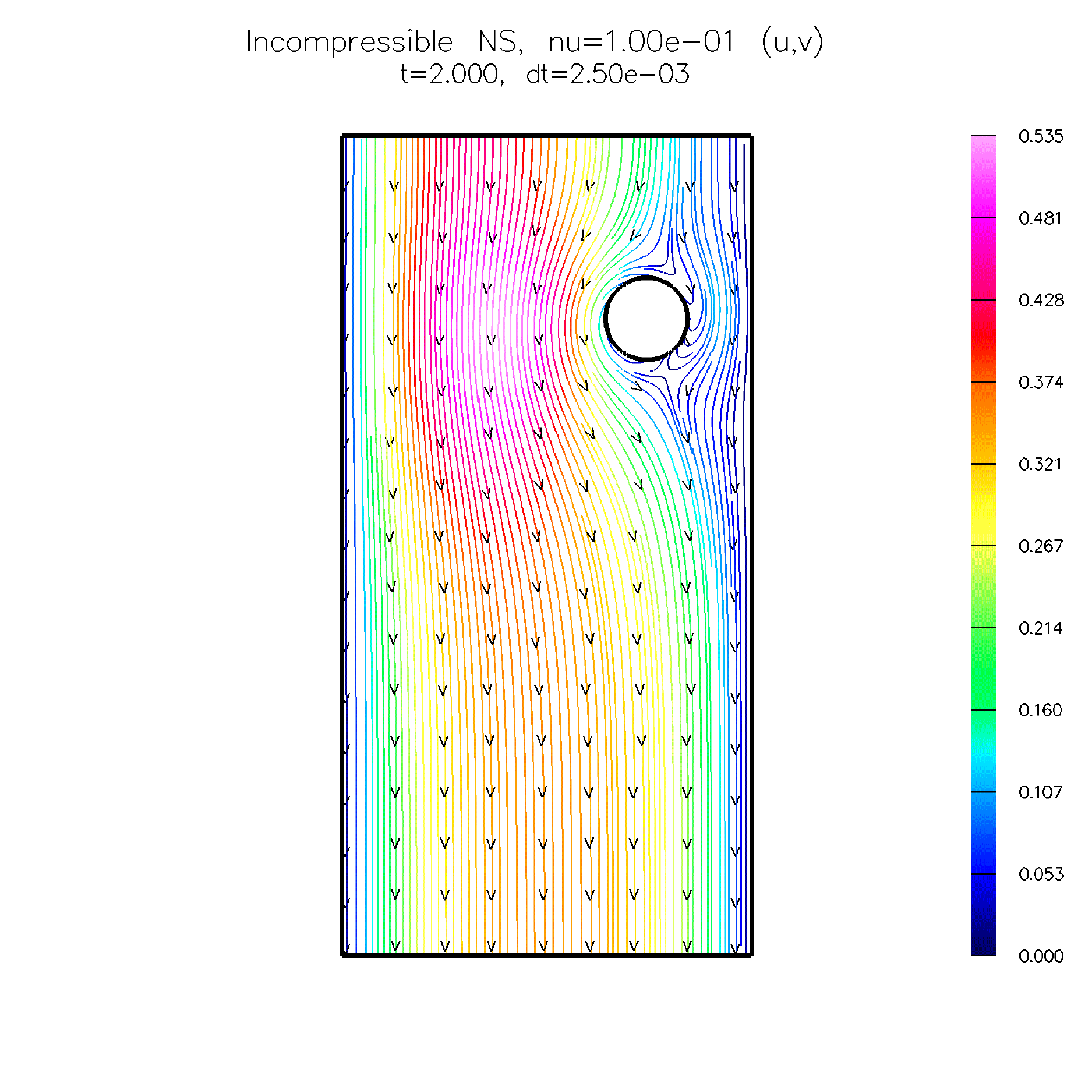}{$p$}{$p$}{$t=2$}{0.0}{0.54};
 \drawContour{xshift= 8.0cm,yshift=0.0cm}{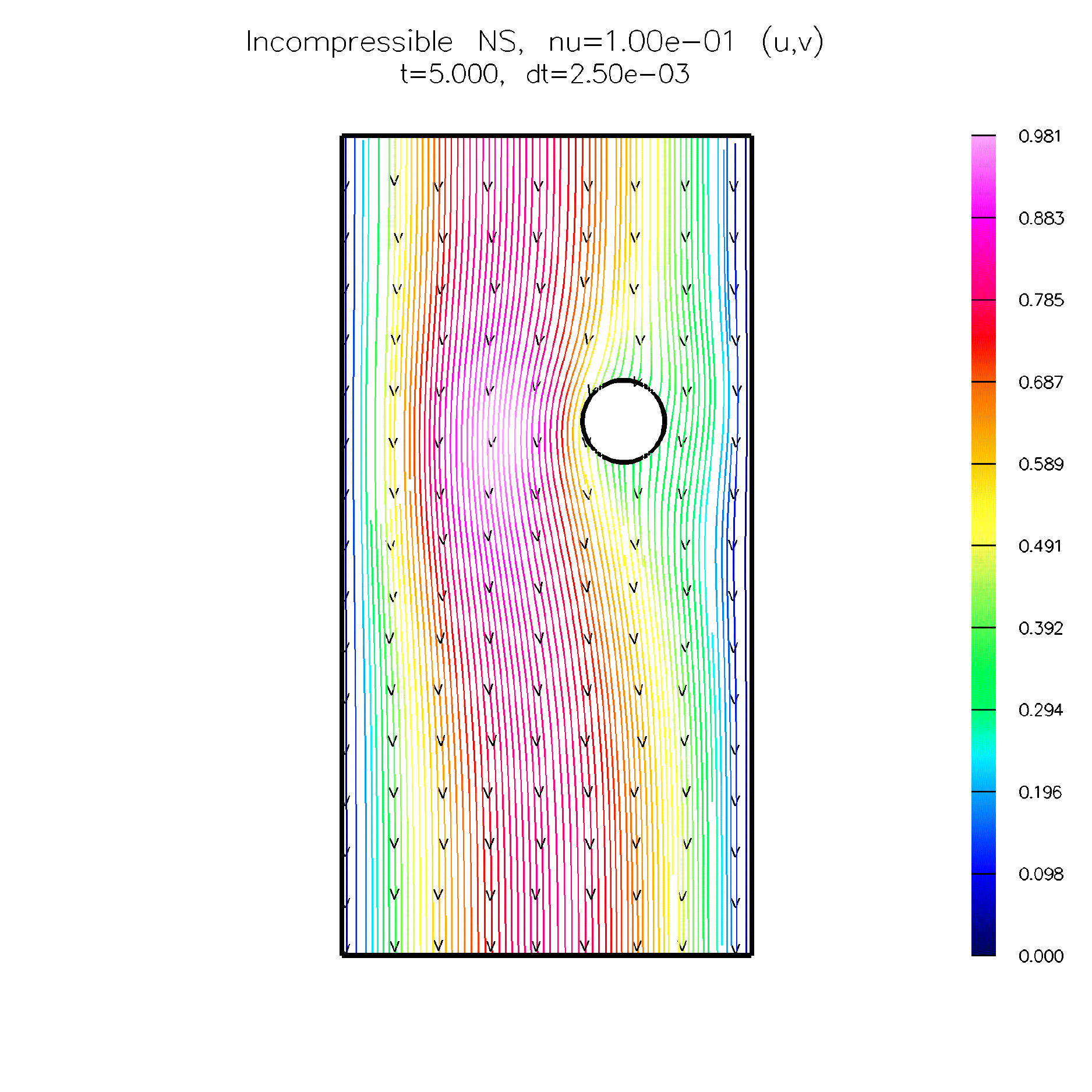}{$p$}{$p$}{$t=5$}{0.0}{0.98};
 \drawContour{xshift=12.0cm,yshift=0.0cm}{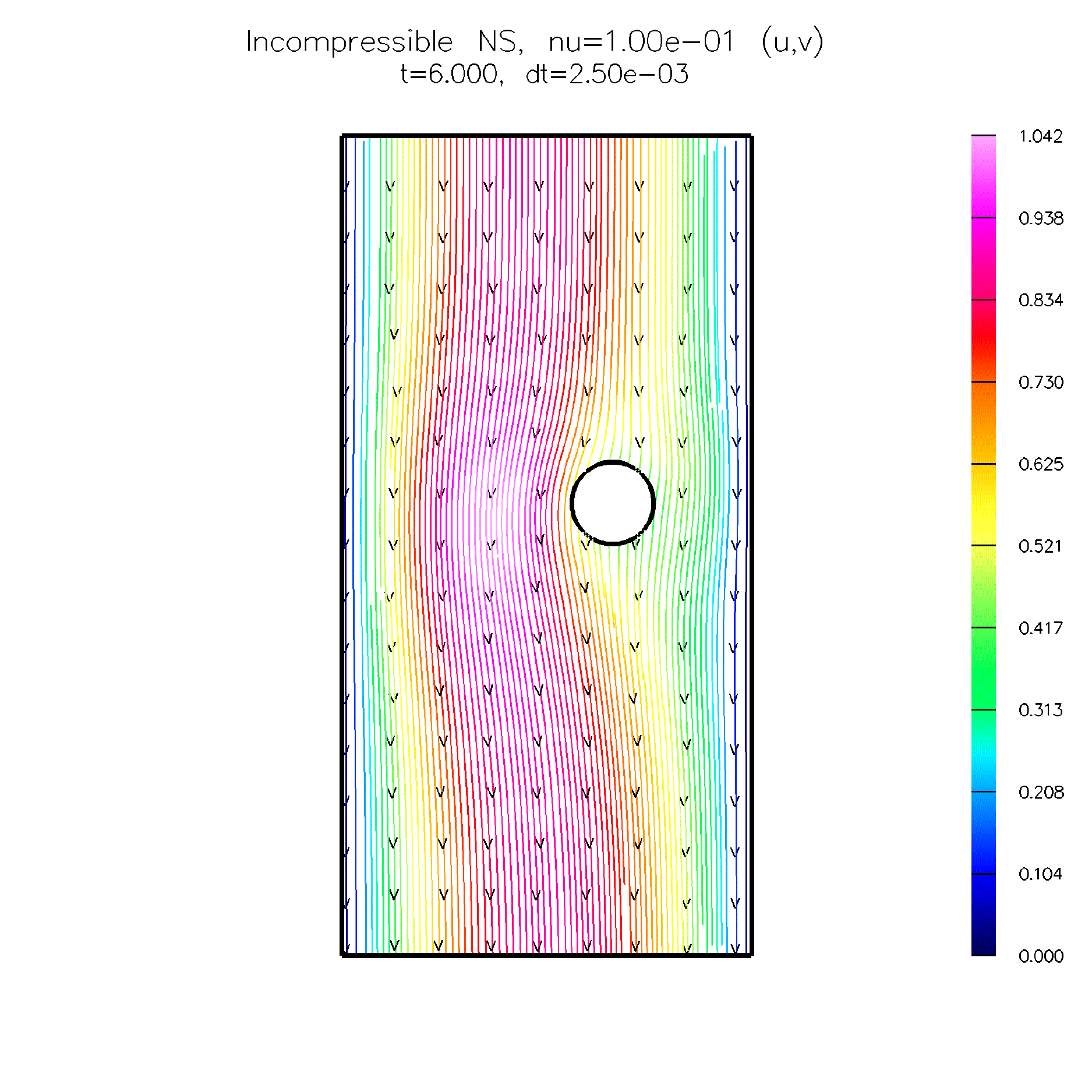}{$p$}{$p$}{$t=6$}{0.0}{1.04};
\end{scope} 
\end{tikzpicture}
}
\end{center}
  \caption{Disk in a counter-flow. Computed streamlines (coloured by the flow speed) at selected times for a moving solid-disk with $\rhob=0$, computed using the composite grid $\Gcrd^{(8)}$.  
}
  \label{fig:cylDropStreamLines}
\end{figure}
}

The composite grid for the problem at $t=0$ is shown in Figure~\ref{fig:cylDropGrid}.
The grid with resolution factor~$j$, and target grid spacing $h^{(j)}=1/(10 j)$,  is denoted by~$\Gcrd^{(j)}$, and 
consists of two component grids. A Cartesian background grid covers the entire domain.
An annular grid of fixed radial width equal to~$0.2$ surrounds the disk.
The grid spacing is stretched in the radial direction to be approximately $h_{BL}^{(j)}=h^{(j)}/3$ near the surface of the disk.
The time-step is taken as $\dt^{(j)}=0.2 \, h^{(j)}$. 

{
\newcommand{\figWidth}{8.cm}
\newcommand{\trimfig}[2]{\trimw{#1}{#2}{.0}{.0}{.0}{.0}}
\newcommand{\figWidthz}{7.5cm}
\newcommand{\trimfigz}[2]{\trimw{#1}{#2}{.0}{.0}{.0}{.0}}
\begin{figure}[htb]
\begin{center}
\resizebox{12cm}{!}{
\begin{tikzpicture}[scale=1]
  \useasboundingbox (0.,1.3) rectangle (16.,21);  
  \draw( 0.0,14.0) node[anchor=south west,xshift=-4pt,yshift=+0pt] {\trimfig{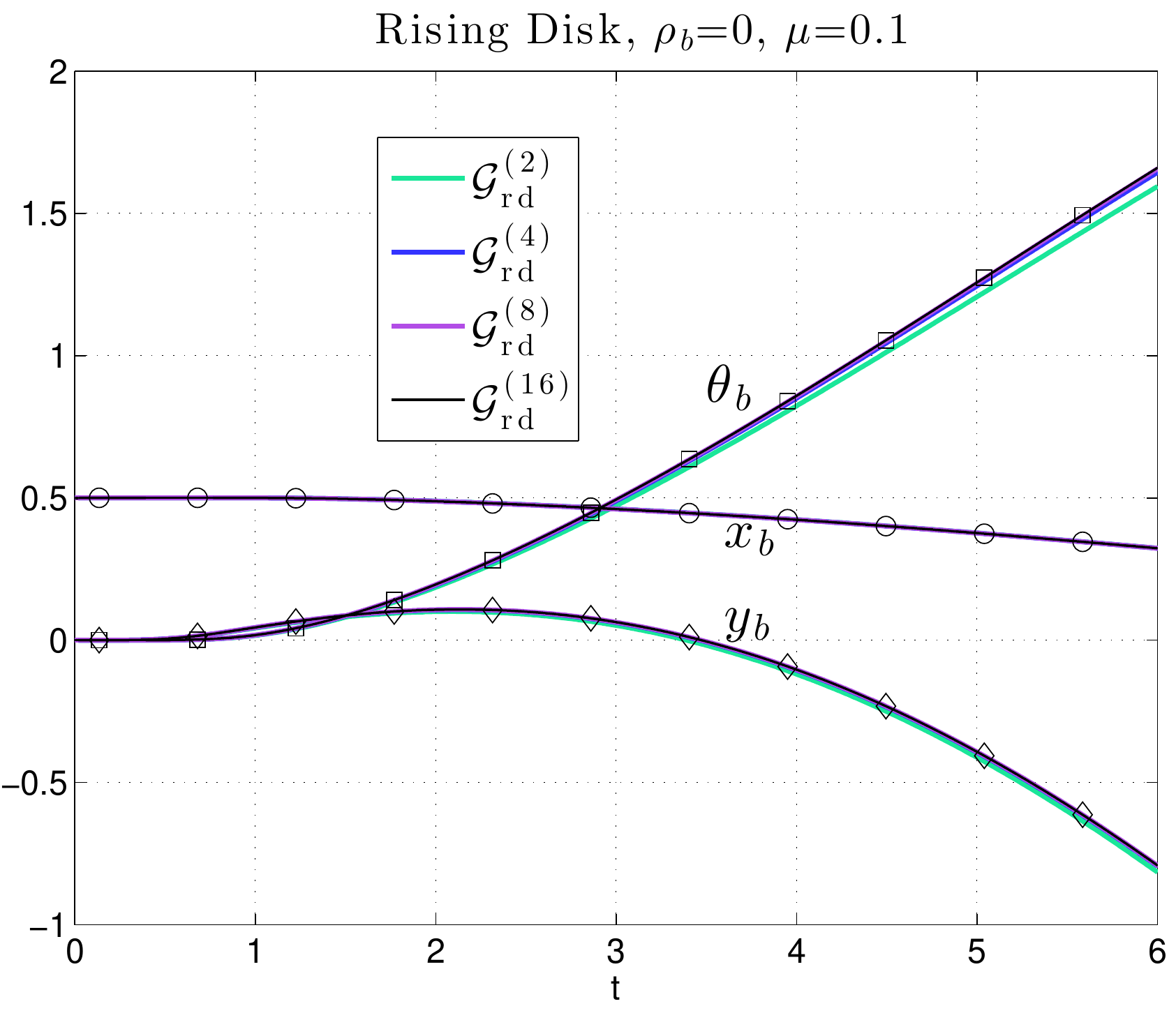}{\figWidth}};
  \draw(8.25,14.25) node[anchor=south west,xshift=-4pt,yshift=+0pt] {\trimfigz{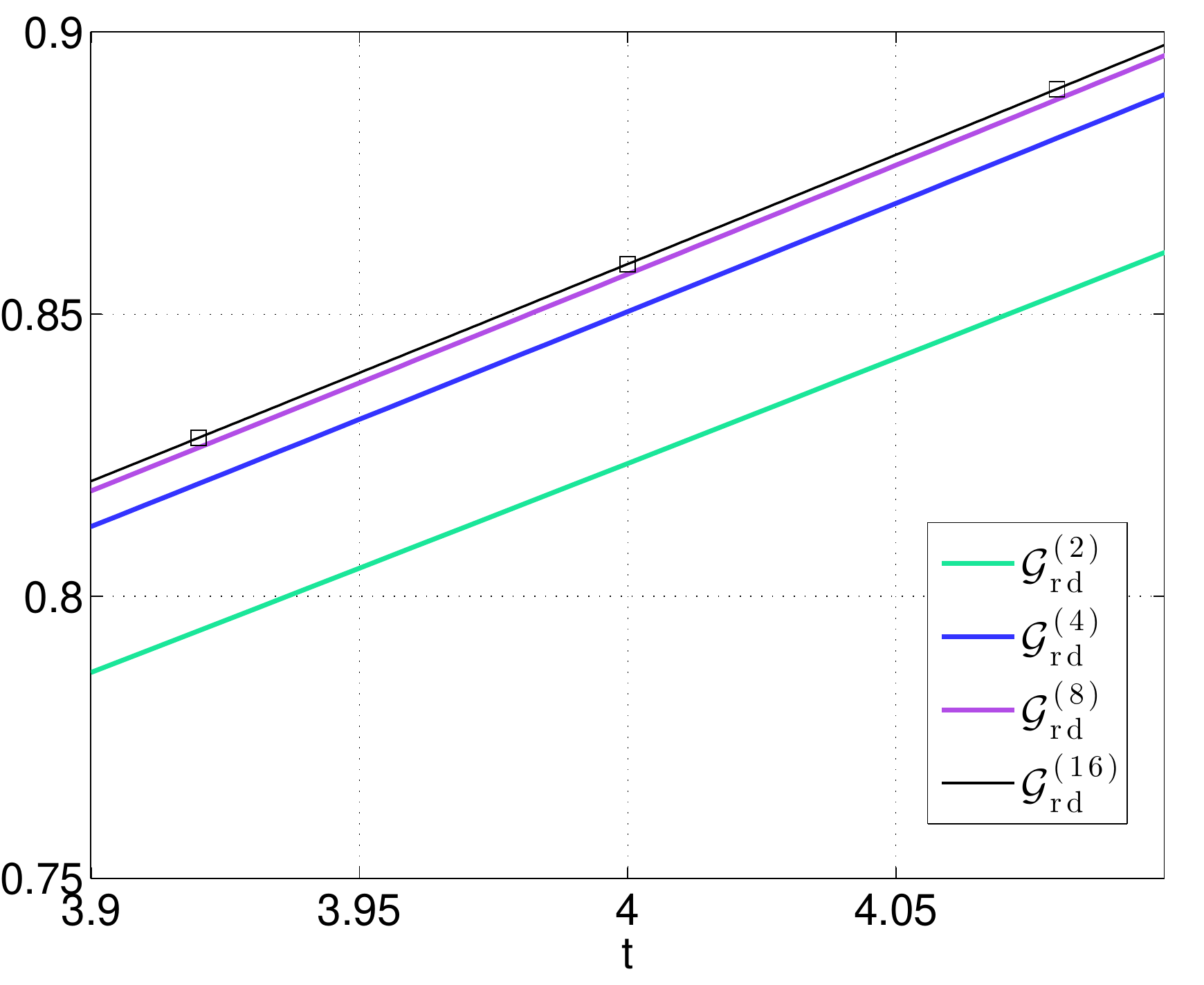}{\figWidthz}};
  \draw( 0.0,7.0) node[anchor=south west,xshift=-4pt,yshift=+0pt] {\trimfig{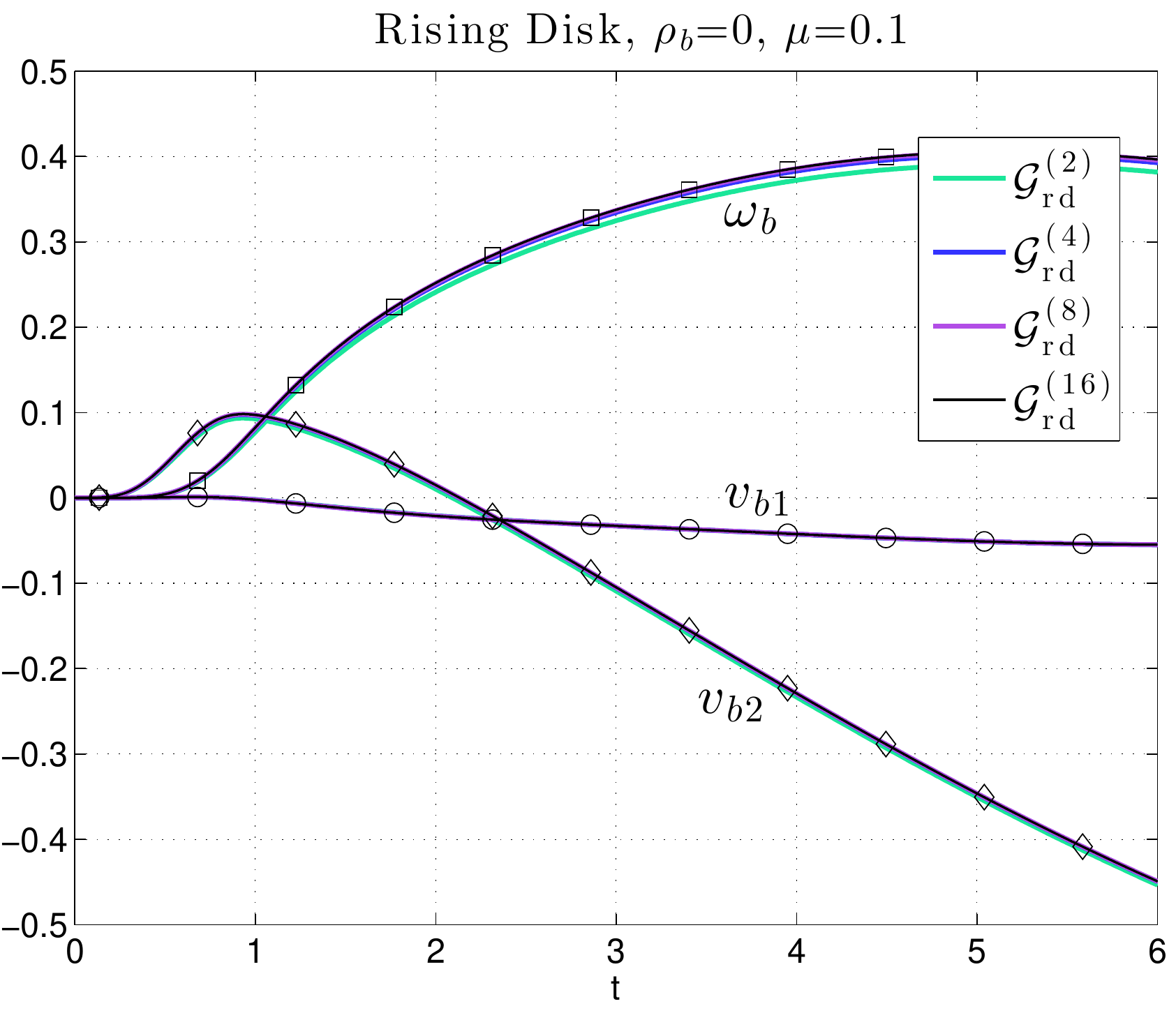}{\figWidth}};
  \draw(8.25,7.25) node[anchor=south west,xshift=-4pt,yshift=+0pt] {\trimfigz{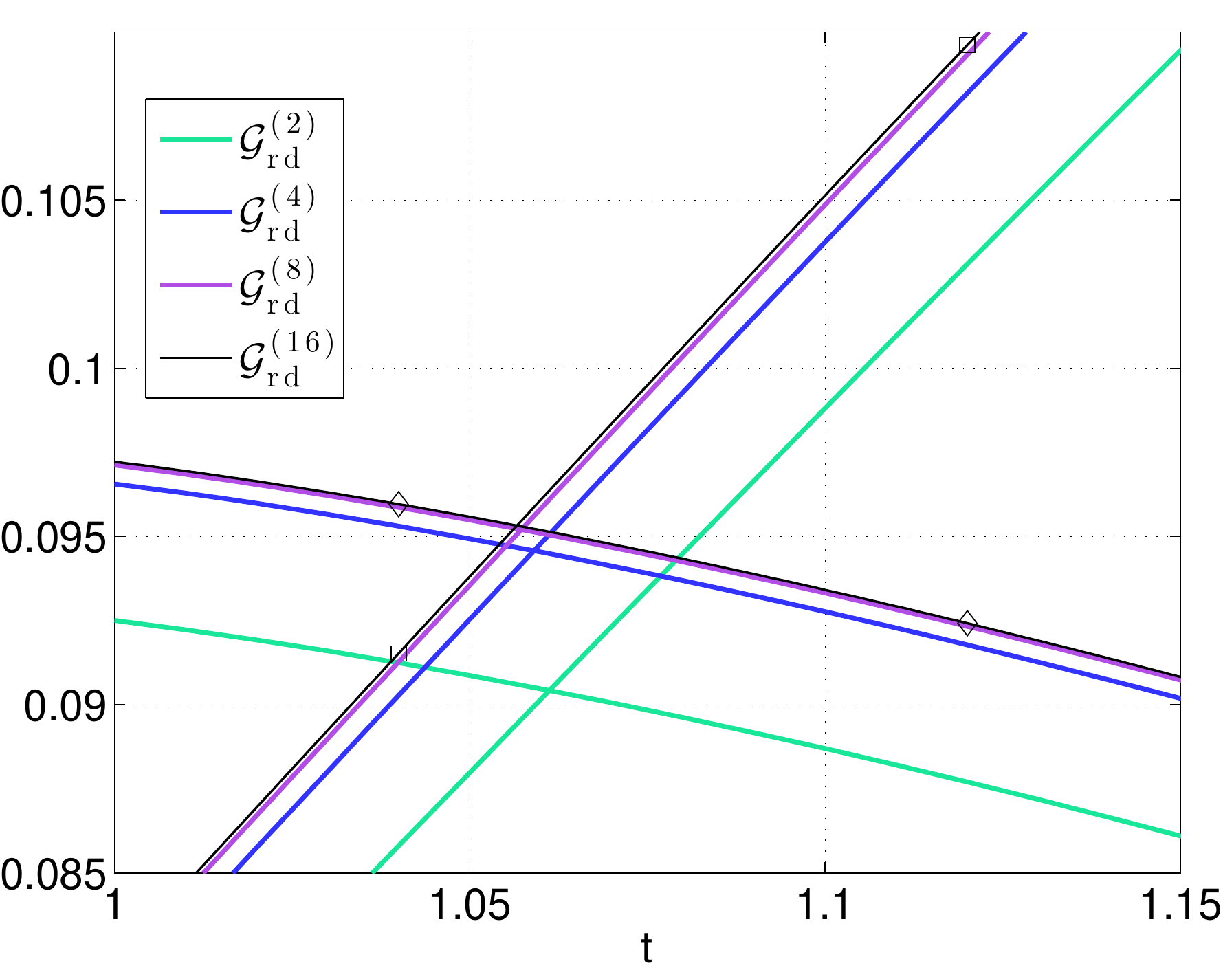}{\figWidthz}};
  \draw(0.0,0.0) node[anchor=south west,xshift=-4pt,yshift=+0pt] {\trimfig{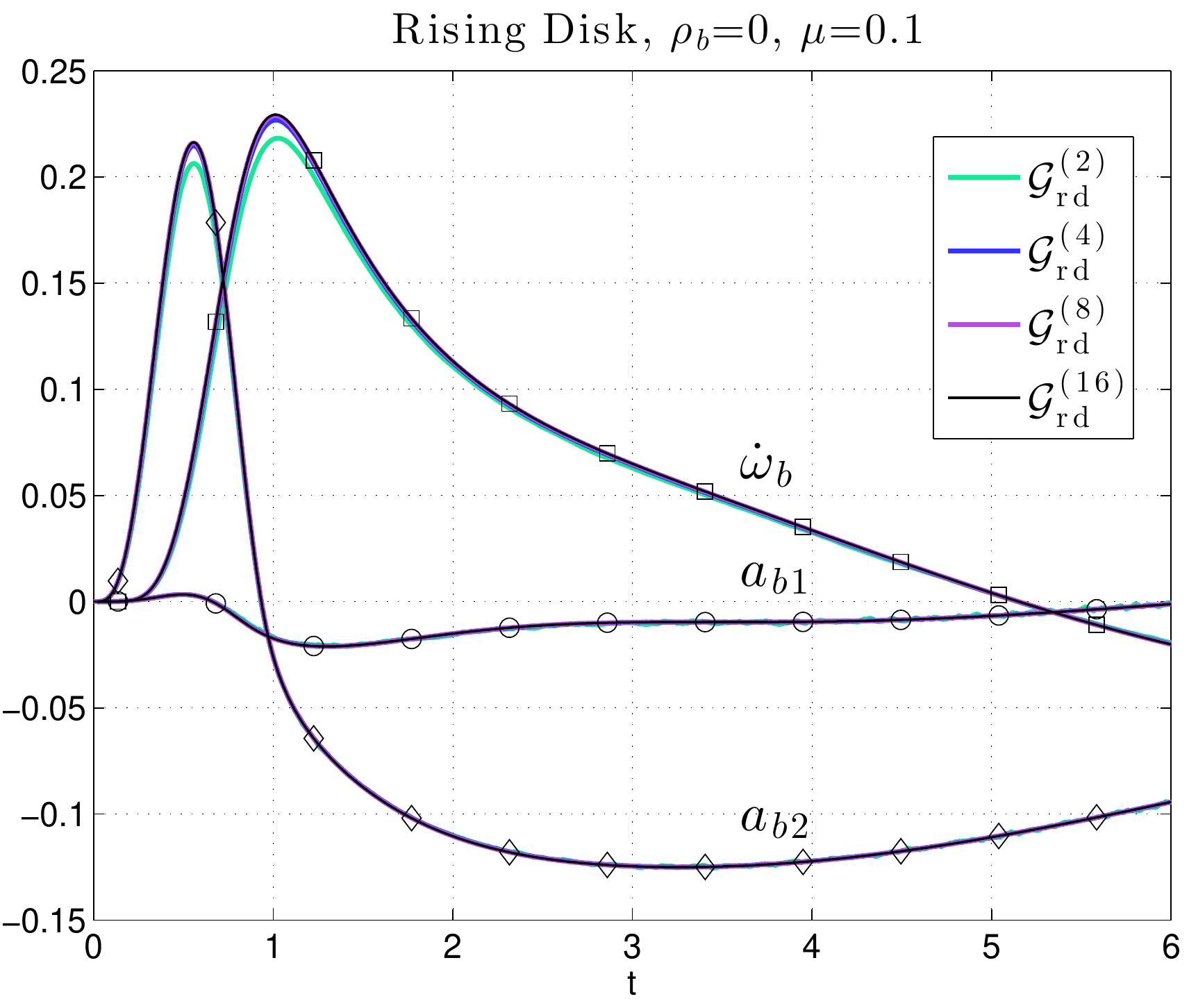}{\figWidth}};
  \draw(8.25,0.25) node[anchor=south west,xshift=-4pt,yshift=+0pt] {\trimfigz{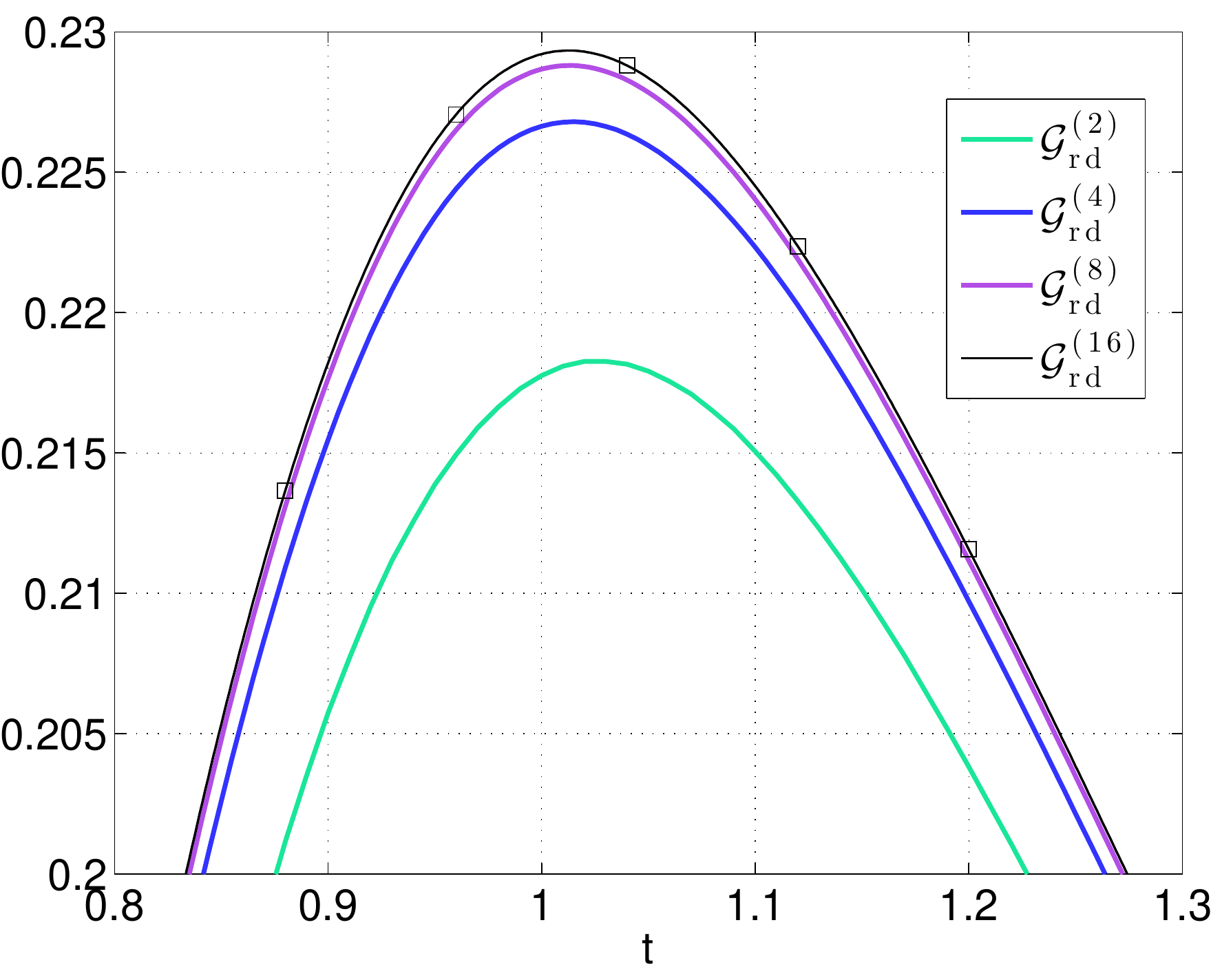}{\figWidthz}};
%
%
\end{tikzpicture}
}
\end{center}
  \caption{Rising disk in a counter-flow. Time history of the position and angular displacement of a light rigid body (top),
  its translational and angular velocities (middle), and its accelerations (bottom).  The plots on the right show selected enlarged views of the corresponding plots on the left.
    }
  \label{fig:risingDropCurves}
\end{figure}
}

\begin{table}[hbt]\tableFont 
\begin{center}
\begin{tabular}{|l|c|c|c|c|c|c|} \hline 
   $h^{(j)}$      & \errFormat{p} &  r   & \errFormat{u} &  r   & \errFormat{v} &  r  \\ \hline
 1/20 & \num{6.1}{-3} &      & \num{6.4}{-4} &      & \num{1.9}{-3} &      \\ \hline
 1/40 & \num{1.3}{-3} &  4.6 & \num{1.3}{-4} &  4.8 & \num{4.1}{-4} &  4.7 \\ \hline
 1/80 & \num{2.9}{-4} &  4.6 & \num{2.8}{-5} &  4.8 & \num{8.7}{-5} &  4.7 \\ \hline
 1/160  & \num{6.4}{-5} &  4.6 & \num{5.7}{-6} &  4.8 & \num{1.9}{-5} &  4.7 \\ \hline
 rate   &     2.19      &      &     2.27      &      &     2.22      &     \\ \hline
\end{tabular}
\qquad
\begin{tabular}{|l|c|c|c|} \hline 
         &  \multicolumn{3}{|c|}{Max-norm Rate} \\ \hline 
   t     &      p         &      u        &      v         \\ \hline
   0.5   &     2.19       &     2.27      &     2.22       \\ \hline
   1.0   &     2.40       &     2.26      &     2.40       \\ \hline
   2.0   &     2.61       &     2.42      &     2.57       \\ \hline
   3.0   &     3.26       &     3.20      &     3.00       \\ \hline
\end{tabular}
\caption{Estimated max-norm errors and convergence rates from a self-convergence grid refinement study.
Left: results for $t=0.5$. Right: estimated convergence rates at different times.
} \label{tab:cylDropMaxNormErrors}
\end{center}
\end{table}

\begin{table}[hbt]\tableFont 
\begin{center}
\begin{tabular}{|l|c|c|c|c|c|c|c|c|c|} \hline 
      & $x_b$ & $y_b$ & $\theta_b$ &  $v_{b1}$ & $v_{b2}$ & $\omega_b$ & $a_{b1}$ & $a_{b2}$ & $\dot\omega_b$  \\ \hline
rate &   3.17  &   3.17  &   2.02  &   1.96  &   3.34  &   1.98  &   2.61  &   2.52  &   1.95   \\ \hline
\end{tabular}
\caption{Rigid-body time-averaged convergence rates from a self-convergence grid refinement study.
     These results correspond to the curves in Figure~\ref{fig:risingDropCurves}.} \label{tab:cylDropCurves}
\end{center}
\end{table}

Figure~\ref{fig:cylDropStreamLines} shows (instantaneous) streamlines of the solution at a sequence
of times. This figure is complemented by the plots in Figure~\ref{fig:risingDropCurves} which show the
position, velocity and acceleration of the body over time.
Initially the disk starts to rise due buoyancy effects but then as the inflow pressure increases
the stronger downward flow drives the disk into a downward motion.
The centre of the disk moves slowly to the left.

Table~\ref{tab:cylDropMaxNormErrors} presents estimated convergence rates from a self-convergence grid-refinement study. 
The rates are computed using a Richardson extrapolation of the results 
from the finest three grids (see~\cite{pog2008a} for more details on this procedure). 
The estimated errors and convergence rates are given at $t=0.5$ and demonstrate close to second-order accurate convergence.
The max-norm rates at different times are also provided and show convergence rates that increase somewhat for later times; 
this is likely an artifact of the coarse grid solutions not being sufficiently resolved at later times and thus the 
asymptotic assumptions made in the Richardson process are not entirely accurate.

Figure~\ref{fig:risingDropCurves} shows the motion of the disk as computed on four grids of increasing resolution.
Magnified views show that the curves are converging nicely in a manner consistent with
second-order accuracy. Table~\ref{tab:cylDropCurves} confirms this observation by providing 
estimated convergence rates computed from a Richardson extrapolation. The values, which are in basic agreement
with the rates in Table~\ref{tab:cylDropMaxNormErrors},  represent the average
rate over time, $t\in[0,6]$, computed using 
\[
   \text{time-averaged rate} = \log_2\left(\f{ \| q_n^{(4)}- q_n^{(2)} \|_{1} }{  \| q_n^{(8)}- q_n ^{(4)} \|_{1} }\right), 
\]
where $\| q_n^{(j)}\|$ denotes the discrete $L_1$-norm of the quantity $q_n^{(j)}$ on grid $\Gcrd^{(j)}$,
over time-steps denoted by $n$.

\def\kgs{\, \rm{kg}/ \rm{s}}
\def\Gcfd{\Gc_{\rm{fd}}}
\def\meter{\, \rm{m}}

\subsection{Solid disk falling to the bottom of a fluid chamber} \label{sec:heavyCylDrop}

An independent check of the \ampRB~scheme can be made by comparing to solutions obtained by other schemes
for a standard test problem available in the literature.  One such problem is considered in~\cite{Robinson-MosherSchroederFedkiw2011, gibou2012efficient} and involves a moderately heavy solid-disk
that falls, due to gravity, in a fluid chamber.  The geometry of the problem is shown in the left plot of Figure~\ref{fig:heavyCylDropGrid}.  A solid disk of radius
$\diskRadius=0.5$ is positioned with centre initially at $(0,0)$
along the centreline and near the top 
of a rectangular channel with boundaries given by
$\channelWidth=2$, $\channelBottom=-14$ and $\channelTop=2$.  The densities of the fluid and disk are taken
to be $\rho=1$ and $\rhob=2$, respectively.  In contrast to the zero-mass disk problem studied in Section~\ref{sec:cylDrop},
the density of the disk is twice that of the fluid for this test problem, and thus
added-mass and added-damping effects are not as strong.  No-slip boundary conditions are assumed for the two
vertical and bottom horizontal walls of the channel.  At the top of
the channel, the pressure and the horizontal component of the fluid velocity are set to zero (while the vertical
component of the velocity is extrapolated to the ghost point).  
The fluid and the solid disk are at rest initially, and a downward
motion of the disk is initiated by an instantaneous application of a body force given by
\[
\fvbe(t)=\pi\diskRadius\sp2(\rho_b-\rho) \, \gv\, H(t),
\]
where $\gv=[0,\,-g]$ is the acceleration due to gravity and $H(t)$ is the Heaviside function.

{
\newcommand{\drawContour}[5]{%
\begin{scope}[#1]
\draw(0.0,0) node[anchor=south west,xshift=-4pt,yshift=+0pt] {\trimfig{fig/#2}{\figWidth}};
\draw(1.1,9) node[draw,fill=white,anchor=west,xshift=2pt,yshift=0pt] {\scriptsize #3};
\end{scope}
}
\newcommand{\cbWidth}{.2cm}
\newcommand{\cbHeight}{2cm}
\newcommand{\xcb}{.5cm}
\newcommand{\ycb}{-.2cm}
\setlength{\ycbTop}{\ycb+\cbHeight}
\setlength{\ycbMid}{\ycb+\cbHeight*\real{.5}}
\newcommand{\trimfigcb}[3]{\includegraphics[width=#2, height=#3, clip, trim=17cm 2.35cm 1.65cm 2.35cm]{#1}}
\def\rad{.25}
\newcommand{\plotDisk}{
\fill[fill=red!20,draw=red,line width=2pt] 
      plot[samples=100, domain=0.:360] ( {\rad*cos(\x)} , {\rad*sin(\x)} ) -- cycle ;
}
\newcommand{\figWidth}{9.35cm}
\newcommand{\trimfig}[2]{\trimh{#1}{#2}{.34}{.35}{.025}{.1}}
\newcommand{\trimfiga}[2]{\trimh{#1}{#2}{.27}{.27}{.025}{.1}}
\newcommand{\labelsize}{\small}
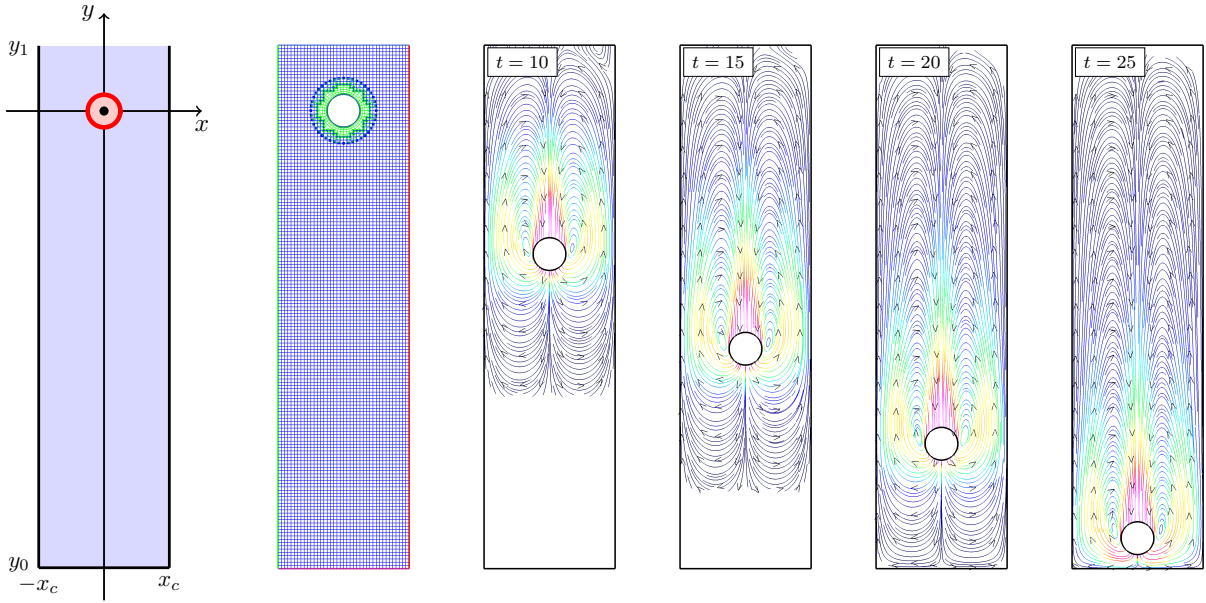
\begin{figure}[htb]
\begin{center}
\resizebox{14cm}{!}{
\begin{tikzpicture}[scale=1]
  \useasboundingbox (-.5,1) rectangle (15.5,9);  
  \begin{scope}[xshift=-.2cm,yshift=7.7cm]
      \fill[fill=blue!15] (-1.,-7) -- (1.,-7) -- (1,1) -- (-1,1) -- cycle; 
      \draw[->,thick] (-1.5,0) -- (1.5,0) node [anchor=north] {$x$}; 
      \draw[->,thick] (0,-7.5) -- (0,1.5) node [anchor=east] {$y$}; 
      \draw[-,very thick] (-1,-7) -- (-1,1);
      \draw[-,very thick] (1,-7) -- (1,1);
      \draw[-,very thick] (-1,-7) -- (1,-7);
    \begin{scope}[xshift=0cm,yshift=0cm]  
      \plotDisk
    \end{scope}
    \fill[black] (0,.0) circle (2 pt); 
	\draw(-1.3,1.2) node[anchor=north] {\labelsize$\channelTop$};
	\draw(-1.3,-6.7) node[anchor=north] {\labelsize$\channelBottom$};
	\draw(-1,-7) node[anchor=north] {\labelsize$-\channelWidth$};
	\draw(1,-7) node[anchor=north] {\labelsize$\channelWidth$};
  \end{scope}
  \draw(1.2cm,0) node[anchor=south west,xshift=0pt,yshift=-.55cm] {\trimfig{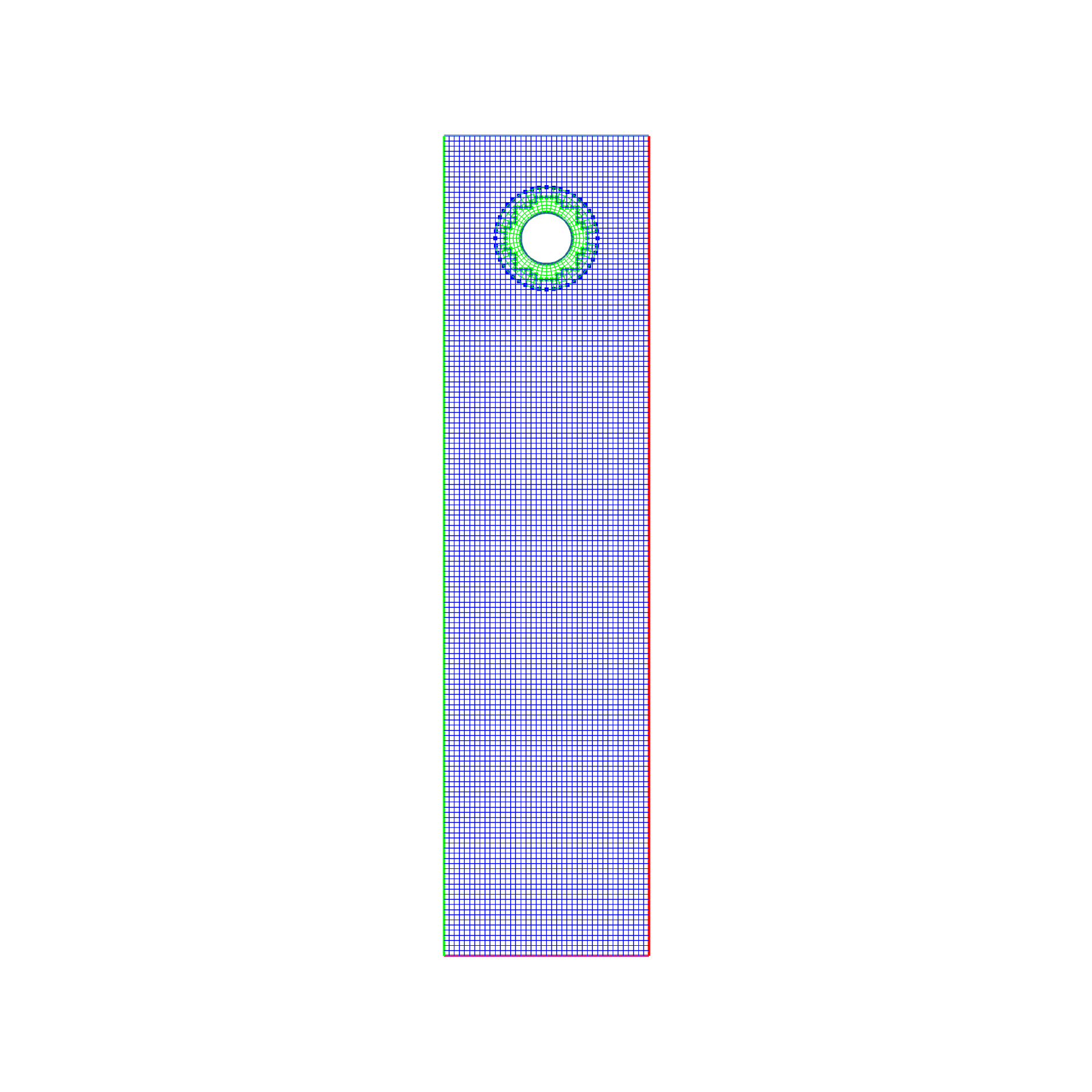}{\figWidth}};
  \begin{scope}[xshift=4.5cm,yshift=-.55cm]
  \drawContour{xshift=0cm,yshift=0pt}{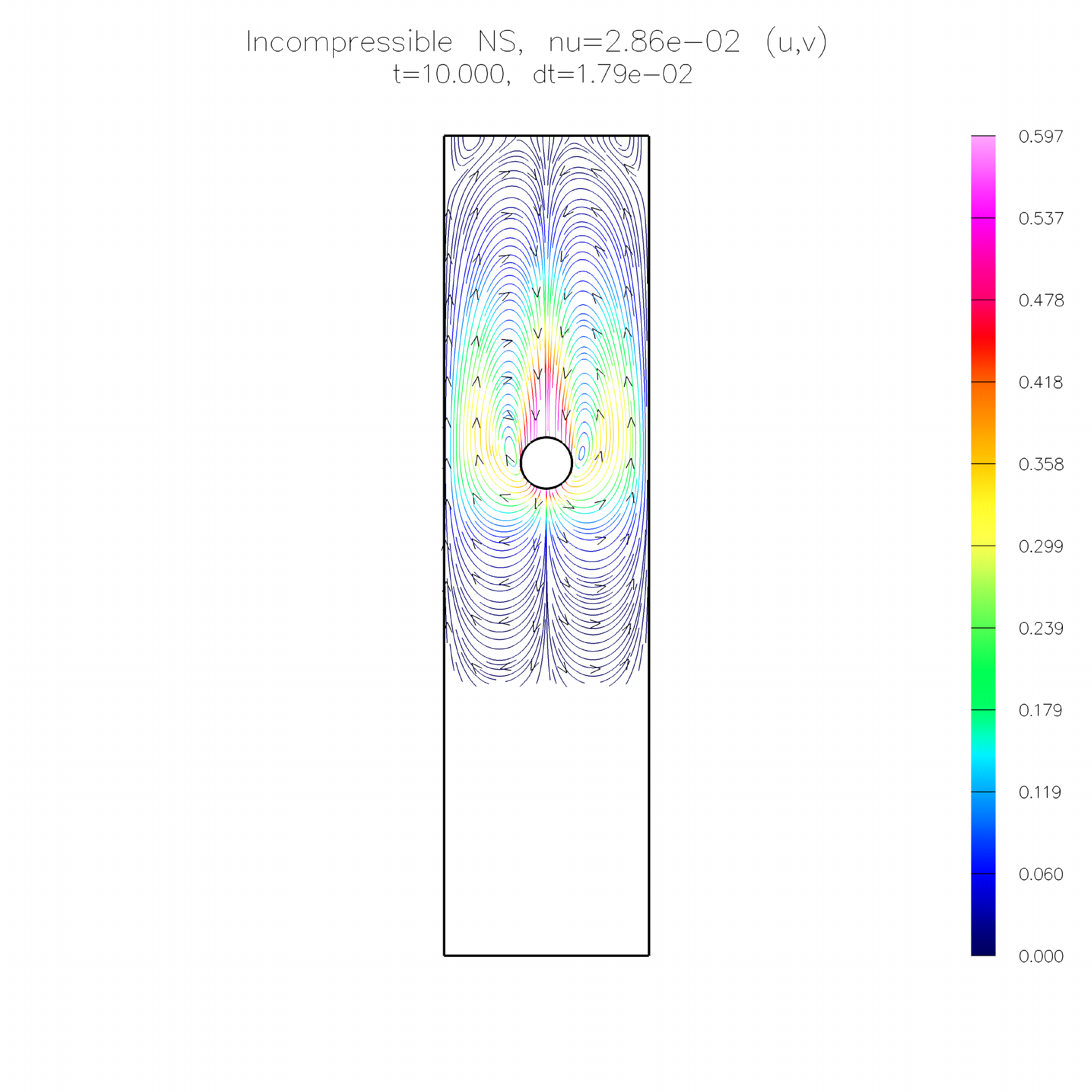}{$t=10$}{0}{0.60};
  \drawContour{xshift=3cm,yshift=0pt}{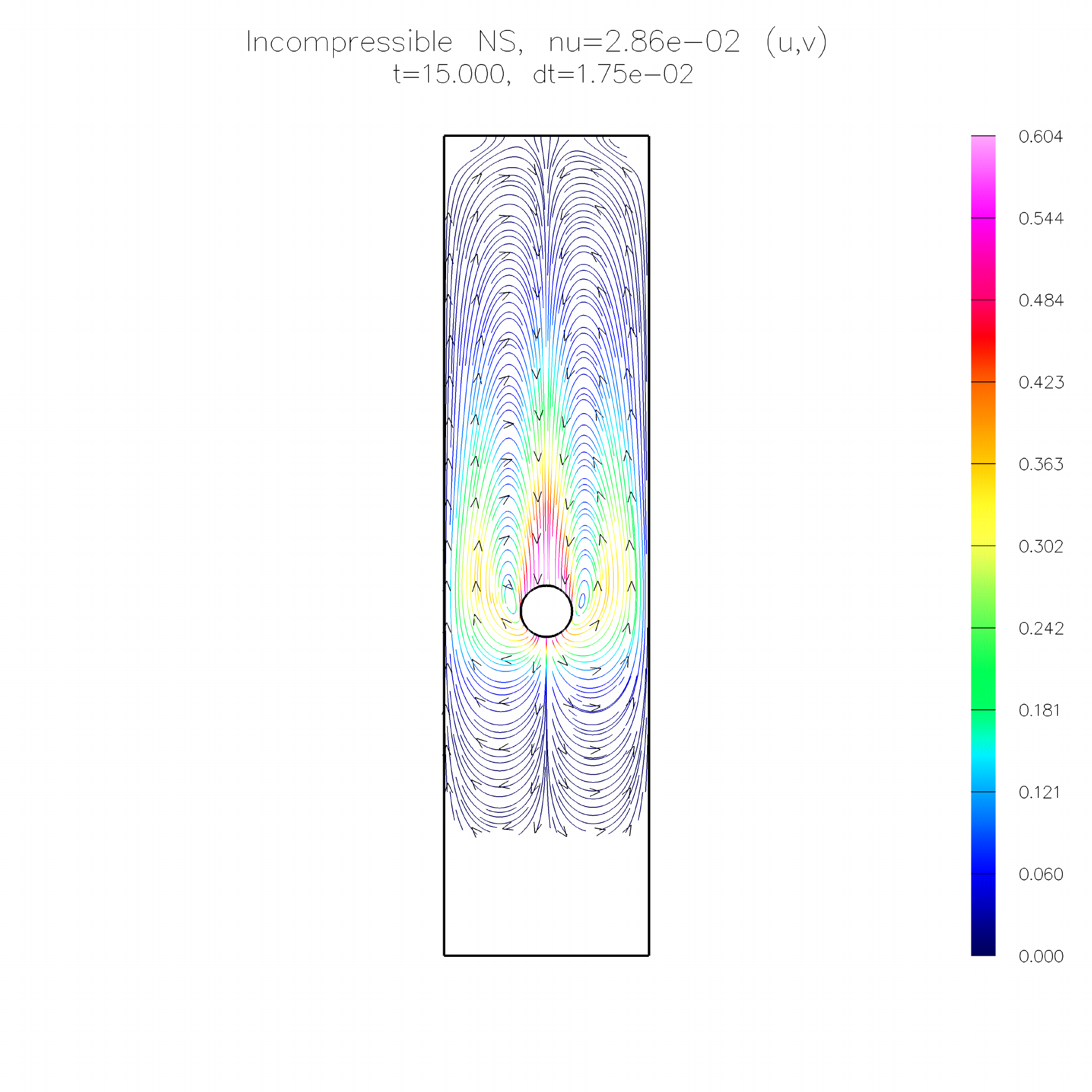}{$t=15$}{0}{0.60};
  \drawContour{xshift=6cm,yshift=0pt}{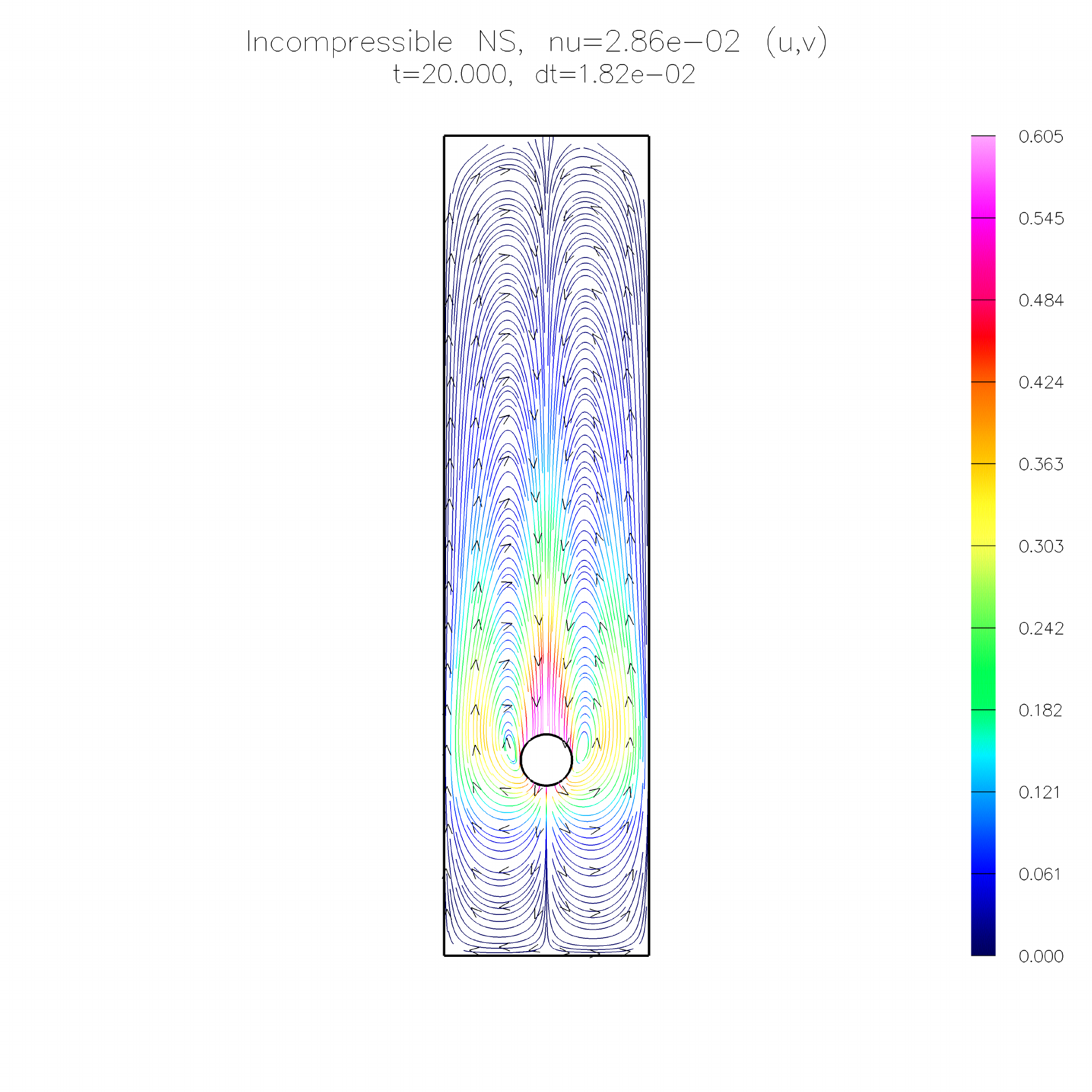}{$t=20$}{0}{0.61};
  \drawContour{xshift=9cm,yshift=0pt}{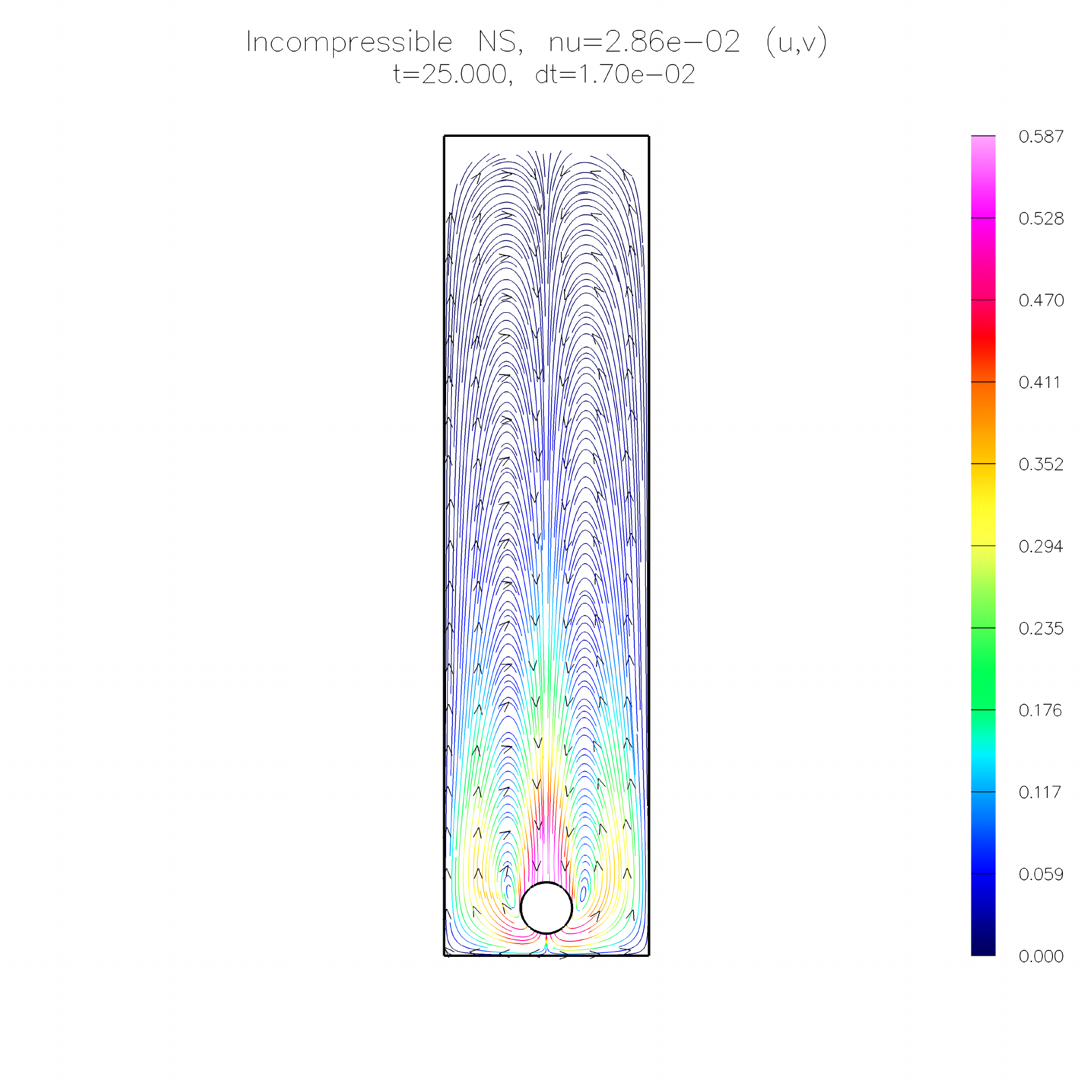}{$t=25$}{0}{0.59};
  \end{scope}
\end{tikzpicture}
}
\end{center}
  \caption{Falling disk in a fluid chamber. Left: Geometrical configuration at $t=0$.
  Middle: Composite grid $\Gcfd^{(1)}$ at $t=0$ (coarse grid). 
  Right: Computed streamlines at selected times using grid $\Gcfd^{(4)}$.
}
  \label{fig:heavyCylDropGrid}
\end{figure}
}

The results in~\cite{Robinson-MosherSchroederFedkiw2011, gibou2012efficient} are described in terms of dimensional
quantities, whereas we have chosen to work with dimensionless variables as was done for the previous test problems
in this section.  The results here can be converted to dimensional quantities using a reference scale for length
given by the diameter of the disk, $\tilde d=0.01\;{\rm m}$, and a reference scale for density given by the
fluid density, $\tilde\rho=10\sp3\;{\rm kg}/{\rm m}\sp3$. 
A reference scale for velocity is taken to be an estimate for the
terminal velocity of a disk falling 
in a channel, assuming Stokes flow, given by $\tilde v_\text{term}=0.35011\;{\rm m}/{\rm s}$ for
the case of a fluid with viscosity given by $\tilde\mu=0.1\;{\rm kg}/({\rm m}\,{\rm s})$, 
see~\cite{Robinson-MosherSchroederFedkiw2011, gibou2012efficient}.  The reference time is thus given by $\tilde d/\tilde v_\text{term}=0.028562\,{\rm s}$.  With these reference scales, the dimensionless fluid viscosity is $\mu=\tilde\mu/(\tilde\rho\, \tilde d\, \tilde v_\text{term})=0.028562$, which is the reciprocal of the Reynolds number (based on the diameter of the solid disk), and the dimensionless
gravitational constant is $g=\tilde g\,\tilde d/(\tilde v_\text{term})\sp2=0.79949$, where $\tilde g=9.8\;{\rm m}/{\rm s}\sp2$.

{
\newcommand{\figWidth}{8.cm}
\newcommand{\trimfig}[2]{\trimw{#1}{#2}{.0}{.0}{.0}{.0}}
\newcommand{\trimfiga}[2]{\trimw{#1}{#2}{.06}{.145}{.02}{.0}}
\newcommand{\figWidthz}{7.5cm}
\newcommand{\trimfigz}[2]{\trimw{#1}{#2}{.0}{.0}{.0}{.0}}
\begin{figure}[htb]
\begin{center}
\resizebox{12cm}{!}{
\begin{tikzpicture}[scale=1]
  \useasboundingbox (0.,1.3) rectangle (16.,21);  
  \draw( 0.0,14.0) node[anchor=south west,xshift=-4pt,yshift=+0pt] {\trimfig{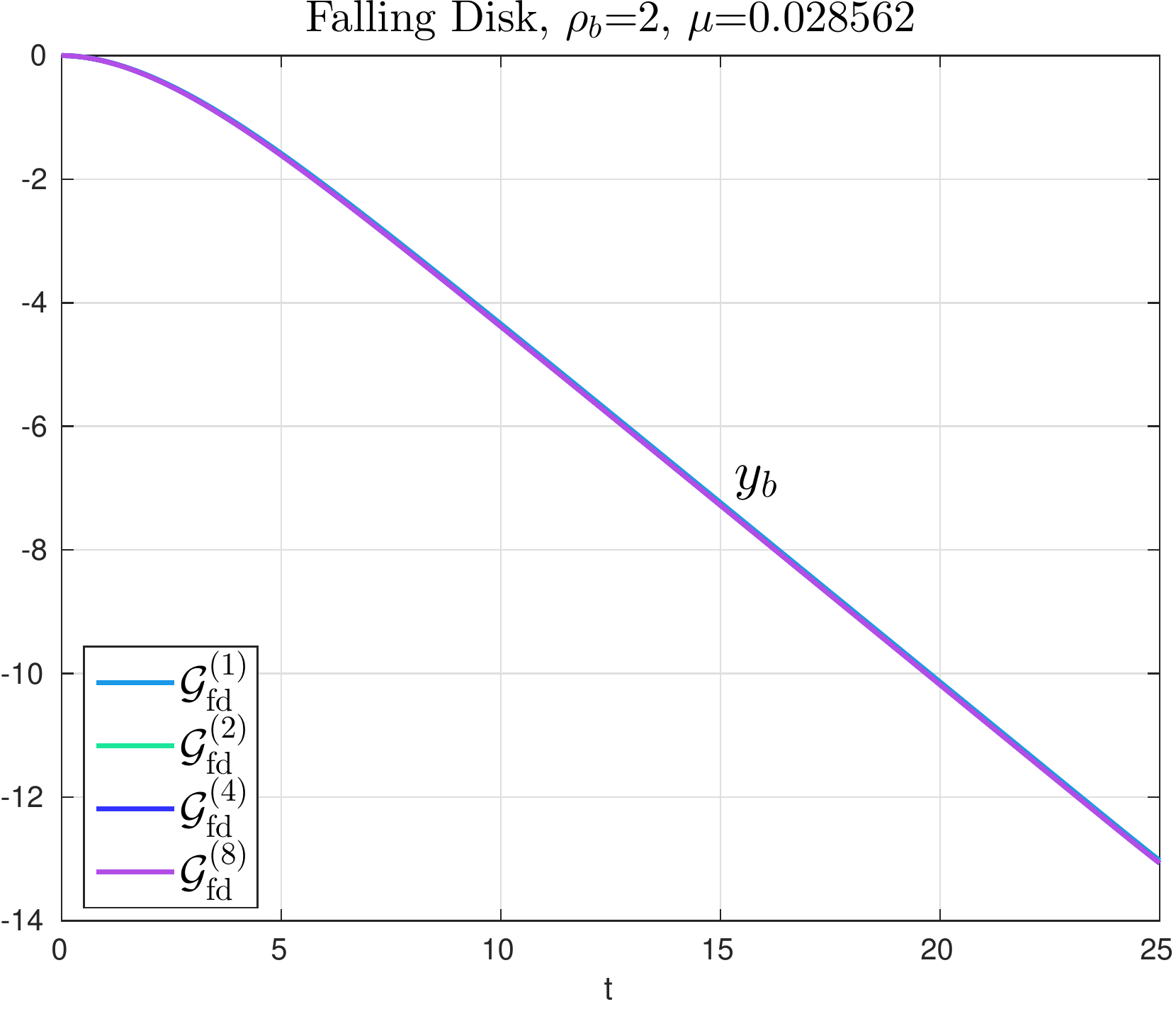}{\figWidth}};
  \draw(8.25,14.25) node[anchor=south west,xshift=-4pt,yshift=+0pt]{\trimfigz{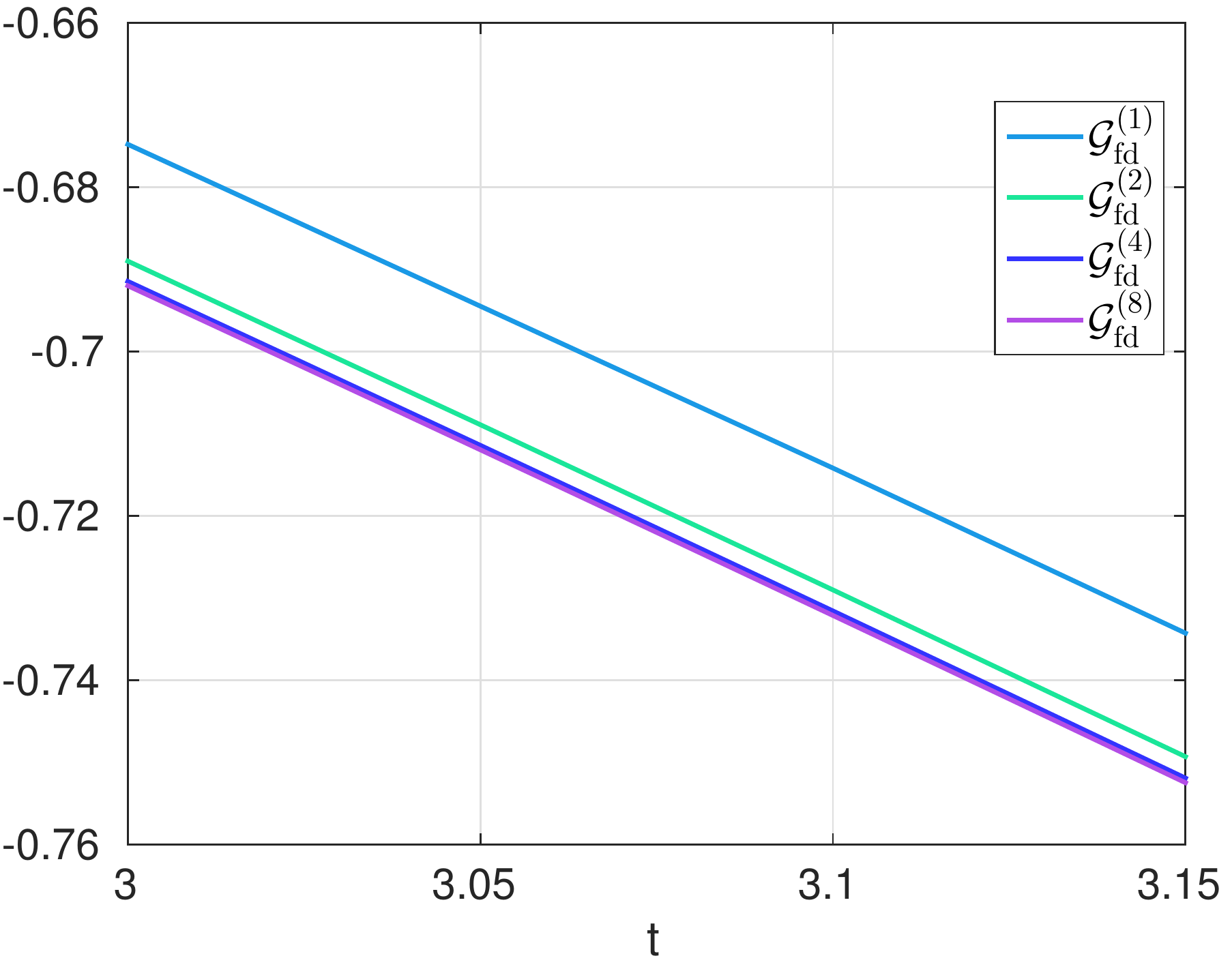}{\figWidthz}};
  \draw( 0.0,7.0) node[anchor=south west,xshift=-4pt,yshift=+0pt] {\trimfig{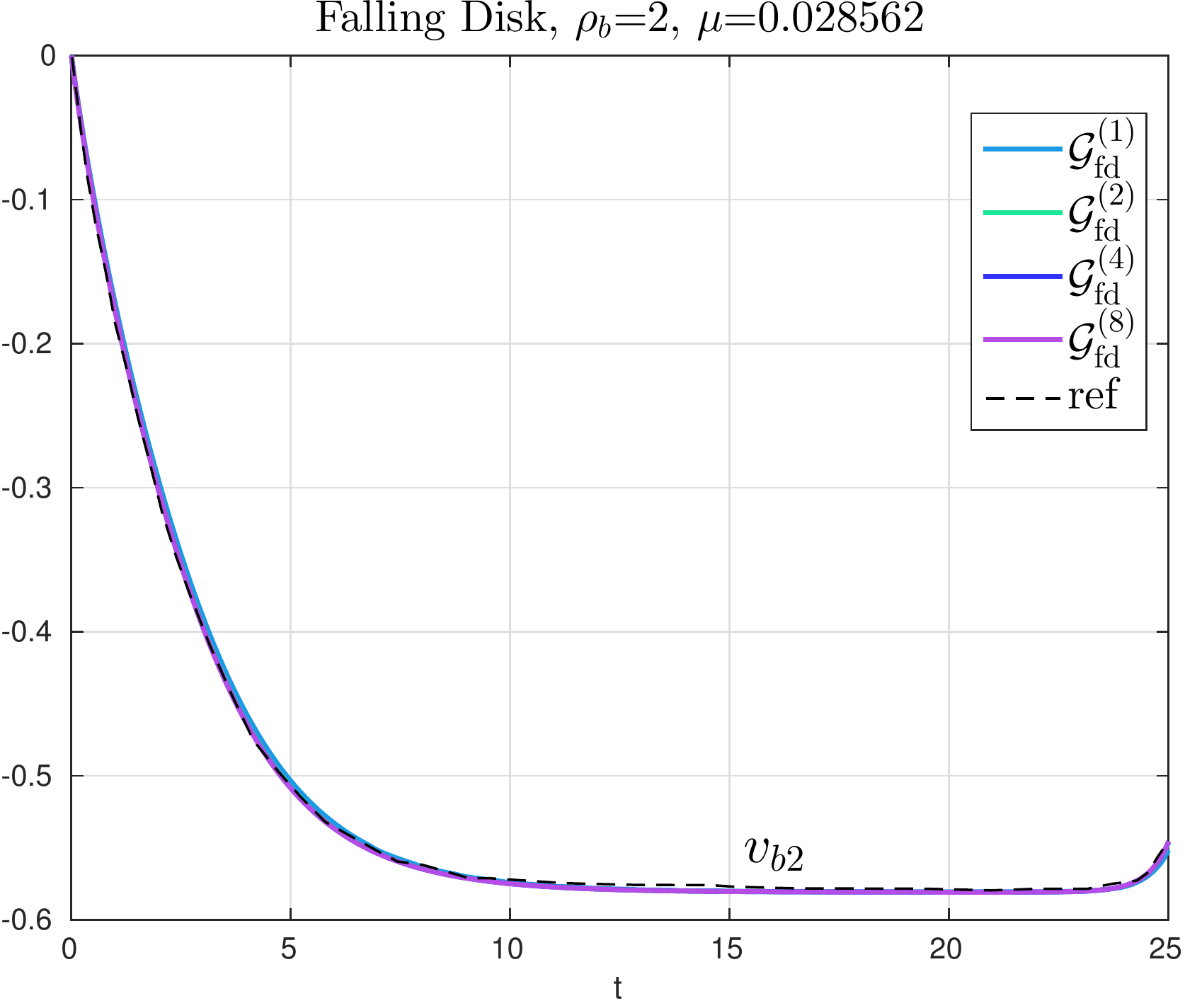}{\figWidth}};
  \draw(8.25,7.25) node[anchor=south west,xshift=-4pt,yshift=+0pt]{\trimfigz{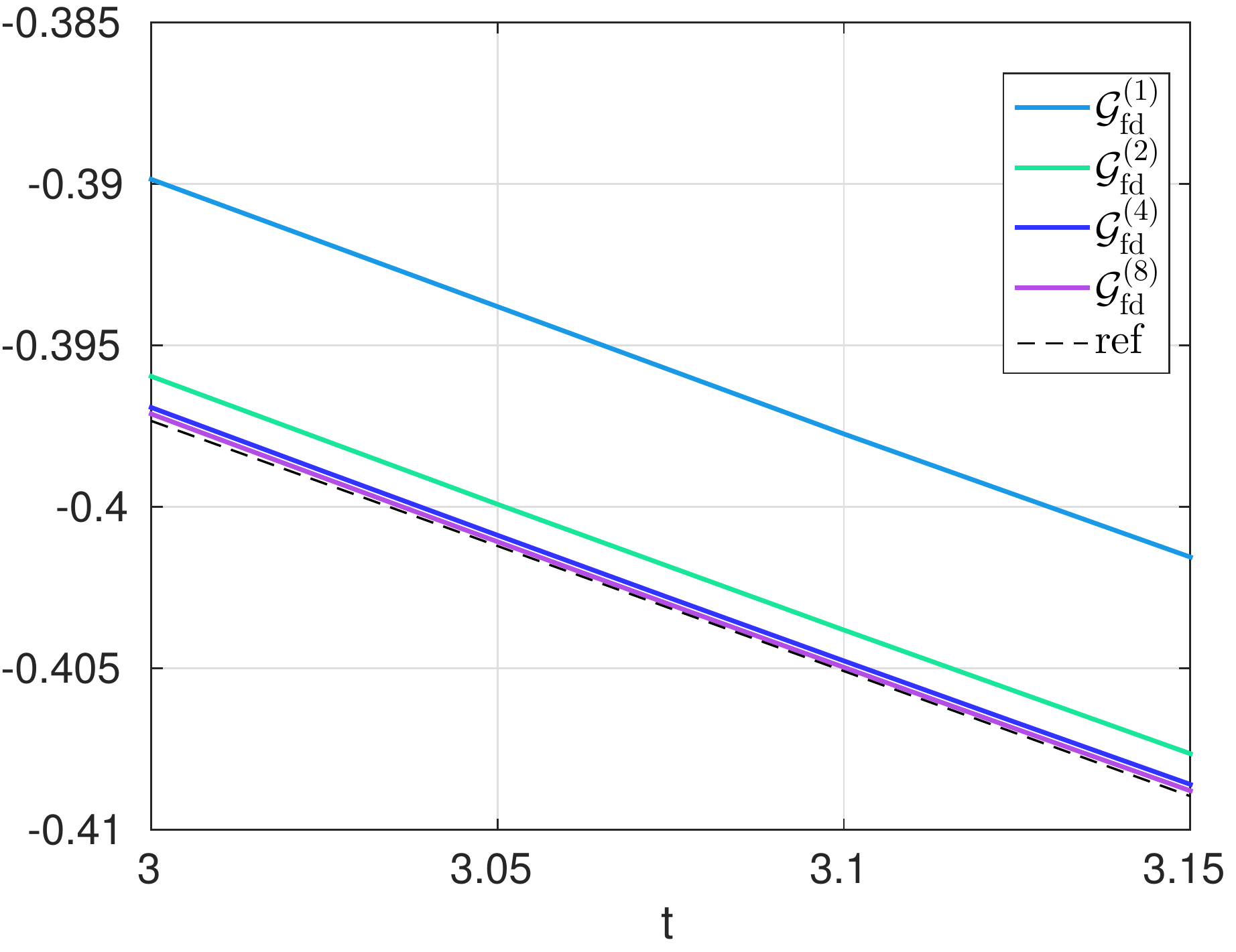}{\figWidthz}};
  \draw( 0.0,0.0) node[anchor=south west,xshift=-4pt,yshift=+0pt]  {\trimfig{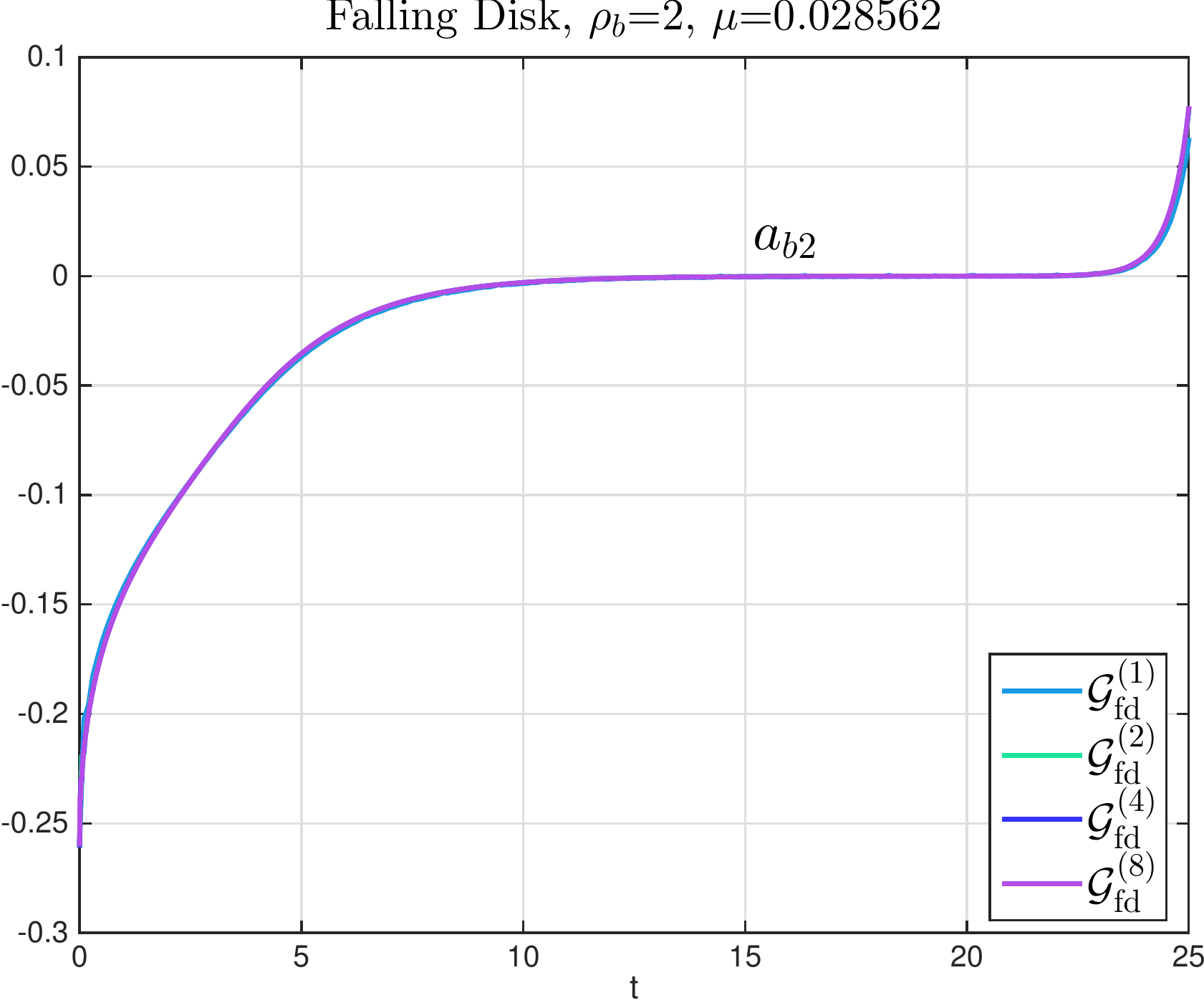}{\figWidth}};
  \draw(8.25,0.25) node[anchor=south west,xshift=-4pt,yshift=+0pt] {\trimfigz{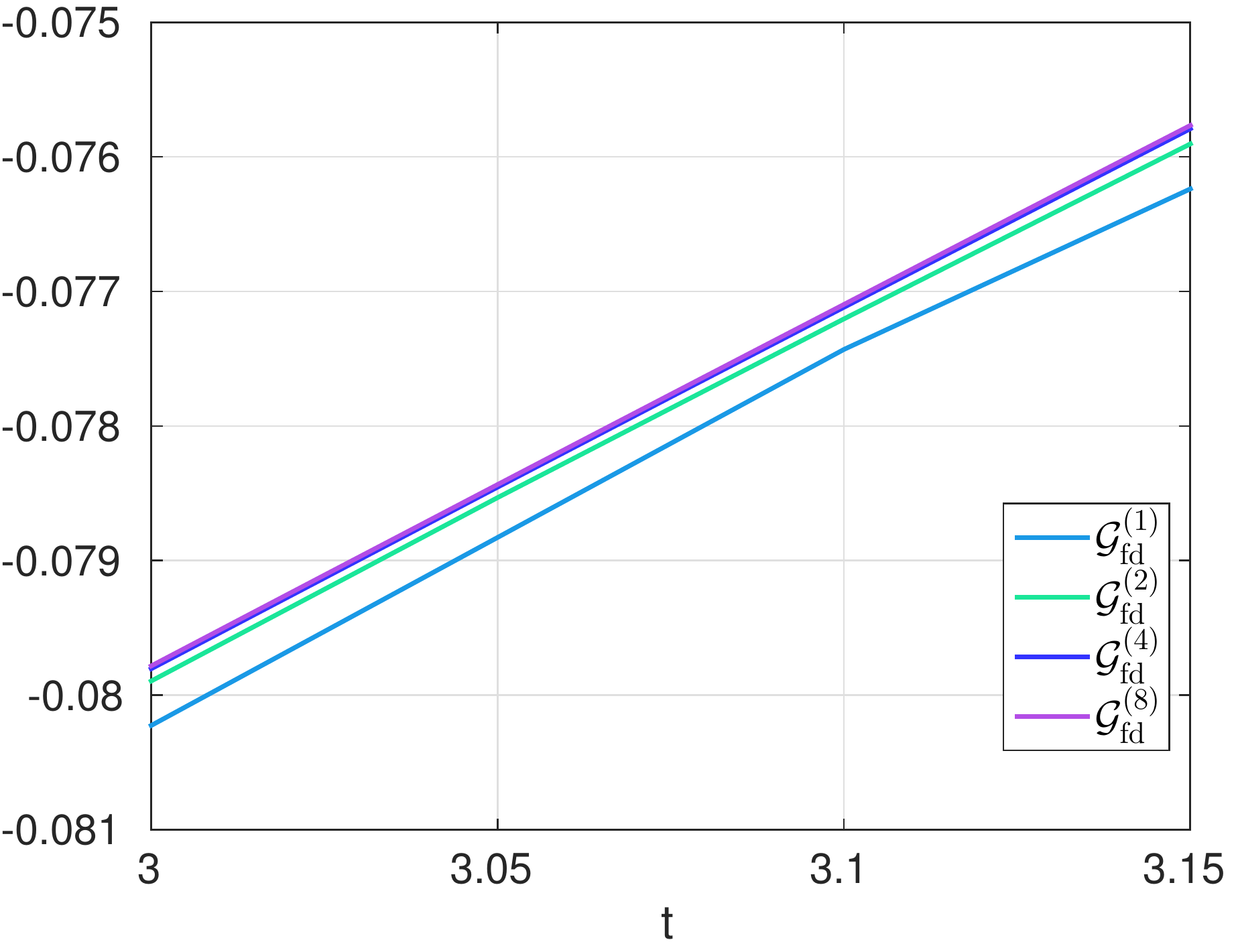}{\figWidthz}};
%
\end{tikzpicture}
}
\end{center}
  \caption{Falling disk in a fluid chamber.
	  Time history of the position (top), velocity (middle) and acceleration (bottom) of the disk in the vertical direction.  
	  The plots on the right show selected enlarged views of the corresponding plots on the left.
          The black dashed-line (``ref'') is taken from Fig.~7.2(b) in~\cite{gibou2012efficient}.
    }
  \label{fig:heavyCylDropCurvesicase1}
\end{figure}
}

Numerical solutions are computed using the \ampRB~scheme with composite grids
of varying resolution, one such grid is shown in Figure~\ref{fig:heavyCylDropGrid}
for the setup at $t=0$.  The fluid domain in the channel is represented by a background Cartesian grid (shown in blue) and a moving annular body-fitted grid (shown in green).  The composite grid has a target grid spacing of $h^{(j)}=1/(10 j)$, with grid resolution factor $j$, and is denoted by~$\Gcfd^{(j)}$.  The annular grid has a fixed radial width
equal to~$0.4$ for all grid resolutions.  The time-step is taken as $\dt^{(j)}= h^{(j)}$.

The streamlines of the solutions at a sequence of times are presented in Figure~\ref{fig:heavyCylDropGrid}, 
complemented by Figure~\ref{fig:heavyCylDropCurvesicase1} which shows the vertical components of the position, $y_b$, velocity, $v_{b2}$, and acceleration, $a_{b2}$, of the falling disk as a function of time, $t$.  The solutions
are computed using the composite grid, $\Gcfd^{(j)}$, with resolution factors $j=1$, 2, 4 and~8.  The plots in the left column of Figure~\ref{fig:heavyCylDropCurvesicase1} show the general behaviour of $y_b$, $v_{b2}$ and $a_{b2}$, while the plots in the right column show enlarged views which indicate the convergence of the various solution components.
The dashed curves in the velocity plots are taken from the numerical results given in Fig.~7.2(b) from~\cite{gibou2012efficient} for their finest grid resolution which corresponds to our grid~$\Gcfd^{(4)}$ approximately.
We observe that the results computed using the \ampRB~scheme are in excellent agreement with the result from~\cite{gibou2012efficient}, and show a convergence consistent with second-order accuracy.  In fact, Richardson-extrapolation estimates for the time-averaged convergence rate are computed as in Section~\ref{sec:cylDrop}, and these rates are found to be $2.11$ for the position,  $2.16$ for the velocity and~$2.43$ for the acceleration, respectively, in agreement with second-order accuracy.
The remaining cases in~\cite{Robinson-MosherSchroederFedkiw2011, gibou2012efficient} using other choices of the viscosity have been computed using the \ampRB~scheme, and we have found similar agreement and convergence behaviour.  
Finally, we note that the \tpRB~scheme is also stable without sub-iterations for this problem, and produces results 
nearly identical to the \ampRB~scheme.

\subsection{Rectangular body rising to the top of fluid chamber} \label{sec:fallingBody}

\newcommand{\Gcrb}{\Gc_{\rm rb}} 

Consider a light and rectangular-shaped body which rises under buoyancy forces in a heavier fluid and settles near the
top of a closed fluid chamber as shown in Figure~\ref{fig:risingBodyGrid}.
The fluid with density $\rho=1$ and viscosity $\mu=0.025$ occupies the square domain $\OmegaF=[-1,1]\times[-1,1]$
and is at rest initially.  No-slip boundary conditions are taken along the four walls of the fluid chamber.
The rigid-body is a rectangular solid of width $\solidWidth=1$, height $\solidHeight=0.5$
and uniform density $\rhob=0.001$, with initial centre-of-mass located at the origin.
The corners of the solid are rounded in order to avoid singular behaviour in the solution 
that can occur near sharp convex corners. 
The buoyant solid rises under the influence of gravity according to the body force in~\eqref{eq:linearAccelerationEquation} given by
\begin{equation}
\fvbe(t)=\solidArea(\rho_b-\rho)\, \gv \, \rampFunction(t),
\label{eq:bodyForce}
\end{equation}
where $\solidArea\approx0.5$ is the computed area of the solid, $\gv=[0,\,-1]\sp{T}$ is an acceleration due to gravity, and
the ramp function $\rampFunction$ is given by~\eqref{eq:ramp}. 
{
\newcommand{\figWidth}{5cm}
\newcommand{\trimfig}[2]{\trimFig{#1}{#2}{.15}{.375}{.25}{.25}}
\begin{figure}[htb]
\begin{center}
\begin{tikzpicture}[scale=1]
  \useasboundingbox (0.0,.75) rectangle (16.,5.25);  
  \draw(0.0, 0) node[anchor=south west,xshift=-4pt,yshift=+0pt] {\trimfig{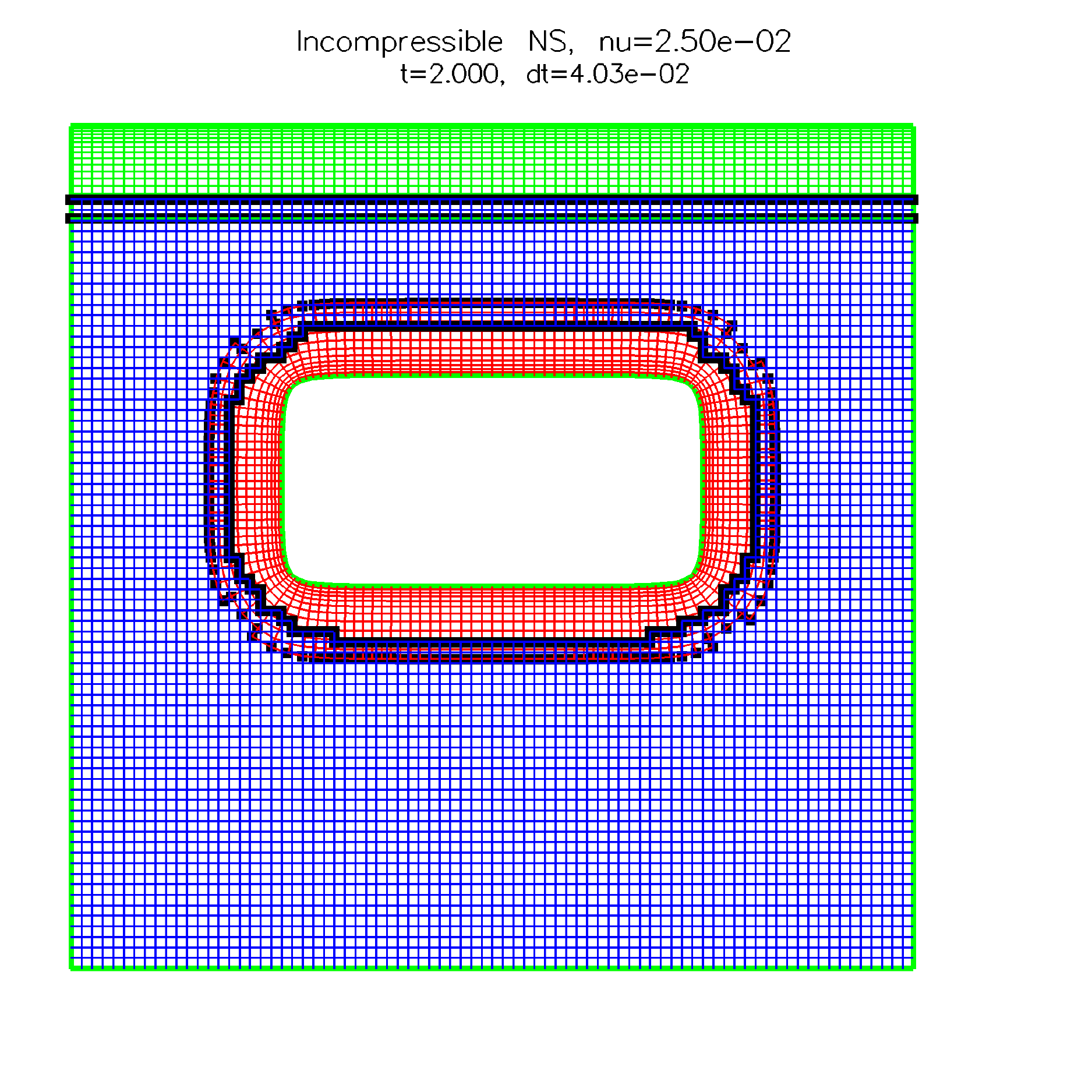}{\figWidth}};
  \draw(5.5, 0) node[anchor=south west,xshift=-4pt,yshift=+0pt] {\trimfig{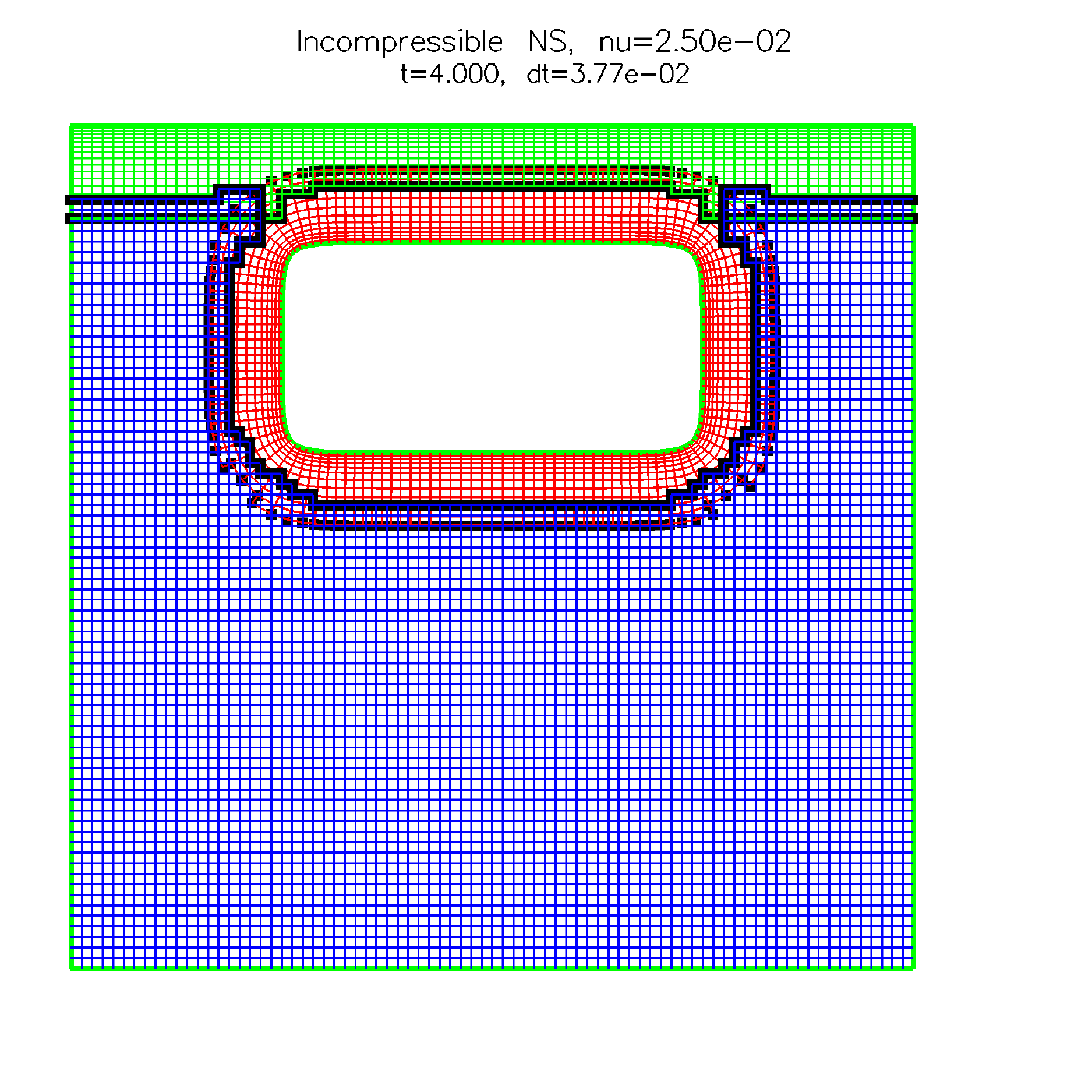}{\figWidth}};
  \draw(11.0,0) node[anchor=south west,xshift=-4pt,yshift=+0pt] {\trimfig{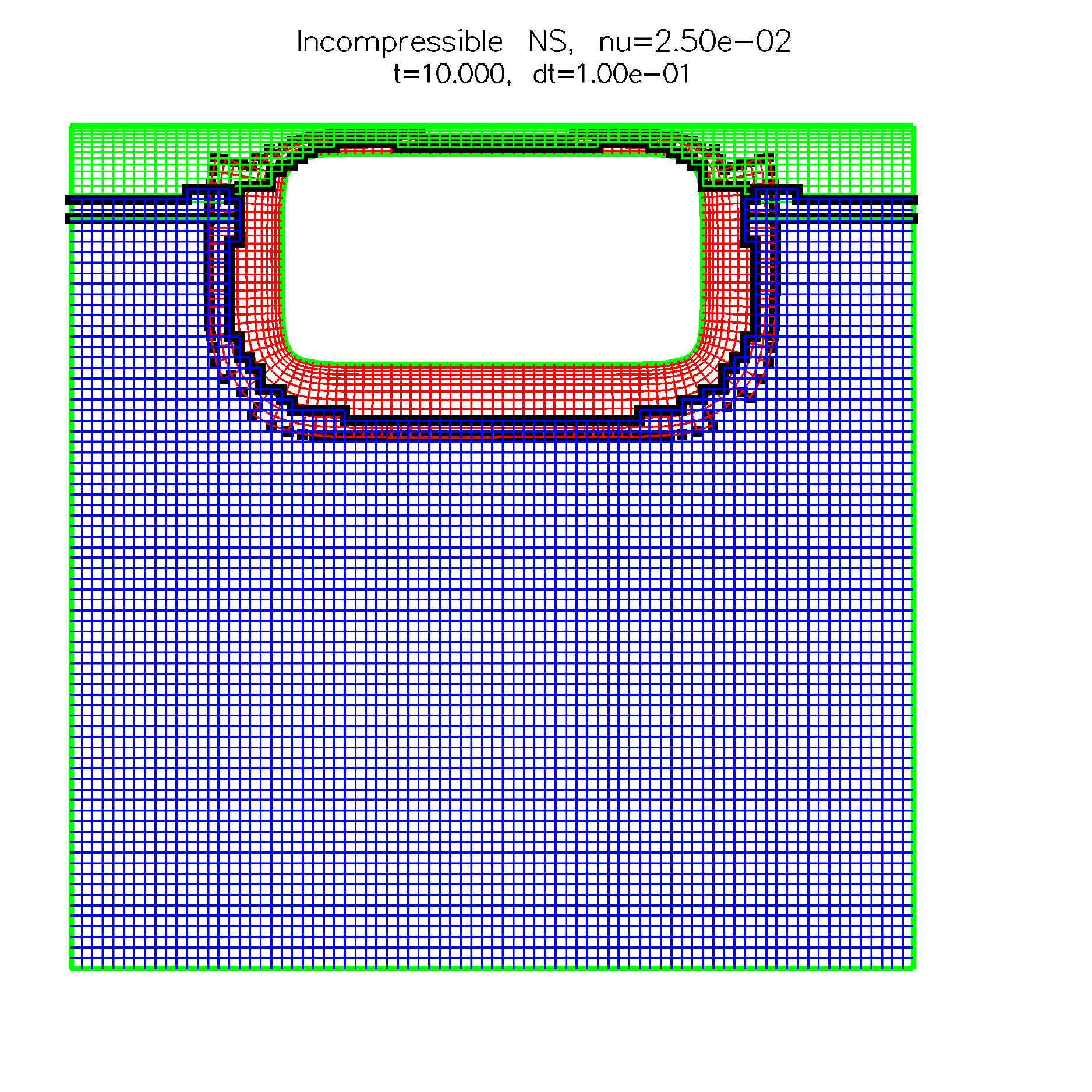}{\figWidth}};
%
\end{tikzpicture}
\end{center}
  \caption{Rising body: Coarse composite grid, $\Gcrb^{(4)}$, at times $t=2$, $4$ and~$10$.  The composite grid consists of 
a background Cartesian component grid, a boundary-fitted grid attached to the surface of the
rising solid, and a
stretched boundary-fitted grid along the top wall of the fluid chamber with a finer grid resolution to accommodate a close encounter with the rising solid. 
}
  \label{fig:risingBodyGrid}
\end{figure}
}

The composite grid used for this problem, $\Gcrb^{(j)}$, is shown in Figure~\ref{fig:risingBodyGrid}, and
consists of a background Cartesian grid (blue), a boundary-fitted grid 
attached to the surface of the rising solid (red), and a
stretched boundary-fitted grid (green) along the top wall of the fluid chamber to accommodate a close encounter with the rising solid. As before, the grid with resolution factor is indicated by the superscript $j$.
The target background grid spacing is $h^{(j)}= 1/(10 j)$, while the target 
boundary-layer spacing is four times finer, $h_{BL}^{(j)}=h^{(j)}/4$, near the top wall.

{
\newcommand{\drawContour}[7]{%
\begin{scope}[#1]
\draw(0.0,0) node[anchor=south west,xshift=-4pt,yshift=+0pt] {\trimfig{fig/#2}{\figWidth}};
  \draw(.5,.5) node[draw,fill=white,anchor=west,xshift=2pt,yshift=1pt] {\scriptsize #3};
  \draw(1.,.5) node[draw,fill=white,anchor=west,xshift=2pt,yshift=1pt] {\scriptsize #5};
\begin{scope}[xshift=9pt]
  \draw (\xcb,\ycb) node[anchor=south west,xshift=0.cm,yshift=.5cm,rotate=-90] {\trimfigcb{fig/colourBarLines}{\cbWidth}{\cbHeight}};
  \draw (.8,0) node[anchor=north,xshift=+3pt,yshift=+2pt] {\scriptsize $#6$};
  \draw (4.4,0) node[anchor=north,xshift=+0pt,yshift=+2pt] {\scriptsize $#7$};
\end{scope}
\end{scope}
}
\newcommand{\cbWidth}{.2cm}
\newcommand{\cbHeight}{4cm}
\newcommand{\xcb}{.5cm}
\newcommand{\ycb}{-.2cm}
\setlength{\ycbTop}{\ycb+\cbHeight}
\setlength{\ycbMid}{\ycb+\cbHeight*\real{.5}}
\newcommand{\trimfigcb}[3]{\includegraphics[width=#2, height=#3, clip, trim=17cm 2.35cm 1.65cm 2.35cm]{#1}}
%
%
\newcommand{\figWidth}{5cm}
\newcommand{\trimfig}[2]{\trimw{#1}{#2}{.05}{.15}{.1}{.1}}
\begin{figure}[htb]
\begin{center}
\resizebox{14cm}{!}{
\begin{tikzpicture}[scale=1]
  \useasboundingbox (0.4,1) rectangle (16.,11);  
 \drawContour{xshift= 0.0cm,yshift=5.75cm}{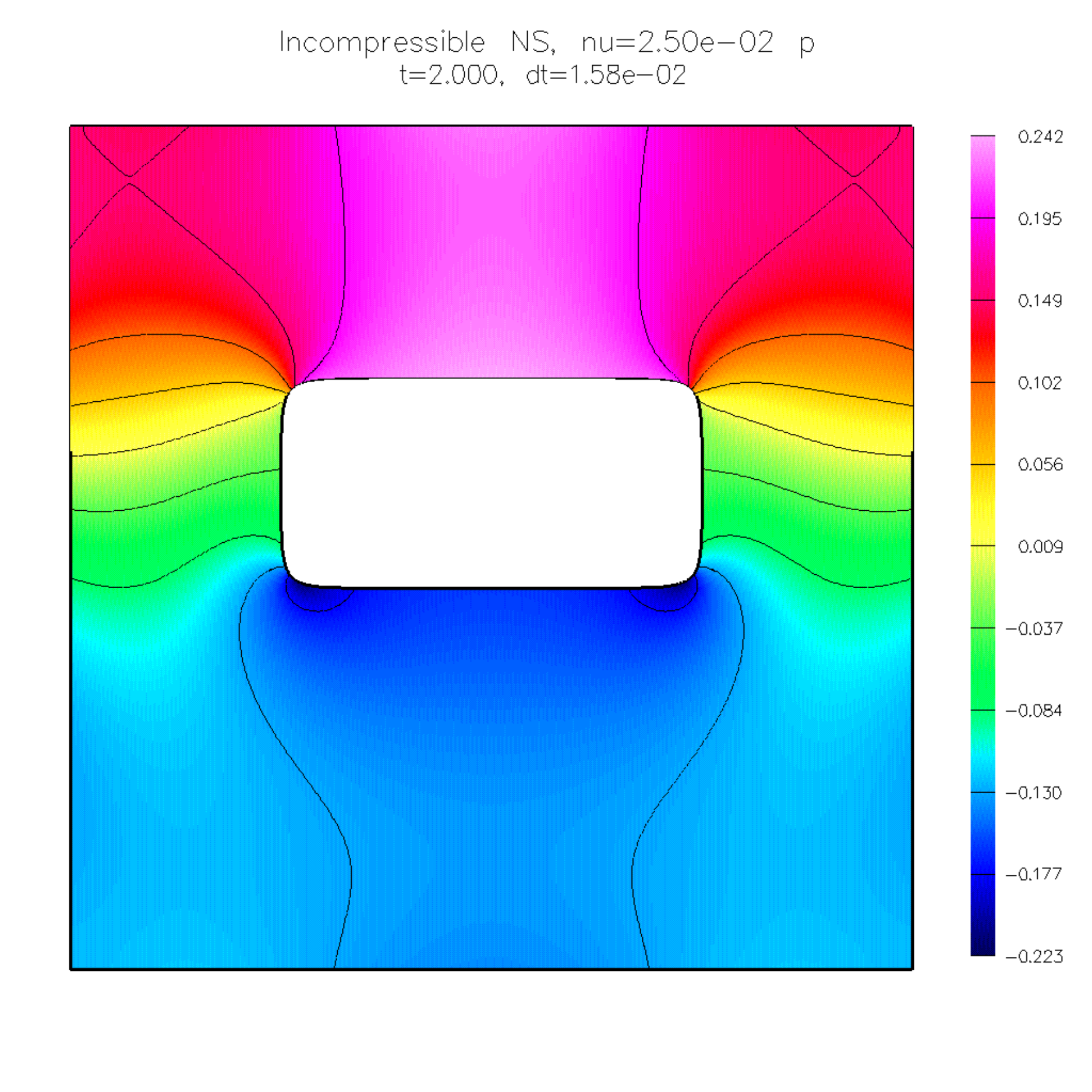}{$p$}{$p$}{$t=2.0$}{-.223}{.242};
 \drawContour{xshift=5.25cm,yshift=5.75cm}{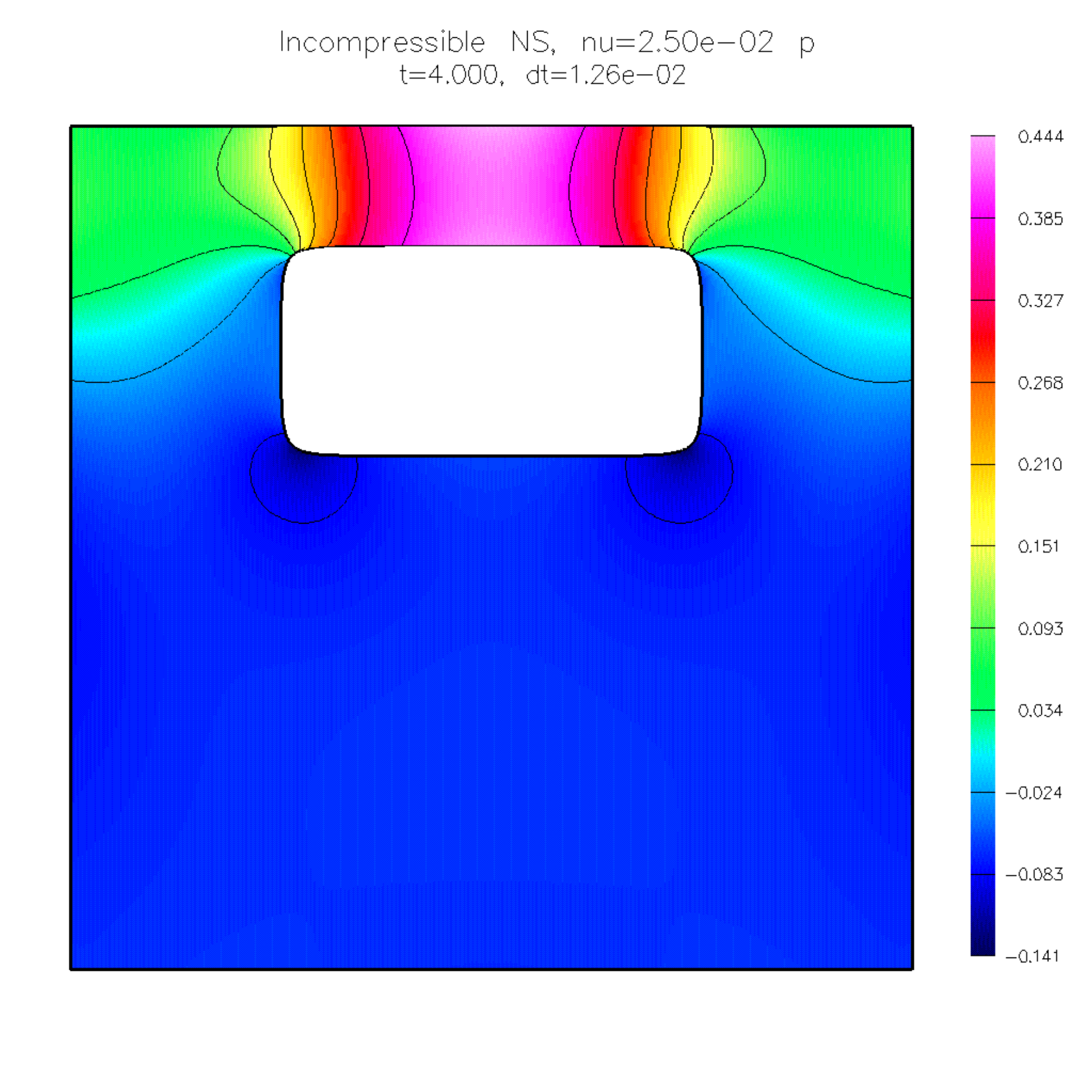}{$p$}{$p$}{$t=4.0$}{-.141}{.444};
 \drawContour{xshift=10.5cm,yshift=5.75cm}{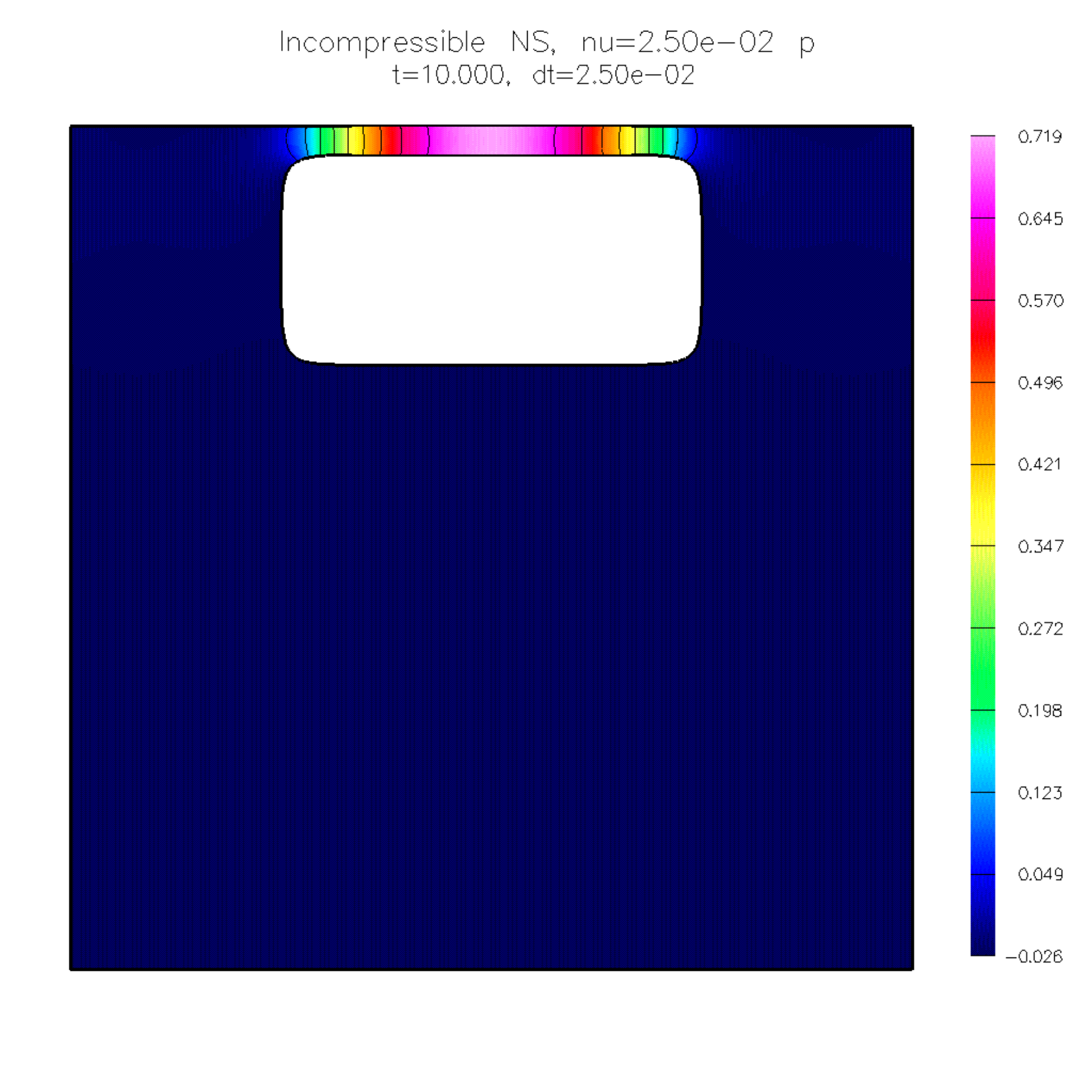}{$p$}{$p$}{$t=10.0$}{-.026}{.719};
  \draw(0.0, 0) node[anchor=south west,xshift=-4pt,yshift=+0pt] {\trimfig{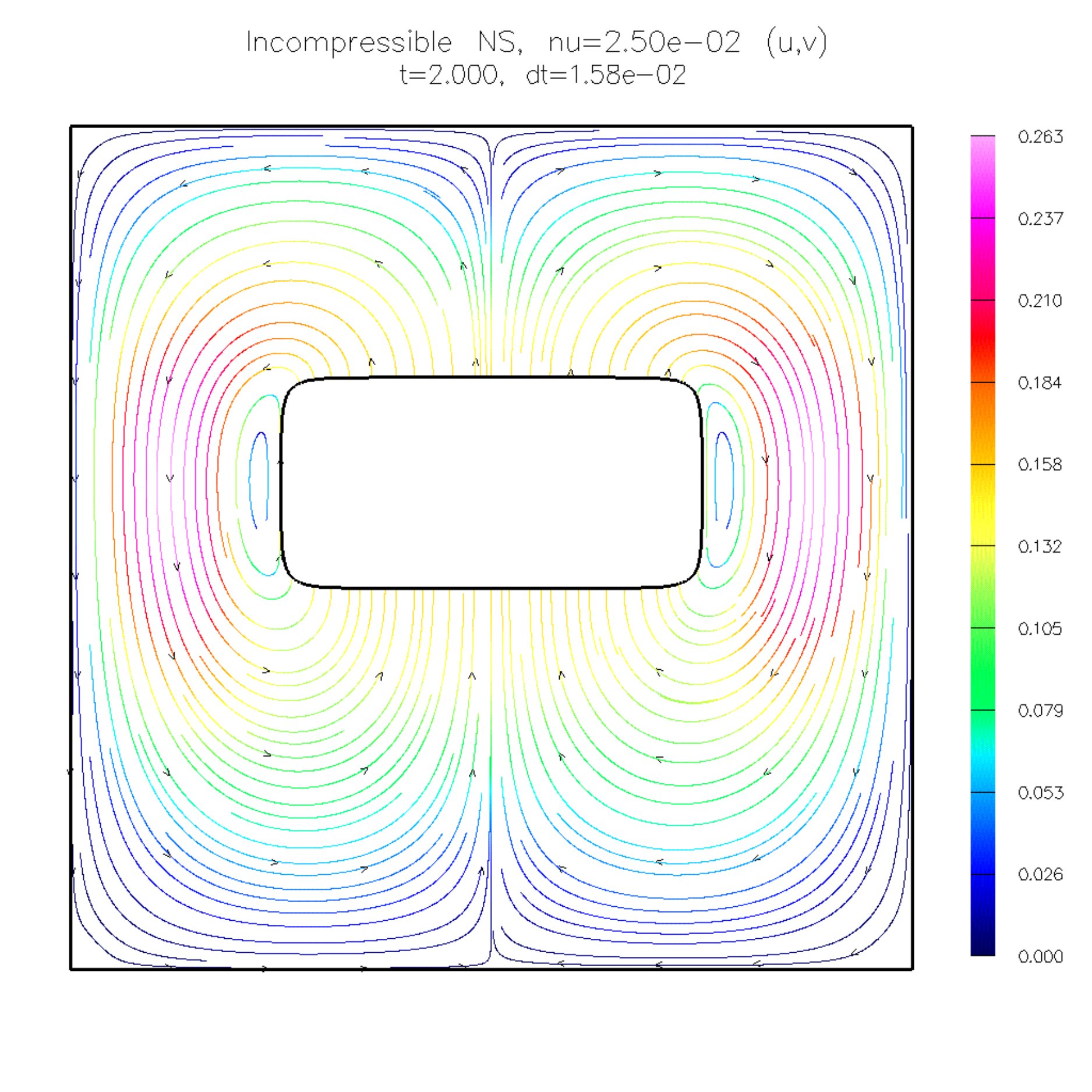}{\figWidth}};
  \draw(5.25, 0) node[anchor=south west,xshift=-4pt,yshift=+0pt] {\trimfig{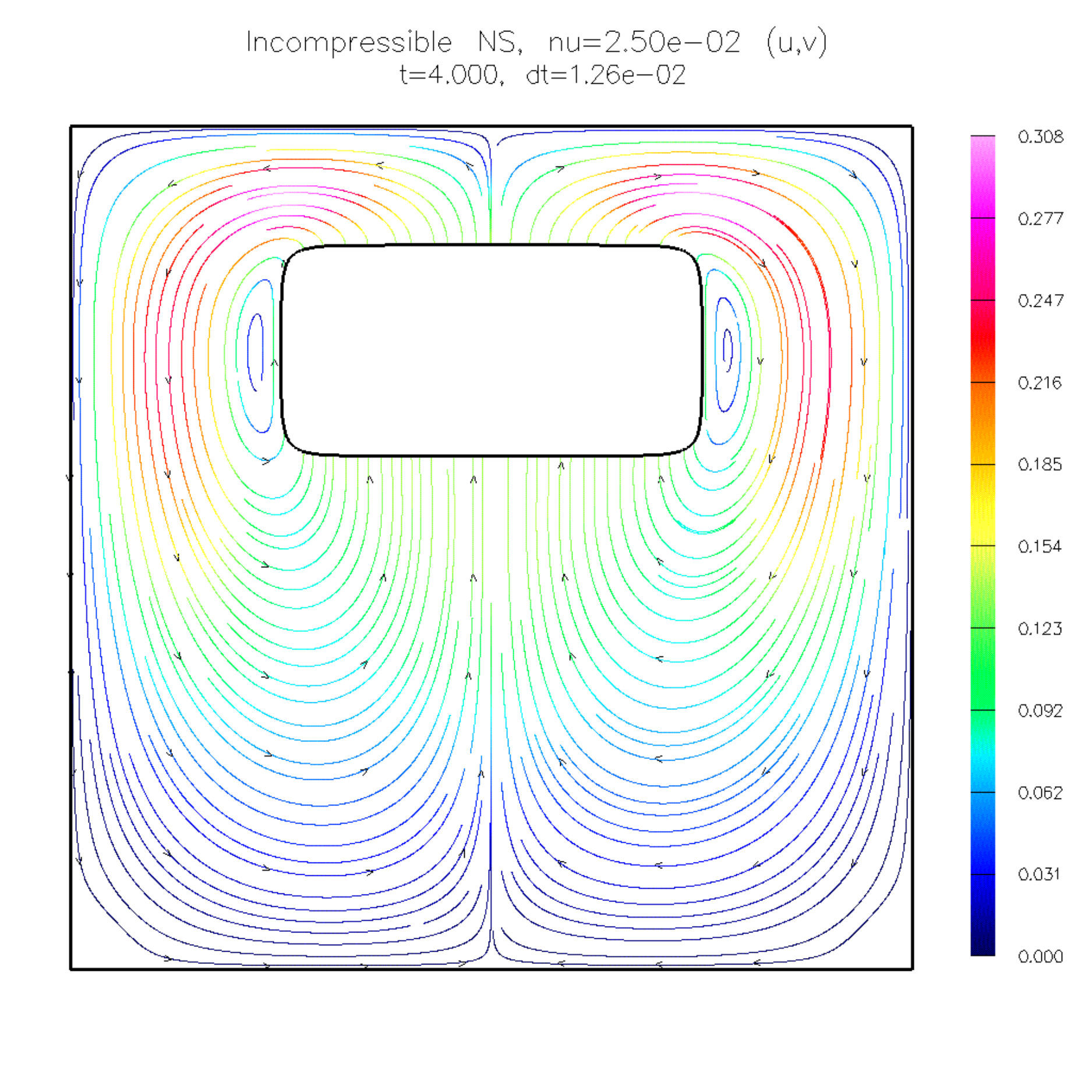}{\figWidth}};
  \draw(10.5, 0) node[anchor=south west,xshift=-4pt,yshift=+0pt] {\trimfig{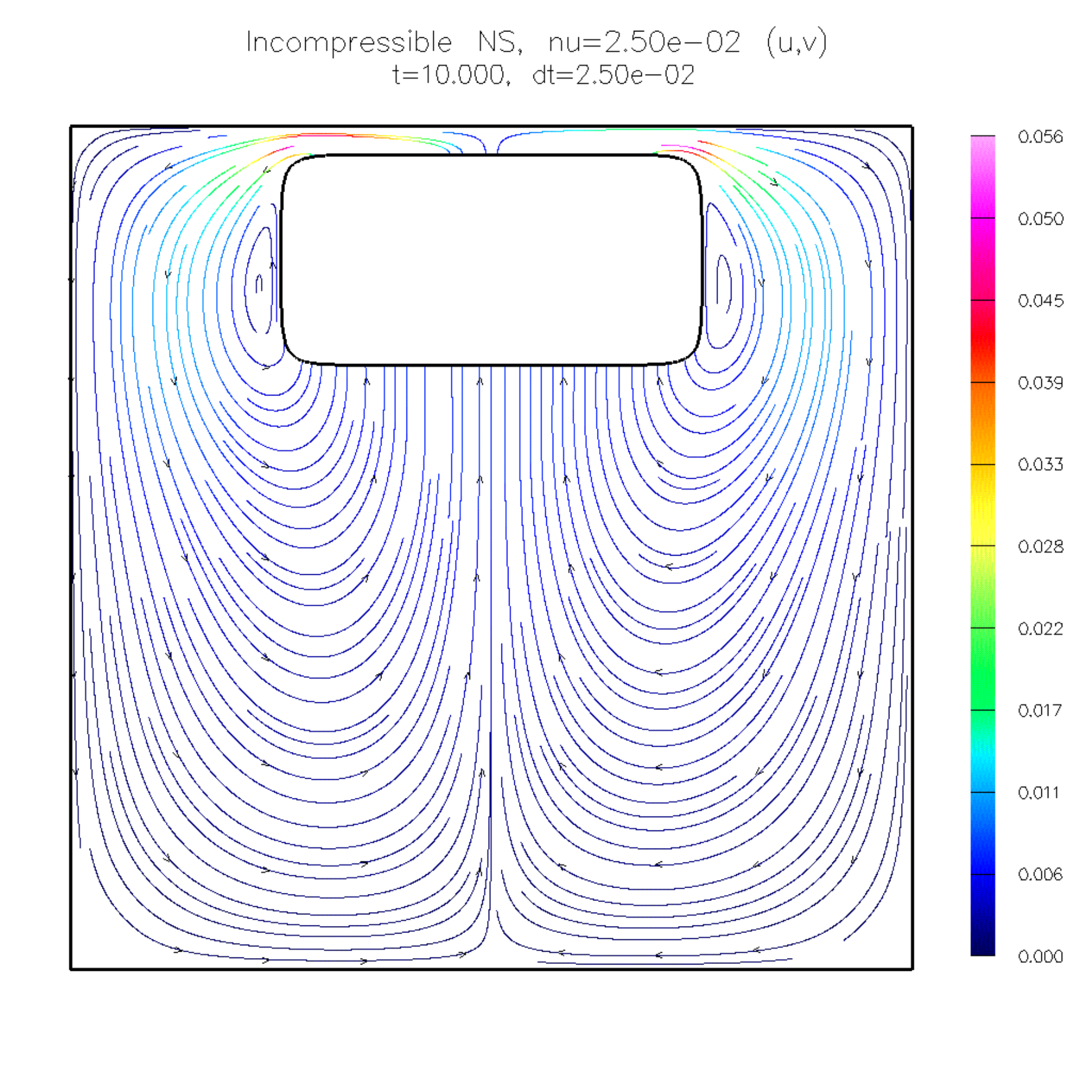}{\figWidth}};
%
\end{tikzpicture}
}
\end{center}
  \caption{Rising body: Contours of the pressure (top row) and instantaneous streamlines (bottom row) for $\rhob=0.001$ at times $t=2$, $4$ and~$10$ computed
    using the composite grid $\Gcrb^{(16)}$. }
  \label{fig:risingBodyContours}
\end{figure}
}

Figure~\ref{fig:risingBodyContours} shows contours
of the pressure and (instantaneous) streamlines from the computed
solution using the composite grid $\Gcrb^{(16)}$ for three different times as the solid rises.
The time history of the vertical position, $y_b(t)$, and vertical 
velocity, $v_b(t)$, 
Figure~\ref{fig:risingBodyCurvesLight}.  The rectangular block initially accelerates upwards due to
buoyant forces before slowing down as it approaches the top wall. At later times, the solid
continues to approach the wall slowly as fluid is expelled from the gap.  By $t=10$ the behaviour of
the fluid in the gap between the solid and the top wall is similar to that of classical lubrication
theory with the pressure showing negligible variation in the vertical direction.  

The behaviours of $y_b(t)$ and $v_b(t)$ are presented from simulations on grids, $\Gcrb^{(j)}$, of
increasing resolution, $j=4$, 8, 16 and 32. 
The self-convergence rates for this problem are not as clean as for the previous cases we have
considered for a variety of reasons. One reason is that the time-step changes discontinuously at
certain times which likely causes a local first-order error since the added-damping coefficients
depend on $\dt$. Another reason is that the boundary-fitted grids we have constructed become
narrower as the grid is refined, to allow a closer approach between the body and wall, but this
usually leads to poor estimates for the convergence rates since the region covered by a given grid changes
as the grid is refined.  However despite these issues, and the
difficulty of the problem, the curves for different grid resolutions appear to be converging reasonably well.
To give a quantitative sense of the convergence rate, a Richardson extrapolation estimate for the time-averaged convergence rate is
computed (as in Section~\ref{sec:cylDrop}) using the finest three simulations for the rigid body position and velocity as a function of
time. 
At the early time, $t=1$, the rates for the position and velocity are estimated as $1.7$ and $1.8$, respectively,
while for a later time, $t=5$, the corresponding rates are  $1.5$ and $2.1$ respectively.

{
\newcommand{\figWidth}{8.cm}
\newcommand{\trimfig}[2]{\trimFig{#1}{#2}{.0}{.0}{.0}{.0}}
\newcommand{\figWidthz}{4.cm}
\newcommand{\trimfigz}[2]{\trimFig{#1}{#2}{.0}{.0}{.0}{.18}}
\begin{figure}[htb]
\begin{center}
\resizebox{13cm}{!}{
\begin{tikzpicture}[scale=1]
  \useasboundingbox (0.0,1) rectangle (16.5,6.75);  
  \draw(0.0,0) node[anchor=south west,xshift=-4pt,yshift=+0pt] {\trimfig{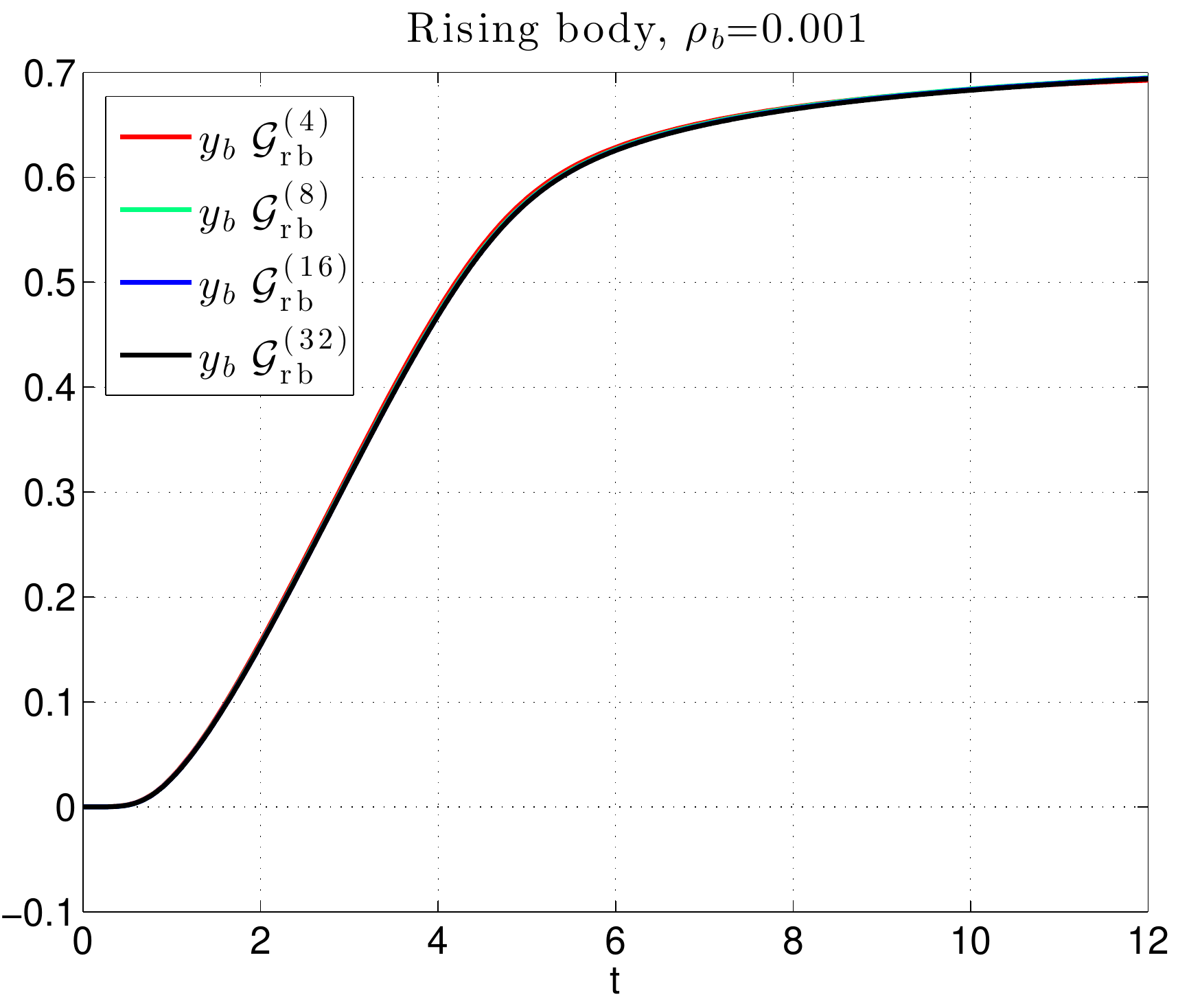}{\figWidth}};
  \draw(3.6,1.5) node[anchor=south west,xshift=-4pt,yshift=+0pt] {\trimfigz{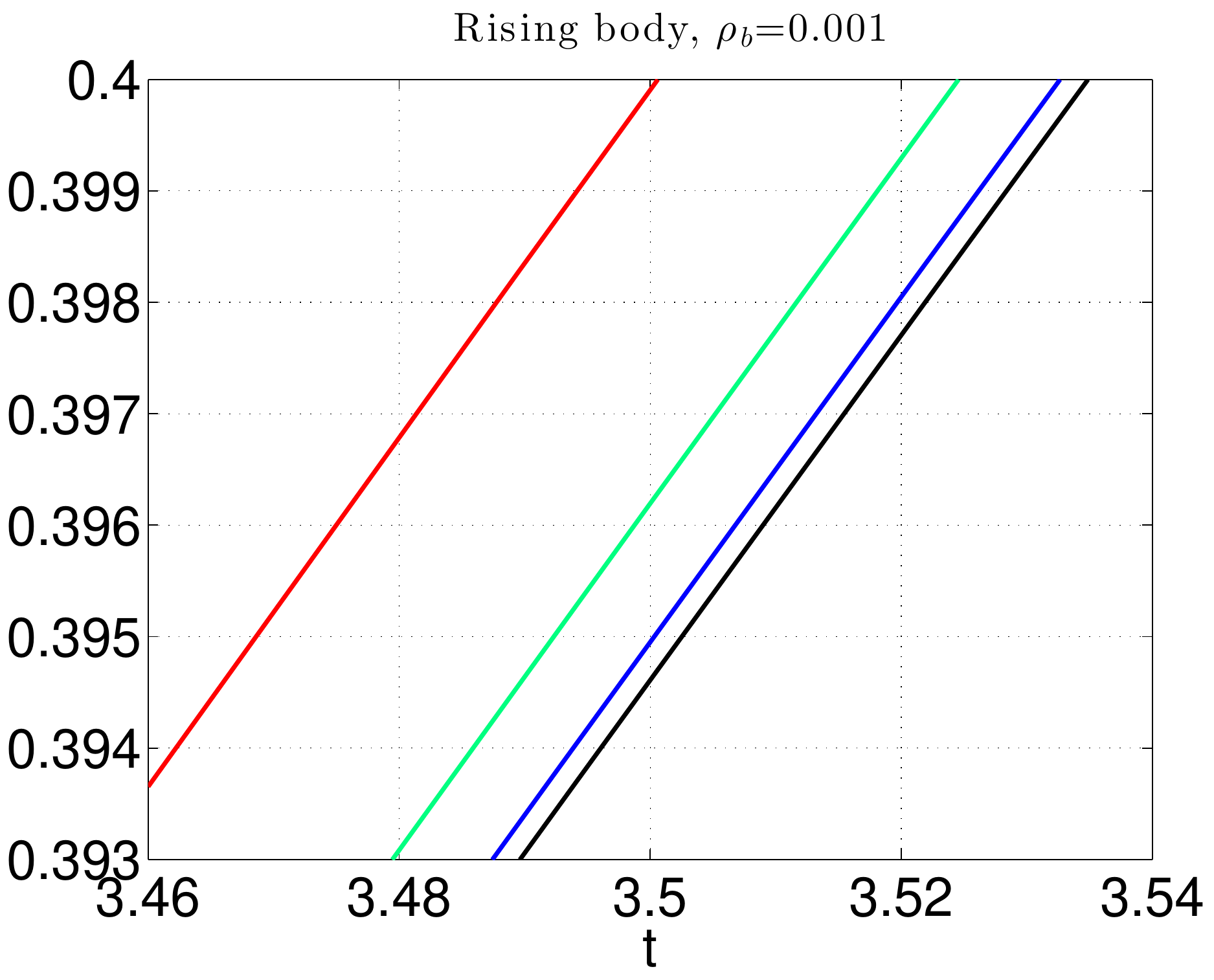}{\figWidthz}};
  \begin{scope}[xshift=.5cm]
  \draw(8.0,0) node[anchor=south west,xshift=-4pt,yshift=+0pt] {\trimfig{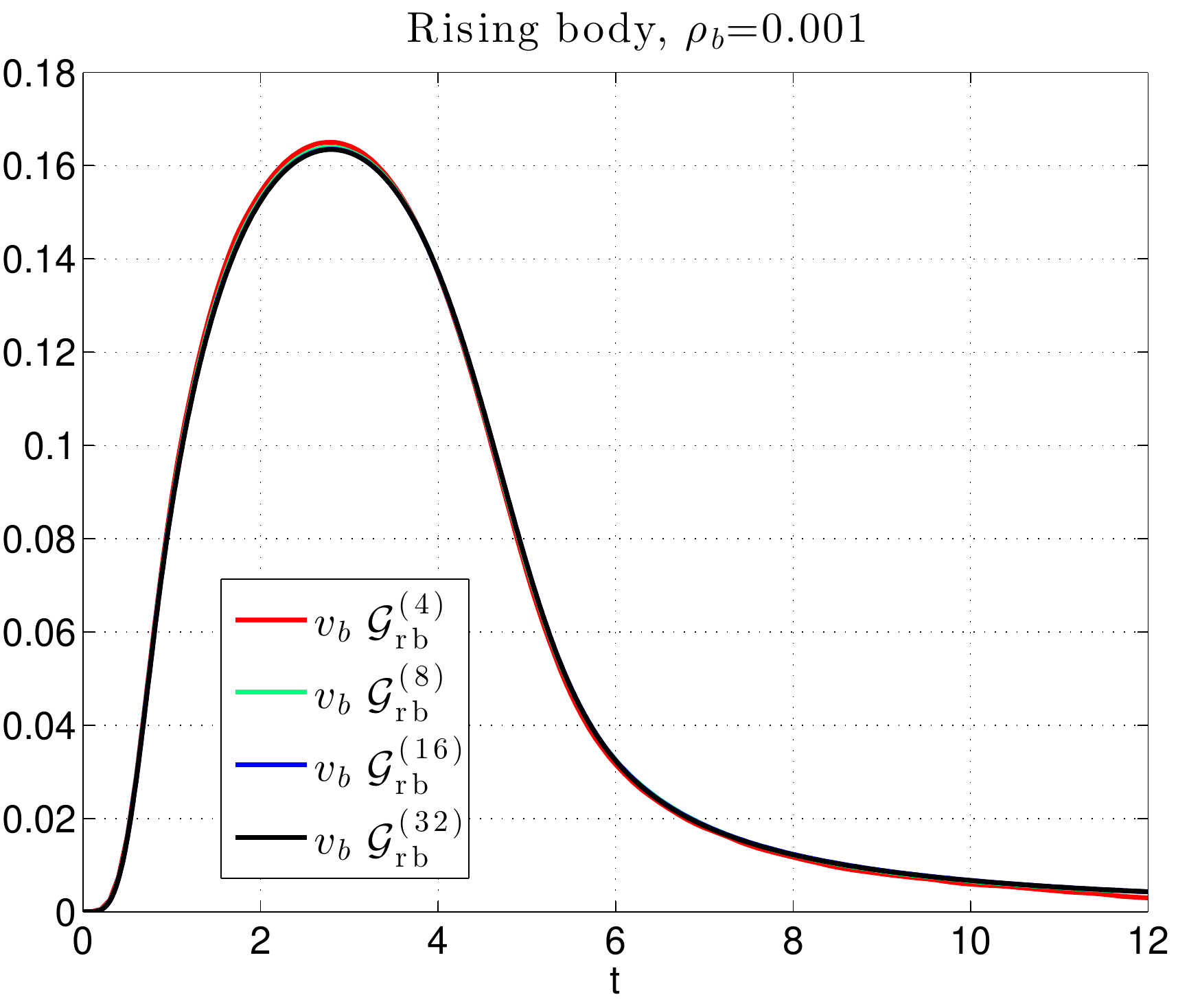}{\figWidth}};
  \draw(11.75,2.5) node[anchor=south west,xshift=-4pt,yshift=+0pt] {\trimfigz{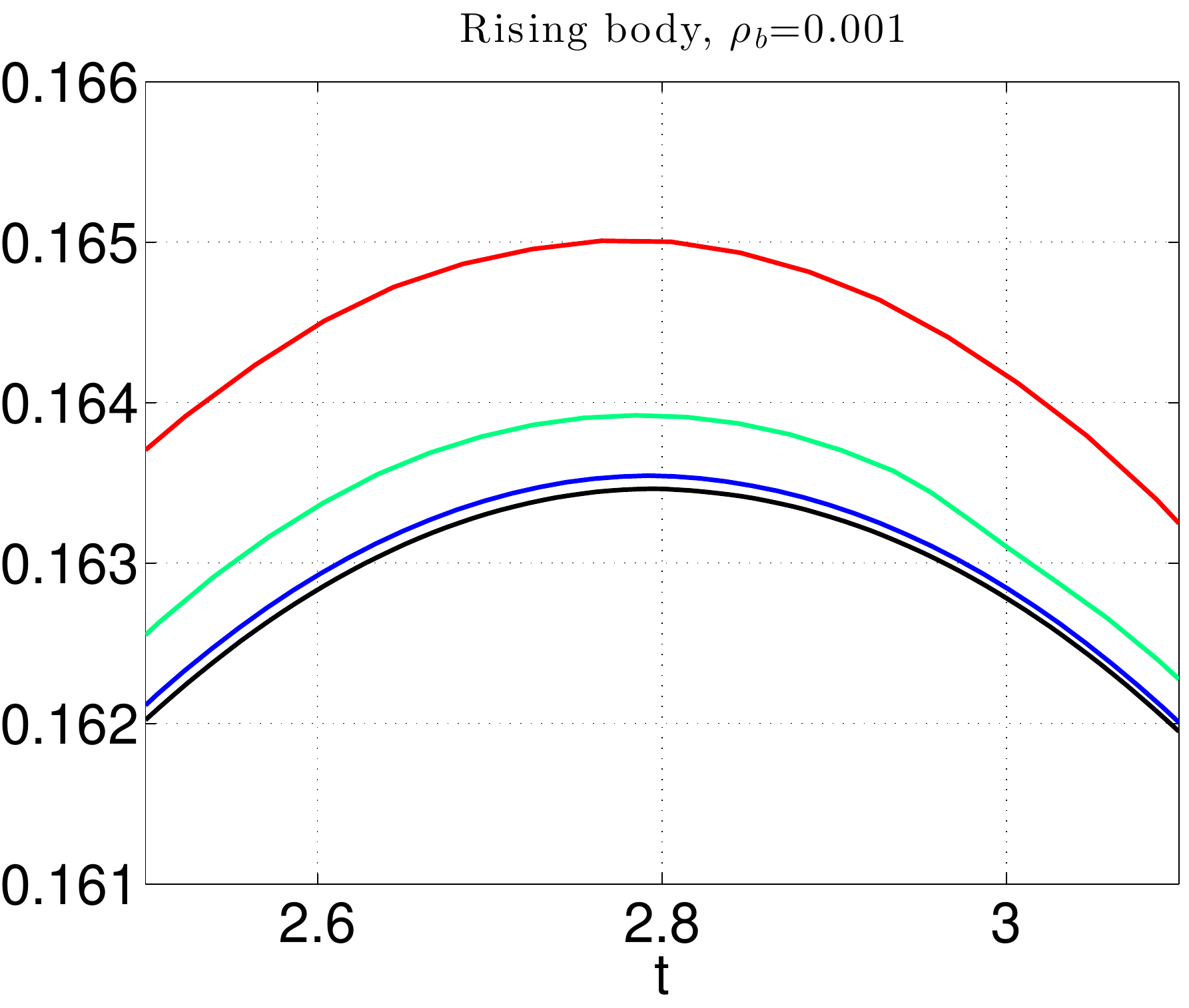}{\figWidthz}};
  \end{scope}
%
\end{tikzpicture}
} 
\end{center}
  \caption{Light rising body: Vertical body position, $y_b(t)$ (left) and vertical velocity of the body, $v_b(t)$ (right), versus time.  Results are shown from calculations using the composite grid, $\Gcrb^{(j)}$, with $j=4$, 8, 16 and~32.  The zoomed views indicate grid convergence.
 }
  \label{fig:risingBodyCurvesLight}
\end{figure}
}

Figure~\ref{fig:risingBodyCompareAcceleration} shows the behaviour of the vertical acceleration of the rigid-body
versus time for different grid resolutions. The acceleration
of the body is not as smooth as the vertical position or velocity, but this is not unexpected. A closer look at the
results for the coarsest grid shows small blips in the acceleration.  These blips correspond to small 
perturbations that are traced to changes
in the overlapping grid connectivity as the body-fitted component grid attached to the solid moves.
The motion creates changes in the overlap as new interpolation points are
computed and as the classification of grid points change between active and inactive (see~\cite{mog2006} for a detailed discussion).
These changes in the grid result in small perturbations to the pressure and velocity fields, which then are reflected
in the fluid forces on the rigid body. The scheme is quite robust to these perturbations which 
are hardly noticeable in the integrated quantities such as the body velocity and position.
In addition, the size of the perturbations is seen to become smaller as the grid is refined. 
In fact, the time averaged convergence rates using Richardson extrapolation for 
the acceleration are $1.7$ and $1.6$ at $t=1$ and $t=5$, respectively.

{
\newcommand{\figWidth}{8.cm}
\newcommand{\trimfig}[2]{\trimFig{#1}{#2}{.0}{.0}{.0}{.0}}
\newcommand{\figWidthz}{4.cm}
\newcommand{\trimfigz}[2]{\trimFig{#1}{#2}{.0}{.0}{.0}{.18}}\begin{figure}[htb]
\begin{center}
\resizebox{14cm}{!}{
\begin{tikzpicture}[scale=1]
  \useasboundingbox (0.5,1) rectangle (16.,6.75);  
  \draw(0,.0) node[anchor=south west,xshift=-4pt,yshift=+0pt] {\trimfig{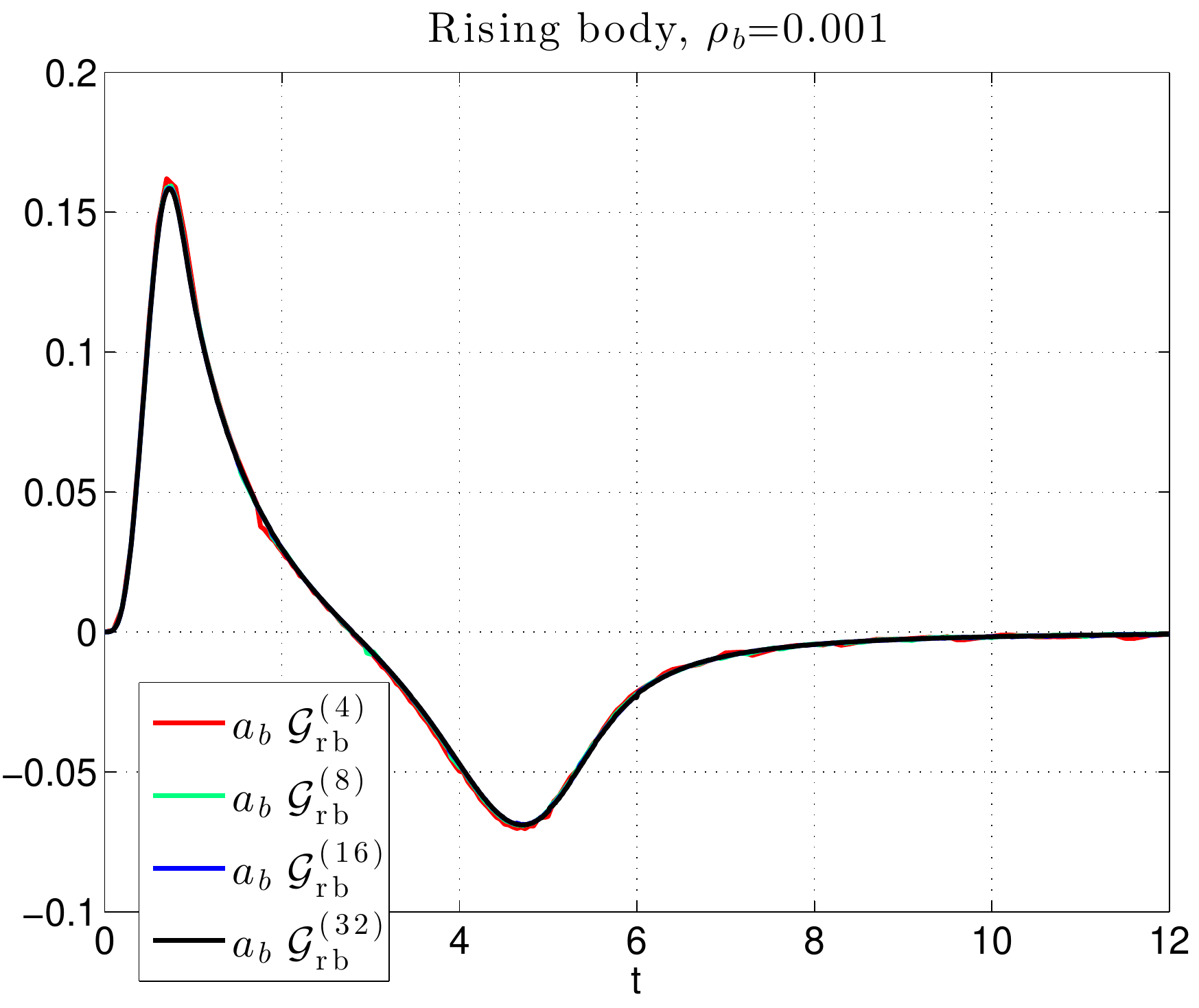}{\figWidth}};
  \draw(2.75,2.75) node[anchor=south west,xshift=-4pt,yshift=+0pt] {\trimfigz{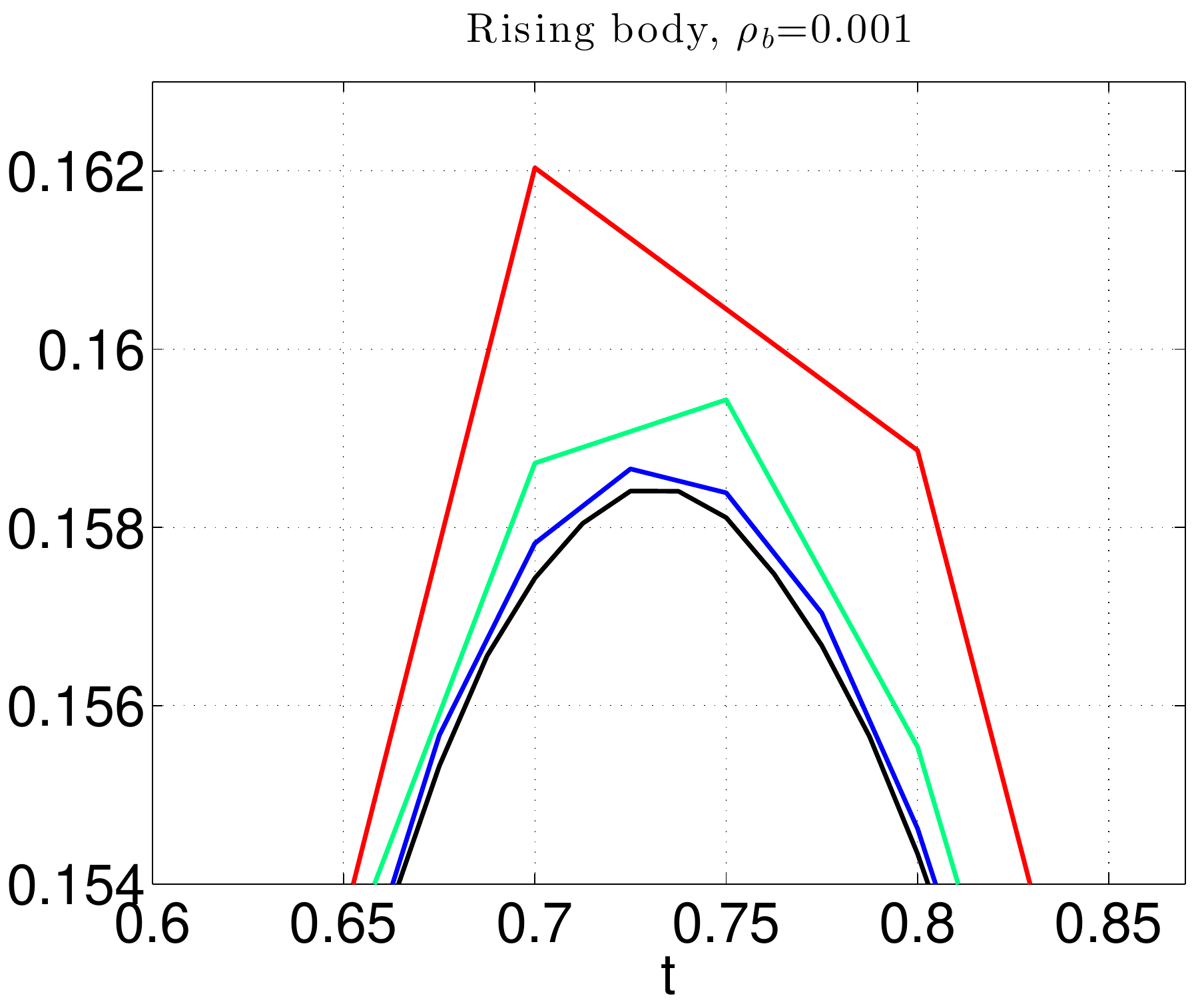}{\figWidthz}};
  \draw(8.5,0) node[anchor=south west,xshift=-4pt,yshift=+0pt] {\trimfig{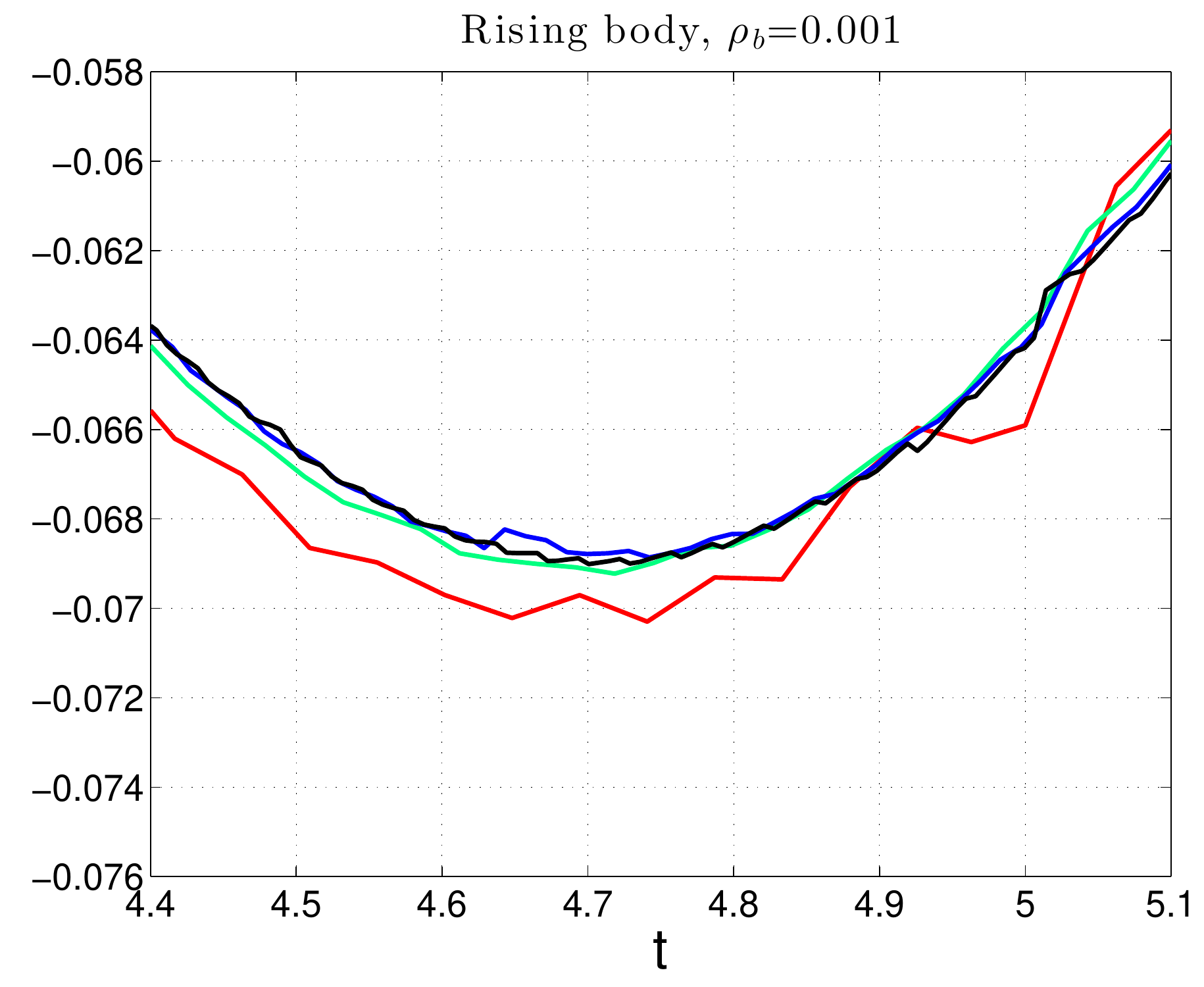}{\figWidth}};
%
\end{tikzpicture}
}
\end{center}
  \caption{Light rising body. Left: vertical acceleration, $a_b(t)$, of the light rising body. Right: zoom near $t=4.75$ showing the perturbations to the
    accelerations caused by the changes in overlapping grid interpolation points as the component grid fitted to the body moves. 
 }
  \label{fig:risingBodyCompareAcceleration}
\end{figure}
}

It is interesting to discover that instabilities
due to added-damping effects can be important for this problem even though the dominant motion is translational.  The added-damping terms in the AMP interface condition~\eqref{eq:AMPDPrigidBodyAcceleration} involve the factor $\adc\dt$ and the added-damping tensor~$\ADtensor$.  The components of this tensor are computed from discrete approximations to the surface integrals in~\eqref{eq:ADvv}--\eqref{eq:ADww}
using the discrete geometry defining the rounded rectangular solid. 
For the grid $\Gcrb^{(4)}$, for example, the main contributions to~$\ADtensor$ are found to be
\[
\Dvva_{h,11}\approx 1.917\,\frac{\mu}{\dn} , \qquad \Dvva_{h,22}\approx 0.917\,\frac{\mu}{\dn},
   \qquad \Dwwa_{h,22}\approx 0.379\,\frac{\mu}{\dn} ,
\]
with other entries in the added-damping tensor being approximately zero.  These contributions are in good agreement with
those obtained using the exact surface integrals for a square-cornered rectangular
solid which give
\[
  \Dvva_{e,11}=2\solidWidth\,\frac{\mu}{\dn}=2\,\frac{\mu}{\dn} ,\qquad  \Dvva_{e,22}=2\solidHeight\,\frac{\mu}{\dn}=\frac{\mu}{\dn},
   \qquad \Dwwa_{e,22}=\frac{\solidWidth\solidHeight(\solidWidth+\solidHeight)}{2}\,\frac{\mu}{\dn}=\frac{3}{8}\,\frac{\mu}{\dn},
\]
(see equation~\eqref{eq:addedDampingRectangle} in \ref{sec:addedDampingTensorExamples}).
Observe that the added-damping terms are proportional to $\mu\dt/\dn$, where $\dn$ is defined in~\eqref{eq:dn}.
For coarser grids where the time-step is larger, $\dn\approx\sqrt{\nu\dt/2}$ so that the added-damping terms scale
as $\sqrt{\nu\dt}$.  This indicates that added-damping effects increase for larger values of $\dt$, in agreement
with the analysis in Part~I.
As an example, Figure~\ref{fig:risingBodyUnstable} shows an added-damping instability for a simulation on the grid $\Gcrb^{(4)}$ when the added-damping coefficient $\adc$ is intentionally and artificially chosen too small ($\adc=0.5$ in this case).
 The instability
appears primarily as a rotational oscillation of the body while it rises.  The time history of the angular acceleration, $\dot\omega_b$,
shows this high-frequency oscillation.  It is interesting to note that the amplitude of this instability saturates due to
counter-acting pressure forces caused by added-mass effects.

{
\newcommand{\figWidth}{8.cm}
\newcommand{\figWidthA}{12.cm}
\newcommand{\trimfig}[2]{\trimFig{#1}{#2}{.0}{.0}{.0}{.0}}
\newcommand{\figWidtha}{6cm}
\newcommand{\trimfiga}[2]{\trimFig{#1}{#2}{.1}{.3}{.2}{.2}}
\begin{figure}[htb]
\begin{center}
\resizebox{13cm}{!}{
\begin{tikzpicture}[scale=1]
  \useasboundingbox (0.0,1) rectangle (15.,6.75);  
  \draw(0.0,0) node[anchor=south west,xshift=-4pt,yshift=+0pt] {\trimfig{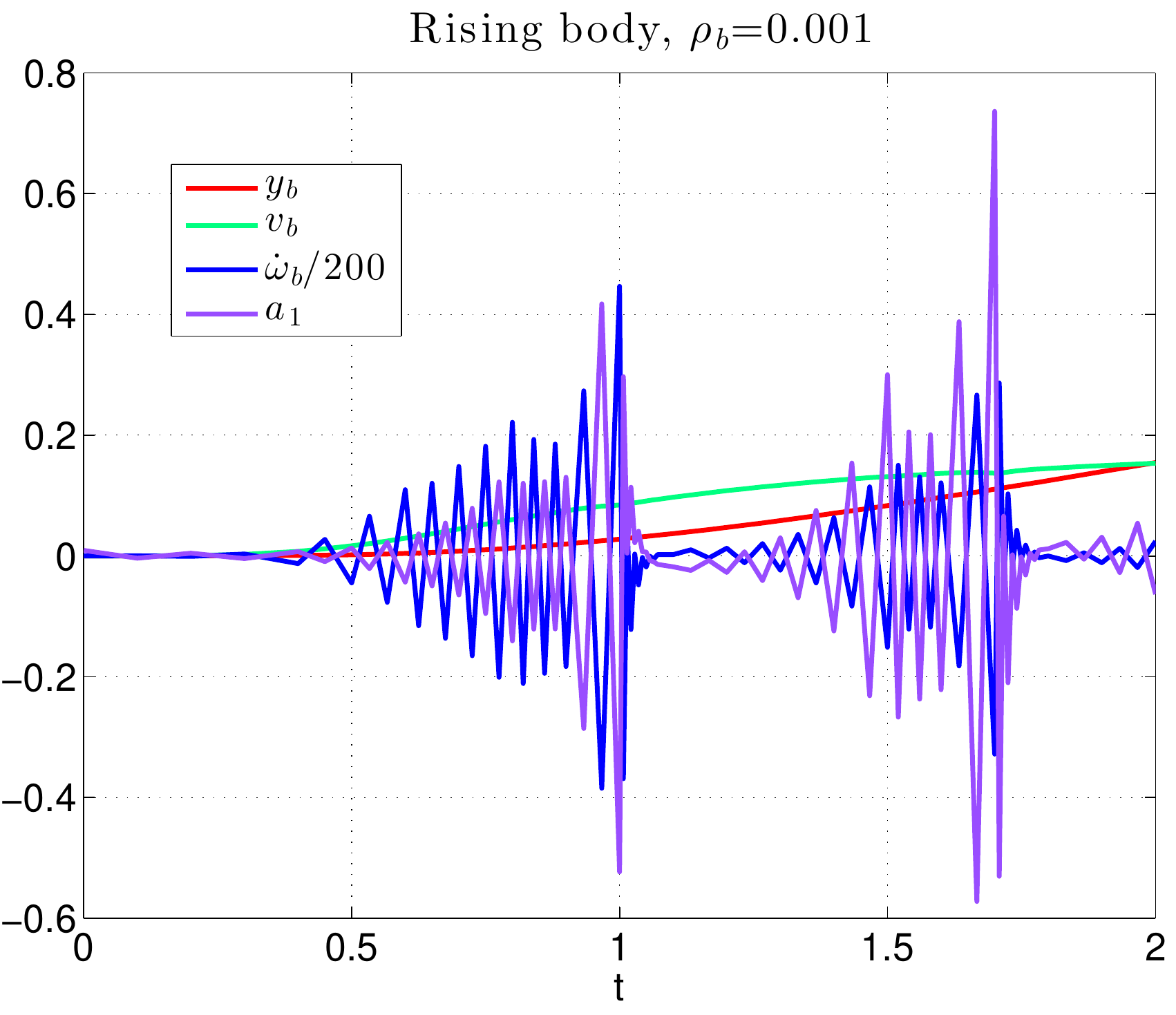}{\figWidth}};
  \draw(8.5,.35) node[anchor=south west,xshift=-4pt,yshift=+0pt] {\trimfiga{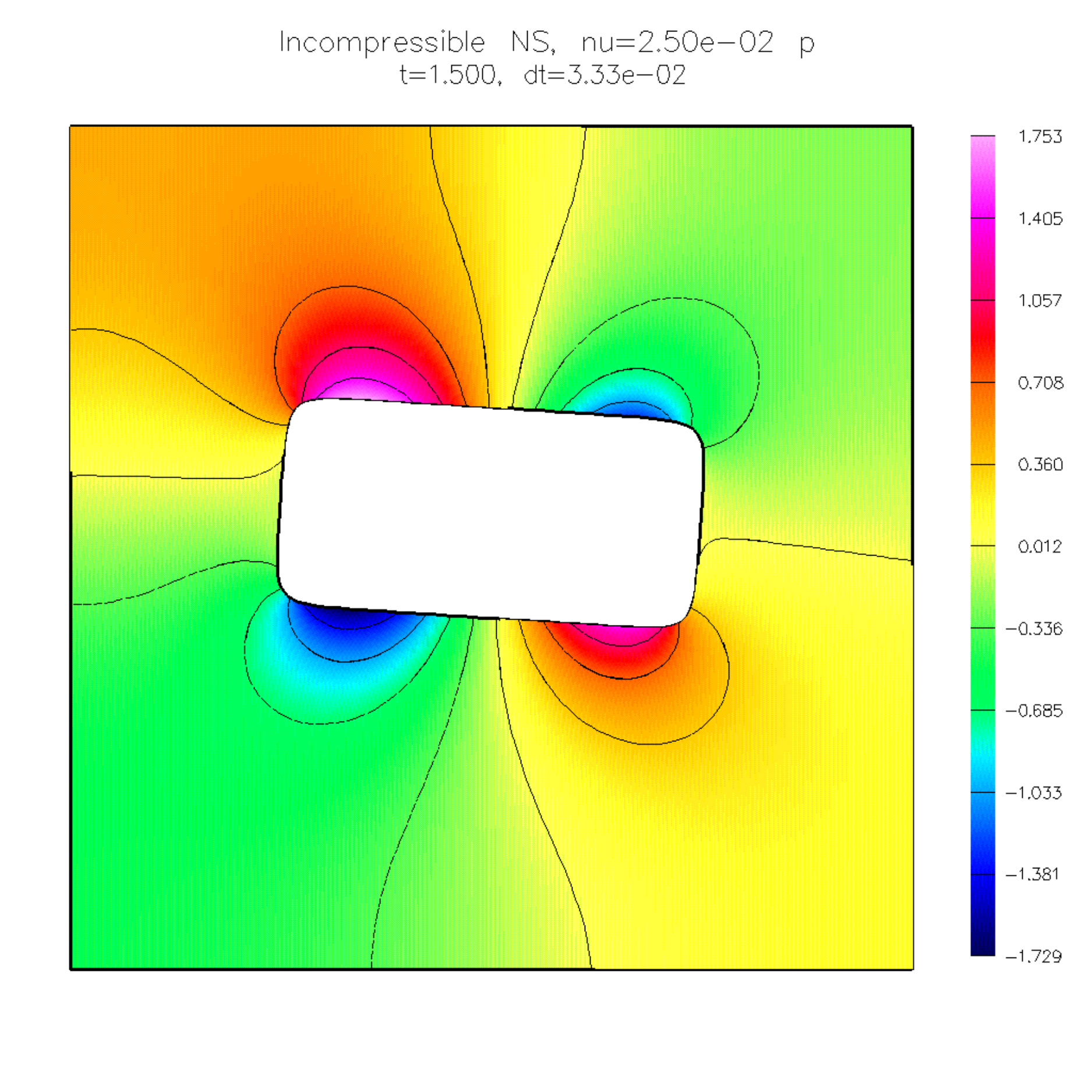}{\figWidtha}};

%
\end{tikzpicture}
}
\end{center}
  \caption{Illustration of the added-damping instability for a light rising-body when the added-damping coefficient $\adc$ is 
    intentionally chosen to be too small.
    The simulation becomes unstable without sufficient added-damping.
   The body undergoes unphysical rotations (as reflected in $\dot\omega_b$) which are partially stabilized by counter-acting
    pressure forces (pressure at $t=1.4$ is shown on right). 
    Results are shown for grid $\Gcrb^{(4)}$, with added-damping coefficient $\adc=0.5$.
 }
  \label{fig:risingBodyUnstable}
\end{figure}
}

Finally, note that the results from the traditional scheme for this problem give nearly identical 
results to those from the {\ampRB} scheme. 
The 
TP-RB scheme, however, requires many sub-time-step iterations at each time step, $85$ on average with a maximum of about 200
at early times, for a simulation using the composite grid $\Gcrb^{(4)}$.

\newcommand{\Gcrf}{\Gc_{\rm rf}}
\subsection{Interaction between rising and falling bodies in a fluid chamber} \label{sec:twoFallingBodies}

{
\newcommand{\figWidth}{5cm}
\newcommand{\trimfig}[2]{\trimh{#1}{#2}{.21}{.21}{.115}{.115}}
\begin{figure}[htb]
\begin{center}
\resizebox{14cm}{!}{
\begin{tikzpicture}[scale=1]
  \useasboundingbox (0.5,1) rectangle (16,5.25);  
  \draw(0.0,0) node[anchor=south west,xshift=-4pt,yshift=+0pt] {\trimfig{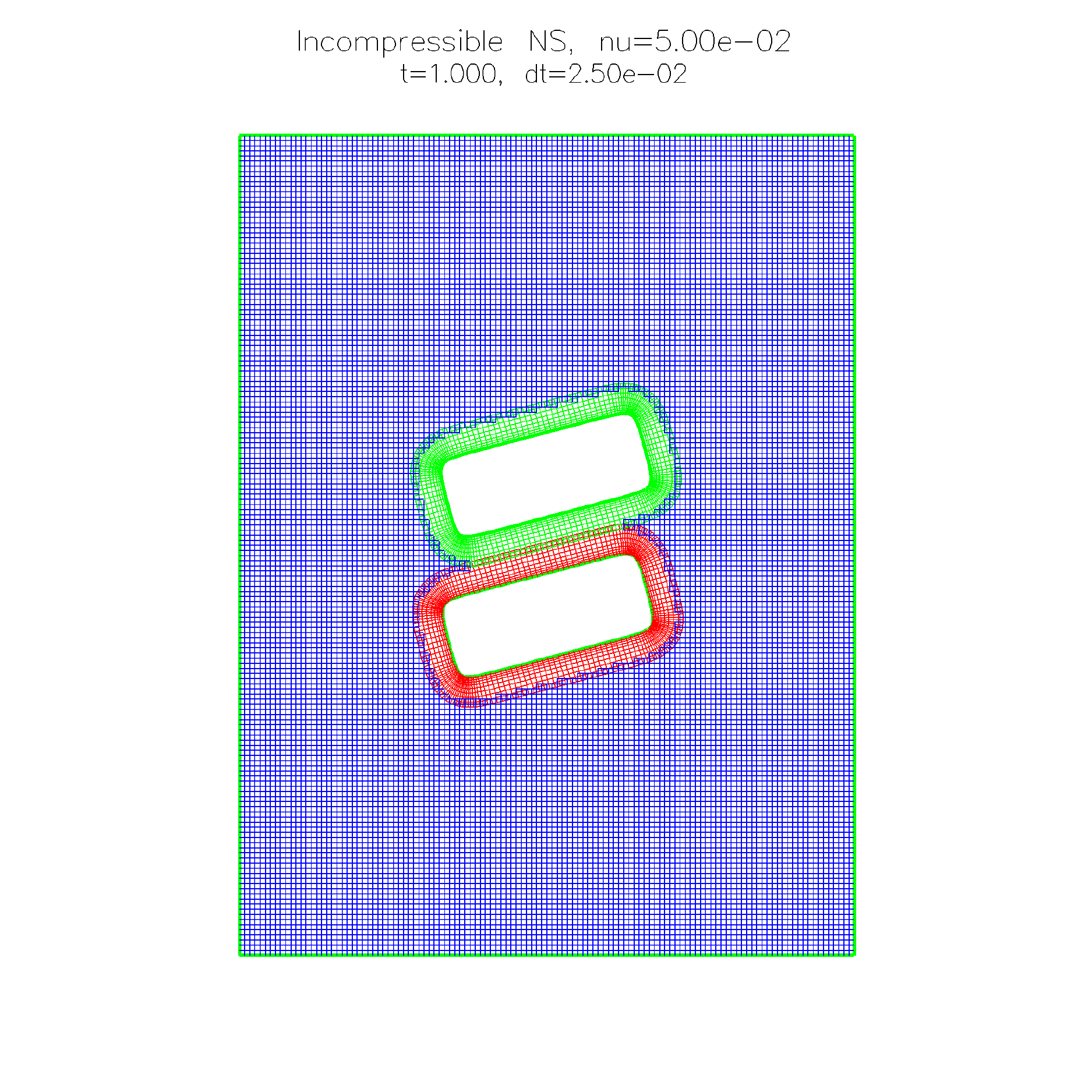}{\figWidth}};
  \draw(4.0,0) node[anchor=south west,xshift=-4pt,yshift=+0pt] {\trimfig{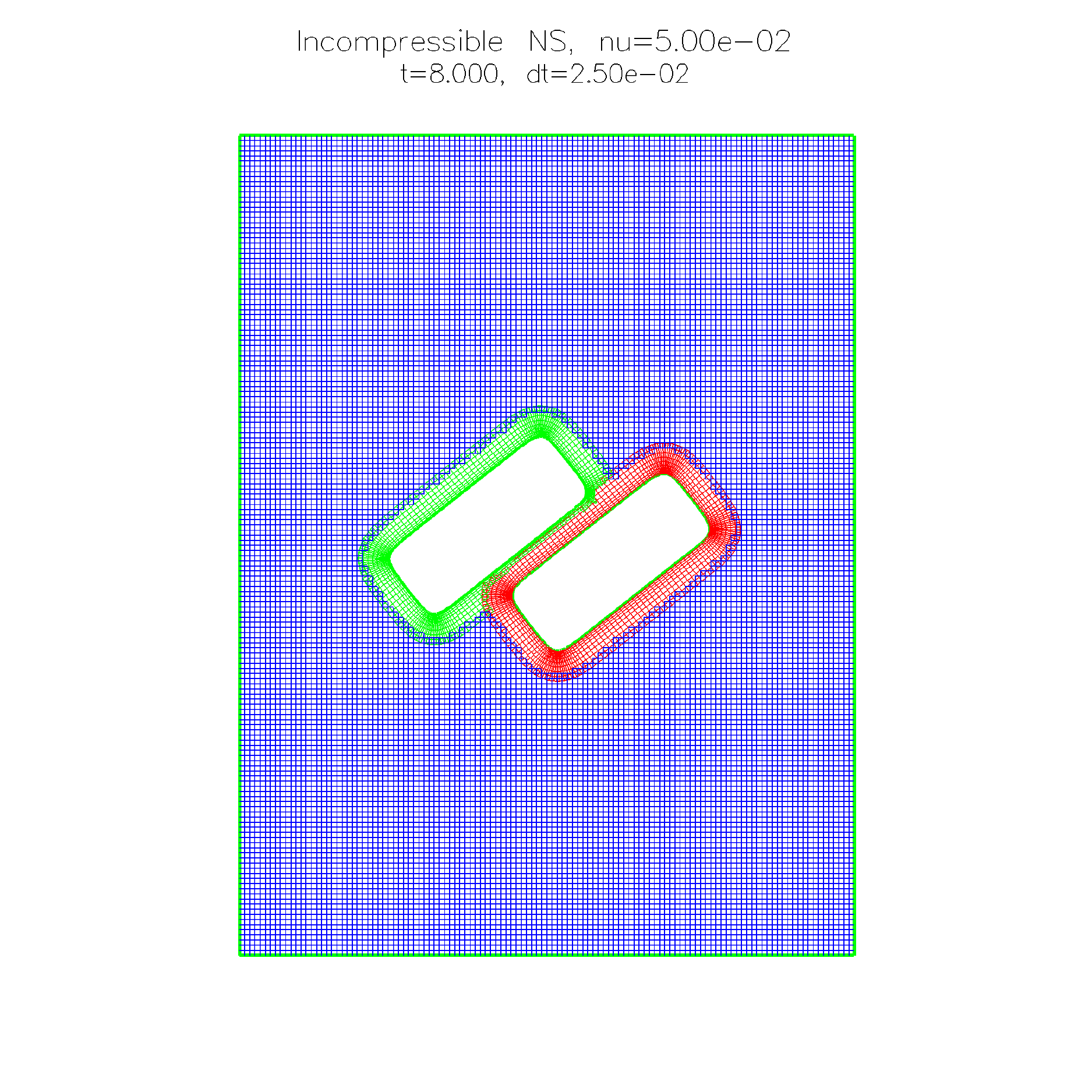}{\figWidth}};
  \draw(8.0,0) node[anchor=south west,xshift=-4pt,yshift=+0pt] {\trimfig{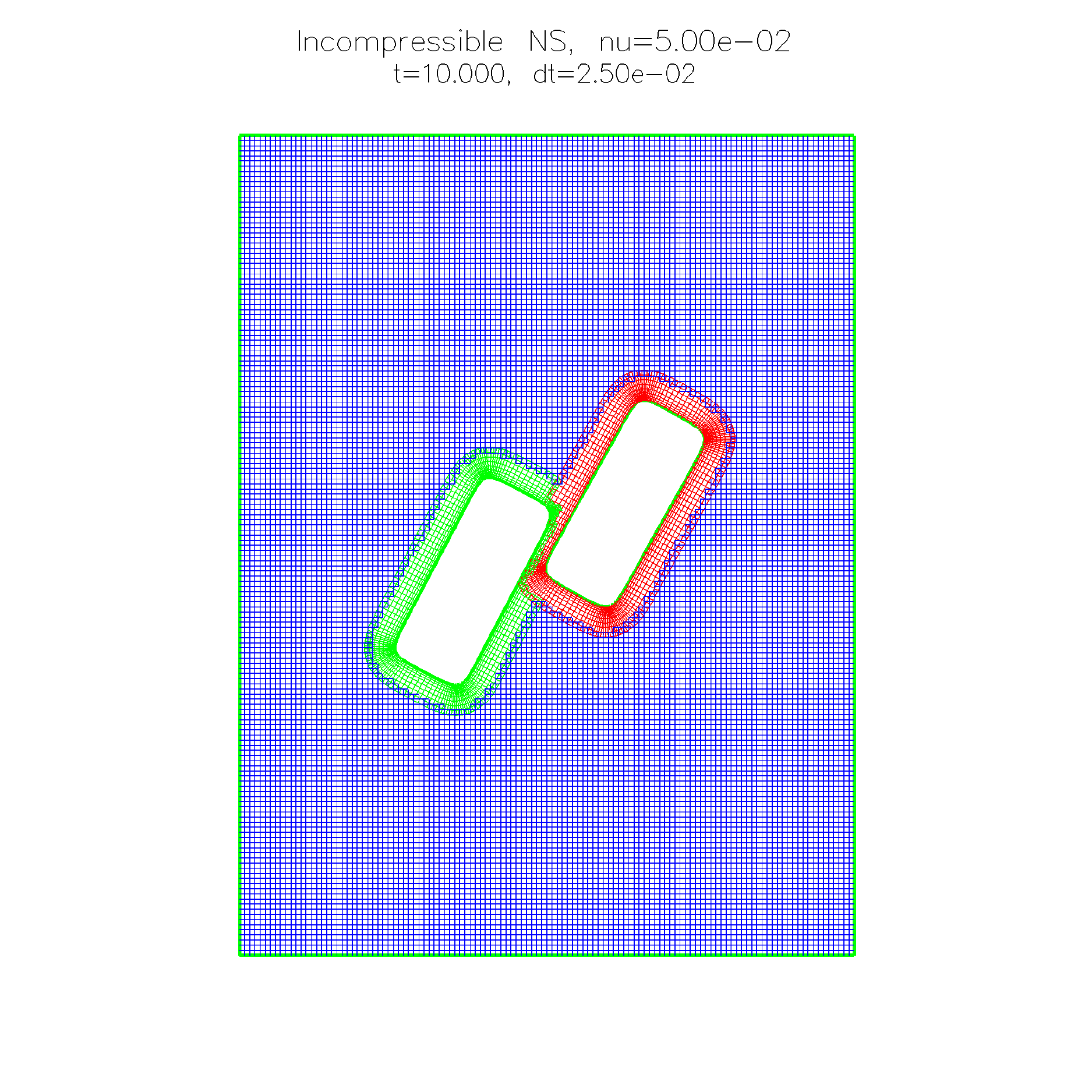}{\figWidth}};
  \draw(12.0,0) node[anchor=south west,xshift=-4pt,yshift=+0pt] {\trimfig{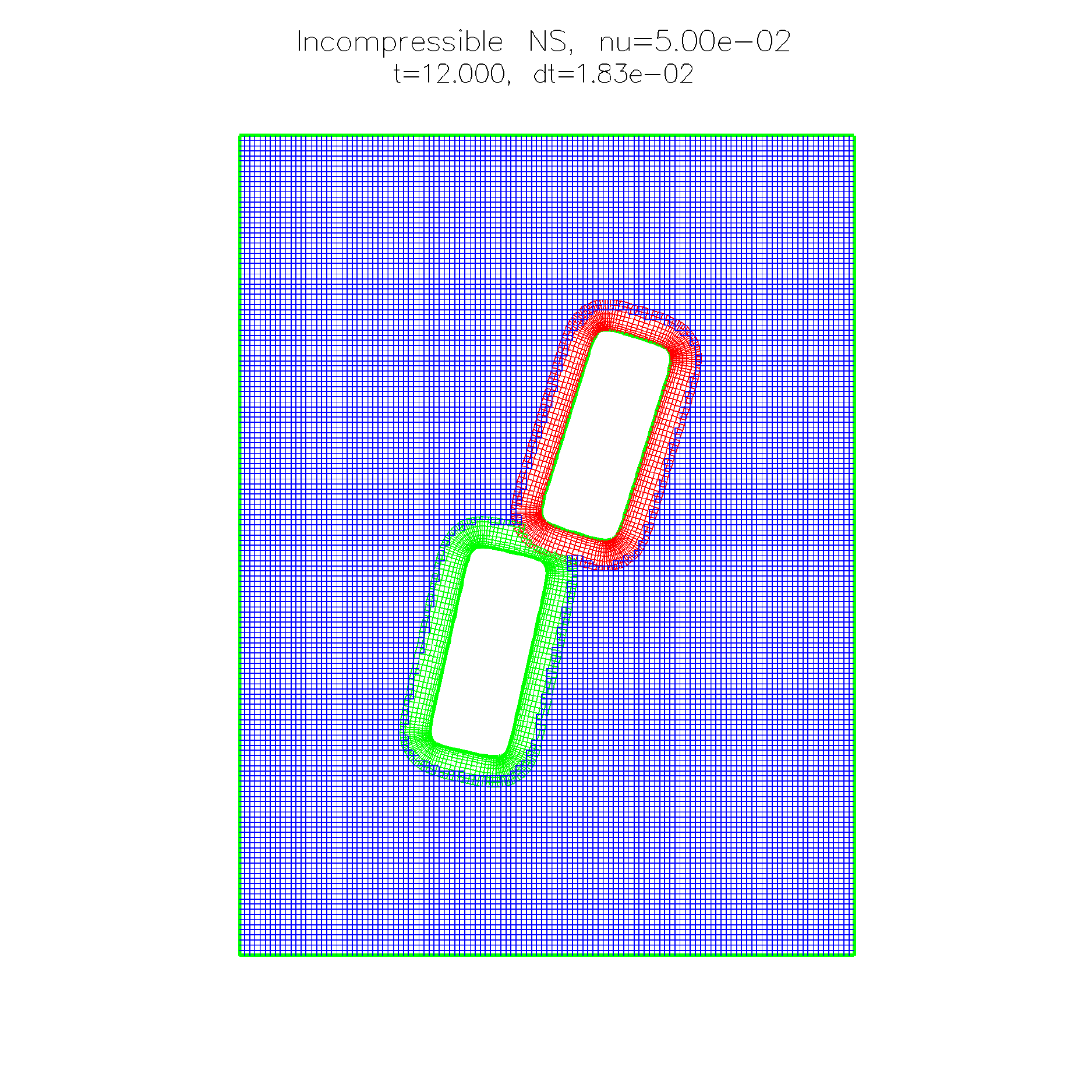}{\figWidth}};
%
\end{tikzpicture}
}
\end{center}
  \caption{
    Rising and falling bodies.  Composite grid $\Gcrf^{(4)}$ at times $t=1$, $8$, $10$ and $12$ (left to right).
}
  \label{fig:fallingBodyGrid}
\end{figure}
}

As a final illustration of the \ampRB~scheme, we consider the buoyancy-driven interaction of two rectangular-shaped bodies, initially located 
one above the other, as shown in Figure~\ref{fig:fallingBodyGrid}.
The lower and lighter body (referred henceforth as the {\em bottom} body) 
of density $\rhob=0.5$ attempts to rise in a fluid of density $\rho=1$, 
but is impeded initially by the heavier body above it (referred henceforth as the {\em top} body)
of density $\rhob=1.5$ that is attempting to fall. 
This example demonstrates the properties of the~\ampRB~scheme for the case of multiple translating
and rotating bodies.  The rigid bodies are both considered to be ``light'' relative the fluid density
in that added-mass and added-damping effects are important, and these effects vary over time
due to the thin gap between the solids at early times and their relative motion.


{
\newcommand{\drawContour}[7]{%
\begin{scope}[#1]
\draw(0.0,0) node[anchor=south west,xshift=-4pt,yshift=+0pt] {\trimfig{fig/#2}{\figWidth}};
\draw(.5,5.75) node[draw,fill=white,anchor=west,xshift=2pt,yshift=0pt] {\scriptsize #3};
\draw(1.,5.75) node[draw,fill=white,anchor=west,xshift=2pt,yshift=0pt] {\scriptsize #5};
\begin{scope}[xshift=2pt]
  \draw (\xcb,\ycb) node[anchor=south west,xshift=0.cm,yshift=.5cm,rotate=-90] {\trimfigcb{fig/colourBarLines}{\cbWidth}{\cbHeight}};
  \draw (.8,0) node[anchor=north,xshift=+3pt,yshift=+2pt] {\scriptsize $#6$};
  \draw (4.4,0) node[anchor=north,xshift=+0pt,yshift=+2pt] {\scriptsize $#7$};
\end{scope}
\end{scope}
}
\newcommand{\cbWidth}{.2cm}
\newcommand{\cbHeight}{4cm}
\newcommand{\xcb}{.5cm}
\newcommand{\ycb}{-.2cm}
\setlength{\ycbTop}{\ycb+\cbHeight}
\setlength{\ycbMid}{\ycb+\cbHeight*\real{.5}}
\newcommand{\trimfigcb}[3]{\includegraphics[width=#2, height=#3, clip, trim=17cm 2.35cm 1.65cm 2.35cm]{#1}}
%
\newcommand{\figWidth}{4.5cm}
\newcommand{\yb}{6.75}
\newcommand{\trimfig}[2]{\trimw{#1}{#2}{.15}{.25}{.1}{.1}}
\begin{figure}[htb]
\begin{center}
\resizebox{11cm}{!}{
\begin{tikzpicture}[scale=1]
  \useasboundingbox (0.0,.75) rectangle (16.,12.75);  
 \drawContour{xshift= 0.0cm,yshift=6.75cm}{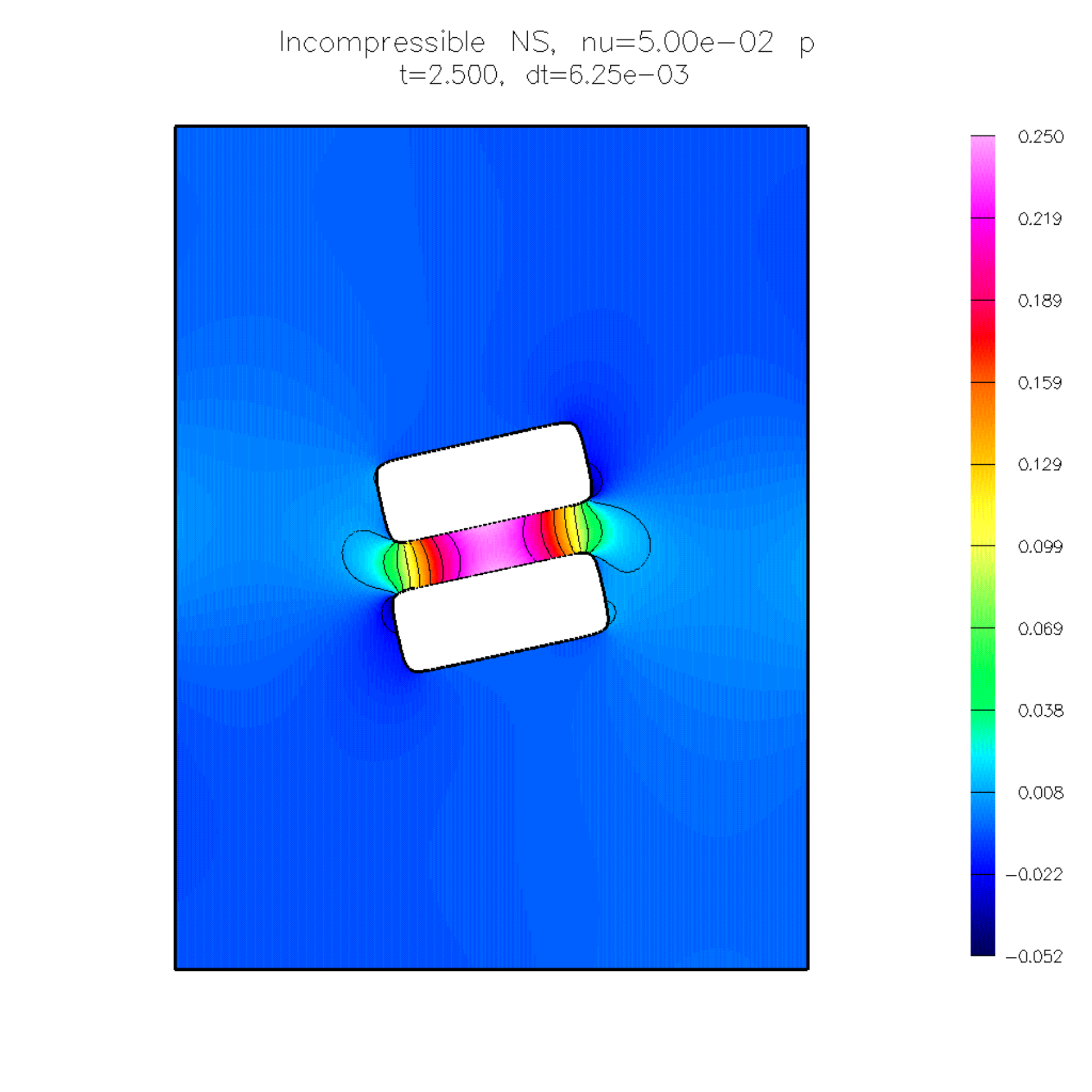}{$p$}{$p$}{$t=2.5$}{-.052}{.250};
 \drawContour{xshift=5.25cm,yshift=6.75cm}{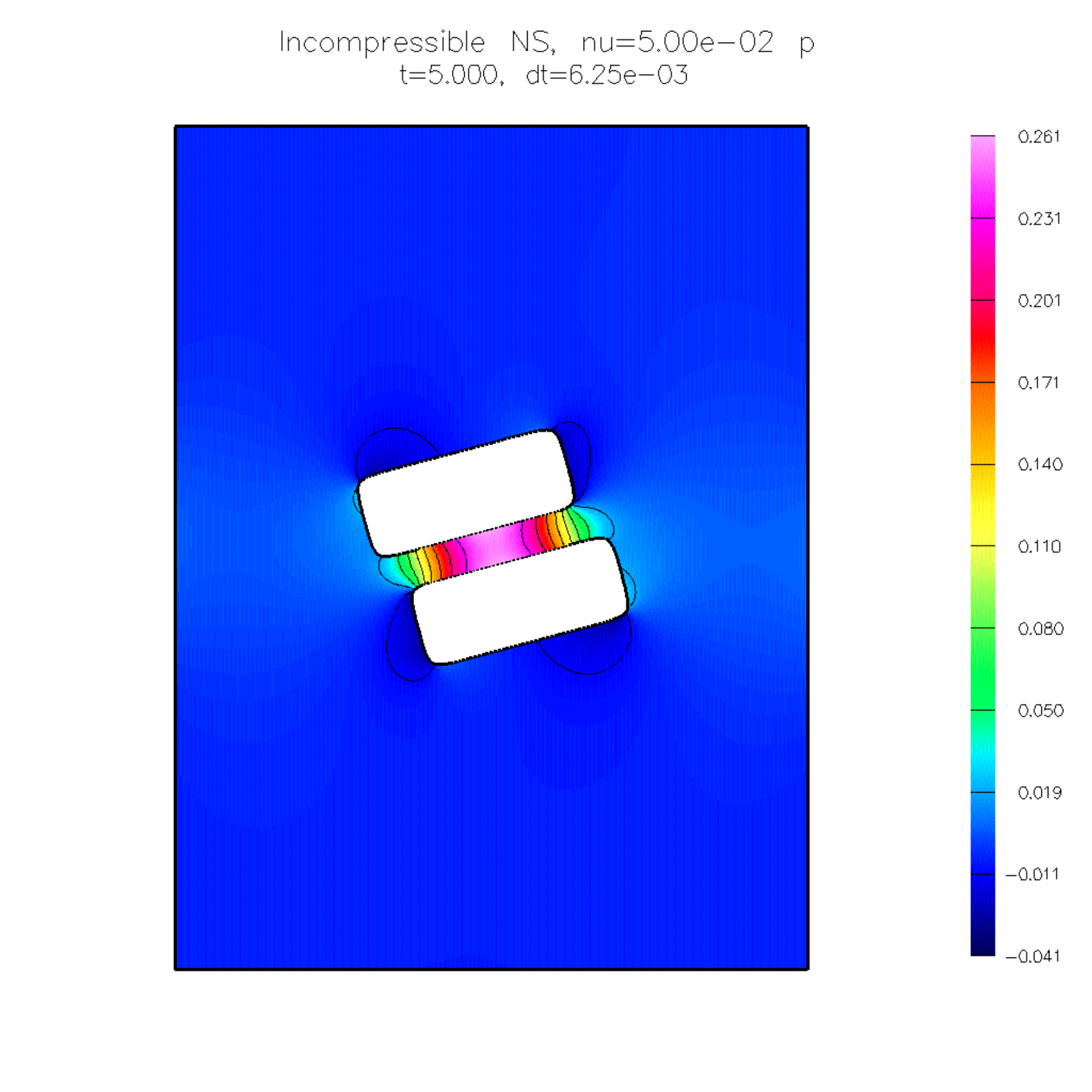}{$p$}{$p$}{$t=5.0$}{-.041}{.261};
 \drawContour{xshift=10.5cm,yshift=6.75cm}{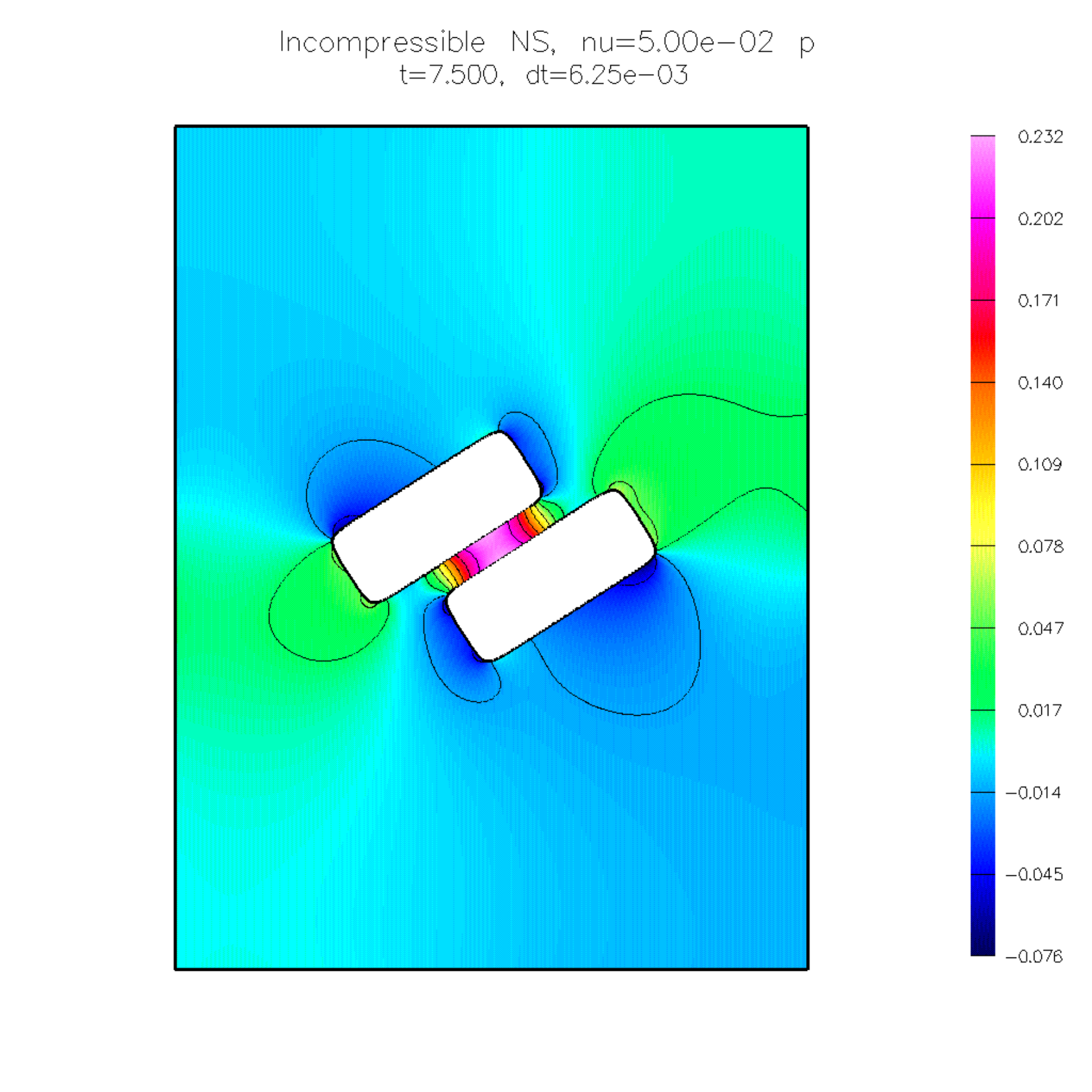}{$p$}{$p$}{$t=7.5$}{-.076}{.232};
 \drawContour{xshift= 0.0cm,yshift=0.cm}{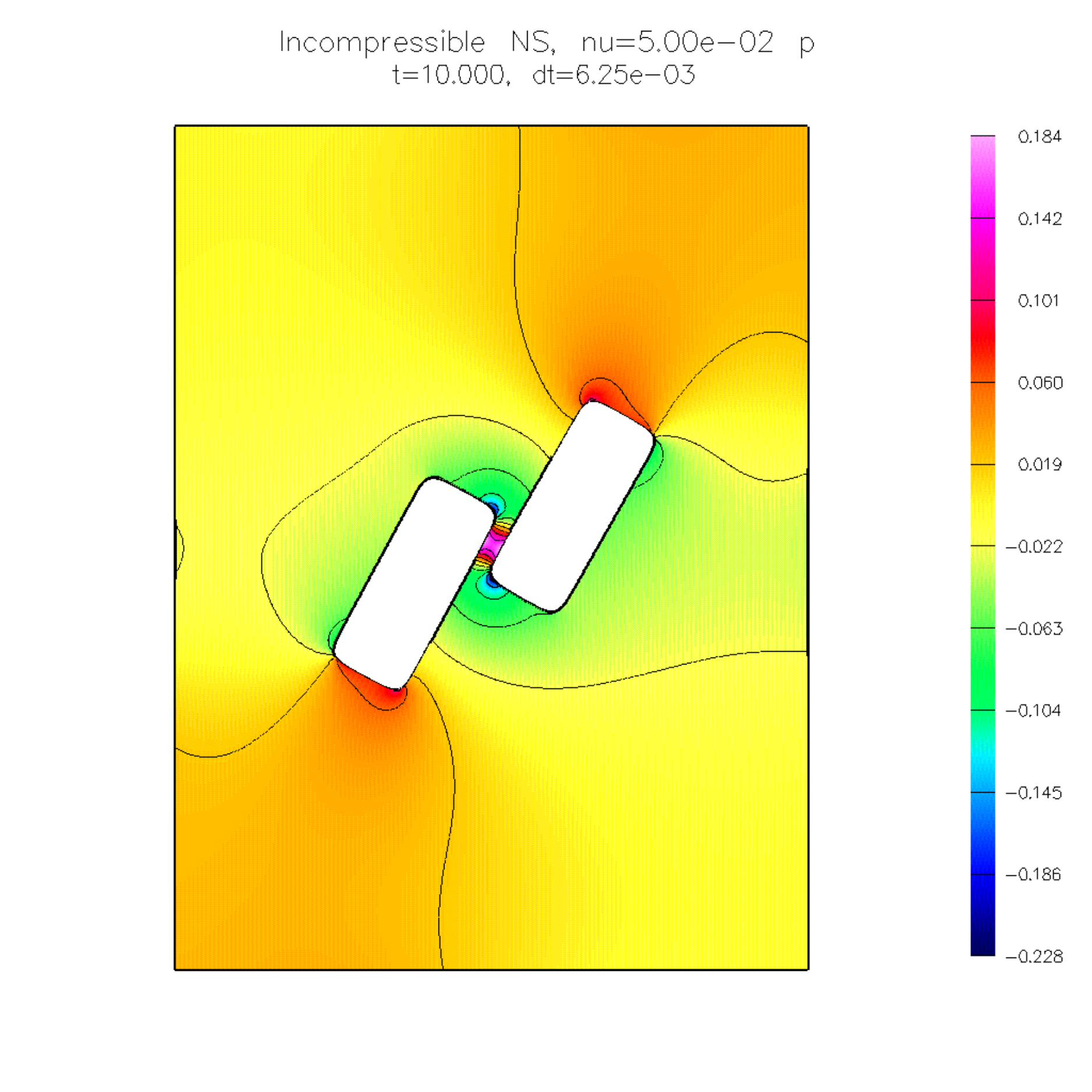}{$p$}{$p$}{$t=10.0$}{-.228}{.184};
 \drawContour{xshift=5.25cm,yshift=0.cm}{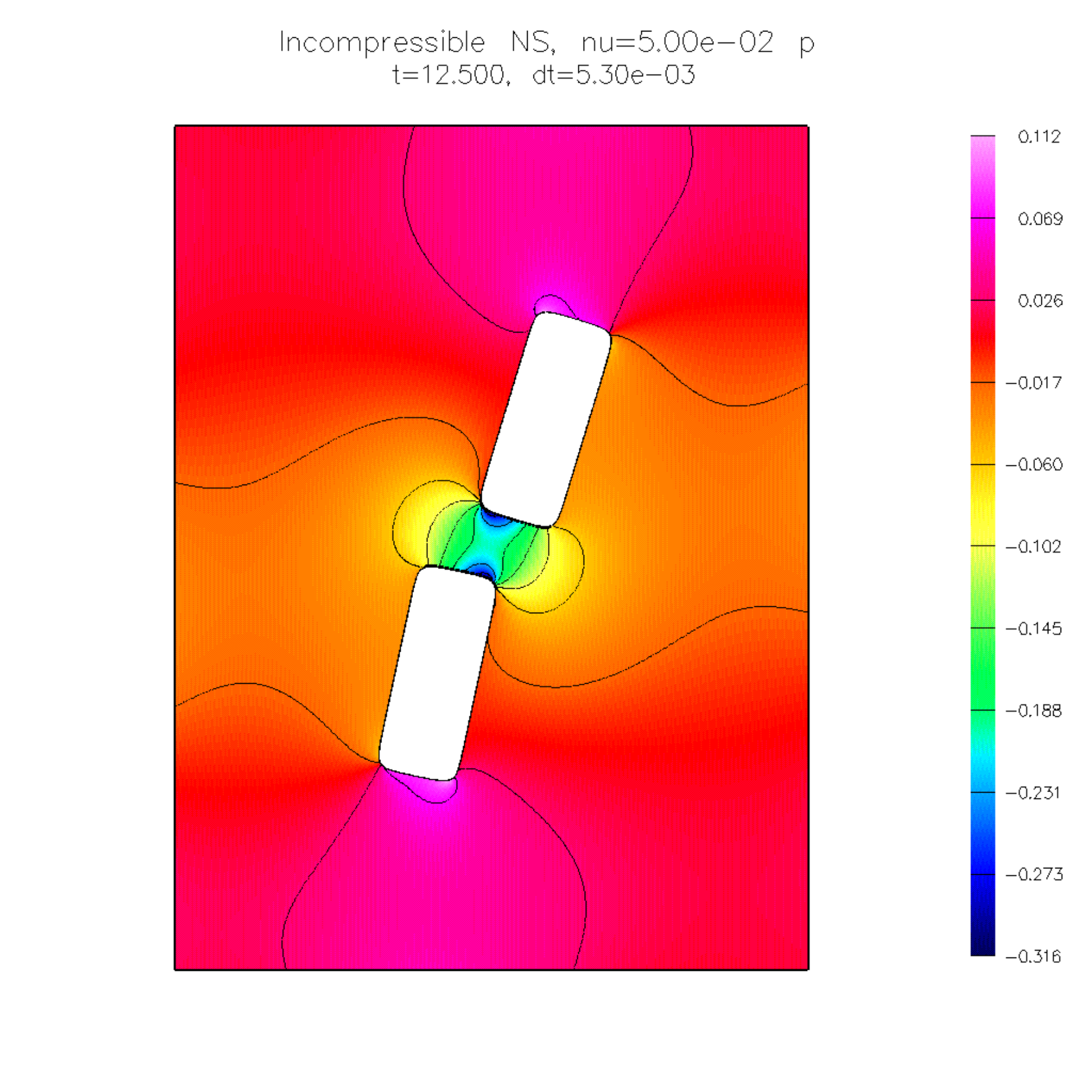}{$p$}{$p$}{$t=12.5$}{-.316}{.112};
 \drawContour{xshift=10.5cm,yshift=0.cm}{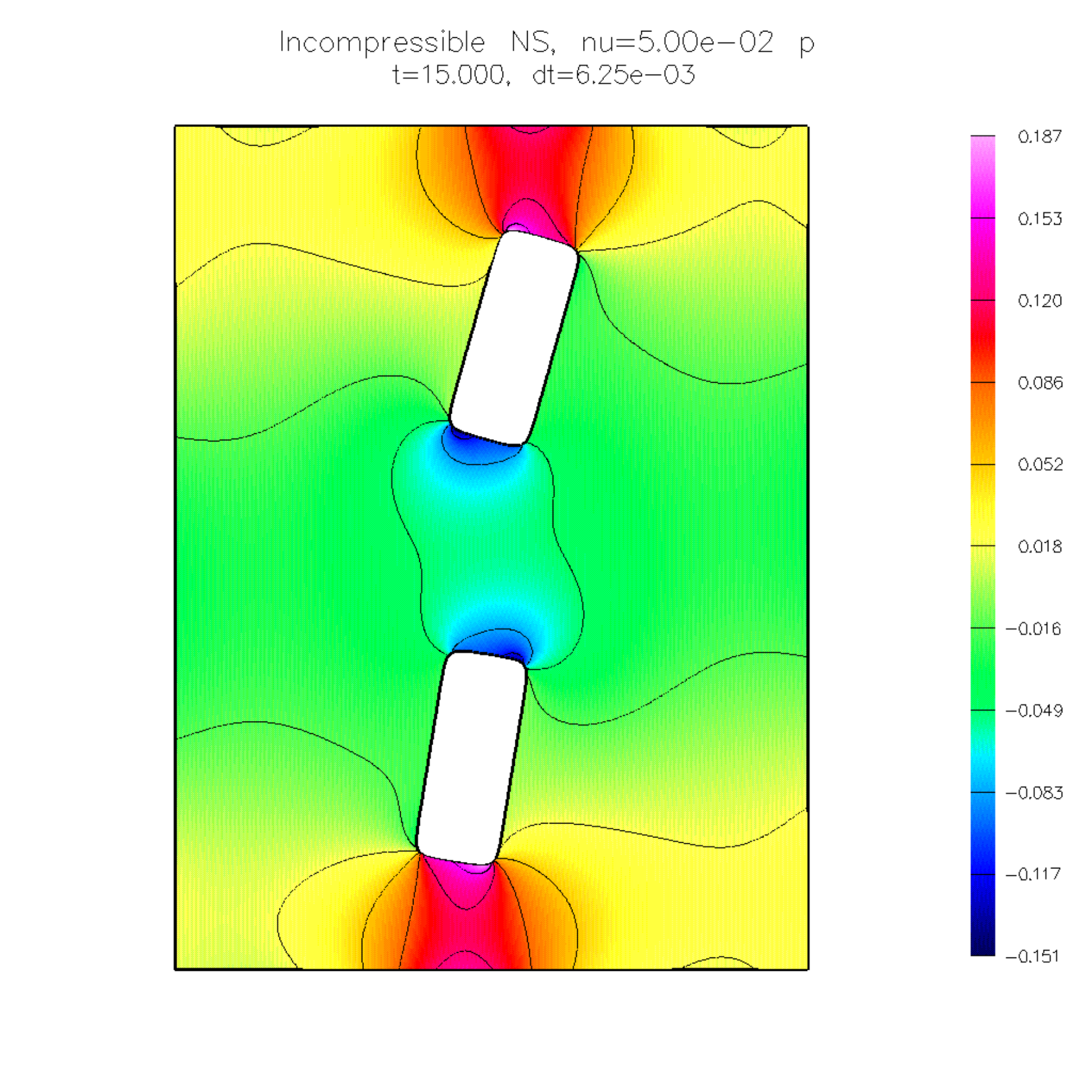}{$p$}{$p$}{$t=15.0$}{-.151}{.187};
%
\end{tikzpicture}
}
\end{center}
  \caption{Rising and falling bodies.  Contours of the pressure at various times.  The results use $\Gcrf^{(16)}$. }
  \label{fig:twoFallingBodiesPressure}
\end{figure}
}

{
\newcommand{\figWidth}{7.75cm}
\newcommand{\trimfig}[2]{\trimFig{#1}{#2}{.0}{.0}{.0}{.0}}
\begin{figure}[htb]
\begin{center}
\resizebox{13cm}{!}{
\begin{tikzpicture}[scale=1]
  \useasboundingbox (0.0,1) rectangle (16.,12.25);  
  \draw(0.0,6.1) node[anchor=south west,xshift=-4pt,yshift=+0pt] {\trimfig{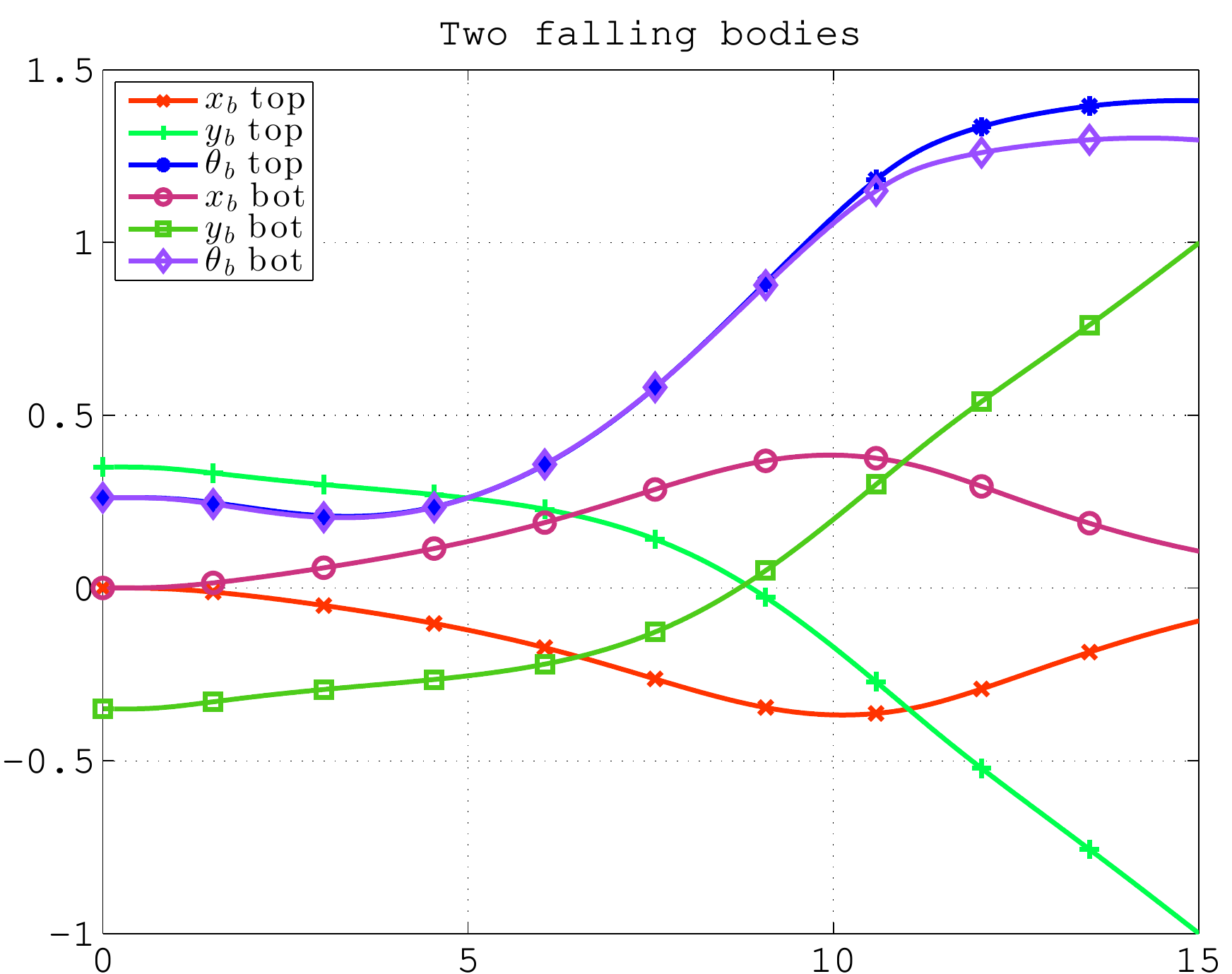}{\figWidth}};
  \draw(0.0,0) node[anchor=south west,xshift=-4pt,yshift=+0pt] {\trimfig{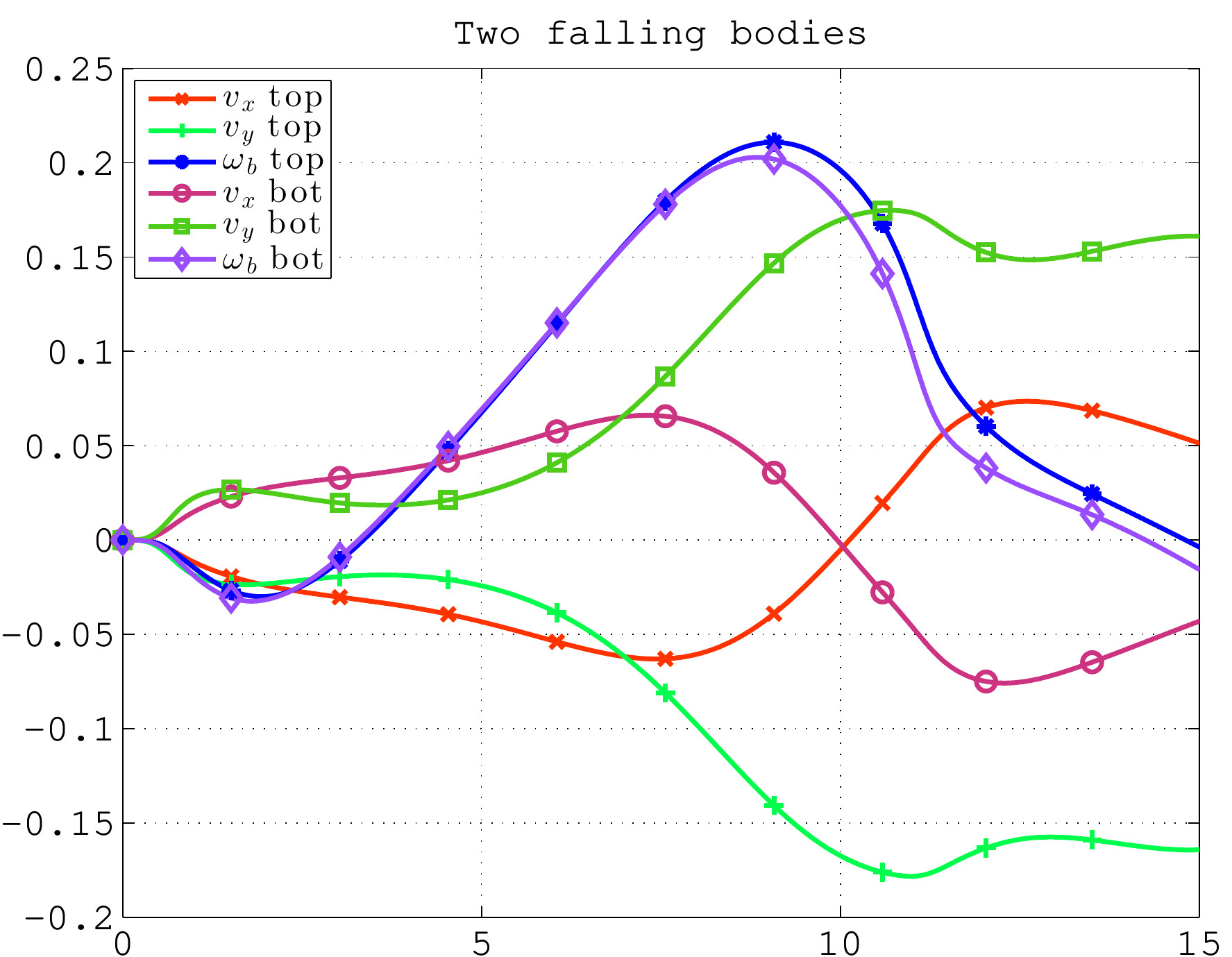}{\figWidth}};
  \draw(8.0,6.1) node[anchor=south west,xshift=-4pt,yshift=+0pt] {\trimfig{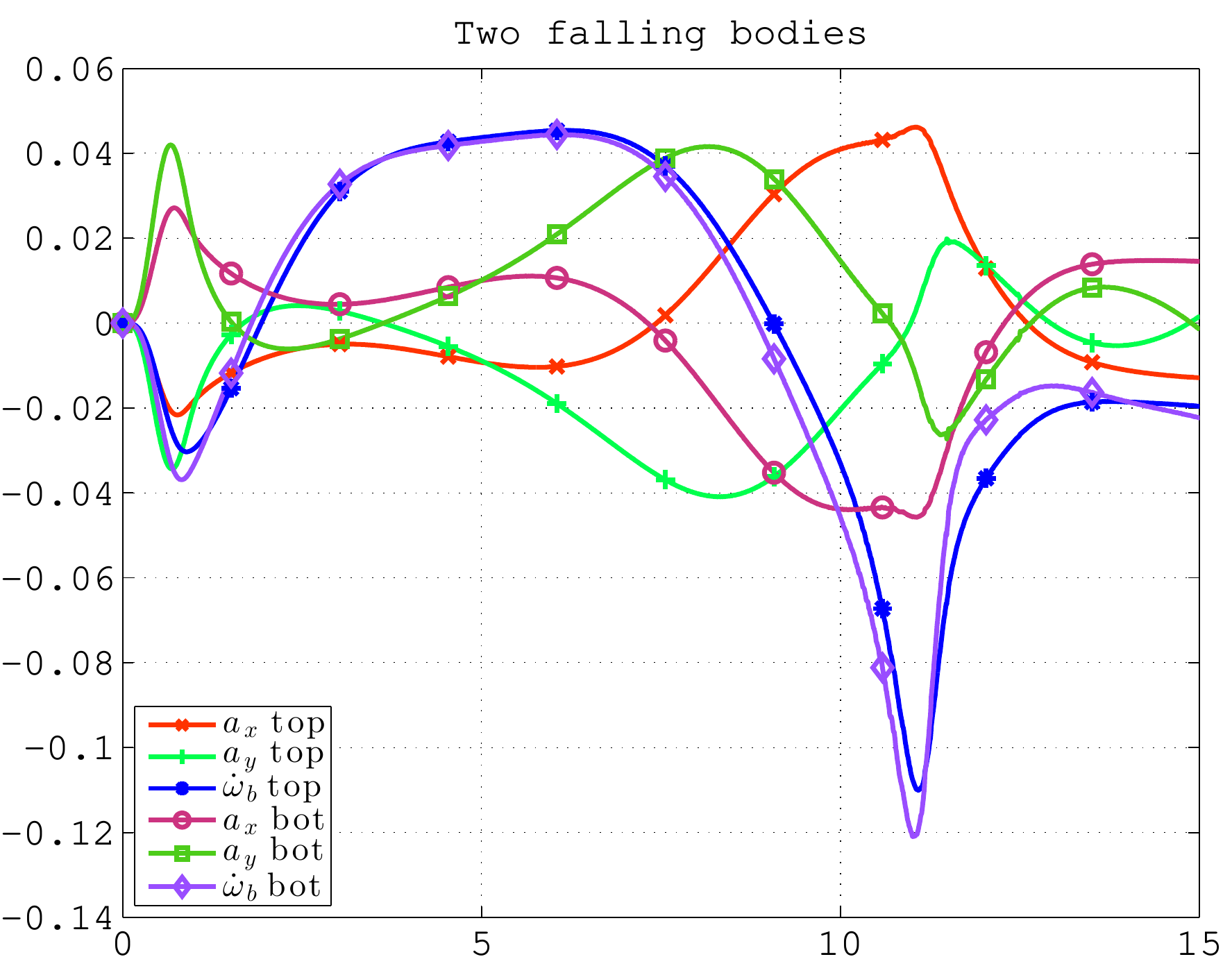}{\figWidth}};
  \draw(8.0,-.25) node[anchor=south west,xshift=-4pt,yshift=+0pt] {\trimfig{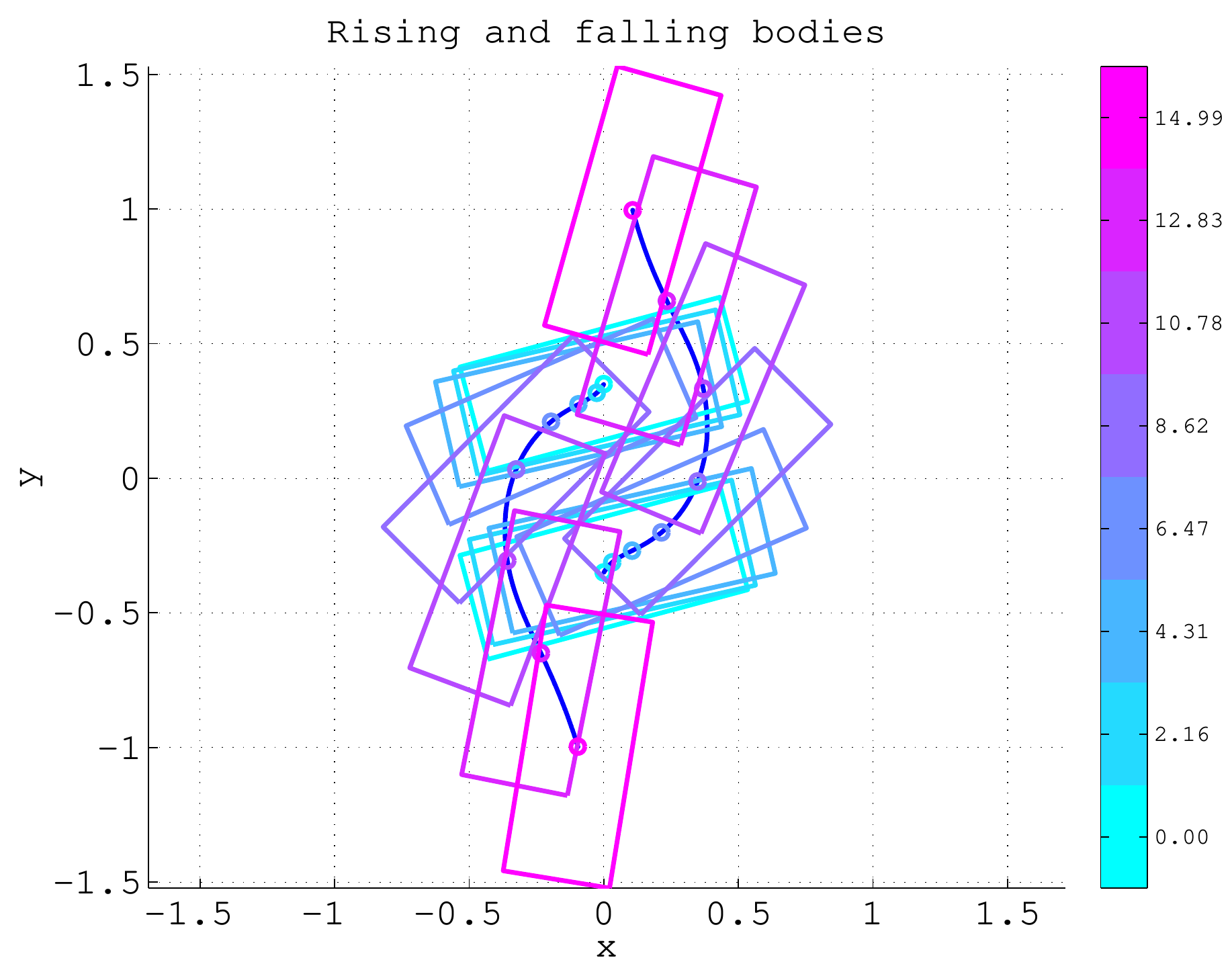}{\figWidth}};
  \draw(15.25,3.25) node [xshift=-3pt,yshift=0pt,rotate=-90] {\footnotesize\textsf{time}};
%
\end{tikzpicture}
}
\end{center}
  \caption{Rising and falling bodies. Time history of the position (top left), velocity (bottom left) and acceleration (top right) of the top and bottom rigid bodies.  The bottom-right plot gives a time animation of the locations 
of the two bodies with the colour indicating time.  The results use $\Gcrf^{(16)}$. 
     }
  \label{fig:twoFallingBodiesMotion}
\end{figure}
}

{
\newcommand{\figWidth}{8.cm}
\newcommand{\trimfig}[2]{\trimw{#1}{#2}{.0}{.0}{.0}{.0}}
\newcommand{\figWidthz}{8.5cm}
\newcommand{\trimfigz}[2]{\trimw{#1}{#2}{.0}{.0}{.0}{.0}}
\begin{figure}[htb]
\begin{center}
\resizebox{12cm}{!}{
\begin{tikzpicture}[scale=1]
  \useasboundingbox (0.0,.95) rectangle (16.,19.1);  
  \begin{scope}[yshift=4pt]
   \draw(0.0,12.8) node[anchor=south west,xshift=-4pt,yshift=+0pt] {\trimfig{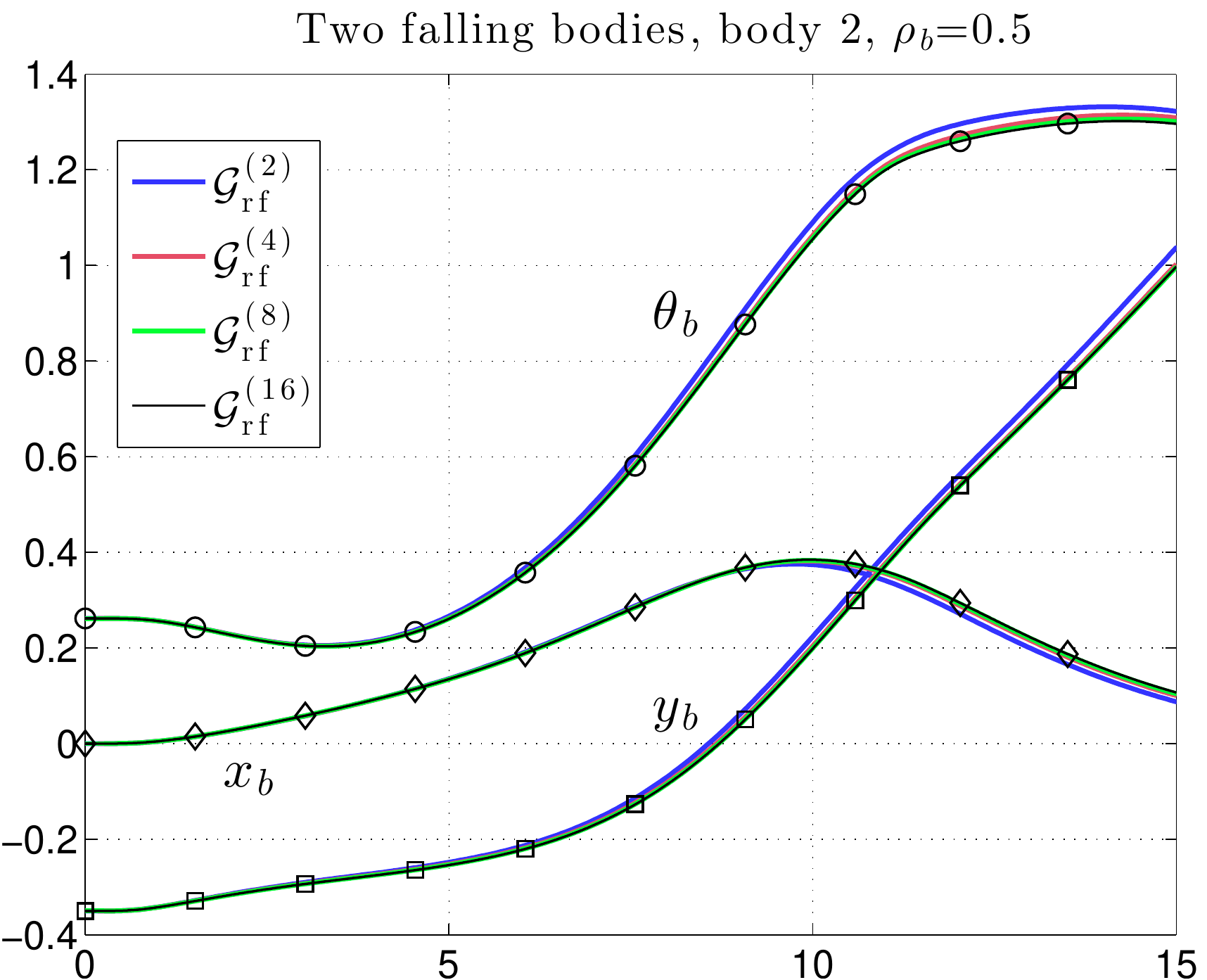}{\figWidth}};
   \draw(8.25,12.55) node[anchor=south west,xshift=-4pt,yshift=+0pt] {\trimfigz{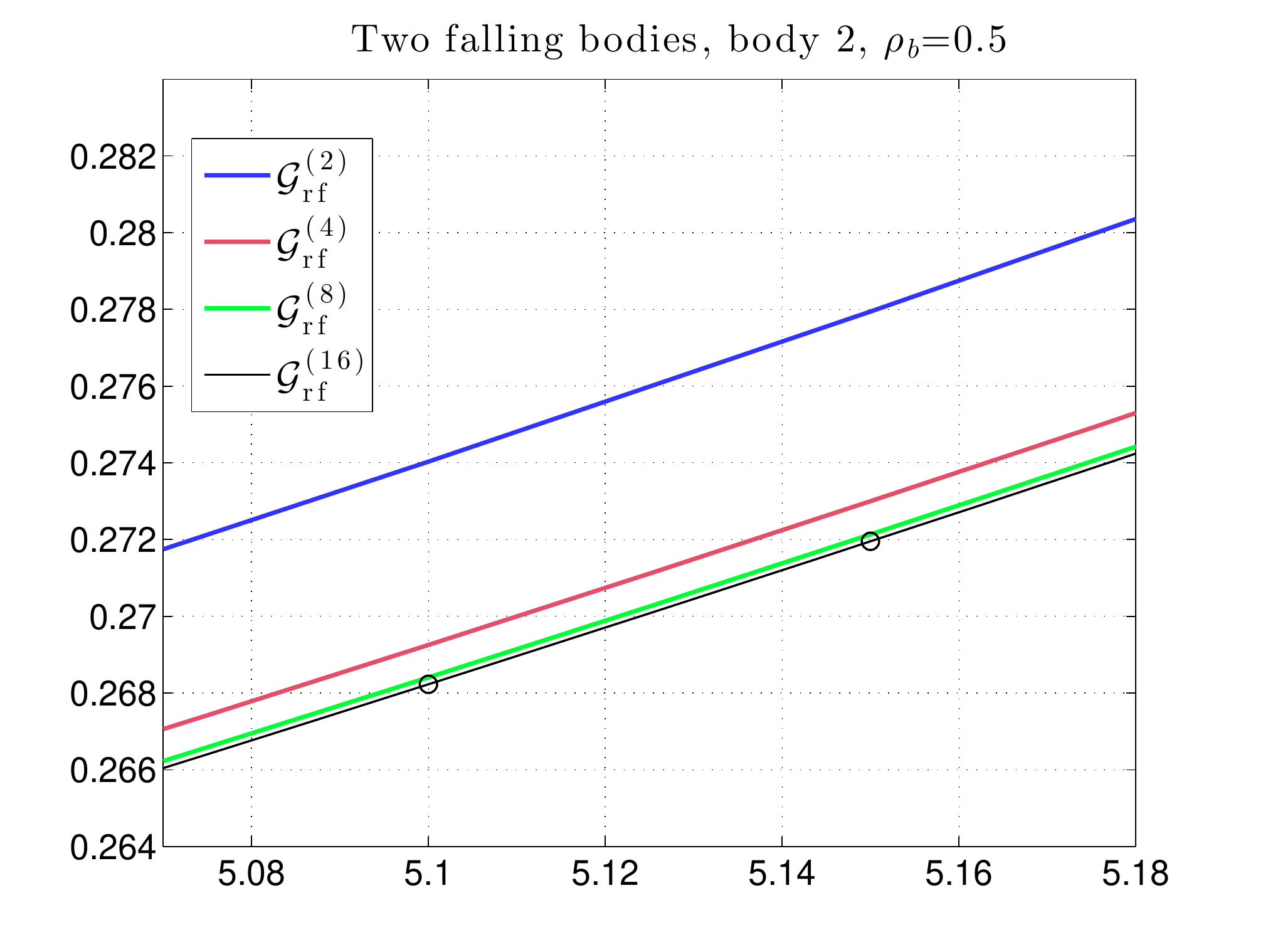}{\figWidthz}};
  \end{scope}
  \draw(0.0,6.4) node[anchor=south west,xshift=-4pt,yshift=+0pt] {\trimfig{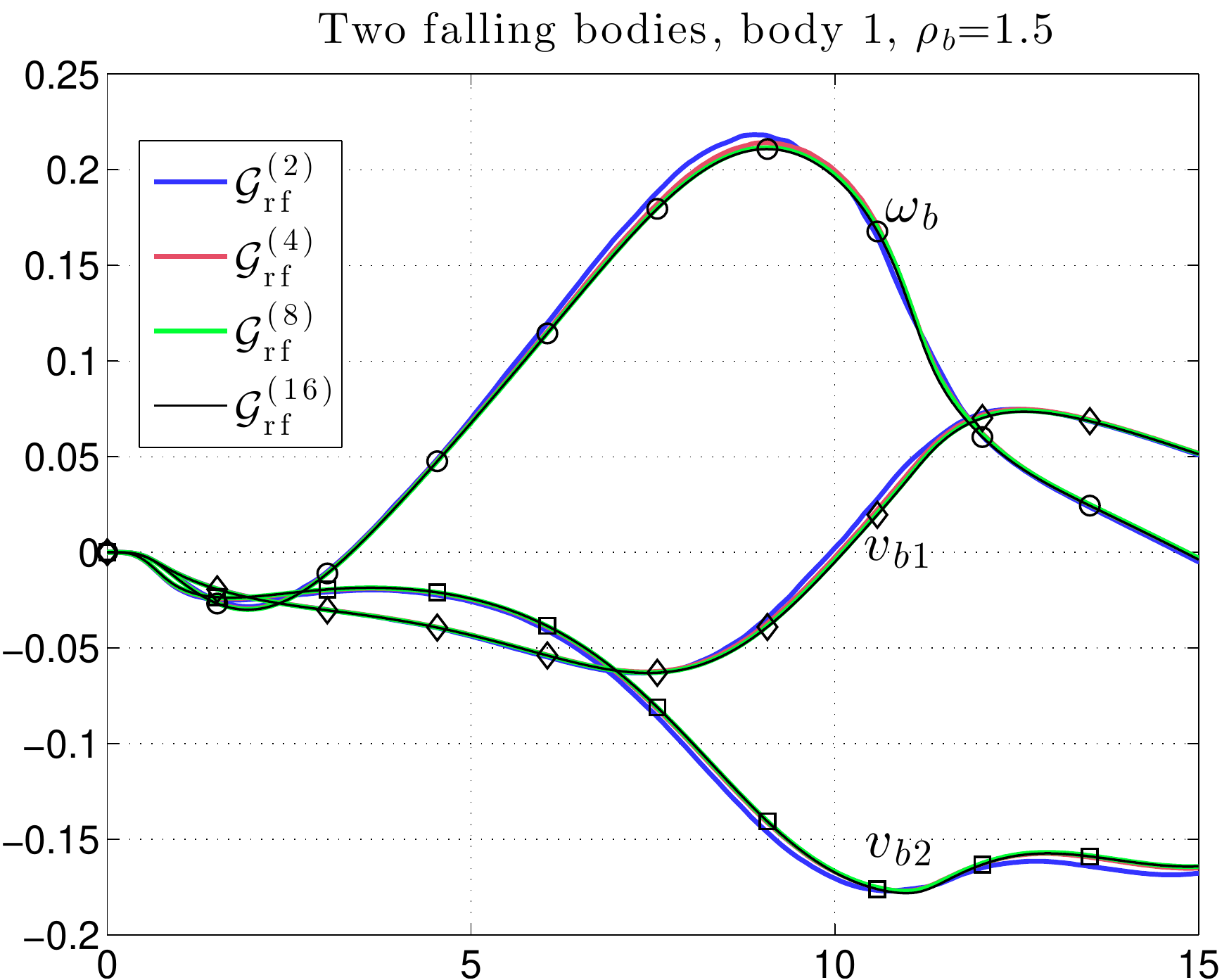}{\figWidth}};
  \draw(8.25,6.05) node[anchor=south west,xshift=-4pt,yshift=+0pt] {\trimfigz{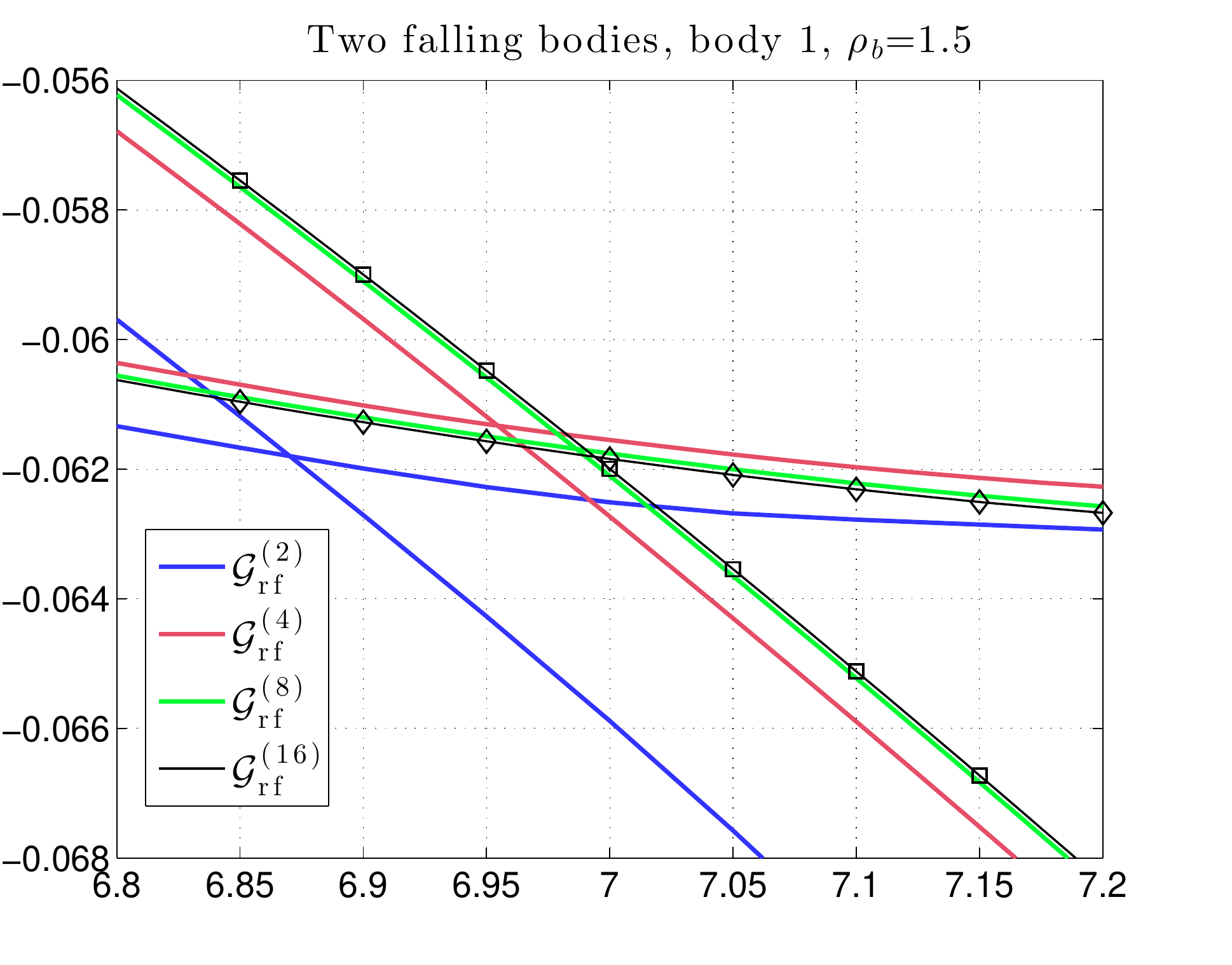}{\figWidthz}};
  \draw(0.0,0) node[anchor=south west,xshift=-4pt,yshift=+0pt] {\trimfig{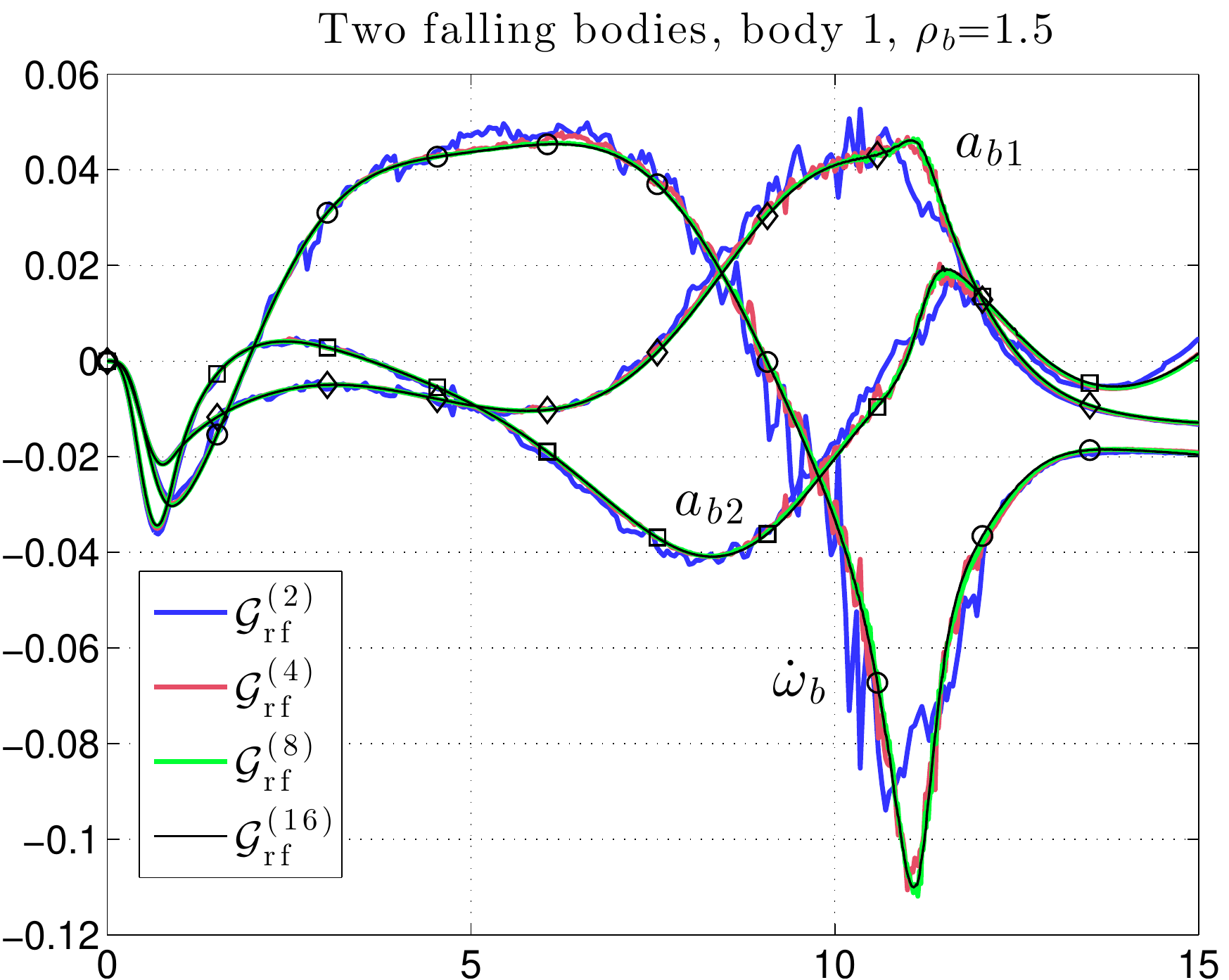}{\figWidth}};
  \draw(8.25,-.4) node[anchor=south west,xshift=-4pt,yshift=+0pt] {\trimfigz{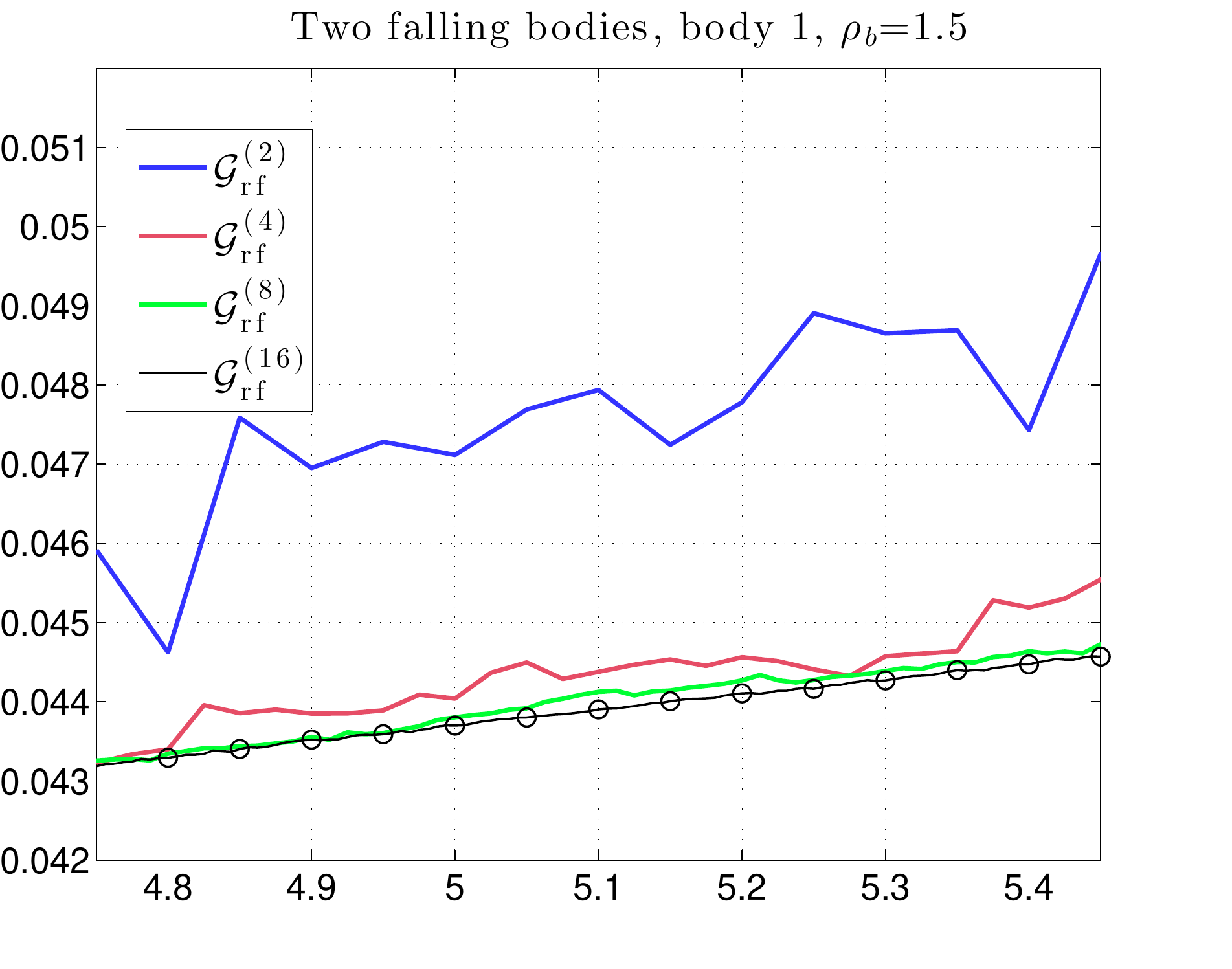}{\figWidthz}};
%
%
\end{tikzpicture}
}
\end{center}
  \caption{Rising and falling bodies. Grid convergence of representative variables for the rigid bodies.
    Top: positions $\xb$, $y_b$ and rotation angle $\theta_b$ for body 2.
    Middle: velocities $v_{b1}$, $v_{b2}$ and angular velocity $\omega_b$ for body 1.
    Bottom: accelerations $a_{b1}$, $a_{b2}$ and angular acceleration $\dot\omega_b$ for body 1.  Corresponding enlarged views are shown on the right.
     }
  \label{fig:twoFallingBodiesCompareZoom}
\end{figure}
}

The fluid domain covers the rectangle $[-1.5,1.5]\times[-2,2]$.
The solids are rectangles with rounded corners, both of the same shape and
with approximate width $\solidWidth=1$ and height $\solidHeight=0.4$. 
The top solid with density $\rhos=1.5$ is centred at $(x,y)=(0,\,0.35)$ initially and is rotated counter-clockwise by $15^\circ$ 
about its centre. 
The bottom solid with density $\rhos=0.5$ is centred at $(x,y)=(0,-0.35)$ initially and is 
rotated counter-clockwise by $15^\circ$ about its centre.
The composite grid for the domain, denoted by $\Gcrf^{(j)}$, where $j$ denotes the 
resolution factor, consists of three component grids as shown in Figure~\ref{fig:fallingBodyGrid}.
A Cartesian background grid covers most of the fluid domain, and boundary-fitted grids
are attached to the surface of the two bodies.  The grid spacing for the Cartesian grid
is $\ds^{(j)}=1/(10 j)$.  The two boundary-fitted grids are stretched in the direction normal
to the surface of the bodies so that $\ds^{(j)}\approx1/(10 j)$ for the grid lines away from
the body, whereas the grids lines near the body surfaces are clustered by a factor of approximately~$4$.

The boundary conditions for the problem are taken as no-slip walls on the surface of the bodies and
all sides of the rectangular fluid domain.  The fluid and the solid bodies are at rest initially.
The viscosity of the fluid is $\mu=.05$, and a body force given by~\eqref{eq:bodyForce} is taken for
each solid with $\gv=[0,\,-1]$.  As before, the body force turns on smoothly to avoid an impulsive start.

Shaded contours of the pressure 
are shown at selected times in Figure~\ref{fig:twoFallingBodiesPressure}, 
while the time-history of the rigid body variables are shown in Figure~\ref{fig:twoFallingBodiesMotion}.
The latter figure also gives an animation of the translation and rotation of the two bodies over time. 
The bodies start from rest and the body force due to gravity is turned on smoothly over the time interval $t\in[0,1]$.
The bodies slowly approach one another at early times. 
During this time, a large pressure develops in the gap between the bodies which
reducing the rate at which the gap closes.
The top body begins to move to the left while the bottom body moves to the right. 
At the same time the bodies begin to rotate in a counter-clockwise direction. 
The bodies continue to slide by each other and rotate. By time $t=10$ the bodies 
have nearly separated from one another, and the gap between the bodies
becomes quite small (see Figure~\ref{fig:twoFallingBodiesPressure}). The angular acceleration of the bodies is the largest in magnitude at about this time.
By $t=12$ the bodies have separated and there is a negative pressure in the remaining gap
that keeps the bodies from separating too quickly. The bodies continue to separate and by $t=15$
they are approaching the lower and upper walls of the fluid container.

Figure~\ref {fig:twoFallingBodiesCompareZoom} shows a grid convergence study for selected rigid body variables; 
other variables show similar behaviour.
As the grid is refined the curves are seen to converge 
with the separation between curves decreasing rapidly as the grids are refined.
The positions and velocities of the bodies are reasonably represented even on the coarsest grid $\Gcrf^{(2)}$, 
although the accelerations on this grid show small oscillations.  The amplitude of these
oscillations in the acceleration decrease as the grid is refined as we have noted in previous example calculations.

\section{Conclusions} \label{sec:conclusions}

In this second part of a two-part series we have developed the general formulation  of the
{\ampRB} scheme for the partitioned solution of rigid bodies moving in an incompressible fluid.
This scheme remains stable, without sub-iterations, for light and even zero-mass bodies.
The primary extension to multiple space dimensions involved the development of the added-damping tensors; these
are incorporated into the equations of motion of the rigid body to address added-damping instabilities.
Linearization of the force and torque terms on the rigid-body led to analytic forms
for these added-damping tensors whose entries depend on solutions to two vector
Helmholtz equations. Approximate forms of these tensors, convenient for use in numerical
solutions were developed. The entries in the approximate added-damping tensors,
which involving surface integrals over the body,  can be pre-computed during a preprocessing step.
The {\ampRB} scheme was implemented in two dimensions for general domains with one or more rigid bodies by using composite
overlapping grids. This is an efficient approach for treating moving body problems and provides
smooth boundary fitted grids on the body for accurate treatment of boundary layers.

Numerical results were used to demonstrate the stability and second-order accuracy of the
scheme. Some challenging problems involving very light bodies demonstrated the effectiveness
of the scheme. The {\ampRB} scheme was shown to produce nearly identical results to the
traditional scheme with sub-iterations, although for hard cases the traditional scheme required tens or hundreds
of sub-iterations per time-step. 
Added-damping effects arising from viscous shear stresses are clearly important for rotationally
symmetric bodies such as a disk in two dimensions or a sphere in three dimensions.
The case of a rising buoyant rectangular body, however, demonstrated that it is also important to account for added-damping effects
for body shapes that are far from rotationally symmetric.

In future work we plan to extend the \ampRB~algorithm to three dimensions, which in principal should
not present any significant issues.
In practice, there are various details that need to be treated properly such as the accurate 
evaluation of surface integrals on bodies covered by multiple overset grids as needed for computing forces, torques
and entries in the added-damping tensors.
We will also consider extending the \ampRB~scheme to higher-order accuracy in space and time.
For this extension, we note that the spatial discretization 
of the Navier-Stokes equations to fourth-order accuracy on overset grids has been
developed previously in~\cite{ICNS}, and this could be used as a basis for a fourth-order accurate scheme. 
Higher-order accurate schemes in time, however,
will need to be developed and these will need to be stable within the fractional-step framework. 
The approximate added-damping tensors developed in this article
could be used in the high-order accurate scheme, likely without change, provided high-order
accurate values for the predicted body accelerations are used, but again the stability of such a procedure
would need to be evaluated.

\appendix
\newcommand{\p}{\partial}

\section{Shear stress of an incompressible fluid on the surface of a rigid body}
\label{sec:shearStressLemma}

In this section we prove Theorem~\ref{theorem:shearStressRigidBody} from Section~\ref{sec:addedDampingTensors}. 
Let $\xv=\gv(\rv,t)$ be  an orthogonal coordinate transformation for the fluid region adjacent to some portion of the body 
whose surface located at $r_1=0$. Thus points $\xv\in \GammaB$ are of the form $\xv = \gv(0,r_2,r_3,t)$, where $r_2$ and $r_3$ parameterize
a portion of the surface. 
The unit normal $\nv=[n_1,\, n_2,\, n_3]^T$ to the surface can be written as 
\begin{align}
    \nv = \frac{1}{\alpha_1} \grad_\xv r_1 , \quad n_i = \frac{1}{\alpha_1} \f{\p r_1}{\p x_i}, 
      \qquad \alpha_1=\| \grad_\xv r_1 \|. 
  \label{eq:normalDef}
\end{align}
The unit tangent vectors can be written in the form,
\begin{alignat}{3}
	\tv_m &= \frac{1}{\alpha_{m+1}} \grad_\xv r_{m+1} , \qquad&& \alpha_{m+1}=\| \grad_\xv r_{m+1} \|, \quad m=1,2.    \label{eq:tangentDef}
\end{alignat}
The normal and tangent vectors form an orthonormal set.
We first derive an expression for the viscous shear stress 
in terms of the normal and tangential  derivatives of the velocity, given by
\begin{align}
  &\vv_n =  \frac{\partial \vv}{\partial n} =(\nv\cdot\grad) \vv = \alpha_1 \frac{\partial \vv}{\partial r_1}, \label{eq:normalDeriv}\\
  &\vv_{t_m} =  \frac{\partial \vv}{\partial t_m} = (\tv_m\cdot\grad) \vv = \alpha_{m+1} \frac{\partial \vv}{\partial r_{m+1}} , \qquad m=1,2, \label{eq:tangentialDeriv}
\end{align}
which is valid for an incompressible fluid in general unconnected
to any rigid body.
\begin{lemma} \label{theorem:shearStressRigidBodyLemma2}
The viscous shear stress for an incompressible fluid, written in terms of the
normal and tangential derivatives of the velocity, is 
\begin{align}
   \tauv\nv = \mu \sum_{m=1}^2  
   \left\{\tv_m \tv_m^T\, \f{\p\vv}{\p n}  + (\tv_m\nv^T -2\nv\,\tv_m^T) \f{\p \vv}{\p t_m}\right\} .  \label{eq:stressForIncompressible}
\end{align}
\end{lemma}
\begin{proof}
The tensor $\grad\vv$ with components $\partial v_i/\partial x_j$ can be written as 
\begin{align}
& \grad\vv  = \frac{\partial\vv}{\partial n} \, \nv^T 
   + \sum_{m=1}^2 \frac{\partial\vv}{\partial t_m} \,\tv_m^T  .   \label{eq:gradvTensor}
\end{align}
Whence the expressions $(\grad\vv)\nv$ and $(\grad\vv)^T\nv$ can be written as 
\begin{align}
& (\grad\vv)\, \nv = \vv_n= \nv\, \nv^T \vv_n +\sum_{m=1}^2 \tv_m\,\tv_m^T \vv_n,    \label{eq:gradv} \\ 
& (\grad\vv)^T \, \nv = \nv\nv^T \vv_n + \sum_{m=1}^2 \tv_m\nv^T \vv_{t_m} , \label{eq:gradvTranspose}
\end{align} 
where the final relation in~\eqref{eq:gradv} comes from the decomposition of the vector $\vv_n$ into its normal and
tangential components. 
Using~\eqref{eq:gradv} and~\eqref{eq:gradvTranspose} 
gives
\begin{align}
  \frac{1}{\mu} \tauv\nv = (\grad\vv + (\grad\vv)^T) \nv 
      = 2 \nv\nv^T \vv_n +  \sum_{m=1}^2 \tv_m\tv_m^T \vv_n + \tv_m\nv^T \vv_{t_m} . \label{eq:stressII}
\end{align}
Note that $\grad\cdot\vv=0$ implies
\begin{align}
&   \nv^T\frac{\partial\vv}{\partial n}  = - \sum_{m=1}^2  \tv_m^T\frac{\partial\vv}{\partial t_m} , 
    \label{eq:divConversionII}
\end{align}
since, using~\eqref{eq:normalDef} for $\nv$ and~\eqref{eq:tangentDef} for $\tv_m$ 
together with~\eqref{eq:normalDeriv} and~\eqref{eq:tangentialDeriv}, 
\begin{align*}
  \nv^T\frac{\partial\vv}{\partial n}  + \sum_{m=1}^2 \tv_m^T\frac{\partial\vv}{\partial t_m} &=
        \sum_{j=1}^3 \left\{  \f{1}{\alpha_1}\f{\p r_1}{\p x_j} \, \alpha_1 \f{\p v_j}{\p r_1} 
		+  \sum_{m=1}^2   \f{1}{\alpha_{m+1}}\f{\p r_{m+1}}{\p x_j} \, \alpha_{m+1} \f{\p v_j}{\p r_{m+1} } \right\} ,\\
       &= \sum_{j=1}^3\sum_{k=1}^3  \f{\p r_k}{\p x_j} \f{\p v_j}{\p r_k }  = \sum_{j=1}^3 \f{\p v_j}{\p x_j} = \grad\cdot\vv=0. 
\end{align*}
Using~\eqref{eq:divConversionII} in~\eqref{eq:stressII} gives 
\begin{align*}
	\frac{1}{\mu} \tauv\nv =  \sum_{m=1}^2 \left\{\tv_m\tv_m^T \vv_n  + (\tv_m\nv^T -2\nv\, \tv_m^T) \vv_{t_m}\right\} ,
\end{align*}
which proves the lemma.
\end{proof}

\begin{lemma} \label{theorem:shearStressRigidBodyLemma3}
The vectors in the orthonormal set $\{\nv, \tv_1, \tv_2\}$ satisfy
\begin{align}
& \nv \nv^T+\tv_1\tv_1^T+\tv_2\tv_2^T = \Iv,  \label{eq:normalSum} \\
& \tv_2 \tv_1^T -\tv_1 \tv_2^T= \pm [\nv]_\times ,  \label{eq:tangentSum}
\end{align}
where the sign of the right-hand-side is positive if the set is right-handed and negative otherwise.
\end{lemma}
\begin{proof}
  Equation~\eqref{eq:normalSum} holds for any orthonormal set since $\xv=\nv \nv^T\xv+\tv_1\tv_1^T\xv+\tv_2\tv_2^T\xv$
for any vector $\xv$. For a right-handed set, equation~\eqref{eq:tangentSum}, with the plus sign,
follows since $\nv\times\nv=0$, $\nv\times\tv_1=\tv_2$ and
$\nv\times\tv_2=-\tv_1$. For a left-handed set the sign is reversed.
\end{proof}
Condition~\eqref{eq:stressForIncompressible} for $\tauv\nv$ can now be specialized to the surface
of a rigid body. 
By taking tangential derivatives of the interface condition
\[
    \vv(\xv,t) = \vvcm(t) + \omegavb(t)\times(\xv-\xvb(t)), \qquad \xv\in\Gamma, 
\]
it follows that the tangential derivatives of $\vv$ on the rigid body are  given by 
\begin{align*}
   \frac{\partial\vv}{\partial t_m} = \omegavb(t)\times \frac{\partial\xv}{\partial t_m}
                =  \omegavb(t)\times\tv_m , \qquad  \xv\in\Gamma , 
\end{align*} 
using ${\partial\xv}/{\partial t_m}=\tv_m$. 
Thus
\begin{align*}
	\frac{1}{\mu} \tauv\nv &= \sum_{m=1}^2\left\{ \tv_m\tv_m^T \vv_n  + (\tv_m\nv^T -2\nv\,\tv_m^T)\, (\omegavb(t)\times\tv_m)\right\} , \\
	&= (\Iv-\nv \nv^T) \vv_n + \sum_{m=1}^2 \left\{ \tv_m\nv^T(\omegavb(t)\times\tv_m)\right\} = (\Iv-\nv \nv^T) \vv_n + \sum_{m=1}^2 \tv_m \left\{ \omegavb(t)^T(\tv_m\times \nv) \right\}, \\
	&= (\Iv-\nv \nv^T) \vv_n \pm (-\tv_1 \tv_2^T + \tv_2 \tv_1^T) \omegavb(t) = (\Iv-\nv \nv^T) \vv_n + [\nv]_\times \omegavb(t),
\end{align*}
where the $\pm$ sign depends on whether the system is right-handed or left-handed and
where we have used the fact that $\tv_m$ and $\omegavb(t)\times \tv_m$ are orthogonal.
This completes the proof of Theorem~\ref{theorem:shearStressRigidBody}. 

\section{Example added-damping tensors for bodies of different shapes} \label{sec:addedDampingTensorExamples}

In this section, examples of the approximate 
added-damping tensors~\eqref{eq:ADvv}--\eqref{eq:ADww} are presented for rigid bodies of different shapes.
The entries in the tensors are computed analytically 
by assuming that $\dn=\dn_\iv$ is constant for $\iv\in\Gamma_h$, and by replacing the surface quadratures in~\eqref{eq:ADvv}--\eqref{eq:ADww} by surface integrals.
These examples can provide some intuition for relating the geometry of a body to 
the form and magnitude of the added-damping as the body translates or rotates. 
Recall that as the bodies rotate the added-damping tensors transform according to~\eqref{eq:addedDampingInTime}.

\subsection{Cylinder and disk}

For a solid three-dimensional circular cylinder of radius $a$ and depth $d$ aligned with the $z$-axis
the moment of inertia about the centre of the disk is
\begin{align*}
  \Ib = \begin{bmatrix} I_x & 0 & 0 \\ 0 & I_y & 0 \\ 0 & 0 & I_z \end{bmatrix}, \qquad
    I_x=I_y = \frac{\ms}{12}( 3 a^2 + d^2), \quad I_z=\frac{\ms}{2} a^2 . 
\end{align*}
For a cylinder restricted to rotate about the $z$-axis, 
the approximate added-damping tensors~\eqref{eq:ADvv}--\eqref{eq:ADww}
are (ignoring the top and bottom circular faces so that the result can be restricted to a two-dimensional disk)
\begin{align}
 & \Dvva = \frac{\mu}{\dn}\, \pi a d \, \begin{bmatrix} 1 & 0 & 0 \\ 0 & 1 & 0 \\ 0 & 0 & 0 \end{bmatrix},  \quad
   \Dvwa = \Dwva = \zerov, \quad
   \Dwwa = \frac{\mu}{\dn}\, (2 \pi a d) a^2\begin{bmatrix} 0 & 0 & 0 \\ 0 & 0 & 0 \\ 0 & 0 & 1 \end{bmatrix}. \label{eq:addedDampingDisk}
\end{align}
Notice that the entries $\Dvva_{11}$ and $\Dvva_{22}$ in $\Dvva$ are 
the product of $\mu/\dn$ and one half of the surface area of the disk, $S=2\pi a d$. 
These values indicate the contribution to the added damping from motions in the $x$- and $y$-directions, respectively. 
The entry $\Dwwa_{22}$ in $\Dwwa$ is the product of the disk surface area $S$ and the square of the radius. 
The surface area $S$ is the area over which a shear force could act and this is multiplied by the
radius of the disk squared. 
The effective moment of inertial for rotations about the $z$-axis is thus
\[
    I_z + \dt \Dwwa_{22} =  \frac{\ms}{2} a^2 + \dt \frac{\mu}{\dn}\, (2 \pi a d)  a^2  
\]

\subsection{Rectangle}

Consider a solid rectangle of constant density $\rhos$ with width $w$ in the $x$-direction and height $h$ in the $y$-direction.
The moment of inertia about the centre of mass for rotations in the plane about the $z$-axis is
\[
   I_z = \frac{\ms}{12}\big( w^2 + h^2 \big) = \frac{\rhos w h}{12} \big( w^2 + h^2 \big) .
\]
The added-damping tensors~\eqref{eq:ADvv}--\eqref{eq:ADww}
for the rectangle 
are 
\begin{align}
 & \Dvva = \frac{\mu}{\dn}\, \begin{bmatrix} 2 w & 0 & 0 \\ 0 & 2 h & 0 \\ 0 & 0 & 0 \end{bmatrix},  \quad
   \Dvwa = \Dwva = \zerov, \quad
   \Dwwa = \frac{\mu}{\dn}\,  2(w+h)\, \frac{w}{2}\frac{h}{2}  \begin{bmatrix} 0 & 0 & 0 \\ 0 & 0 & 0 \\ 0 & 0 & 1 \end{bmatrix}.
     \label{eq:addedDampingRectangle}
\end{align}
The entries $\Dvva_{11}$ and $\Dvva_{22}$ in $\Dvva$ are 
the product $\mu/\dn$ and the surface areas in the $x$- and $y$-directions, respectively. 
These values indicate the contribution to the added damping from motions in the $x$- and $y$-directions, respectively. 
The entry $\Dwwa_{22}$ in $\Dwwa$ is the product of the surface area of the rectangle $S=2(w+h)$ and the quantity $\frac{w}{2}\frac{h}{2}$ 
which corresponds to the radius-squared term for the the disk in~\eqref{eq:addedDampingDisk}. 

\newcommand{\wt}{w}
\newcommand{\at}{a}
\newcommand{\bt}{b}
\begin{figure}[hbt]
\begin{center}
\resizebox{6cm}{!}{
\begin{tikzpicture}
 \useasboundingbox (0,0.25) rectangle (8,4);  
\begin{scope}[xshift=1cm,yshift=2cm,scale=1]
\draw[-,very thick,fill=blue!30] (0,-1) -- (6,-2) -- (6,2) -- (0,1) -- (0,-1);
\fill[red] (3.5,0) circle (3 pt) node[anchor=west,black,xshift=3pt] {$(\frac{w}{6}\frac{b-a}{b+a},0)$};
\fill[blue] (6,-2) circle (3 pt) node[anchor=north west,black] {$(\frac{\wt}{2},-\frac{\bt}{2})$};
\fill[blue] (6, 2) circle (3 pt) node[anchor=north west,black] {$(\frac{\wt}{2}, \frac{\bt}{2})$};
\fill[blue] (0,-1) circle (3 pt) node[anchor=north east,black] {$(-\frac{\wt}{2},-\frac{\at}{2})$};
\fill[blue] (0, 1) circle (3 pt) node[anchor=north east,black] {$(-\frac{\wt}{2}, \frac{\at}{2})$};
\end{scope}
%
%
\end{tikzpicture}
} 
\end{center}
\caption{Geometry of the trapezoid used for computation of the added-damping tenors. The centre
of mass is marked 
as a solid red circle.}    
\label{fig:addedDampingTrapezoid}
\end{figure}
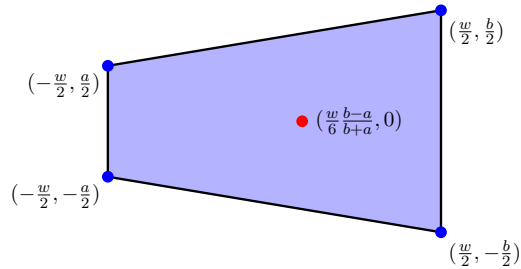

\subsection{Trapezoid}

Consider a solid trapezoid of constant density $\rhos$
as shown in Figure~\ref{fig:addedDampingTrapezoid}.
The moment of inertia about the centre of mass for rotations about the $z$-axis is
\begin{align*}
  I_z = \frac{\rhos w}{12}\frac{(a + b)}{2} \left( w^2\,\frac{2}{3}\frac{(a^2+4ab+b^2)}{(a+b)^2} + \frac{a^2 + b^2}{2} \right) .
\end{align*}
The approximate added-damping tensors are 
\newcommand{\deltaTrap}{\Delta}
\newcommand{\Dvvaa}{\frac{4 w^2}{\deltaTrap}}
\newcommand{\Dvvbb}{a + b  + \frac{(a-b)^2}{\deltaTrap}}
\newcommand{\Dvwbc}{\frac{w}{2}(b-a)\Big( 1 + \frac{a+b}{\deltaTrap}\Big)}
\newcommand{\Dvwbcnew}{\frac{w}{3}(b-a)\Big( 1 -\frac2 \deltaTrap\frac{a^2+ab+b^2}{a+b}\Big)}
\begin{align*}
 & \Dvva = \frac{\mu}{\dn}\, 
      \begin{bmatrix} \Dvvaa & 0 & 0 \\ 0 & \Dvvbb & 0 \\ 0 & 0 & 0 \end{bmatrix},  \quad
 \Dvwa = (\Dwva)^T = \frac{\mu}{\dn}\, \Dvwbcnew\, 
      \begin{bmatrix} 0  & 0 & 0 \\ 0 & 0 & 1 \\ 0 & 0 & 0 \end{bmatrix},  \quad \\
	  & \Dwwa = \frac{\mu}{\dn}\, \frac{w^2}{9}\Big( \frac{a^2+7ab+b^2}{a + b}  + \frac{4}{\deltaTrap}\frac{(a^2+ab+b^2)^2}{(a + b)^2}\Big)
         \begin{bmatrix} 0 & 0 & 0 \\ 0 & 0 & 0 \\ 0 & 0 & 1 \end{bmatrix}, 
\end{align*}
where
\[
  \deltaTrap = \sqrt{ 4 w^2 + (a-b)^2}. 
\]
In this example there can be non-zero entries in $\Dvwa$ and $\Dwva$, which can
be interpreted as follows. The entry $\Dvwa_{2,3}$ represents a force in the $y$-direction, in the equation for $\vs_2$,
due to unequal shear forces arising from rotation about the centre of mass  (e.g.~the shear
force on the left face under a rotation is different from the shear force on the right face if $\at \ne \bt$). 
The entry $\Dwva_{3,2}$ represents a torque about the centre of mass due to motions in the $y$-direction,
again due to the torque that arises from unequal shear forces.

\subsection{L-shaped domain}

\newcommand{\al}{a}
\newcommand{\bl}{b}

\begin{figure}[hbt]
\begin{center}
\resizebox{6cm}{!}{
\begin{tikzpicture}
 \useasboundingbox (0,0.25) rectangle (8,4.5);  
\begin{scope}[xshift=1cm,yshift=0cm,scale=1]
\draw[-,very thick,fill=blue!30] (0,0) -- (6,0) -- (6,2) -- (3,2) -- (3,4) -- (0,4) -- (0,0);
\fill[red] (2.5,1.67) circle (4 pt) node[anchor=north east,black] {$(\frac{5\bl}{6},\frac{5\al}{6})$};
\fill[blue] (0,0) circle (3 pt) node[anchor=north east,black] {$(0,0)$};
\fill[blue] (6,0) circle (3 pt) node[anchor=north west,black] {$(2 \al,0)$};
\fill[blue] (3,2) circle (3 pt) node[anchor=south west,black] {$(\al,\bl)$};
\fill[blue] (0,4) circle (3 pt) node[anchor=south east,black] {$(0,2\bl)$};
\fill[blue] (6,2) circle (3 pt);
\fill[blue] (3,4) circle (3 pt);
\end{scope}
%
%
\end{tikzpicture}
} 
\end{center}
\caption{An L-shaped domain use in the computation of added-damping tensors. The centre of mass is at $(5b/6,5a/6)$. }
\label{fig:addedDampingLshape}
\end{figure}
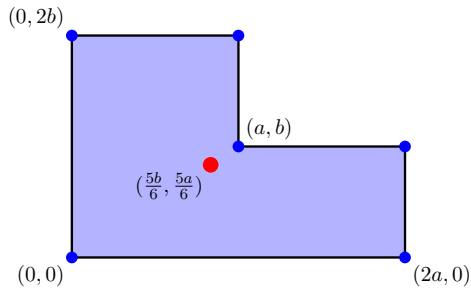

Consider an L-shaped domain of constant density $\rhos$
as shown in Figure~\ref{fig:addedDampingLshape}. 
The moment of inertia about the centre of mass for rotations about the $z$-axis is
\begin{align*}
   I_z = \frac{11}{12} {\rho_b} a b (a^2+b^2).
\end{align*}
The approximate added-damping tensors are 
\begin{align*}
 & \Dvva = \frac{\mu}{\dn}\, 
      \begin{bmatrix} 4b & 0 & 0 \\ 0 & 4a & 0 \\ 0 & 0 & 0 \end{bmatrix},  \quad 
   \Dvwa = (\Dwva)^T = \frac{\mu}{\dn}\, \frac{1}{3} a b \, 
      \begin{bmatrix} 0  & 0 & 1 \\ 0 & 0 & -1 \\ 0 & 0 & 0 \end{bmatrix},  \quad 
   \Dwwa = \frac{\mu}{\dn}\, \frac{25}{9} a b(a+b) 
         \begin{bmatrix} 0 & 0 & 0 \\ 0 & 0 & 0 \\ 0 & 0 & 1 \end{bmatrix}. 
\end{align*}
The non-zero entries in $\Dvwa$ and $\Dwva$ can be interpreted following
a similar explanation as for the trapezoid, i.e.~translations of the body can cause rotational added-damping effects, while
rotations can cause translational added-damping effects. 

\subsection{Ellipse}

The approximate added-damping tensors for a ellipse with boundary $(x/a)^2 + (y/b)^2 =1$,
rotating about the $z$-axis are 
\begin{align*}
 & \Dvva = \frac{\mu}{\dn}\, \begin{bmatrix} \Dvv_{11} & 0 & 0 \\ 0 & \Dvv_{22} & 0 \\ 0 & 0 & 0 \end{bmatrix},  \quad
   \Dvwa = \Dwva = \zerov, \quad
   \Dwwa = \frac{\mu}{\dn}\, \begin{bmatrix} 0 & 0 & 0 \\ 0 & 0 & 0 \\ 0 & 0 & \Dww_{33} \end{bmatrix},\\
&  \Dvva_{11} =-4 a \, \frac{ b^2 K(\eta) - a^2 E(\eta)}{a^2-b^2}, \quad
   \Dvva_{22} = 4 b^2 a \frac{ K(\eta)-E(\eta)}{a^2-b^2}, \quad 
   \Dwwa_{33} = 4 a b^2 K(\eta) , \\
& \eta = \frac{(a^2-b^2)^{1/2}}{a}, \quad K(\eta)= \int_0^1 \frac{1}{\sqrt{ (1-t^2)(1-\eta^2 t^2)}} ~dt , \quad
         E(\eta) = \int_0^1 \frac{\sqrt{1-\eta^2 t^2}}{\sqrt{ 1-t^2}} ~dt
\end{align*}
where $K$ is the complete elliptic integral of the first kind
and  $E$ is the complete elliptic integral of the second kind.

\subsection{Sphere}

The approximate added-damping tensors for a sphere of radius $a$ rotating about its centre are 
\begin{align*}
 & \Dvva = \frac{\mu}{\dn}\, \frac{8}{3} \pi a^2 \, 
      \begin{bmatrix} 1 & 0 & 0 \\ 0 & 1 & 0 \\ 0 & 0 & 1\end{bmatrix},  \quad 
   \Dvwa = \zerov,  \quad 
   \Dwwa =   \frac{\mu}{\dn}\, \frac{8}{3} \pi a^4 \, 
         \begin{bmatrix} 1 & 0 & 0 \\ 0 & 1 & 0 \\ 0 & 0 & 1 \end{bmatrix}. 
\end{align*}

\bibliographystyle{elsart-num}
\bibliography{journal-ISI,jwb,henshaw,henshawPapers,fsi}

\end{document}